\documentclass[graybox,envcountchap,sectrefs]{svmono}

\usepackage{type1cm}         

\usepackage{makeidx}         % allows index generation
\usepackage{graphicx}        % standard LaTeX graphics tool
\usepackage{multicol}        % used for the two-column index
\usepackage[bottom]{footmisc}% places footnotes at page bottom

\usepackage{newtxtext}       % 
\usepackage{newtxmath}       % selects Times Roman as basic font

\usepackage[dvipsnames]{xcolor}
\usepackage{tcolorbox}
\usepackage{thmtools}

\colorlet{LightGray}{White!90!Black}
\colorlet{LightRed}{BrickRed!15}
\colorlet{LightGreen}{Green!15}
\colorlet{LightBlue}{NavyBlue!15}
\colorlet{LightPurple}{Orchid!15}
\colorlet{LightYellow}{Goldenrod!15}

\tcolorboxenvironment{theorem}{colback=LightGray}

\tcolorboxenvironment{lemma}{colback=LightRed}

\tcolorboxenvironment{corollary}{colback=LightPurple}

\tcolorboxenvironment{definition}{colback=LightGreen}

\tcolorboxenvironment{example}{colback=LightBlue}

\tcolorboxenvironment{proposition}{colback=LightYellow}

\usepackage[style=phys]{biblatex}
\usepackage{enumitem}

\makeindex             % used for the subject index
\graphicspath{{./}{images/}}

\begin{document}

\author{Eugene Tan, David Walker, Michael Small, Braden Thorne}
\title{Dynamics, Complexity and Time Series Analysis}
\subtitle{First steps with theory and application}
\maketitle

\frontmatter%%%%%%%%%%%%%%%%%%%%%%%%%%%%%%%%%%%%%%%%%%%%%%%%%%%%%%

%%%%%%%%%%%%%%%%%%%%%%% dedic.tex %%%%%%%%%%%%%%%%%%%%%%%%%%%%%%%%%
%
% sample dedication
%
% Use this file as a template for your own input.
%
%%%%%%%%%%%%%%%%%%%%%%%% Springer %%%%%%%%%%%%%%%%%%%%%%%%%%

\begin{dedication}
Use the template \emph{dedic.tex} together with the Springer document class SVMono for monograph-type books or SVMult for contributed volumes to style a quotation or a dedication\index{dedication} at the very beginning of your book
\end{dedication}

%%%%%%%%%%%%%%%%%%%%%%foreword.tex%%%%%%%%%%%%%%%%%%%%%%%%%%%%%%%%%
% sample foreword
%
% Use this file as a template for your own input.
%
%%%%%%%%%%%%%%%%%%%%%%%% Springer %%%%%%%%%%%%%%%%%%%%%%%%%%

\foreword

%% Please have the foreword written here
Use the template \textit{foreword.tex} together with the document class SVMono (monograph-type books) or SVMult (edited books) to style your foreword\index{foreword}. 

The foreword covers introductory remarks preceding the text of a book that are written by a \textit{person other than the author or editor} of the book. If applicable, the foreword precedes the preface which is written by the author or editor of the book.

\vspace{\baselineskip}
\begin{flushright}\noindent
Place, month year\hfill {\it Firstname  Surname}\\
\end{flushright}

%%%%%%%%%%%%%%%%%%%%%%preface.tex%%%%%%%%%%%%%%%%%%%%%%%%%%%%%%%%%%%%%%%%%
% sample preface
%
% Use this file as a template for your own input.
%
%%%%%%%%%%%%%%%%%%%%%%%% Springer %%%%%%%%%%%%%%%%%%%%%%%%%%

\preface

%% Please write your preface here
% Use the template \emph{preface.tex} together with the document class SVMono (monograph-type books) or SVMult (edited books) to style your preface.

This book has been adapted from a series of lecture notes by the same authors for an undergraduate course in Dynamical Systems at the University of Western Australia. Whilst much of the content has remained unchanged, several accommodations and extensions have been made to suit the form of the text and improve clarity. This has been periodically updated to accommodate variations in special topics in subsequent years both for variety and to reflect current state-of-the-art mathematical findings and methods.

The aim of this text is to provide a linguistically accessible, but comprehensive introduction into a variety of topics in dynamical systems and its applications. Whilst preliminary knowledge of dynamical systems is useful, it is not essential and readers are only assumed to have familiarity with foundational undergraduate mathematics topics of calculus, linear algebra and rudimentary statistics. A variety of extended topics on recent publications and research activities in the field have been included in the last four chapters, which the interested reader may use as an introduction into further reading. A collection of exercises and questions both theoretical and computational are also included in this text. These questions have been curated over the years during which original course material was used.

\vspace{1em}

\textbf{Version History}
\begin{enumerate}
    \item v1.0 - First edition (29/4/2024) E. Tan
    \item v1.1 - Added chapter on reservoir computers, and recurrence plots. (29/5/2024) B. Thorne
    \item v1.2 - Added chapter on dynamical networks (12/5/2025) E. Tan
    \item v1.3 - Added chapter on ordinal time series analysis (17/5/2025) E. Tan
    \item v1.3.1 - Updated foreword (17/5/2025) E. Tan
    \item v1.4 - Updated hand drawn figures to digital (25/6/2025) E. Tan
    \item v1.5 - First round of edits and clean up (1/7/2025) E. Tan
    \item v1.6 - Added exercises at end (2/7/2025) E. Tan
    \item v1.7 - Compiled in book form (2/7/2025) E. Tan
\end{enumerate}

\vspace{\baselineskip}
\begin{flushright}\noindent
Place(s),\hfill {\it Firstname  Surname}\\
month year\hfill {\it Firstname  Surname}\\
\end{flushright}

%%%%%%%%%%%%%%%%%%%%%%acknow.tex%%%%%%%%%%%%%%%%%%%%%%%%%%%%%%%%%%%%%%%%%
% sample acknowledgement chapter
%
% Use this file as a template for your own input.
%
%%%%%%%%%%%%%%%%%%%%%%%% Springer %%%%%%%%%%%%%%%%%%%%%%%%%%

\extrachap{Acknowledgements}

Use the template \emph{acknow.tex} together with the document class SVMono (monograph-type books) or SVMult (edited books) if you prefer to set your acknowledgement section as a separate chapter instead of including it as last part of your preface.

\tableofcontents

\mainmatter%%%%%%%%%%%%%%%%%%%%%%%%%%%%%%%%%%%%%%%%%%%%%%%%%%%%%%%

\chapter{Dynamical Systems Theory}\label{chap:DynamicalSystemsTheory}

Dynamics refers the occurrence of change both large and small over time. This term, owing to its breadth, is highly relevant and descriptive of a wide range of phenomena that is observed in the physical world. From the foundational principles of classical Newtonian physics to chaos, complex networks and chimeras, the presence of dynamics is characteristically present and has been the topic of persistent inquiry. If one wishes to understand the governing rules of observed dynamics, it is inevitable then that the inquisitive mind will encounter the field of dynamical systems theory in the hopes that a degree of logical structure may be established.

Naturally, the study of dynamical systems theory aims to understand and codify the mathematical rules that govern dynamics. This is done with the aim of providing guidance and tools for better understanding the observed dynamics that one encounters in the application of their work. As we begin our discussion on this topic, it is no doubt necessary that we must first provide a definition for the object of our study.

\section{Dynamical Systems}
\begin{definition} \textbf{-- Dynamical system}\\
    A dynamical system is a triple $(T,M,\Phi)$ where:
    \begin{itemize}
        \setlength\itemsep{0em}
        \item $T$: time index
        \item $M$: state space
        \item $\Phi$: evolution operator
    \end{itemize}
    such that
    $$\Phi: U\subset(T\times M) \to M$$
\end{definition}

In other words, dynamical systems are mathematical objects consisting some state located in state space whose position evolves across time (or any ordered variable) based on the rules given by the evolution operator $\Phi$. Given an initial condition $x(0) = x_0$, the evolution operator is given by
\begin{subequations}
    \begin{gather}
        \Phi_{t}{(x(0))}= x(t)\\
        \Phi_{t+s}(x) = \Phi_{t}\circ \Phi_{s} (x)= \Phi_{s}\circ \Phi_{t}(x)
    \end{gather} 
\end{subequations}

Here, $x(t)\in M$ represents a trajectory (orbit) of the dynamical system in phase space as a function of time. An alternative representation of a dynamical system can be expressed in as an ordinary differential equation (ODE):
\begin{equation}
    \dot{x} = F(x,t)
\end{equation}
A dynamical system is said to be autonomous if $F$ is time-independent (i.e. $F(x,t) = F(x)$). Therefore, realisations of trajectories are only dependent on the current state the evolution operator. Solutions to the above ODE initial value problem yields trajectories $x(t)$. There are two main ways in which a dynamical system can be visualised:
\begin{enumerate}
    \setlength\itemsep{0em}
    \item As a vector field $\dot{x} = F(x)$
    \item As a trajectory in phase space $x(t) \in M$ (alternatively $(x,t)\in M\times T$)
\end{enumerate}

\begin{figure}
    \centering
    \includegraphics[width = 0.85\textwidth]{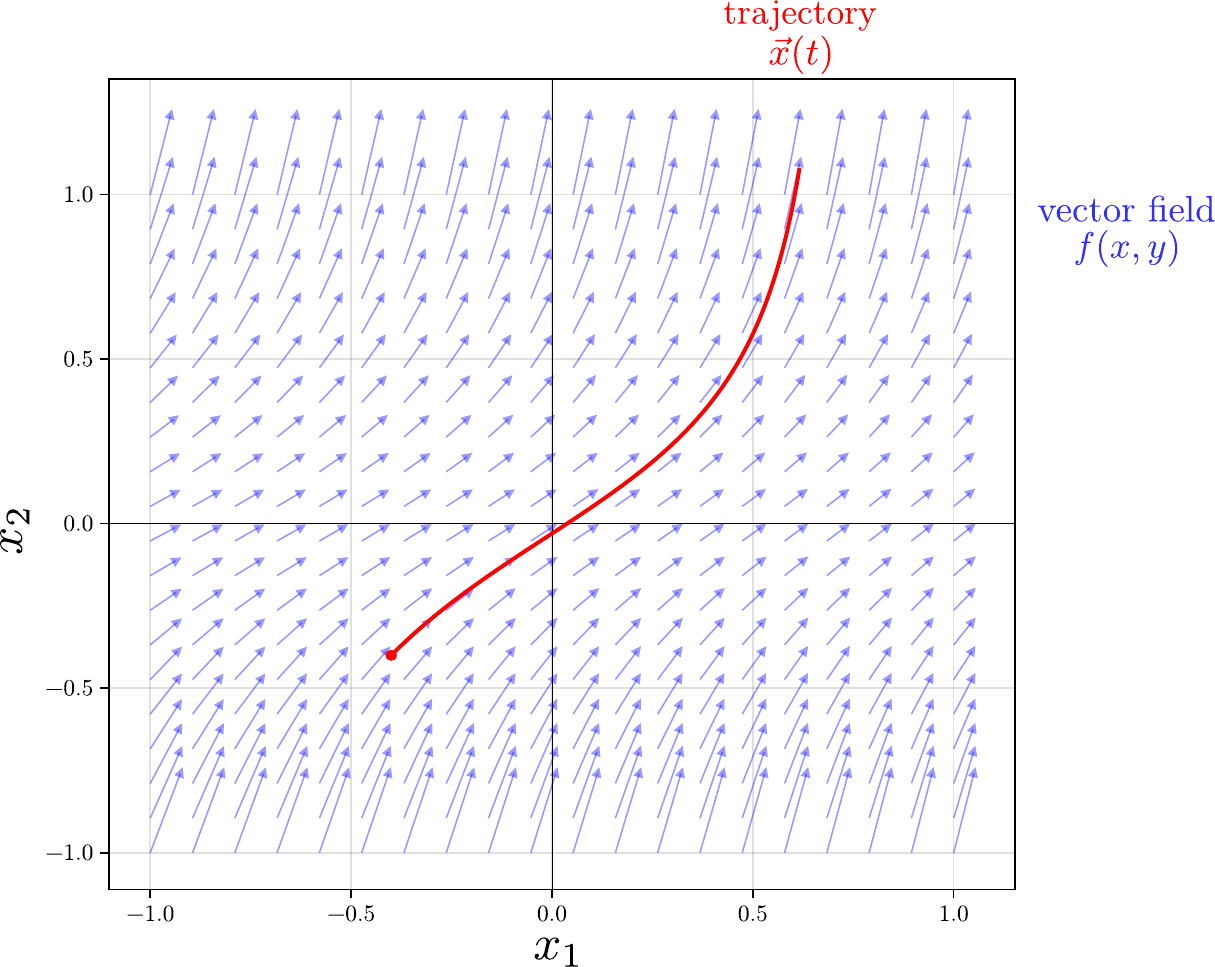}
    \caption{Vector field and trajectory of dynamical system}
    \label{fig:vector_field}
\end{figure}

\section{Flows and Linearisation}
In most cases, one is usually interested in the long-term behaviour of typical solutions of the differential equation. Information about this behaviour can be obtained by analysing the fixed points of the differential equation and its iterates. Specifically, the long term behaviour can be determined by the looking at the stability of the fixed points.

\begin{definition} \textbf{-- Fixed points}\\
    Given a general ODE $\dot{x}=f(x)$ defined in the space $\mathbb{R}^{n}$. Any point $x^{*}\in \mathbb{R}^n$ which satisifies $f(x^*)=0$ is a fixed point.
\end{definition}

At a fixed point $x^*$, the property that $f(x^*)=0$ means that trajectories have zero velocity and hence cannot depart from the fixed point for all times $t$. Note that for continuous and differentible ($C^1$) systems, this also applies for past histories $t$ as well. The above definition allows us to easily calculate the location of fixed points. However, it does not provide any information about the stability and behaviour in the neighbourhood around the fixed points.

\begin{figure}
    \centering
    \includegraphics[width = 0.75\textwidth]{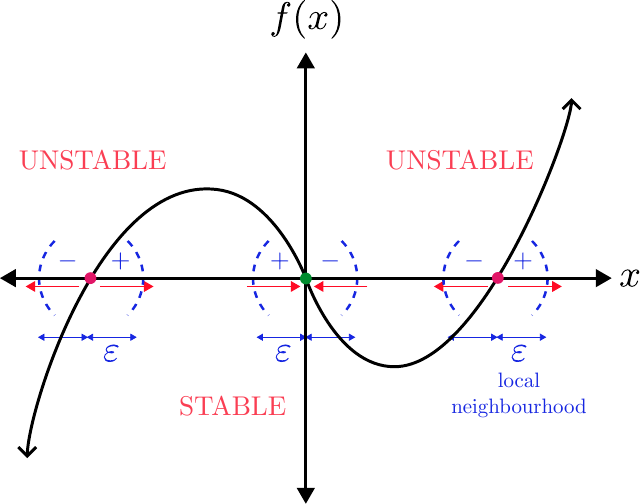}
    \caption{Phase space diagram of $\dot{x} = x(x+1)(x-1)$. Interesting behaviour occurs in the local region near fixed points. }
    \label{fig:phase_space}
\end{figure}

Consider the simplest case of a 1D flow (i.e. $x(t)\in \mathbb{R}$), and let the corresponding ODE be
\begin{equation}\label{eq1:cubic_ode}
    \dot{x} = f'(x) = x(x+1)(x-1).
\end{equation}
Plotting the phase plot of this system reveals three fixed points at $x^*=-1,0,1$ with the following behaviour (see Figure \ref{fig:phase_space})
\begin{itemize}
    \setlength\itemsep{0em}
    \item At $x = -1$, trajectories that start left of the fixed point accelerate to the left. Similarly trajectories that start right of the fixed point accelerate to the right. $\to$ \textbf{UNSTABLE}
    \item At $x = 0$, trajectories that start left of the fixed point move to right. Similarly, trajectories that start right of the fixed point move left. $\to$ \textbf{STABLE}
    \item At $x = 1$, behaviour is same as $x = -1$ $\to$ \textbf{UNSTABLE}
\end{itemize}

\vspace{1em}
In general, for differentiable ODEs, the stability of fixed points can be determined using the derivative evaluated at the fixed point.
\begin{itemize}
    \item $f'(x^*)>0$ $\implies$ \textbf{UNSTABLE}
    \item $f'(x^*)<0$ $\implies$ \textbf{STABLE}
\end{itemize}

\vspace{1em}
In the case of Eq \ref{eq1:cubic_ode}, $f$ is a simple function whose derivative is easy to evaluate. Note that $f'(x)$ is also nonlinear (a quadratic). The method of using derivatives to infer stability can be simplified by way of a Taylor series of $f$ centred around each fixed point $x^*$,
\begin{equation}
    \dot{x} = f'(x^*)(x-x^*)+R_{1}(x,x^*).
\end{equation}

If $f$ is $C^1$ (i.e once differentiable function), the size of the higher order terms $\vert\vert R_{1}(x,x^{*})\vert\vert\to 0$ faster than $||x-x^{*}||$ as $x\to x^{*}$. In other words, the dynamics around the neighbourhood of the fixed points are dominated by the first derivative $f'$. Therefore, one can use a change of coordinates $u = x-x^*$ to \textit{linearise} the nonlinear ODE to yield a locally linear approximation of the dynamics,
\begin{equation}
    \dot{u}=Df(x^*)u,
\end{equation}
where $Df(x) = f'(x)$. Letting $Df(x^{*}) = A$, the corresponding solutions of the approximate linear system for a trajectory with initial conditions $u(0) = u_0$ are be given exactly as,
\begin{equation}\label{eq1:exp_linearisation}
    u(t) = u_{0}e^{tA}.
\end{equation}
In Eq. \ref{eq1:exp_linearisation}, the stability of a the fixed point is entirely determined by the sign of $A$. This result is also true in general for $x\in\mathbb{R}^n$, the final linear ODE is expressed in terms of eigenvectors and eigenvalues of the matrix $Df(x^*)$
\begin{equation}
    \vec{u}(t) = \vec{v}_{1}e^{t\lambda_{1}} + \vec{v}_{2}e^{t\lambda_2}...+\vec{v}_{1}e^{t\lambda_n}.
\end{equation}

\begin{theorem} \textbf{-- Linearisation theorem}\\
    Let $f: \mathbb{R}^{n}\to \mathbb{R}^n$ be a $C^1$ vector field, and $x^*$ is a fixed point. If $A = Df(x^{*})$ has all eigenvalues with non-zero real part, then $\exists \epsilon>0$ such that $\forall x$, $||x-x^*||<\epsilon$, $\dot{x}=f(x)$ is topologically equivalent to $\dot{y}=Ay$
\end{theorem}

You might find it a useful a exercise to derive this expression on your own. Note that this formulation is only useful if $\lambda_{i}\neq 0$ for all $i$.  If this is true, the corresponding fixed point is said to be a \textbf{hyperbolic}. 

\begin{definition} \textbf{-- Stable and unstable manifolds}\\
    Consider $\dot{x} = f(x)$ and $x\in \mathbb{R}^n$, where $x^*$ is a hyperbolic fixed point centred at the origin, and $A = Df(x^*)$ is the Jacobian matrix of $f$.\\
    
    Let $\lambda_{1},...,\lambda_{n}$ and $\vec{v}_{1},...,\vec{v}_{n}$ be the ordered eigenvalues and eigenvectors of $A$ such that 
    $$\lambda_{1}<\lambda_{2}<...<\lambda_{i}<0<\lambda_{j},...\lambda_{n}.$$
    Then

    \begin{itemize}
        \item $E_{s}= \text{span}(\vec{1},...,\vec{v}_i)$ is the stable manifold
        \item $E_{u}= \text{span}(\vec{v}_j,...,\vec{v}_n)$ is the unstable manifold
    \end{itemize}
\end{definition}

\begin{figure}
    \centering
    \includegraphics[width = \textwidth]{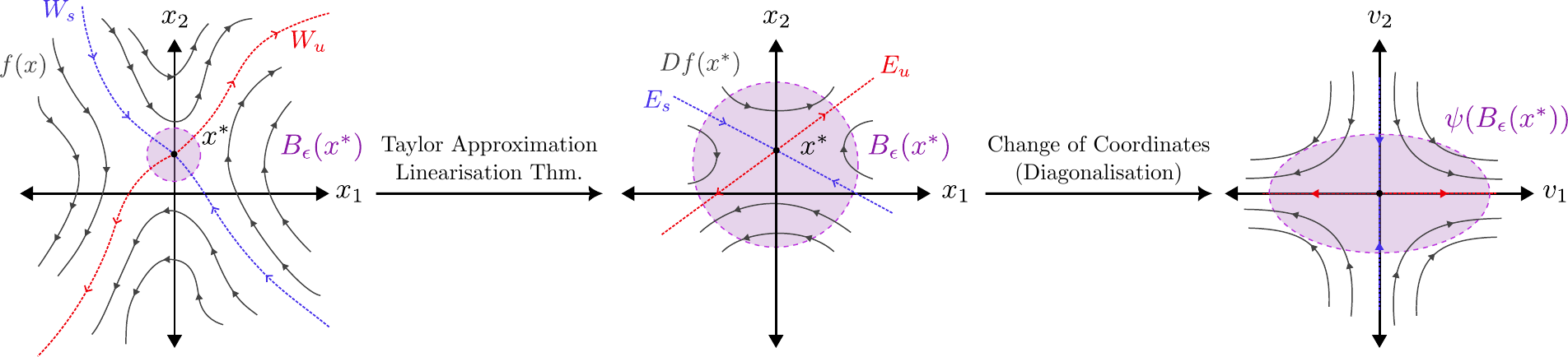}
    \caption{Linearisation about a fixed point}
    \label{fig:linearisation}
\end{figure}

The stable and unstable manifolds define the tangent directions in which trajectories approach the fixed point. Hyperbolic fixed points can be largely categorised into the following cases:
\begin{itemize}
    \item Stable nodes $\lambda_{1} <...<\lambda_n < 0$
    \item Unstable nodes $0<\lambda_{1} <...<\lambda_n$
    \item Saddle nodes $\lambda_{1}<...<\lambda_{p}<0<\lambda_{q},...\lambda_{n}$
    \item Stable focus $\text{Re}(\lambda_{1}) <...<\text{Re}(\lambda_n) < 0$ and $\text{Im}(\lambda_i)\neq0$
    \item Unstable focus $\text{Re}(\lambda_{1}) >...>\text{Re}(\lambda_n) > 0$ and $\text{Im}(\lambda_i)\neq0$
\end{itemize}

\begin{figure}
    \centering
    \includegraphics[width = 0.7\textwidth]{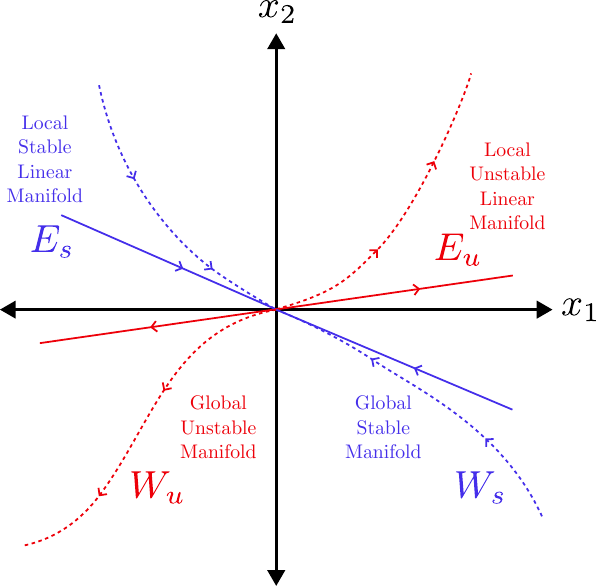}
    \caption{Global and local manifolds}
    \label{fig:global_local}
\end{figure}

\vspace{1em}
Accompanied with the above classifications are some theoretical guarantees for dynamics of ODEs pertaining to existence and uniqueness.

\begin{theorem} \textbf{-- Existence and uniqueness for flows}\\
    If $f: \mathbb{R}^{n}\to \mathbb{R}^n$ is a vector field of class $C^1$, then for all $x_{0} \in \mathbb{R}^{n}$, $\exists T>0$ such that $\dot{x}=f(x)$ has a unique solution $\psi(t,x_0)$ for $-T<t<T$ where $\psi(0,x_{0})=x_{0}$. Hence:
    \begin{enumerate}
        \item Two different states cannot evolve to the same state at the same time
        \item It two different states do evolve to the same state, they are on the same orbit
        \item Trajectories never intersect, cross, nor merge
    \end{enumerate}
\end{theorem}

This theorem guarantees that as long as $f$ is well behaved (once differentiable), there will always exist a unique trajectory solutions that passes through each point locally. This condition does not hold for cases where $f$ is not $C^1$. To illustrate this, we consider a pathological example that violates the conditions for the uniqueness of solutions.

Consider the following initial value problem,
\begin{equation}
    \dot{x} = \sqrt{|x|},\quad x(t_0) = x_0.
\end{equation}
Clearly $f(x)$ is not differentiable at the point $x=0$ and thus we would expect that the existence and uniqueness guarantees fail to apply. To demonstrate this, we consider a trajectory that starts at the fixed point $x_0 = 0$. Explicit solutions for the trajectory can be evaluated by separating the variables as follows,
\begin{align*}
    \frac{1}{\sqrt{x}} \, dx &= dt\\
    x(t) &= \left( \frac{t+\mathcal{C}}{2} \right) ^2.
\end{align*}
Substituting the initial condition $x(t_0) = 0$ to solve for the integration constant $\mathcal{C}$ yields the following trajectory solution
\begin{equation}
    x(t) = \frac{(t-t_0)^2}{4}.
\end{equation}
The above results show a contradiction where the same initial condition $x(t_0) = 0$ evolve into two different trajectories. In the first case $\dot{x} = 0$ is a fixed a point and thus should remain $x = 0$ for all $t>t_0$. However, the explicit trajectory solution show that for $t>t_0$ the trajectory $x(t)$ spontaneously leaves the fixed point $x = 0$ and follows a quadratic trajectory. This non-uniqueness of solutions is a direct consequence of the non-differentiability of the $f$ at the fixed point. 

\begin{equation}
    \dot{x} = \frac{1}{x^2}
\end{equation}

\begin{theorem} \textbf{-- Continuity of solutions}\\
    \label{thm1:continuity}
    If $f: \mathbb{R}^{n}\to \mathbb{R}^n$ is a vector field of class $C^1$, then for all $x_{0} \in \mathbb{R}^{n}$, $\exists \epsilon>0$ and $T>0$ such that $\dot{x}=f(x)$ has solutions $\phi(t,x)$ for $-T<t<T$ and $||x-x_0||<\epsilon$ where $\phi(0,x)=x$, and $\psi(t,x)$ is continuous in $t$ and $x$. 
\end{theorem}
Theorem \ref{thm1:continuity} provides guarantees that within local neighbourhoods, states evolve smoothly and deform continuously (i.e. no jumps in position). It is important to emphasise that these are only local in space and time and do not necessarily guarantee the solutions exist for all $T>0$ or $\epsilon>0$. This can be seen with the following example system
\begin{equation}
    \dot{x} = x^2, \quad x(t_0) = x_0.
\end{equation}
As the ODE is $C^1$, we can guarantee that trajectories exist for all initial conditions $x_0$ and are unique. These solution trajectories can be explicitly solved for by separating the variables.
\begin{equation}
    x(t) = \frac{1}{\frac{1}{x_0}+(t_0-t)}.
\end{equation}
We can see that there is a singularity at $t^*=\frac{1}{x_0}+t_0$ and trajectories move to infinity in finite time (\textit{finite time blow-up}). As a result, the solutions for $t<t^*$ are not continuous with $t>t^*$.

\begin{definition} \textbf{-- Topological conjugacy (equivalence)}\\
    Let there be two dynamical systems
    \begin{itemize}
        \item $\dot{x}=f(x), x\in \mathbb{R}^n$ with solutions $q(t,x)$
        \item $\dot{y}=g(y), y\in \mathbb{R}^n$ with solutions $r(t,y)$
    \end{itemize}
    These two dynamical systems are \textbf{equivalent (topologically conjugate)} on a domain $U\subseteq \mathbb{R}^n$ if there exists a change of coordinates $\Psi(x)$ that is a differentiable one-to-one mapping $\Psi:U \to \mathbb{R}^n$, such that $r(t, \Psi(x))=\Psi(q(t,x))$ for all $x$ and $t$ where $q(t,x)$ exists in $U$.
\end{definition}

\begin{theorem} \textbf{-- Fundamental theorem of flows}\\
    If $f: \mathbb{R}^{n}\to \mathbb{R}^n$ is a vector field of class $C^1$, $x_{0} \in \mathbb{R}^{n}, \epsilon >0, f(x) \neq 0$ for $||x-x_0||<\epsilon, v \in \mathbb{R}^{n}, v\neq 0$, then there exists a differentiable change of coordinates $y=\Psi(x)$ for $||x-x_0||<\epsilon$ such that $\dot{x}=f(x)$ is equivalent to $\dot{y} = v$. Hence:
    \begin{enumerate}
        \item Trajectories far away from fixed points all look like parallel flows with constant speed. (i.e. straight lines).
        \item Interesting behaviour only occurs near the fixed points.
    \end{enumerate}
\end{theorem}

\section{Centre Manifolds}
The linearisation approach for determining stability only works for fixed points that are hyperbolic (i.e. non-zero eigenvalues). In order to determine the stability of non-hyperbolic fixed points, we require a little bit more sophistication.

\begin{definition} \textbf{-- Centre manifolds (eigenspaces)}\\
    Let $\dot{x}=f(x)$ be a nonlinear dynamical system in $\mathbb{R}^{s+u+c}$, with an associated linearised flow given by $\dot{y} = Ay$. Let $\Lambda_{-} = \{\lambda_{1}^{-},...,\lambda_{s}^{-}\}$, $\Lambda_{+} = \{\lambda_{1}^{+},...,\lambda_{u}^{+}\}$, $\Lambda_{0} = \{\lambda_{1}^{0},...,\lambda_{c}^{0}\}$ be the ordered set of eigenvalues with positive, negative and zero real parts respectively. Consequently, let $V_{-} = \{ \vec{v}_{1}^{-}, ..., \vec{v}_{s}^{-} \}$, $V_{+} = \{ \vec{v}_{1}^{+}, ..., \vec{v}_{u}^{+} \}$ and $V_{0} = \{ \vec{v}_{1}^{0}, ..., \vec{v}_{c}^{0} \}$ be the corresponding eigenvectors. The stability manifolds are thus defined as 
    \begin{itemize}
        \item Stable manifold $E^{s}= \text{span}(V_{-})$
        \item Unstable manifold $E^{u}= \text{span}(V_{+})$
        \item Centre manifold $E^{c}= \text{span}(V_0)$
    \end{itemize}
\end{definition}

From linearisation theory, we know that the stability of hyperbolic fixed points can be determined by the stable and unstable manifolds. However, for non-hyperbolic fixed points, linearisation is too coarse a reduction to be able to make any useful or accurate conclusions. Similar to how linearisation attempts to reduce dynamics to be invariant on the fixed point, centre manifold theory attempts to reduce dynamics to be localised with respect to some invariant nonlinear \textbf{centre manifold} passing through the fixed point.\\

\begin{figure}
    \centering
    \includegraphics[width = 0.85\textwidth]{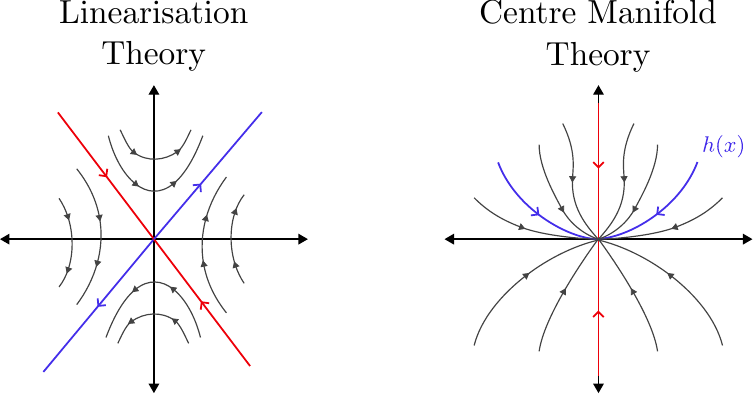}
    \caption{Linearisation vs. centre manifold theory}
    \label{fig:centre_manifold}
\end{figure}

Consider a general dynamical system defined with the following vector fields
\begin{subequations}
    \label{eq1:centre_manifold_ODE1}
    \begin{align}
        \dot{x} &= Ax + f(x,y),\\
        \dot{y} &= By + g(x,y)
    \end{align}
\end{subequations}
where $(x,y) \in \mathbb{R}^{c}\times \mathbb{R}^s$ with the following boundary conditions:
\begin{subequations}
    \label{eq1:centre_manifold_ODE2}
    \begin{gather}
        f(0,0) = 0, \qquad Df(0,0) = 0,\\
        g(0,0) = 0. \qquad Dg(0,0) = 0.
    \end{gather}
\end{subequations}
This condition focuses our attention on a fixed point with zero first derivatives located at the origin (Note that any fixed point in a dynamical system can be written in terms of the form of Eq \ref{eq1:centre_manifold_ODE1} and \ref{eq1:centre_manifold_ODE2} by applying a linear change of coordinates. In the above system, the first $c$ components corresponding to matrix $A$ describe the centre manifold (i.e. eigenvalues with zero real part), Matrix $B$ describes the components with eigenvalues of negative real parts. Without loss of generality, one could also generalise the results to have all non-zero eigenvalue components. The result would correspond to a saddly point, which is not of particular interest in the current discussion.

\begin{definition} \textbf{-- Centre manifold}\\
    An invariant manifold with respect to the above system is a centre manifold if it can be locally represented as follows:
    $$W^{c}(0)= \{ (x,y) \in \mathbb{R}^{c} \times \mathbb{R}^{s} | y = h(x), |x| <\delta, h(0) = 0, Dh(0) = 0 \}$$
    for $\delta$ sufficiently small. See Figure \ref{fig:centre_manifold}
\end{definition}

To find the stability of the fixed point, we need to find the stability of the centre manifold. Conveniently, \cite{carr2012applications} provides several useful theorems for this.

\begin{theorem} \textbf{-- Existence of centre manifolds}\\
    For the system given by Eqs \ref{eq1:centre_manifold_ODE1} and \ref{eq1:centre_manifold_ODE2}, there exists a $C^r$ centre manifold. The dynamics of the system restricted to the centre manifold is, for $u$ sufficiently small, given by the following $c$-dimensional vector field,
    $$\dot{u} = Au + f(u,h(u)), \quad u \in \mathbb{R}^c.$$
\end{theorem}

\begin{theorem} \textbf{-- Stability of centre manifolds}\\
    For the dynamical system restricted to the centre manifold.
    \begin{enumerate}
        \item Suppose the zero solution (fixed point) is stable (likewise asymptotically stable/unstable), then the zero solution of the original dynamical system is also stable (asymptotically stable/unstable)
        \item Suppose the zero solution of the restricted system is stable. Then if $(x(t), y(t))$ is a solution of the original system with $(x(0), y(0))$ sufficiently small, there is a solution $u(t)$ of the restricted system such that as $t\to \infty$

        \begin{align*}
            x(t) &= u(t) + \mathcal{O}(e^{-\gamma t}),\\
            y(t) &= h(u(t)) + \mathcal{O}(e^{-\gamma t})
        \end{align*}
        where $\gamma >0$ is a constant.
    \end{enumerate}
\end{theorem}

The first statement of the theorem states that it is possible to infer the stability of the original system by determining the stability of the centre manifold. Namely, the centre manifold approximation is sufficient for representing the dynamics. The second statement states that this result holds as long as one restricts the analysis to a small enough neighbourhood of the fixed point. If this is the case, the our dynamics is well represented by the centre manifold, with an exponentially decaying error as $t\to \infty$.

\section{Calculating Centre Manifolds}

We are now well equipped to begin the calculation of the centre manifold $h(x)$. Once this is done, we can use the results to subsequently infer the stability of non-hyperbolic fixed points. \\

Suppose we have an invariant centre manifold
\begin{equation}
    W^{c}(0) = \{ (x,y) \in \mathbb{R}^{c} \times \mathbb{R}^{s} | y = h(x), |x| <\delta, h(0) = 0, Dh(0) = 0 \}
\end{equation}

\begin{enumerate}
    \item The $(x,y)$ coordinates of any point on $W^{c}(0)$ must satisfy
    $$y=h(x)$$
    \item Taking the time derivative implies that the $(\dot{x}, \dot{y})$ coordinates of any point on $W^{c}(0)$ must satisfy
    $$\dot{y} = Dh(x)\dot{x}$$
    \item For the original dynamical system,
    \begin{align*}
        \dot{x} &= Ax + f(x,h(x)),\\
        \dot{y} &= Bh(x) + g(x,h(x))
    \end{align*}
    substituting the time derivative $\dot{y}$ yields the following condition
    \begin{equation}
        \label{eq:centre_manifold_condition}
        \mathcal{N}(h(x)) \equiv Dh(x)[ Ax +f(x,h(x)) ] - Bh(x) - g(x,h(x)) = 0.
    \end{equation}
    This expression $\mathcal{N}$ can be solved by assuming a power series expansion for $h(x)$ \cite{wiggins2003introduction}.
\end{enumerate}

As an illustrative example of the above fixed point stability calculation for centre manifolds, we consider the following system taken from Wiggins et al. \cite{wiggins2003introduction},
\begin{subequations}
    \label{eq1:centre_manifold_ode_example}
    \begin{align}
        \dot{x} &= \frac{x}{2} + y + x^{2}y,\\
        \dot{y} &= x + 2y + y^2.
    \end{align}
\end{subequations}
For this example, we seek to determine the stability of the fixed point located at the origin $(x,y)=(0,0)$. A first attempt to determine the stability by linearising about $(0,0)$ provides the following Jacobian matrix,
\begin{gather}
        Df(x,y) = \begin{bmatrix}
        \frac{1}{2}+2xy & 1+x^2\\
        1 & 2+2y
    \end{bmatrix},\\
    Df(0,0) = A =\begin{bmatrix}
        \frac{1}{2} & 1\\
        1 & 2
    \end{bmatrix}.
\end{gather}
The eigenvalues and eigenvectors of the matrix $A$ are $(\lambda_1,\lambda_2)=(0,5/2)$ and $v_1 =(-2,1), v_2=(1,2)$ respectively. From this, we can conclude that the fixed point is indeed non-hyperbolic and that there is an unstable manifold in the direction of $v_2$. However, the linear approximation does not shed any light on the stability in the $v_1$ direction. To determine this, we first apply a change of coordinates using the calculated eigenvectors,
\begin{subequations}
    \begin{align}
        x &= -2u+v,\\
        y&=u+2v,
    \end{align}    
\end{subequations}
where $u$ and $v$ correspond to the centre manifold direction and unstable directions respectively. Thus the original ODE in Eq. \ref{eq1:centre_manifold_ode_example} can be expressed in terms of $u$ and $v$,
\begin{subequations}
    \begin{align}
        \dot{u} &= \frac{2}{5}(-2u+v)^{2}(u+2v) + \frac{1}{5}(u+2v)^{2} = g_1(u,v),\\
        \dot{v} &= \frac{5}{2}v + \frac{1}{5}(-2u+v)^{2}(u+2v) + \frac{2}{5}(u+2v)^2 = g_2(u,v).
    \end{align}
\end{subequations}
We choose a polynomial approximation to represent the centre manifold,
\begin{equation}
    h(u) = au^2 + bu^3 + \mathcal{O}(u^4),
\end{equation}
where $\mathcal{O}(u^4)$ are 4\textsuperscript{th} order or higher terms which are negligible for small $u$ (i.e. close to the fixed point). Substituting into the condition given by Eq. \ref{eq:centre_manifold_condition} yields the following equation,
\begin{align}
    (2au+bu^2+...)[ 0 &+ \frac{2}{5}(-2u + au^2 + bu^3 +...)^{2}(u + 2au……2 + 2bu^3 + ...) \notag\\
    &+ \frac{1}{5}(u+2au^2+2bu^3+...)^2 ] =0.
\end{align}
Removing all higher order terms $\mathcal{O}(u^4)$ and equating coefficients for $u^2$ and $u^3$ terms yields the values $a = -\frac{4}{25}$ and $b = -\frac{116}{625}$, and thus the centre manifold is approximated as
\begin{equation}
    h(u) = -\frac{4}{25}u^2-\frac{116}{625}u^3 + \mathcal{O}(u^4).
\end{equation}
Substituting $h(u)$ into $g_1(u,h(u))$ to restrict our analysis to the dynamics on the centre manifold provides the following simplified equation,
\begin{equation}
    \dot{u} = \frac{1}{5}u^2 - \frac{8}{5}u^3 - \frac{12}{5}au^4+\mathcal{O}(u^5),
\end{equation}
which is stable for $u<0$ and unstable for $u>0$. Hence, we can conclude that the fixed point at the origin is unstable.

\section{Maps}
As previously discussed, flows can be understood as the family of continuous solution trajectories for an associated ordinary differential equation (ODE),
\begin{equation}
    \dot{x} = f(x),\quad x\in M\subseteq \mathbb{R}^n.
\end{equation}

Similarly, maps are a discrete analogue of flows that can be understood as the family of discrete solution sequences for an associated difference equation (DE),
\begin{equation}
    x_{n+1} = f(x_{n}).
\end{equation}
Many of the ideas for flows carry over to maps, with definitions provided below
\begin{definition} \textbf{ - Fixed points and stability of maps}\\
    Let there be a map given by $x_{n+1} = f(x_{n})$, where $x_{n}\in M$. A fixed point $x^*$ is any solution that satisfies,
    $$x^{*}=f(x^*)$$

    The stability of fixed points are determined by the first magnitude of $f'(x^{*})$,
    \begin{itemize}
        \item $|f'(x^{*}|<1\implies$STABLE
        \item $|f'(x^{*}|>1\implies$UNSTABLE
    \end{itemize}
\end{definition}

\subsection{Poincar\'{e} Maps}
Continuous flows may be linked with the dynamics of a discrete map by defining a Poincar\'{e} surface section $S$. To do so:
\begin{enumerate}
    \item Define a trajectory $x(t) \in \mathbb{R}^n$ and a open section $S\in \mathbb{R}^n$ where $S$ is of dimension $n-1$. 
    \item Let $t_{1}< t_{2}<...$ be the times where $x(t)$ intersects with $S$ (i.e. $x(t_{i}) \in S$)
    \item Define a new sequence consisting of values given by
    $$y_{n} = x(t_{n})$$
\end{enumerate}

In essence, Poincar\'{e} sections allow the dynamics from continuous flows to be represented in terms of discrete maps. This, as we shall see in later chapters, is a useful tool that allows results for maps to be easily extneded to the study of continuous flows as well. 

\begin{figure}
    \centering
    \includegraphics[width = 0.6\textwidth]{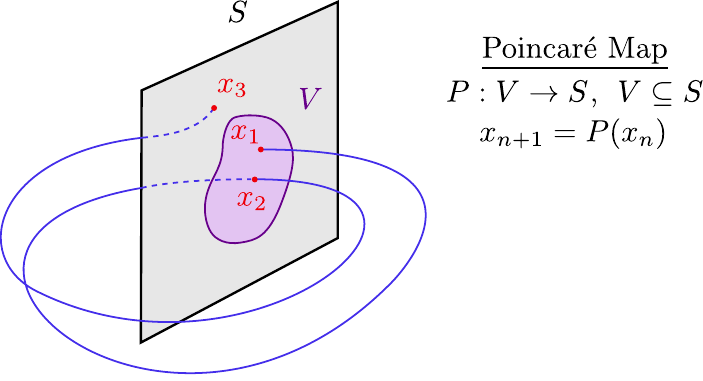}
    \caption{Construction of a Poincar\'{e} map}
    \label{fig:poincare_map}
\end{figure}

\section{Periodicity and Chaos}

\subsection{Periodic Behaviour and Symbolic Dynamics}

Consider the case of the logistic map:
\begin{equation}
    x_{n+1} = rx_{n}(1-x_{n})
\end{equation}
For parameter $r \in [0,4]$. Simple calculations show that there exists two fixed points at $x_{1}= 0$ and $x_{2}=\frac{r-1}{r}$. The stability of these solutions are determined by the magnitude of $f'(x,r)$ that is in turn dependent on $r$,
\begin{equation}
    |f'(x,r)| = |r(1-2x)|.
\end{equation}

\textit{\textbf{Exercise: }Determine the existence and stability of these fixed points for various values of $r$.}\\

From the above, what can we say about the existence of periodic points?

\begin{definition} \textbf{-- Periodic points}\\
    For a map given by function $x_{n+1} = f(x_n)$, a point $a$ is periodic with period $n\in \mathbb{N}$ if,
    $$a = f^{(n)}(a),\quad a \neq f^{(i)}(a) \quad \forall i\neq kn, k\in \mathbb{Z}$$
    Similarly, for a flow with evolution operator $\Phi_t$, a point $a$ is periodic with period $T\in \mathbb{R}, T>0$ if,
    $$a = \Phi_{T}(a), \quad a \neq \Phi_{t(a)} \quad \forall t \neq kT,k \in \mathbb{Z}$$
\end{definition}

\begin{figure}
    \centering
    \includegraphics[width = 0.85\textwidth]{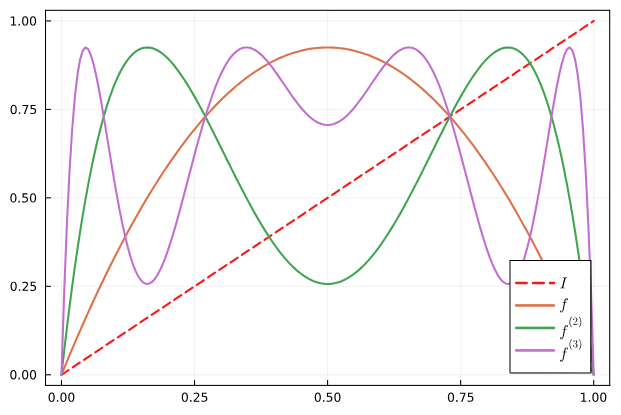}
    \caption{Iterates of the logistic function}
    \label{fig:logistic_iterates}
\end{figure}

Observe that $n$-periodic points correspond to fixed points of $f^{(n)}$. The \textit{trajectory} (sequence) of observations all lie somewhere along the continuous interval $[0,1]$. An alternative way of observing the dynamics would be to see if a given value in the sequence lies on the left half of the interval $x_n<0.5$, or the right half $x_{n}\geq0.5$. By partitioning the sequence in this way, we can construct another sequence $y_{n}$ consisting of symbols (say $L$ and $R$),

\begin{equation}
    y_{n}=\begin{cases} L, \quad x_{n}<0.5 \\
R, \quad x_{n} \geq 0.5\end{cases}
\end{equation}
As it turns out, this representation of the trajectory consisting of discrete symbols (rather than continuous values) retains all the dynamical properties of the logistic map. Periodic dynamics in $x_n$ correspond to recurring patterns in $y_{n}$. The method of partitioning state space to produce symbolic dynamics is a useful tool in reducing the complexity of a trajectory whilst preserving the dynamical behaviour of the system. Symbolic dynamics also forms the foundation for explaining the emergence of chaotic behaviour as we will show later. However, the choice of methof for paritioning state space is not always clear and is further discussed in Chapter \ref{chap:OrdinalTSA}.\\

Returning to original question, what can we say about the existence of periodic points for different values of $r$? To do this, we can plot the bifurcation diagram of the logistic map.\\

\begin{figure}
    \centering
    \includegraphics[width = 0.95\textwidth]{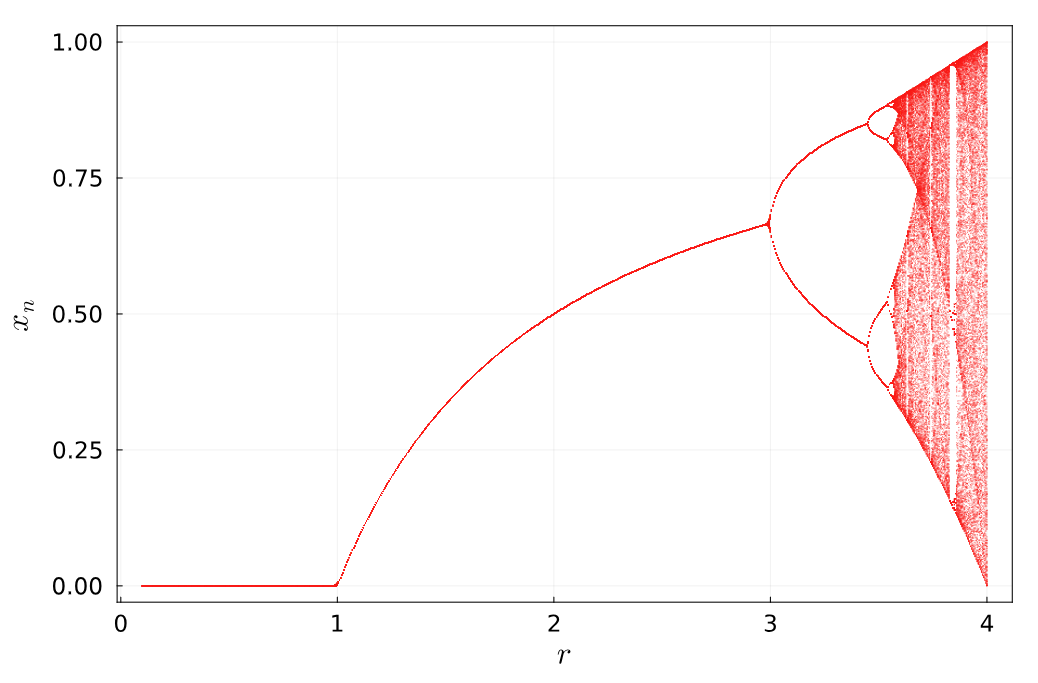}
    \caption{Bifurcation diagram of logistic map}
    \label{fig:logistic_bif}
\end{figure}

From the bifurcation diagram, periodic points only exist for $r>3$. Increasing values of $r$ result in the formation of longer period (doubling in each case) until approximately $r = 3.7$ where trajectories appear to cover an interval of the state space, indicating chaotic behaviour.

\subsection{Period Doubling and Chaos}

Up to now, we have alluded to brief references on the topic of \textbf{chaos} with no clarified definition. There are numerous definitions/properties for chaos in the dynamical systems literature. Below are some that are commonly encountered:
\begin{itemize}
    \item Sensitive dependence on initial conditions
    \item Deterministic but not predictable
    \item Positive Lyapunov exponent
    \item ``Butterfly effect"
    \item Bounded, deterministic and aperiodic
    \item Topological mixing
    \item Existence of fractal invariant sets
\end{itemize}

\vspace{1em}
Most notably, the definition \textit{sensitive dependence on initial conditions} (SDIC) and its more colloquial name of the \textit{butterfly effect} is useful as it describes the behaviour for the exponential divergence of trajectories for close initial conditions. Namely, trajectories with extremely small differences in initial conditions evolve into exponentially diverging trajectories. For flows, chaotic behaviour requires state space of at least 3 dimensions in order to fulfill sufficient topological mixing and aperiodicity, whilst maintaining avoiding crossings of trajectories. However, for discrete maps only 1 dimension is required to produce chaos.\\

For the curious reader, we provide some exanokes if chaotic flows and discrete maps:
\begin{example} \textbf{-- Lorenz (1963) \cite{lorenz1963deterministic} }\\

    \begin{subequations}
        \begin{align} \dot{x} &= \sigma (y-x)\\
            \dot{y} &= -xz +rx-y\\ 
            \dot{z} &= xy-bz\\ 
        \end{align}
    \end{subequations}
    
    For parameters $\sigma = 16, b = \frac{8}{3}, r = 28$
\end{example}

\begin{example} \textbf{-- R\"{o}ssler (1976) \cite{rossler1976equation}}\\

    \begin{subequations}
        \begin{align} 
            \dot{x} &= -y-z\\
            \dot{y} &= x+ay\\ 
            \dot{z} &= b+z(x-c)\\ 
        \end{align}
    \end{subequations}
    
    For parameters $a = 0.1, b = 0.1, c \in [4,20]$
\end{example}

Examples of chaotic maps:
\begin{example} \textbf{-- Henon Map \cite{henon1976two}}\\
    \begin{subequations}
        \begin{align}
            x_{n+1} &= 1 + y_{n} -ax_{n}^2\\
            y_{n} &= bx_{n}
        \end{align}
    \end{subequations}
    
    For parameters $a = 1.4, b = 0.3$
\end{example}

\begin{example} \textbf{-- Ikeda Map \cite{ikeda1979multiple}}\\
    \begin{subequations}
        \begin{align} 
            z_{n} &= x_{n} + iy_{n}\\
            z_{n+1} &= p + Bz_{n}\text{exp}\{ i\kappa - i \frac{\alpha}{1+|z_{n}|^2} \}
        \end{align}
    \end{subequations}
    For parameters $p = 1.0, B = 0.9, \kappa = 0.4, \alpha = 6.0$
\end{example}

One of the primary features of chaotic behaviour is the existence of bounded, but aperiodic trajectories. As seen from the logistic map (and in general the entire family of one hump maps), chaotic behaviour arises from multiple period doubling bifurcations occurring in sequence. For increasing $r$, periodic behaviour increases from length $2, 4, 8,...,3$ following the Sharkovskii ordering of integers until at some critical bifurcation value $r^*$, following which the system enters the fully chaotic regime. This perspective also gives rise to a key theorem linking periodicity and chaos.

\begin{definition} \textbf{-- Chaotic functions}\\
    A function map $f$ is chaotic if for any $n$, it posseses an orbit of period $n$.
\end{definition}

\begin{theorem} \textbf{-- Period 3 implies chaos}\\
    Let there be a measure preserving map given by function $f$. If $f$ contains a point $a$ that is period 3, then $f$ contains points of all periods $n \in \mathbb{N}$. Therefore, $f$ is chaotic.
\end{theorem}

Recall that Poincare sections allow any flow to be represented as a discrete map. Therefore, if the discrete Poincare map of a continuous system exhibits period doubling bifurcations, and eventually period $3$, this implies the existing of all possible periods lengths in the original continuous system. Furthermore, these periodic orbits are dense in the state space (i.e for any point in state space $M$ that is aperiodic, it is arbitrarily close to a periodic orbit), resulting in very complex behaviour.

\chapter{Number Theory and Recurrence}

\section{Random vs. Deterministic Systems}

The triple formulation $(T,M,\Phi)$ of a dynamical system is very broad and does not give explicit rules on the properties of the evolution operator $\Phi$. If $\Phi$ contains some stochastic process (e.g. randomly generating a sequence of next observations from a uniform distribution), it is still a valid dynamical (albeit stochastic) system. Clearly, stochastic processes causes problems for the $C^1$ assumption needed for the many convenient guarantees for dynamical systems. On the other end of the spectrum, we may define a fully deterministic system (e.g. a uniformly increasing sequence of numbers) that possesses all the useful characteristics of continuity, existence and uniqueness of solutions. Such a system, if simple enough, produces very predictable behaviour with none of the more interesting and complex dynamics that one may wish to study.\\

In the above context, chaotic systems can be understood as cases which lie on the boundary between deterministic and random. They can be fully described analytically (deterministic), but at the same time exhibit topological mixing and diverging trajectories such that temporally distant observations are highly decorrelated and appear almost random. More importantly, this apparently 'randomness' arises from deterministic rules rather than in built randomness in the initial conditionsor the evolution operator. In this chapter, we explore this distinction by way of an example taken from number theory.

\section{The Kronecker System: Dynamics of Fractional Numbers}

Let there be two functions $\text{frac}(\alpha)$ and $\text{int}(\alpha)$ that takes in a positive real number $\alpha$ and returns the fractional part (i.e. $\alpha\,\, \text mod\,\, 1$) and the integer part (i.e. $\lfloor \alpha \rfloor$) respectively. Given this, consider a shift map:

\begin{subequations}
    \begin{align}
        \tau_{\alpha}(x) &= \text{frac}(x+\alpha)\\
        x_{n+1} &= \tau_\alpha(x_n)
    \end{align}
\end{subequations}
for some predefined constant $\alpha$. This dynamical system is called the Kronecker system and can be expressed as a dynamical system with the triple $(\mathbb{Z}, [0,1], \tau_{\alpha})$.

\begin{lemma} \textbf{-- Kronecker System}\\
    Let there be a discrete map defined as
    $$x_{n+1} = \tau_\alpha(x_n),$$
    where $\tau_{\alpha}$ is the shift map with respect to $\alpha \in [0,1]$.
    \begin{enumerate}
        \item $\tau_{\alpha}(x)$ is aperiodic $\forall \alpha \in [0,1]\setminus \mathbb{Q}$
        \item $\tau_{\alpha}(x)$ is eventually periodic $\forall \alpha \in [0,1]\cap \mathbb{Q}$
    \end{enumerate}
\end{lemma}

The behaviour of the Kronecker system is entirely dependent on the choice of $\alpha$. The shift map $\tau_{\alpha}$ essentially produces a sequence of fractional numbers that procedurally truncates the first $n$ decimals of $\alpha$ and returns the remainder (the initial condition $x_0$ merely adds a shift). Recall that irrational numbers have decimal ($n$-ary) fractional parts that never repeat. Therefore, irrational values of $\alpha$ will yield a sequence $x_n$ that never repeats. The converse is also true, if $\alpha$ is rational, its decimal ($n$-ary) representation contains a periodic sequence, and thus $x_{n}$ will also be periodic, excepting an initial transient period.\\

The Kronecker system is quite interesting and leads us to ask several questions:
\begin{enumerate}
    \item[\textbf{Q1}.] What are the possible values that $x_n$ can take?
    \item[\textbf{Q2}.] What can we say about the distribution of $x_n$ on the domain $[0,1]$?
    \item[\textbf{Q3}.] Is the Kronecker system considered chaotic?
\end{enumerate}

\vspace{1em}
The answer to the first question is given by Kronecker's theorem.

\begin{theorem} \textbf{-- Kronecker's Theorem}\\
    Let $\alpha\in \mathbb{R}\setminus \mathbb{Q}$. Then, $\forall U \subseteq [0,1]$, $U \neq \emptyset$, $\exists m \in \mathbb{N}$ such that $\text{frac}(m \alpha) \in U$
\end{theorem}

Simply put, Kronecker's theorem states that for any subinterval in $[0,1]$ and for any irrational $\alpha$, there exists a multiple of $\alpha$ whose fractional part is contained within the subinterval. Therefore, we can define $\epsilon$ width subintervals on $[0,1]$ and show that the Kronecker system is guaranteed to visit any number within some distance $\epsilon$. There are several approaches to prove Kronecker's theorem. An intuitive one is as follows:

\begin{enumerate}
    \item Observe that for $\alpha < 0.5$, the shift map moves an interval $I_n$ right on the number line by $\alpha$. Likewise, $\alpha > 0.5$ moves $x_n$ left of the number line by $1-\alpha$.
    \item Without loss of generality, let us take $\alpha$ to be small. For small number of iterations of the shift map, $I_n$ increases linearly and $\text{frac}()$ does not alter the trajectory.
    \item Let $k$ be the first instance where part of $I_n$ exceeds the boundary value of 1. Therefore, the shift map will result in a some portion of $I_n$ to jump to the left of the number line taking values given by a decimal shift and truncation of $\alpha$.
    \item If $\alpha$ is irrational, truncation and successive additions should contain every possible subsequence of decimals. Therefore, by repeating process on the interval $I_n$, the whole unit interval will eventually be fully covered by images of future intervals $I_n$ vis-\'{a}-vis the pigeon hole principle.
    \item Repeat the above for the case of $1-\alpha$.
\end{enumerate}

\vspace{1em}
We have established that the Kronecker system with irrational $\alpha$ will eventually visit every number in $[0,1]$. What then can we say about \textbf{Q2} -- the distribution of points?

\begin{theorem} \textbf{-- Weyl's Theorem}\\
    Let $\alpha\in\mathbb{R}\setminus \mathbb{Q}$ . Then the sequence $X=\{\text{frac}(n\alpha)\}_{n=1}^N$ is distributed in $[0,1)$ such that $\forall J \subseteq [0,1)$,
    
    $$|X\cap J|\to |J|$$
    
    Alternatively, the proportion of the the first $n$ terms that belong to $J$ has a limit as $n\to \infty$, and is equal to $\mu (J)$, where $\mu$ is a measure of length. See Nillsen for an excellent proof and discussion \cite{nillsen2010randomness}.
\end{theorem}

In essence, Weyl's theorem implies that for the Kronecker system, the sequence of points is uniformly distributed across the unit interval (i.e. each unique subinterval of width $\epsilon$ is visited equally often).\\

The above discussions leads us to address \textbf{Q3}: is the Kronecker system chaotic? One way to tackle this problem is from the perspective of information loss. Take $x_{0}= 0$ and small irrational value $\alpha$ (say $\frac{\pi}{10}$). For short periods of successive iterations of the shift map, the position $x_{n+1}$ is ``predictable" from $x_n$. It is a linear growth by $\alpha$ at each step. Alternatively, knowing the position $x_n$ provides perfect information on the previous values $x_{n-1}, x_{n-2}...$ as long as we stay in the window where the sequence grows linearly. In this respect, the Kronecker system certainly contains chaotic-like properties (i.e. topological mixing), which we define below:

\begin{definition} \textbf{-- Topological transitivity}\\
    Let $S$ be an interval, and $f: S\to S$. The system $(S,f)$ is topologically transitive if $\forall x,y \in S$ and $\epsilon>0$, $\exists z\in S$ and $n\in \mathbb{N}$ such that,
    
    $$|z-x|<\epsilon,\quad |f^{(n)}(z)-y|<\epsilon$$
    Conversely, if a system is not transitive, then $\exists S_{1},S_{2} \subseteq S$ such that
    \begin{itemize}
        \item $S = S_{1} \cup S_{2}$, $S_{1} \cap S_{2} =\emptyset$
        \item $f^{(n)}(S_{1})\subseteq S_1$, $f^{(n)}(S_{2}) \subseteq S_2$ $\forall n$
    \end{itemize}
\end{definition}

The topological transitivity (topological mixing) property states that one can always find trajectories joining the neighbourhoods of two points in state space $S$ within an arbitrary level of accuracy $\epsilon$. Systems that are not topologically transitive contain subsets that do not have trajectories joining them.\\

However, a key feature of chaotic systems is sensitivity to initial conditions (SDIC). To conclusively determine if the Kronecker system is chaotic, we need to provide definitions for SDIC.

\begin{definition} \textbf{-- Sensitive dependence on initial conditions (SDIC)}\\
    A system $(S,f)$ is sensitive to initial conditions (SDIC) if $\exists \delta >0$ such that $\forall x\in S$ and $\epsilon >0$, $\exists z\in S$ and $n\in \mathbb{N}$ where
    
        $$|x-z|<\epsilon,\quad |f^{(n)}(x)-f^{(n)}(z)|>\delta$$
\end{definition}

Is the Kronecker system chaotic based on SDIC? \textbf{\textit{(Yes)}}

\section{Recurrence}

Another feature of dynamical systems is recurrence. Recall that any continuous flow can be reduced to a discrete map by identifying the successive intersections of trajectories with a defined Poincare section. From this, it makes sense to ask several questions:
\begin{enumerate}
    \item How does one define recurrence?
    \item What can be said about recurrence of states in dynamical systems?
\end{enumerate}

\vspace{1em}
We begin by addressing the first question with a mathematical definition of recurrence of states in dynamical systems.

\begin{definition} \textbf{-- Recurrent points and intervals}\\
    Let $(S,f)$ be a dynamical system. A point $x\in S$ is recurrent within a bound $\epsilon$ if $\exists n\in \mathbb{N}$ such that 
    $$|x-f^{(n)}(x)|<\epsilon.$$
    Similarly, a dynamical system is recurrent if for $J\subset S$ where $\mu(J)>0$, $\exists x\in J$ and $n\in \mathbb{N}$ such that
    $$f^{(n)}(x)\in J.$$
    Alternatively, $\exists n\in \mathbb{N}$ such that
    $$f^{(n)}(J)\cap J \neq \emptyset.$$
\end{definition}

Keep in mind that the above definition does not distinguish whether the number of recurrent points or $n$ is finite. It is possible to define dynamical systems where recurrence only occurs for a finite number of times. \\

The second question pertaining to recurrence in dynamical systems is given in general by the Poincar\'{e} Recurrence Theorem.

\begin{theorem} \textbf{-- Poincar\'{e} Recurrence Theorem}\\
    Let $S$ be a bounded interval, let $U$ be a basic subset of $S$ of positive length, and let $f:S\to S$ be a length preserving transformation. Then for almost all $x\in U$, there is $n\in \mathbb{N}$ such that $f^{(n)}(x) \in U$.
\end{theorem}

As with many mathematical proofs, the Poincare recurrence theorem requires unpacking several key words.

\begin{definition} \textbf{-- Basic sets}\\
    A basic set $S$ is a set that can be defined as a finite union of intervals.
    
    $$S = S_{1}\cup S_{2}\cup ... \cup S_{n}$$
    
    It follows that the measure (length) $\mu$ of $S$ is the sum of the lengths of intervals in a finite family of pairwise disjoint intervals whose union is $S$,
    
    $$\mu(S) = \mu (S_{1})+ \mu (S_{2}) + ... + \mu (S_{n}),\quad S_{i} \cap S_{j} = \emptyset$$
\end{definition}

\begin{definition} \textbf{-- Length preserving function}\\
    A function $f$ defined on domain $S$ is length preserving with respect to a measure $\mu$ if $\forall J\subseteq S$ where $J$ is basic,
    
    $$\mu(f^{-1}(J))=\mu(J)$$
    
    and $f^{-1}(J)$ is also basic.
\end{definition}

\begin{definition} \textbf{-- Almost all}\\
    For a space $S$ and measure $\mu$, $Z \subset S$ is a set of measure zero if $\forall \epsilon >0$, there is a collection of connected subsets $A_{n}\subset S$ such that
    
    $$Z \subseteq \bigcup_{n=1}^{\infty}A_{n,}, \quad \sum_{n=1}^{\infty}\mu(A_n)<\epsilon.$$
    
    A statement is true for almost all points in $S$ if it is true for all points in $S\setminus Z$ for some $Z$.
\end{definition}

The Poincar\'{e} Recurrence Theorem states that any bounded dynamical system with dynamics given by a length preserving function $f$ is recurrent. Furthermore, points are recurrent infinitely often. This result presents several key conclusions for chaotic dynamical systems:

\begin{enumerate}
    \item Unstable periodic orbits (i.e. sets of measure zero) are dense.
    \item Topological mixing occurs (bouncing between unstable periodic orbits)
    \item States return within arbitrarily close infinitely often
\end{enumerate}

\vspace{1em}
Furthermore, for chaotic systems, short-term dynamics are governed by $f$ and diverge due to SDIC and topological mixing, but long-term behaviour is determined by properties of the phase space rather than initial conditions (i.e. settling on an attractor manifold.)

\chapter{Embedding Theory and Reconstruction}

Recall that any general dynamical system may be written as a triple $(T,M,\Phi)$ consisting of a time ordering $T$, state space manifold $M$ and evolution operator/dynamics $\Phi$. If $\mathbf{s}(t) \in M$ is the system state, an autonomous dynamical system can be represented as,
\begin{equation}
    \dot{\mathbf{s}} = f(\mathbf{s})
\end{equation}
How does this relate to real world problems? Almost always, it is impossible to observe the full state $\mathbf{s}$. Even if this were possible, the observations would be discretised in time due to measurement and computational constraints. This limitation can be described using a measurement function $h: M \to \mathbb{R}$,
\begin{equation}
    x(t) = h(\mathbf{s}(t)).
\end{equation}
where $x(t)$ is a lower dimensional, incomplete representation of the dynamics. Typically, when we are sampling (with $h$), the observation is applied at discrete times. Thus, we may also write observations as a time series $\{x_{t} \}_{t = 1}^{T}$.\\

Within the above context, what are some interesting and useful things we can do with the observed time series? Generally we can do two type of tasks
\begin{enumerate}
    \item \textbf{System characterisation:} \textit{What are the underlying dynamics and operating mechanisms of the underlying system?} While this is a qualitative question, it can be quantified using descriptive measures like Lyapunov exponents, correlation dimension etc.\
    \item \textbf{State prediction:} \textit{Given some state history $\mathbf{s}(t_0),...,\mathbf{s}(t)$, what is the predicted state $\mathbf{s}(t+\tau)$ for some time $\tau$ in the future?}
\end{enumerate}
For state prediction, we often do not or cannot know $\mathbf{s}$, and so we must settle for an alternative question: \textit{Given some observed history $x(t_{0}), ... x(t)$ where $x(t) = h(\mathbf{s}(t))$, what are the best estimates of future observations $x(t+\tau)$ for some $\tau$ in the future?}

These two problems have several challenges to tackle
\begin{itemize}
    \item Prediction requires the characterisation of dynamics $f$ (or the evolution operator).
    \item If $f$ is known, this would only be useful if we have knowledge of the full state $\mathbf{s}(t)$.
    \item We don't have direct access to the full state space $\mathbf{s}(t)$
\end{itemize}

\vspace{1em}
One can view $f$ as a description of how the components of $\mathbf{s}$ interact and feedback into each other. Therefore, it would be sensible to think that information from one component would inevitably propagate into other components given enough time has passed. Therefore, there is an equivalence between a high dimensional state at single time step and a low dimension state across multiple time steps. This provides a hint to how we may tackle the above problems. Fortunately, embedding theory provides a way for us to address this interesting challenge.

\section{Embedding Theorems}
\begin{definition} \textbf{- Embedding}\\
    Consider a dynamical system with state $\mathbf{s}(t) \in \mathcal{S} \subseteq \mathbb{R}^m$ with the evolution operator $f_T$ such that,
    $$\mathbf{s}(t+T) = f_T(\mathbf{s}(t)),$$
    and a corresponding observation function $h: \mathbb{R}^m \to \mathbb{R}$ that produces a time series,
    $$x(t) = h(\mathbf{s}(t)).$$
    An embedding is a transformation $\Psi: \mathbb{R} \to X \subseteq \mathbb{R}^d$, such that there exists a diffeomorphism $\Phi : \mathcal{S} \to X$ where,
    $$f_T = \Phi^{-1} \circ F_T \circ \Phi,$$
    such that,
    $$x(t+T) = \Psi^{-1} \circ F_T \circ \Psi(x(t)).$$
\end{definition}

We begin the discussion by looking at manifolds and their basic embedding properties. Specifically, we define manifolds as any surface in $M \subseteq \mathbb{R}^n$ that is locally Euclidean in some small neighbourhood. (There are more rigorous descriptions pertaining to mappings of neighbourhood (charts) to produce an \textit{atlas}, whose further discussion we will relegate to some other more authoritative text.)\\

Firstly, we should ask what are the motivations for studying embeddings of manifolds? Recall that in many cases, chaotic systems (and periodic as well) typically do not occupy the entirety of the ambient state space and instead settle on some subspace manifold. Given that we are also restricted in our observations to lower dimensions, the study of embedding theorems of manifolds and their dimensions may prove useful.

\begin{figure}
    \centering
    \includegraphics[width = 0.9\textwidth]{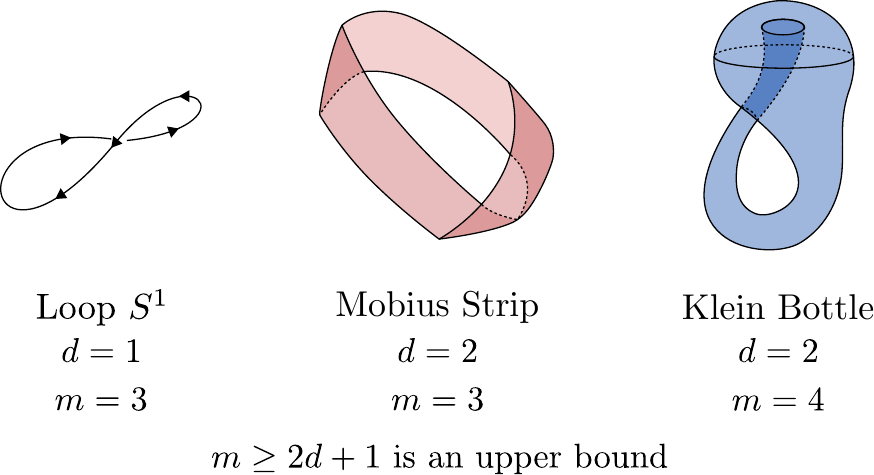}
    \caption{Examples of Whitney embedding theorem. Note that $m=2d+1$ is a minimum needed to guarantee an embedding. Embeddings in lower dimensions can sometimes be possible.}
    \label{fig:whitney}
\end{figure}

\begin{theorem} \textbf{-- Whitney's Embedding Theorem}\\
    Any continuous function mapping a $d$-dimensional manifold to $m$-dimensional manifold($d,m \in \mathbb{Z}^{+}$) can be approximated by a smooth embedding provided $m \geq 2d+1$. See Figure \ref{fig:whitney} for examples.\\
    
    Examples:
    \begin{itemize}
        \item Klein bottles: $d=2 \implies \mathbb{R}^4$
        \item Circle with loop: $d=1 \implies \mathbb{R}^3$
    \end{itemize}
\end{theorem}

Whitney's embedding theorem provides a useful upper bound for the requisite embedding dimension. However, this was further refined by Sauer et al. in the seminal work ``\textbf{Embedology}'' that provided a tighter bound using the more generalised quantity of \textbf{correlation dimension}.

\begin{corollary} \textbf{-- Whitney's Embedding Theorem (Sauer et al., 1991)}\\
    Every $d$-dimensional manifold ($d \in \mathbb{R}^+$) can be embedding in $\lceil 2d \rceil$ where $d$ is the correlation dimension. 
\end{corollary}

Back to the problem of state space estimation. Suppose that we have a trajectory $\psi (t, \mathbf{s}_{0})\subset M \subset \mathbb{R}$ and the measurement function $h(\psi(k \Delta t, \mathbf{s}_{0}))= x_k$ for $n\in \mathbb{Z}$. For $n=1$, we can define the following:

\begin{equation}
    x_{k+1} = \int_{0}^{\Delta t} f(\mathbf{s})\, dt = \phi_{\Delta t}(\mathbf{s}).
\end{equation}
For $n >1$, we would need to estimate multiple derivatives ($f$) to calculate $x_{k+1}$. Using the reverse argument, could we use $(x_{k}, x_{k+1}, x_{k+2},..., x_{k+(n-1)})$ successive measurements to represent these derivatives? (i.e. can we represent high dimensional state space as a vector of successive observations?)

\begin{theorem} \textbf{-- Takens' Embedding Theorem} \cite{takens1981dynamical}\\
    Let $\mathcal{S}$ be a compact manifold of dimension $d$. For pairs $(\phi, h)$, with $\phi \in \text{Diff}^{2}(M),\,h \in C^2(\mathcal{S}, \mathbb{R})$, it is a generic property that the map
    
    $\Psi_{(\phi,h)}: \mathcal{S}\to \mathbb{R}^{2d+1}$, defined by
    
    $$\Psi_{(\phi, h)}(s) = (h(s), h(\phi(s)),...,h(\phi^{2d}(s))),$$

    is a valid embedding.
\end{theorem}

Here, $\phi$ is equivalent to the evolution operator. Takens' theorem essentially states that for almost all selections of observation functions $h$ and differentiable evolution dynamics $f_T$, the construction of a vector $\Psi_{(\phi, h)}(s)$ leads to a valid embedding. This result relies heavily on the ideas transversality and linearisation of $C^k$ maps between $C^k$ manifolds. More intuitively, the theorem guarantees that for a given dynamical system, the construction of a delay vector consisting of time lagged scalar observations

\begin{equation}
    \mathbf{x}(t) = (x(t), x(t-\tau),..., x(t-2d\tau)),
\end{equation}
is a valid embedding. Therefore, for any dynamical system with known dimension $m$,

\begin{equation}
    \dot{\mathbf{s}} = f(\mathbf{s}),
\end{equation}
a topologically conjugate dynamical system may be reconstructed using delay vectors of scalar observations where

\begin{equation}
    \dot{\mathbf{x}} = F(\mathbf{x}), \quad \mathbf{x}\in X \subseteq \mathbf{R}^m ,
\end{equation}
as long as $m \geq 2d+1$. As such, Takens' theorem also highlights close links between the evolution of an autonomous dynamical system to one of a general scalar autoregressive process. Numerous embedding approaches have been proposed since, but invariably capitalise on the intricate link between valid embeddings and delayed observations. However, the delay embedding approach has remained the popular method in time series analysis applications due to its simplicity and relative robustness to noise. See \cite{noakes1991takens} proof sketch on Takens' theorem.\\

\begin{figure}
    \centering
    \includegraphics[width = \textwidth]{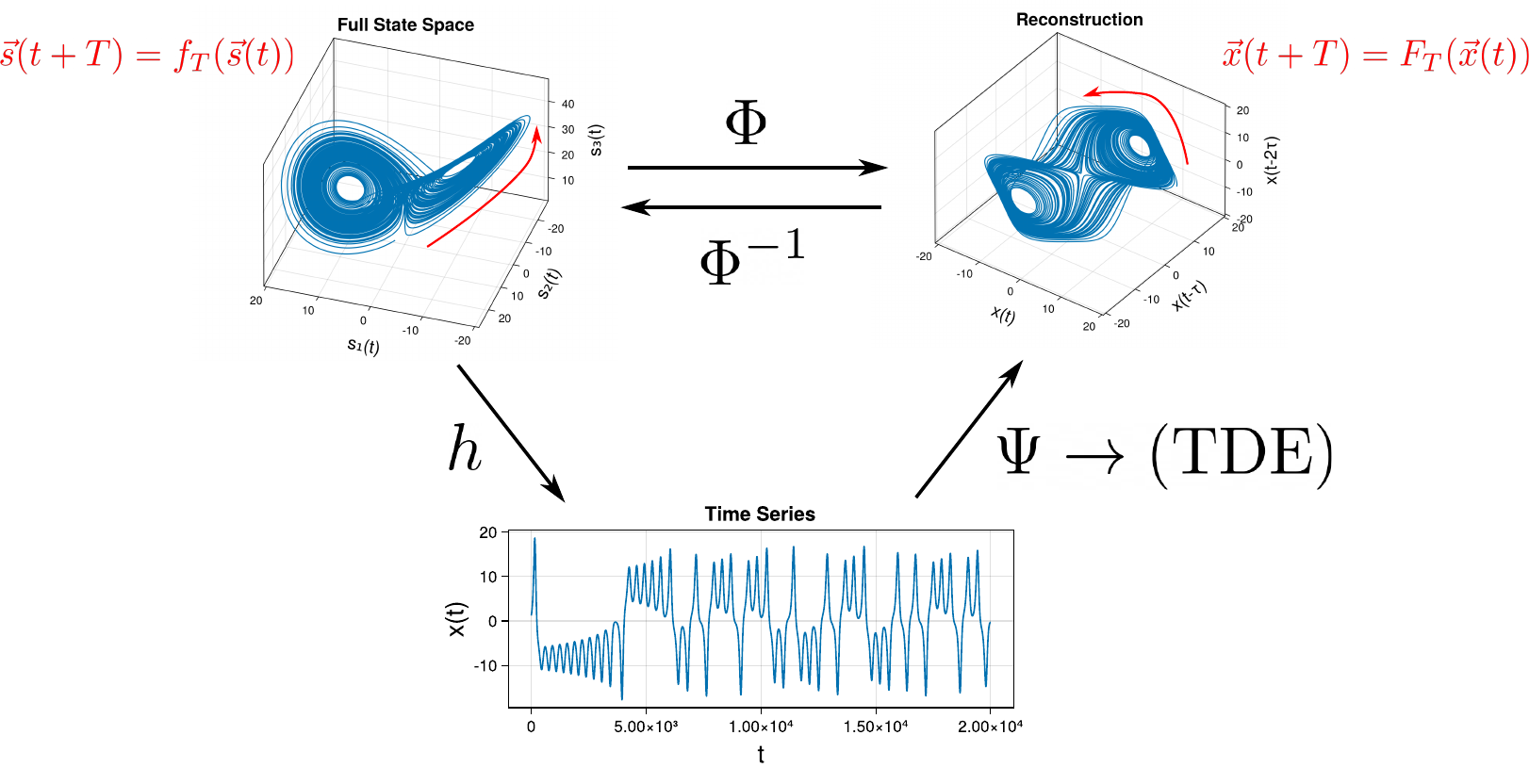}
    \caption[Functional relationships in embedding theory]{Schematic of the embedding process and the relationship between its components.}
    \label{fig:1_Embedding}
\end{figure}

\section{Embedding Methods}

Up until now, our discussion has formed a theoretical basis for validity of embedding approaches when trying to reconstruct phase space. However, the details on how this is achieved in practice has remained sparse. In this section, we present several computational methods that can be used for constructing embeddings.

\subsection{Time Delay Embedding}
First described by \cite{packard1980geometry}, time delay embedding involves the augmentation of a scalar time series $x(t)$ into a higher dimension through the construction of delay vector $\vec{x}(t)$ given as
\begin{equation}
    \vec{x}(t) = (x(t), x(t-\tau),...,x(t-(m-1)\tau)),
\end{equation}
where the embedding parameters to be selected are the delay lag $\tau$ and embedding dimension $m$. According to the guarantees of Takens' theorem, any value of $\tau$ will yield a valid embedding given sufficiently large $m$ and measurement values of infinite precision. However, this is not achievable in practice and different selections of delay lag and embedding dimension can yield varying results.\\

Note that the task of selecting ideal delay lag $\tau$ and embedding dimension $m$ is not unique to time delay embedding. Selecting values of $\tau$ and $m$ are also key decisions in other time series analysis methods such as permutation entropy \cite{bandt2002permutation} and ordinal partition networks \cite{mccullough2015time}. In both of these instances, a delay vector is constructed and represented by an encoding based on the size order each component. The time series may then be viewed as a transitions between different encoding states and used for further analysis.

\subsection{Derivatives Embedding}
The embedding method of derivatives reconstructs an embedding vector using successively increasing order of time derivatives from the observed time series,
\begin{equation}
    \vec{x}(t) = \left( x(t),\frac{dx(t)}{dt},...,\frac{d^mx(t)}{dt^m} \right)
\end{equation}
Derivatives are taken via numerical approximations. The derivatives embedding method is a valid embedding for sufficiently large $m$ if one is able to accurately calculate the required derivatives. 

\subsection{Integral-Differential Embedding}

One weakness of the derivatives embedding approach is the need to evaluate numerical derivatives from data. Whilst this may be acceptable for the first derivative, approximations of successive higher order derivatives are generally inaccurate as the signal to noise ratio tends to be negatively impacted. This is true even if one possesses very clean data sets.\\

An alternative to derivatives embedding is integral-differential embedding \cite{gilmore1998topological}. This approach avoids the calculation of successive higher order derivatives by replacing the second order the derivative with an integral instead. This yields the following embedding construction:
\begin{equation}
    \vec{x}(t) = \left( \int_{-\infty}^t x(t)-\langle x(t)\rangle_t \, dt,x(t),\frac{dx(t)}{dt} \right),
\end{equation}
where the first component is first set to zero mean before integration. The usage of a first order integral and numerical derivative results in a degradation of the signal to noise ratio by only one order each for the first and third embedded components. This is in contrast with the derivatives embedding approach where each successive numerical derivative has a signal to noise ratio that is degraded with increasing orders of magnitude.  However, the integral-differential embedding approach suffers from the same noise effects as the pure derivatives method for higher dimensional embedding. This limits its applicability to systems where system dynamics are presumed to be high dimensional.

\subsection{Global Principal Value Embedding}

\begin{figure}
    \centering
    \includegraphics[width = 0.9\textwidth]{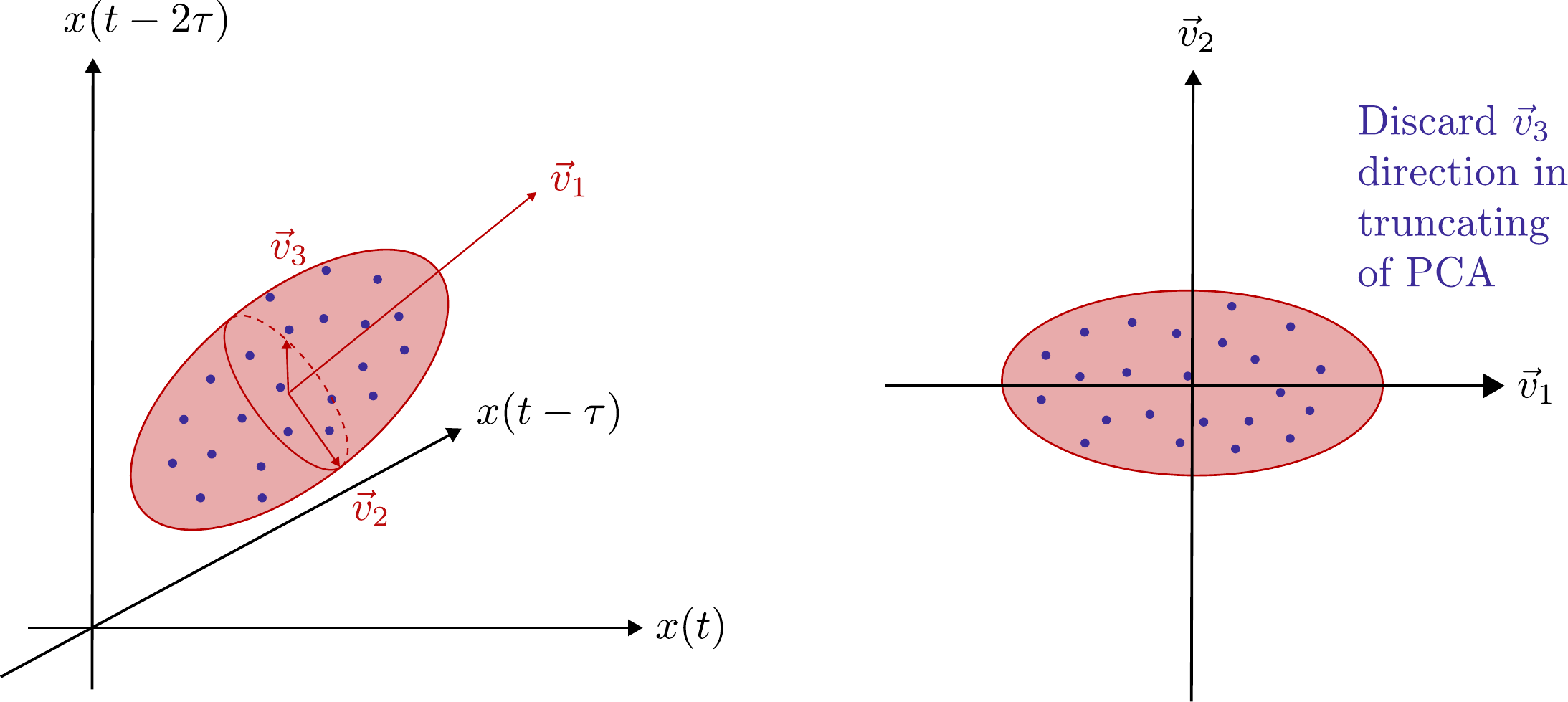}
    \caption{PCA Embedding with dimension reduction from 3 dimensions to 2 (i.e. 2 principal components $v_1$ and $v_2$)}
    \label{fig:PCA}
\end{figure}

The method of principal value embedding was proposed by Broomhead and King \cite{broomhead1986extracting} as a modified alternative to time delay embedding using the theorems by Takens. This method draws upon the ideas of principal component analysis to find an ideal rotation of the time delay embedding with a sufficiently high dimension (see Figure \ref{fig:PCA}). Given a time series $x(t)$ of length $N_T$ and a sliding window of length $M$, we can construct a collection of $N=N_T-(M-1)$ delay vectors,
\begin{equation}
    \mathbf{X}=N^{-1/2} 
    \begin{bmatrix}
    \vec{x}_1\\
    \vec{x}_2\\
    \vdots\\
    \vec{x}_N\\
    \end{bmatrix},
\end{equation}
where $\vec{x}_i$ is the delay vector constructing using the $i^{th}$ value in the time series as the first component,
\begin{equation}
    \vec{x}_i = (x(t_i), x(t_{i-1}),...,x(t_{i-(M-1)})).
\end{equation}
An $M\times M$ covariance matrix $C$ can be calculated from $\mathbf{X}$. The elements $C_{ij}$ of this matrix can be simply given as,
\begin{equation}
    C_{ij} = \langle x(t) \,x(t+(i-j))\rangle_t
\end{equation}
where $\langle ... \rangle_t$ denotes a time average. The principal components of $C$ are then found by calculating its respective eigenvalues and eigenvectors. Taking the first $m$ principal components corresponding to the desired number of embedding dimensions, the eigenvector matrix can be used to calculate a projection of $\mathbf{X}$ corresponding to the final embedded coordinates. Interested readers are advised to refer to \cite{broomhead1986extracting} and \cite{casdagli1991state}.\\

Principal component value embedding essentially aims to distill and simplify a high dimensional delay embedding (usually obtained by taking a large number of lagged components) into a lower dimensional subspace. The remaining subspaces are argued to correspond to component directions with little dynamical variation and importance. One application of this method was as an attempt to simplify the selection of the optimal embedding dimension, where the ideal embedding dimension $m$ corresponds to the number of singular values that are distinctly greater than some `noise floor'.

\section{Uniform vs. Non-uniform Delay Embedding}

Delay embedding occurs in two forms. The first form is uniform delay embedding where a single lag $\tau$ and embedding dimensions $m$ is selected to constructed a delay vector composed of observed values that are equally spaced in time,
\begin{equation}
    \mathbf{x}(t) = (x(t), x(t-\tau), ..., x(t-(m-1)\tau))).
\end{equation}
This approach is the classical form of delay embedding commonly used due to the need to only select two parameters. \\

Whilst the uniform approach works well for systems that exhibit a single dominant periodicity and is easy to implement, its convenience comes at a cost of reduced versatility and limitations when analysing systems with more complex multiscale dynamics \cite{judd1998embedding}.\\

Firstly, the choice to use a single delay limits the ability for the reconstruction to highlight features across multiple disparate time-scales \cite{pecora2007unified}. For example, a fast-slow system with characteristic time scales $\tau_1$ and $\tau_2$ where $\tau_1/ \tau_2 \gg1$ , the choice of selecting $\tau_1$ (i.e. slow dynamics) as the embedding lag can limit the reconstruction's ability to fully unfold attractor topologies corresponding to the fast dynamics. The dynamics the time scale of $\tau_2$ (i.e. fast dynamics) will appear as noisy fluctuations within the reconstructed state space.\\

Secondly, reconstruction from a uniform delay embedding that is sufficient is not necessarily optimal. Here, we must clarify that the definition of optimal presumes some criterion or notion of quality. \cite{casdagli1991state} noted that the quality of an embedding, defined as the reconstruction's robustness to noisy data for prediction, can vary locally throughout different regions of the attractor. This behaviour was also highlighted by \cite{uzal2011optimal} in his extension of Casdagli's noise amplification and distortion methods. Additionally, we should also consider that invariant measures such as the Lyapunov exponent also vary locally. Hence, the selection of a single embedding lag implies that all these variations may be averaged. \\

An alternative proposed form of delay embedding is non-uniform delay embedding. 
\begin{equation}
    \mathbf{x}(t) = (x(t), x(t-\tau_1), ..., x(t-\tau_{m-1}))
\end{equation}
This is a natural extension of uniform embedding that aims to address some of the latter's limitations by allowing multiple time scales to be encoded into the delay vector.

\chapter{Embedding in Practice}

Embedding methods seem like a viable way to address the logistical challenges and limitations when performing time series analysis. Assuming we are using a some form of delay embedding (because it's simple and easy), we still need to answer two questions
\begin{enumerate}
    \item[\textbf{Q1}.] Uniform or non-uniform? - This is application dependent and up for debate. See \cite{tan2023selecting} for a detailed discussion. 
    \item[\textbf{Q2}.] What values for embedding parameters (dimension and lag)?
\end{enumerate}

\section{Do embedding parameters parameter?}

Recall that Takens' theorem made no stipulation regarding the value of $\Delta t$ (or $\tau$) that is required for the delay embedding to be valid. Only that we must have sufficiently high dimension. However, as we will see, this does not mean that all embeddings are made equal and perform equally well for time series analysis. This is discounting the fact that we have not discussed the meaning of ``large enough dimension'' and what is considered ``large enough''.\\

In the original formulation of Takens' theorem, the quoted evolution operator $\phi: x(t) \to x(t+\tau)$ is general and thus should hold for any selection of $\tau$. However, this is only true for the case of infinite sampling resolution and precision, and the absence of noise. In practice, such conditions are rare if not impossible to achieve. \\

To illustrate the effects of noise and precision, consider the time series $s(t) = \cos (t)$, which may also be expressed as a solution to linear ODE,
\begin{equation}
    \ddot{s} +s= 0,
\end{equation}
where the underlying attractor is represented by a loop. This system may also be written as a first order 2 dimensional dynamical system that adheres to the form of Equation
\begin{subequations}
    \begin{align}
        \dot{x} &= y \\
        \dot{y} &= x
    \end{align}
\end{subequations}

\begin{figure}
    \centering
    \includegraphics[width = \textwidth]{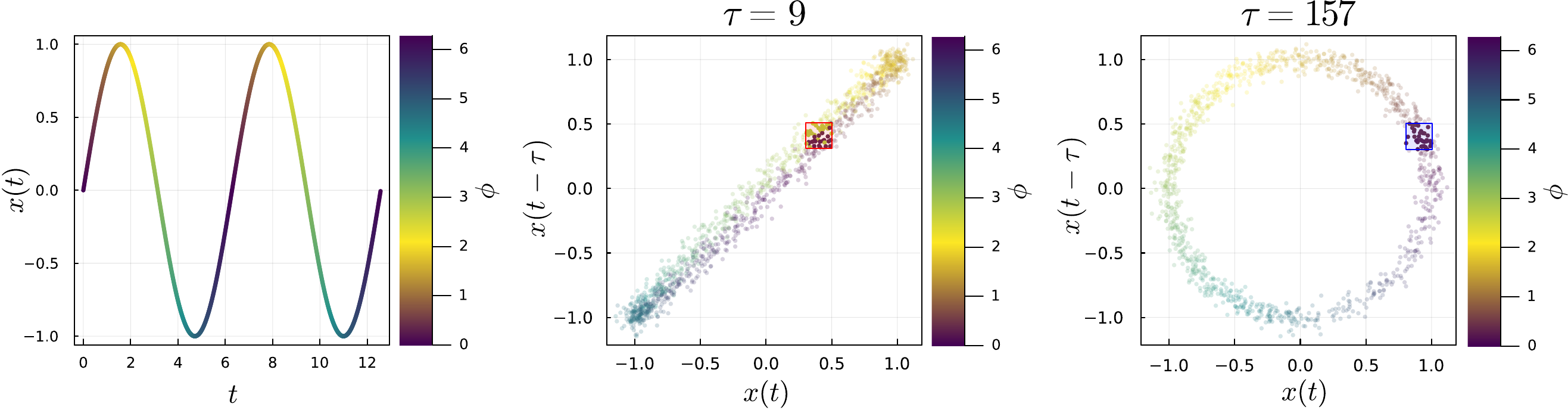}
    \caption[Embedding and state estimation of periodic time series]{Sinusoidal time series \textit{(left)}. Delay embeddings constructed with observational different time lags for time series contaminated with Gaussian noise $\xi \sim \mathcal{N}(0,0.05^2)$. Poor lag selection that is too small can result incorrect identification of the phase \textit{(middle)}. Lag of a quarter period $\tau = 157$ provides better separation of trajectories \textit{(right)}. Bounding boxes correspond to a potential precision error of $\epsilon = 0.05$}
    \label{fig:3_circle_fig}
\end{figure}

A delay vector constructed as $\mathbf{x}(t) = (x(t), x(t-\tau))$ forms an ellipse for all values $\tau >0$ (see Figure \ref{fig:3_circle_fig}). Increasing values of $\tau$ reduces the eccentricity of the ellipse for a minimum when $\tau$ is equal to a quarter of the period. Hence, there is a clear one-to-one mapping between the the position $\mathbf{x}(t)$ and the phase of the oscillator $\phi$. The state of the underlying system can be perfectly determined based on observations $x(t)$.\\

The presence of observational noise injects uncertainty into the estimation of the real state. In the case of finite precision $\epsilon$, the uncertainty forms a square region centred around each observation point with width $\epsilon$. This effect is most detrimental when regions of uncertainty from multiple observations corresponding to different phases overlap creating ambiguity on the real state $s(t)$. For the case of a periodic time series, this effect is minimal at $\tau$ equal to a quarter of the period where the ellipse is maximally unfolded. Whilst the effect of noise and precision appear relatively small for periodic signals, the effect of uncertainty is not insignificant for chaotic systems containing dense orbits and can make estimation of the real underlying state difficult. Furthermore, in the case of systems with dynamics across multiple spatial and time scales there is no guarantee that an ideal $\tau$ exists.\\

Let $\Psi_{\theta}$ be an embedding mapping with respect to the selection of embedding parameters $\theta$. What selection of parameters $\theta$ should be selected to produce the best embedding?\\

In the case for uniform delay embedding, $\theta$ would consist of the embedding lag $\tau$ and embedding dimension $m$. More generally, non-uniform delay embedding would require the selection of a set of lags $\mathcal{T} = \{ \tau_1, ..., \tau_{m-1} \}$ to construct the delay vector,
\begin{equation}
    \mathbf{x}(t) = (x(t), x(t-\tau_1), ..., x(t-\tau_{m-1})).
\end{equation}
However, the embedding problem is ill-posed as the notion of embedding quality is not defined. Therefore, answer the embedding problem would likewise require first addressing a separate problem: ``How does one quantify a good embedding?".

\section{Measuring Embedding Quality}

The main useful feature of an embedding is that it preserves dynamics of the original (partially observed) system. This also means that there is a unique one-to-one correspondence between real states $\mathbf{s}\in M$, and reconstructed states $\mathbf{x}\in \mathbb{R}^m$. As previously demonstrated, this cannot be guaranteed when there is noise and finite resolution in the data, the effect of which is that some embeddings are better at preserving this one-to-one correspondence than others. Therefore, we would like our embedding (delay embedding with parameters $\theta$) to contain as little ambiguity regarding out real state as much as possible. But what does this mean topologically?

\begin{enumerate}
    \item[\textbf{1.}] \textbf{Minimise the number of self-crossings}\\
    From the existence and uniqueness, and fundamental theorems of flows, trajectories cannot cross. Therefore it would also make sense that our trajectories and manifolds should not intersect with themselves. This would suggest that there is an insufficiently large embedding dimension.\\

    In terms of state estimation, self-crossings create ambiguity when performing state inference. (What happens when you try to embed a sum of two sines with different frequencies in 2 dimensions?)
    \vspace{1em}
    \item[\textbf{2.}]  \textbf{Minimise laminar manifolds}\\
    Similar to self-crossings but arguably less egregious, laminar structures in manifolds while not directly in violation of existence and uniqueness, results in ill-posed manifolds for reasons similar to self crossings. Unless observations and trained prediction models are exact, laminar surfaces are regions of uncertainty for state estimation.
\end{enumerate}

\begin{figure}[h]
    \centering
    \includegraphics[width = 0.85\textwidth]{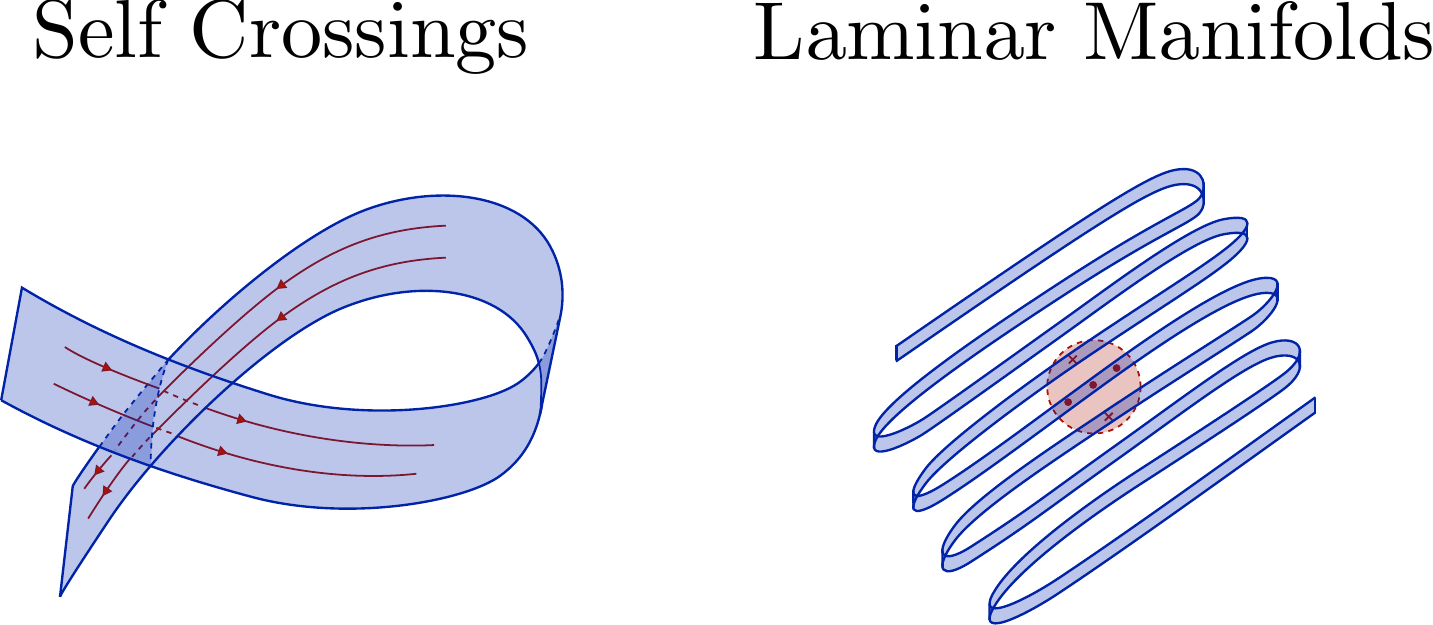}
    \caption{Self crossings and laminar manifolds. Open ball on laminar manifolds reveal true and false nearest neighbours.}
    \label{fig:redundance_irrelevance}
\end{figure}

We can summarise these considerations in a general by considering two concepts when choosing delay coordinates (i.e. $\tau$ and $m$): (1) irrelevance and (2) redundancy. Instead of considering $\tau$ and $m$ separately, we can consider the size of the embedding window $\tau_{w} = m\tau$ (for uniform embedding). Ideally, the embedding window $\tau_w$ should not too short (redundance) as this creates laminarity along the main diagonal since components have high temporal correlation. Similarly, $\tau_w$ should not be too long (irrelevance) to the point where coordinates are no longer correlated at all and embedding manifolds become too folded with large numbers of self crossings. Generally, all embedding parameter selection methods aim to optimise with respect to some notion of quality that can be related to the ideas of redundancy and irrelevance.

\section{Uniform Delay Embedding - Selecting $\tau$ and $m$ simultaneously}

\subsection{Gao \& Zheng - Characteristic Lengths}

The embedding method proposed by Gao and Zheng \cite{gao1993local, gao1994direct} is based on the incidence of false nearest neighbours. False nearest neighbours can be attributed to either redundancy (insufficiently unfolded) and irrelevancy (spurious intersections in the attractor). The method proposed by Gao and Zheng operates on the notion that the separation distance and proportion of false nearest neighbours, should be minimised in an ideal embedding.\\

Consider a pair of points in embedded space $\vec{x}_i, \vec{x}_j$ and their evolution $k$ steps into the future $\vec{x}_{i+k}, \vec{x}_{j+k}$. Points that are false nearest neighbours will tend to separate faster than real neighbours as the attractor unfolds in a time delay embedding. As a result, the ratio between their distances $\lvert \vec{x}_{i+k}- \vec{x}_{j+k}\rvert / \lvert \vec{x}_i- \vec{x}_j\rvert$ will be larger for pairs of false nearest neighbours and approximately equal to 1 for real neighbours. Gao and Zheng then propose the following measure $\Lambda$ to optimise the embedding parameters,
\begin{equation}
    \Lambda(k,m,\tau) = \frac{1}{N_{ref}} \sum_{i,j} \ln \frac{\lvert \vec{x}_{i+k}- \vec{x}_{j+k} \rvert}{\lvert \vec{x}_i- \vec{x}_j \rvert},
\end{equation}
where $N_{ref}$ is the number of randomly sampled point pairs over which the distance ratio is averaged. There are several additional restrictions on the selection of point pairs $\vec{x}_{i}, \vec{x}_{j}$. Firstly, the initial separation of these points should satisfy $\lvert \vec{x}_i - \vec{x}_j \rvert \leq r$ where $r$ is a small selected threshold, i.e. the initial separation of points should be small enough such that the calculation of growing separation is sensible. Secondly, the selection of pairs of points should not have an intersecting Theiler window $\lvert i-j \rvert>l_{\rm Theiler}$, where $l_{\rm Theiler}\in \mathbb{N}^+$. This is done to prevented unwanted correlations between points on the same local trajectory \cite{theiler1986spurious}. Finally, the constant $k$ should not be too large and selected with respect to the natural time scale of the system dynamics.\\

To identify good embedding parameters, profiles of $\Lambda(\tau)$ are calculated for increasing values of embedding dimensions $m$. The value of $m$ that corresponds to the largest decrease across the profile $\Lambda(\tau)$ is selected as the embedding dimension. The embedding lag $\tau$ is then selected as the first minimum of $\Lambda(\tau)$.\\

This approach has a problem in that it requires the selection of an evolution time $k$. Instead of arbitrarily selecting $k$, a characteristic length $J(m,\tau)$ describing the natural spatial scale of the system's attractor can be calculated,
\begin{equation}
    J(m,\tau) = \langle \lvert \vec{x}_i - \vec{x}_j \rvert \rangle,
\end{equation}
where $\langle ... \rangle$ denotes an average over sampled pairs of points of the attractor. The characteristic length is then used to calculate the separation time $T_J(\vec{x_i},\vec{x_j})$ defined as the time taken for pairs of nearest neighbours to diverge by some proportion of the characteristic length $J(m, \tau)$. For real neighbours, $T_J$ will converge to a value related to the Lyapunov exponent of the system with increasing embedding dimension $m$, whilst false nearest neighbours will result in a smaller value $T_J$ as trajectories quickly separate. The new measure that is used to determine the embedding parameters is given by,
\begin{equation}
    C(m, \tau) = \frac{1}{N_{ref}} \sum_{i,j}  T_J(\vec{x_i},\vec{x_j}),
\end{equation}
where $N_{ref}$ is the number of sampled pairs of nearby neighbours. The values for $m$ and $\tau$ that maximise $C(m,\tau)$ are selected as the embedding dimension and lag.

\begin{figure}
    \centering
    \includegraphics[width = 0.6\textwidth]{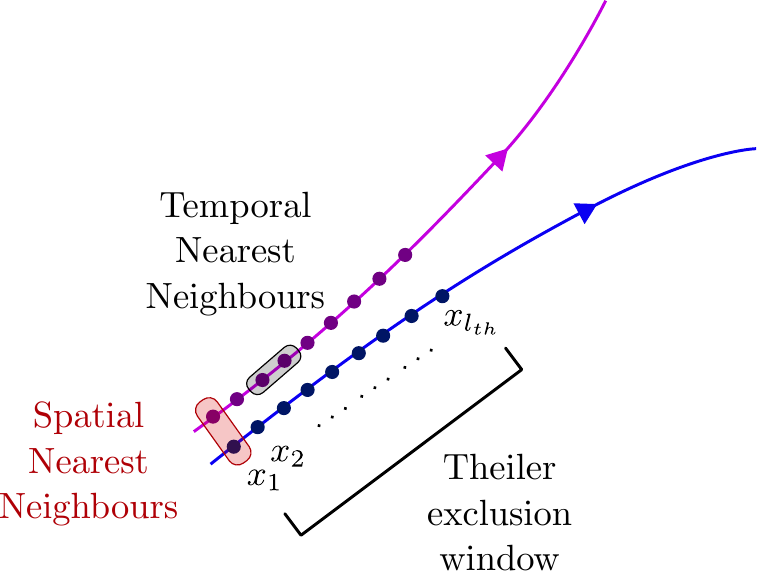}
    \caption{Theiler window used to exclude points that are near neighbours by virtue of being temporally close. This is done to prevent bias in the statistics.}
    \label{fig:Theiler}
\end{figure}

\subsection{Wavering Product}

The wavering product \cite{liebert1991optimal} is similar to that of Gao and Zheng and characteristic lengths in that all are based around the concepts of nearest neighbours. The authors propose that good embeddings should preserve the correspondence between the order of nearby neighbours of a given reference point (i.e. the order of neighbours sorted according to distance from some reference point $\vec{x}_i$ should be preserved). This is done by comparing the order of $N$ nearby neighbours of a point $\vec{x}_i$ between a given embedding $\Phi_k$, whose ordered sequence neighbours are given by,
\begin{equation}
    X_{\Phi_k}=\{ \vec{x}_{i,1},...,\vec{x}_{i,N} \}_{\Phi_{k}},
\end{equation}
and its projection onto its next order embedding $\Phi_{k+1}$ (by increasing $m$ or $\tau$) with the sequence given by,
\begin{equation}
    Z_{\Phi_{k+1}}=\{ \vec{z}_{i,1},...,\vec{z}_{i,N} \}_{\Phi_{k+1}}.
\end{equation}
Here, $\vec{x}_{i,n}$ corresponds to the $n^{th}$ nearest neighbour of the $i^{th}$ reference point $\vec{x}_i$. The projection $\vec{z}_{i,n}$ corresponds to the same neighbour data point $\vec{x}_{i,n}$ whose position is recalculated from the next order embedding $\Phi_{k+1}$.\\

Similar comparisons can also be made into a projection into an embedding of lower order (by decreasing $m$ or $\tau$) giving a new set of ordered points,
\begin{equation}
    V_{\Phi_{k-1}}=\{ \vec{v}_{i,1},...,\vec{v}_{i,N} \}_{\Phi_{k-1}}.
\end{equation}
Ideally, a good embedding should preserve a one to one correspondence in these ordered sequences. This will yield a value equal to 1 for the following ratios,
\begin{equation}
    \frac{\lvert \vec{x}_i-\vec{z}_{i,n} \rvert}{\lvert \vec{x}_i-\vec{x}_{i,n} \rvert}
,\;{\rm and}\;
\frac{\lvert \vec{x}_i-\vec{x}_{i,n} \rvert}{\lvert \vec{x}_i-\vec{v}_{i,n} \rvert}.
\end{equation}

The method presented by Liebert et al. propose the following measure as the product of the above two ratios,
\begin{equation}
    W_i(m,\tau)=\prod_{n=1}^{N}
\left\{
\left(
\frac{\lvert \vec{x}_i-\vec{z}_{i,n} \rvert}{\lvert \vec{x}_i-\vec{x}_{i,n} \rvert}
\right)
\cdot
\left(
\frac{\lvert \vec{x}_i-\vec{x}_{i,n} \rvert}{\lvert \vec{x}_i-\vec{v}_{i,n} \rvert}
\right)
\right\}.
\end{equation}

The measure to be optimised is given by the average over $N_{ref}$ randomly sampled reference reference points,
\begin{equation}
    W(m,\tau)=\ln \left( 
\frac{1}{N_{ref}}\sum_{i=1}^{N_{ref}} W_i(m,\tau)
\right)
\end{equation}
with $m$ being selected as the dimension which achieves the limiting behaviour of $W$ and $\tau$ corresponds to the first minimum of the resulting profile. 

\section{Uniform Delay Embedding - Separate selection}

In practice, simultaneous parameter selection methods are rarely used when performing uniform delay embedding. Instead, the common practice is to select appropriate values of $m$ and $\tau$ separately. In methods where the selection of one relies on a prior selection of the other, one can employ an iterative approach (e.g. select $\tau_{1}\to m_1|\tau_{1} \to  \tau_2|m_1...$).

\subsection{Selecting $\tau$}
There are various methods for guiding the selection of embedding lag $\tau$, the sophistication of which varies from heuristic approaches to more theoretical methods. Below are brief descriptions on several methods.

\textbf{Quarter Period}\\
 This approach allows the natural time scale of the system dynamics to be encoded within the embedding procedure. This heuristic is inspired from the problem of embedding a sine wave in 2D $x(t) = \sin (\omega t)$. In this case, $\tau = \frac{2\pi}{4\omega}$ yields a 2D delay embedding that is maximally circular. However, this heuristic cannot be directly applied to chaotic systems where signals are aperiodic. Instead an estimation of some form of the pseudo-period is required, such as the average time between successive maxima,

\vspace{1em}

\textbf{Autocorrelation}\\
This approach aims to select lags such that components in the delay vector have minimal linear correlation in order to minimise redundancy in the embedding,
\begin{equation}
    R(\tau) = \frac{\mathbb{E}[(x(t)-\bar{x})(x(t+\tau)-\bar{x})]}{\sigma^2}.
\end{equation}
Delay lag $\tau$ can be selected as either the first zero, first minima or the first crossing of $R(\tau)<\frac{1}{e}$. The choice of each criterion is arbitrary and is selected based on the one that produces the best performance for the target application. 

\begin{figure}
    \centering
    \includegraphics[width = 0.95\textwidth]{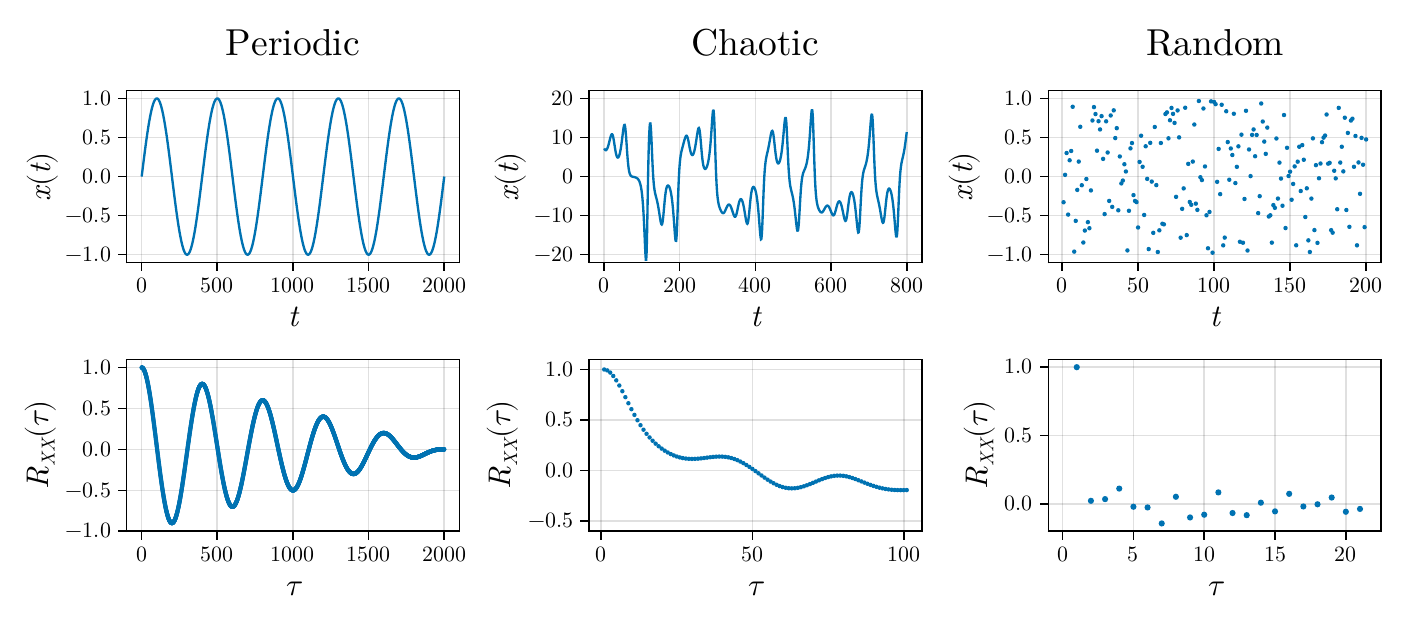}
    \caption{Autocorrelation profiles for periodic, chaotic and random signals}
    \label{fig:autocorrelation}
\end{figure}

\vspace{1em}

\textbf{Minimum Mutual Information}\\
One weakness of autocorrelation is its inability to account for non-stationarities in the time series (e.g. drifts in phase, frequency and magnitude). Additionally, its application only measures linear correlation (Abarbanel, 1993), which might be problematic for nonlinear dynamical systems. In all but the simplest cases, dynamical systems exhibit some level of non-stationary and nonlinear behaviour. \cite{fraser1986independent} proposed that the auto mutual information of the system should be used in place of autocorrelation. In their original paper, Fraser and Swinney first provide a geometrical interpretation to complement the theoretic arguments for mutual information. Namely, consider a set of points whose values in one component $x$ lie within some fixed window. From this set, track their positions $\tau$ steps into the future and calculate the distribution of values $p_\tau(x)$ in the same component for the same set of points. A value of $\tau$ that results in a wider distribution $p_\tau(x)$ should correspond to a good lag, which also corresponds to small values in the mutual information.\\

The mutual information can be interpreted as the nonlinear analogue of the autocorrelation function evaluated using ideas from information theory. Given two random variables $X, Y$, the mutual information between them is given as,
\begin{equation}
    I(X,Y) = \sum_{x\in \mathcal{X}} \sum_{y\in \mathcal{Y}} P_{X,Y}(x,y) \log_{2}\left( \frac{P_{X,Y}(x,y)}{P_{X(x)}P_Y(y)} \right)
\end{equation}
or for continuous random variables,
\begin{equation}
    I(X,Y) = \int_{\mathcal{X}} \int_{\mathcal{Y}} P_{X,Y}(x,y) \log_{2}\left( \frac{P_{X,Y}(x,y)}{P_{X}(x)P_Y(y)} \right) \, dy dx.
\end{equation}

Simply put, the quantity $I(X,Y)$ expressed in bits (if base 2) or nats (base $e$) measures the amount of "information" contained in measurements of $X$ that is also contained in $Y$. Alternatively, how much does an observation of $X$ tell you about an unseen observation of $Y$? Adapting this analogy to time series, we can ask a similar question for $\tau$: How much information is shared between an observation $x(t)$ and at $x(t+\tau)$? If the aim is to minimise redundancy, we need to select $\tau$ such that the mutual information $I(\tau)$ is minimised where,
\begin{align}
    I(\tau) = \int P(x(t),x(t+\tau)) \log_2 \left(
\frac{P(x(t),x(t+\tau))}{P(x(t)) P(x(t)+\tau)} \right) dt\\
    I(\tau) = \sum_{n=1}^{N} P(x(n),x(n+\tau)) \log_2 \left(
\frac{P(x(n),x(n+\tau))}{P(x(n)) P(x(n+\tau))} \right).
\end{align}
where $P(x(t))$ is the probability of observing a state $x(t)$ at any given time and $P(x(t),x(t+\tau))$ is the joint probability defined similarly for both time $t$ and a future time $t+\tau$. Drawing from information theory, mutual information $I(\tau)$ aims to quantify the amount of information about a future state at time $t+\tau$ that is contained in an observation at time $t$. High levels of mutual information for a given lag $\tau$ imply a high degree of correlation between states and will result in higher redundancy for the delay reconstruction. \\

The strengths of the minimum mutual information and autocorrelation lies in its ability to provide reasonable estimates for lag with relatively simple and quick computation. However, there are no guarantees for the existence of a clear minimum for a given mutual information profile $I(\tau)$ \cite{wallot2018calculation}. Additionally, calculating mutual information requires the numerical estimation of probability density functions $P(x(t))$ and $P(x(t), x(t+\tau))$, and thus requires consideration regarding optimal histogram bin size and data length requirements \cite{kraskov2004estimating, papana2009evaluation}.

\begin{figure}
    \centering
    \includegraphics[width = 0.95\textwidth]{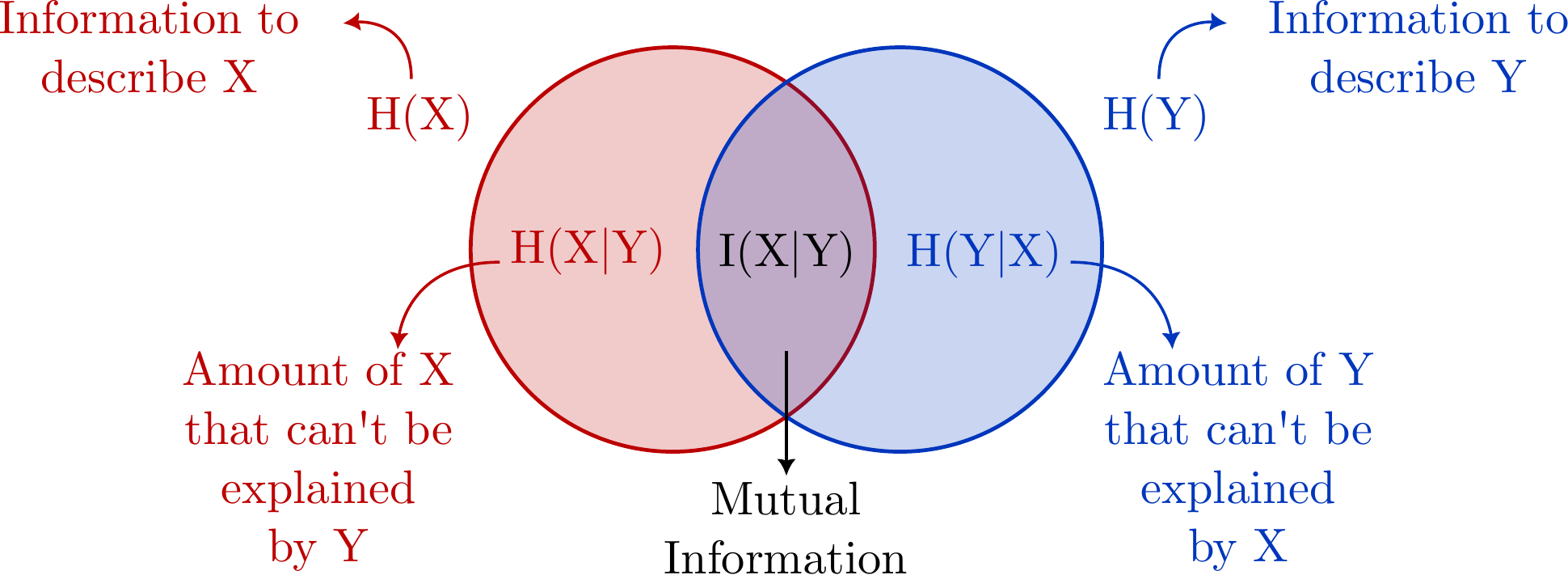}
    \caption{Illustration of mutual information and Shannon entropy relationships}
    \label{fig:mutualinfo}
\end{figure}

\vspace{1em}
\textbf{Fill-factor}\\
The fill-factor approach first proposed by \cite{buzug1992comparison} is an entirely geometrical approach to calculating the quality of a given embedding. This method assumes that an ideal embedding should be able to unfold an attractor and maximise the separation between the trajectories. The authors argue that such an embedding optimally utilises the ambient space and reduces the ambiguity of the true state of the system for different points in the reconstructed state space.\\

The fill-factor is calculated by first sampling $m+1$ random points from an $m$ dimensional delay embedding of the data. A reference point $\vec{x}_r$ is then selected from this collection and the corresponding relative distance vectors can be calculated,
\begin{equation}
    \vec{d}_i (\tau)=
    \begin{bmatrix}
    x_i(t)-x_r(t)\\
    x_i(t-\tau)-x_r(t-\tau)\\
    \vdots\\
    x_i(t-(m-1)\tau)-x_r(t-(m-1)\tau)\\
    \end{bmatrix}.
\end{equation}
The corresponding $m \times m$ matrix can then be expressed as
\begin{equation}
    M(\tau) = (\vec{d}_1, \vec{d}_2,...,\vec{d}_m),
\end{equation}
and the volume of the resulting parallelepiped is given by calculating the determinant of $M$,
\begin{equation}
    V(\tau) = \det(M(\tau)).
\end{equation}

The final expression for the fill-factor is given by calculating the average volume over a collection of randomly sampled parallelepipeds $V_i(\tau)$, normalised by the range of the sampled data points,
\begin{equation}
    f = \log \left( \frac{\frac{1}{N}\sum_{i=1}^N V_i(\tau)}{(\max_{k}x(t_k)-\min_kx(t_k))^m}\right).
\end{equation}
The authors recommend the selection of $\tau$ that maximises the fill-factor $f$ over the interval $\tau \in (0,T_c/2)$, where $T_c$ is the characteristic recurrence time. The value of $T_c$ is given by,
\begin{equation}
    T_c=\frac{1}{\omega_c},
\end{equation}
where $\omega_c$ is the most dominant frequency from the power spectrum of the time series. 

\begin{figure}
    \centering
    \includegraphics[width = 0.95\textwidth]{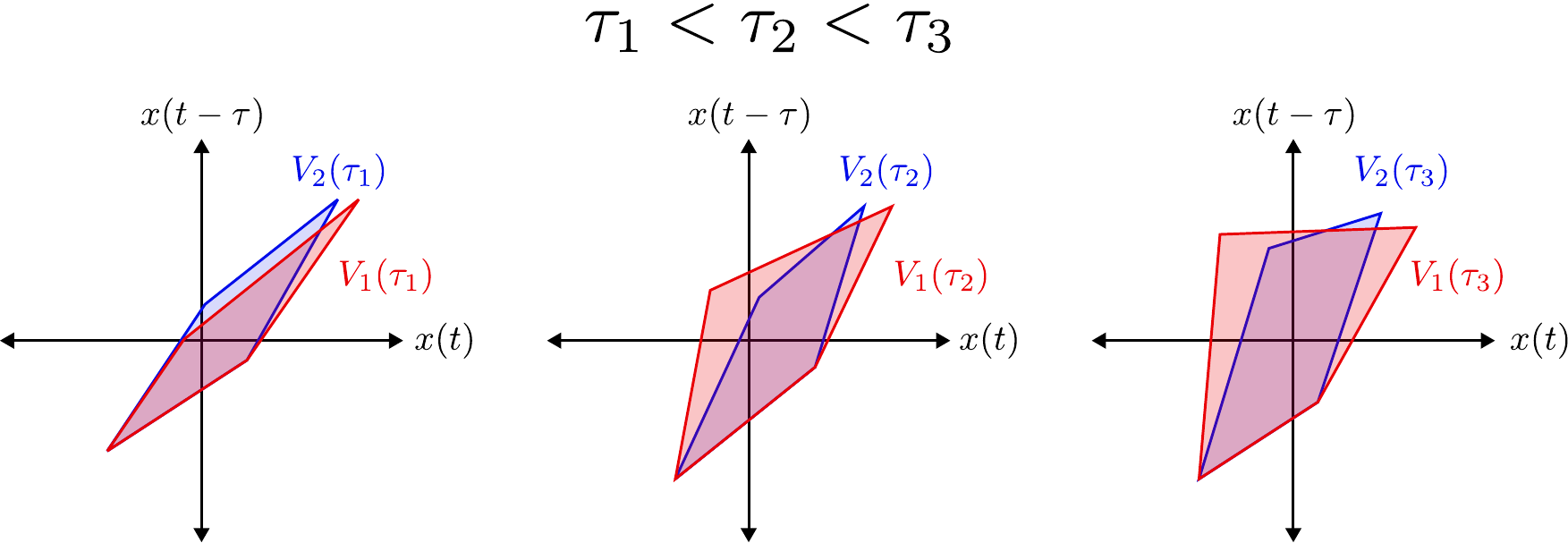}
    \caption{Fill factor method for evaluating $V$}
    \label{fig:fillfactor}
\end{figure}

\vspace{1em}

\textbf{Noise Amplification}\\
Noise amplification was a measure proposed by Casdagli \cite{casdagli1991state} in an attempt to quantify the quality of an embedding. This is supported by the notion that a good embedding should be useful in performing predictions. Additionally, good embeddings should be able to still perform relatively well even in the presence of noise. Noise amplification for a given embedding $\vec{x}(t)$ is defined with respect to predictability of the system $T$ steps into future under the presence of noise. Generally, this is given by:
\begin{equation}
    \sigma(T,\vec{x}) = \lim_{\epsilon \to 0} \sigma_\epsilon (T,\vec{x}),
\end{equation}
where
\begin{equation}
    \sigma_\epsilon (T,\vec{x}) = \frac{1}{\epsilon}\sqrt{\text{Var}[x(T)|B_\epsilon(\vec{x})]}.
\end{equation}

Here, $\text{Var}[x(T)|B_\epsilon(\vec{x})]$ corresponds to the conditional variance of $T$ step predictions into the future in $\mathbb{R}$ from an initial condition $\vec{x}$ in embedding space $\mathbb{R}^m$ contaminated with added small observation noise $\epsilon$. In this case, it is assumed that predictions have no model errors. This condition may be fulfilled by choosing nearby neighbours in the embedding $\mathbb{R}^m$ as a proxy for noisy initial conditions \cite{uzal2011optimal}. \\

Finally, the noise amplification quantity $\sigma(T,\vec{x})$ is averaged over a collection of reference points sampled across the time series in order to calculate the noise amplification value $\sigma$. Embeddings with high noise amplification imply that nearby neighbours in embedded space $\mathbb{R}^m$ tend to have future trajectories that rapidly diverge because they do not correspond to real neighbors in the true manifold $\mathcal{M}$ state space. Therefore, the impact of noise is greatly amplified as small perturbations in the reconstructed space $\mathbb{R}^m$ result in large uncertainties in the true state of the system.

\vspace{1em}

\textbf{L-Statistic}\\
One weakness of the noise amplification measure is its requirement to define $T$, the prediction horizon over which to calculate the noise amplification. This was addressed by Uzal et al. by modifying the definition of noise amplification to the following:
\begin{equation}
    \sigma_\epsilon^2 (\vec{x}) = \frac{1}{T_M}\int_0^{T_M} \sigma_\epsilon^2 (T,\vec{x}) \,dT.
\end{equation}

This definition calculates the noise amplification with respect to a range of prediction horizons up to a maximum value of $T_M$ and is found to be relatively robust for sufficiently large $T_M$. \\

The algorithm used to calculate $\sigma$ relies on using $k$ nearest neighbours from a reference point $\vec{x}_i$ as a proxy. Based on the distribution of points, this can result in effective noise levels $\epsilon$ of different sizes for each point. Therefore, Uzal proposed a normalisation constant $\alpha_k$ accommodate for this variation given by:
\begin{equation}
    \alpha_k^2 = \left[ \sum_i \epsilon_k^{-2}(\vec{x}_i)\right]^{-1}
\end{equation}

Combining these two ideas, the authors propose that noise amplification $\sigma$ measures some notion of redundancy, and $\alpha_k$ measures some notion of irrelevance. The L-statistic is then described as a cost function to minimise both of these values simultaneously,
\begin{equation}
    L = \log(\sigma \alpha_k)=\log(\sigma)+\log(\alpha_k).
\end{equation}

\subsection{Selecting $m$ - Global False Nearest Neighbours (GFNN)}
The guarantees of Takens' theorem requires that the $m>2d+1$, where $d$ is typically unknown. It is possible to use Sauer, York and Casdagli's refinement \cite{sauer1991embedology} such that $m>\lceil c_{d}\rceil$ is sufficient where $c_d$. But again, estimating $c_d$ is non trivial, and also requires data to be in the form of a state space (i.e. prior embedding is required). One common approach is that of global false nearest neighbours (GFNN) similar to the ideas found in Gao \& Zheng/Characteristic Lengths/Wavering product approaches. 

\begin{enumerate}
    \item Select a lag $\tau$ and starting embedding dimension $m$ (e.g. $m = 2$) to produce a reconstructed time series
    $$y(t) = [x(t), x(t+\tau), ..., x(t+(m-1)\tau)]$$
    \item For each point $y(t^*)$ in the reconstructed attractor, find its nearest spatial neighbour, making sure to exclude temporal neighbours using a Theiler window,
    $$y^{NN}(t^*) = [x(t^*_{NN}), x(t^*_{NN}+\tau), ..., x(t^*_{NN}+(m-1)\tau)]$$
    To exclude based on a Theiler window of length $\ell_{Th}$, we make sure that the $|t^*-t^*_{NN}|>\ell_{Th}$. This ensures that $y^{NN}(t)$ is a near neighbour of $y(t)$ due to the dynamics and structure of the manifold rather than temporal correlation. See Chapter \ref{chap:Invariants} for a more detailed discussion. 

    \item Produce another reconstruction with dimension $m+1$ and test to see whether $y(t^*)$ and $y^{NN}(t^*)$ remain close neighbours (keep the same time index) in the increased dimension. If pairs of points cease to be nearest neighbours by increasing the dimension by 1, this means that they were originally close together due to lack of unfolding in the attractor.

    \item Calculate the distance between nearest neighbours in $\mathbb{R}^{m+1}$
    \begin{align*}
        R_{m+1}(t^{*})^{2} &= \sum_{d=1}^{m+1}[x(t^* + (d-1)\tau)-x(t^*_{NN}+(d-1)\tau)]^2\\
        &= \left( \sum_{d=1}^{m}[x(t^* + (d-1)\tau)-x(t^{*}_{NN}+(d-1)\tau)]^{2}\right) \\
        &+ [x(t^{*}+ m\tau)-x(t^*_{NN}+m\tau)]^2\\
        &= R_{m}(t^{*})^{2} + [x(t^{*}+ m\tau)-x(t^*_{NN}+m\tau)]
    \end{align*}

    \item Calculate the ratio of distances $R(t^*)$ between nearest neighbours taking dimension scaling into account
    $$\frac{|x(t^*+m\tau)-x(t^*_{NN}+m\tau)|}{R_{m}(t^{*})}=\sqrt{\frac{R_{m+1}(t^*)^2-R_{m}(t^{*})^2}{R_{m}(t^{*})^2}}$$

    \item If $R(t^{*})>15$ (or some chosen threshold), then the pair of nearest neighbours at time index $t^*$ are false neighbours. Additionally, nearest neighbours whose distance in the dimension $m+1$ are much larger than the average \textit{size} of the attractor are also classified as false neighbours:
    $$\frac{|x(t^*+m\tau)-x(t^*_{NN}+m\tau)|}{R_A}>2$$
    where $R_A$ is a the root mean square of data centred around its mean.

    \item Repeat the above and calculate the percentage of points that are false nearest neighbours for every increase in $m$. Select the embedding dimension $m$ when the $\%FNN$ is close to 0. 
    
\end{enumerate}

\section{Non-uniform Delay Embedding}

Uniform embedding is by far the most common method for embedding time series. This is primarily due to its ease of implementation and optimisation. In a direct application, uniform delays only require the selection of two parameters, $\tau$ and $m$. However, the convenience of such an approach comes at the cost of reduced versatility and limitations, particularly when analysing systems with dynamics occurring on multiple disparate timescales.\\

Firstly, the choice to use a single delay limits the ability for the reconstruction to highlight features across multiple disparate time-scales. For example, a fast-slow system with characteristic time scales $\tau_1$ and $\tau_2$ where $\tau_1/ \tau_2 \gg1$ , the choice of selecting $\tau_1$ (i.e. slow dynamics) as the embedding lag can limit the reconstruction's ability to fully unfold attractor topologies corresponding to the fast dynamics. The dynamics the time scale of $\tau_2$ (i.e. fast dynamics) will appear as noisy fluctuations within the reconstructed state space.\\

\begin{figure}
    \centering
    \includegraphics[width = 0.85\textwidth]{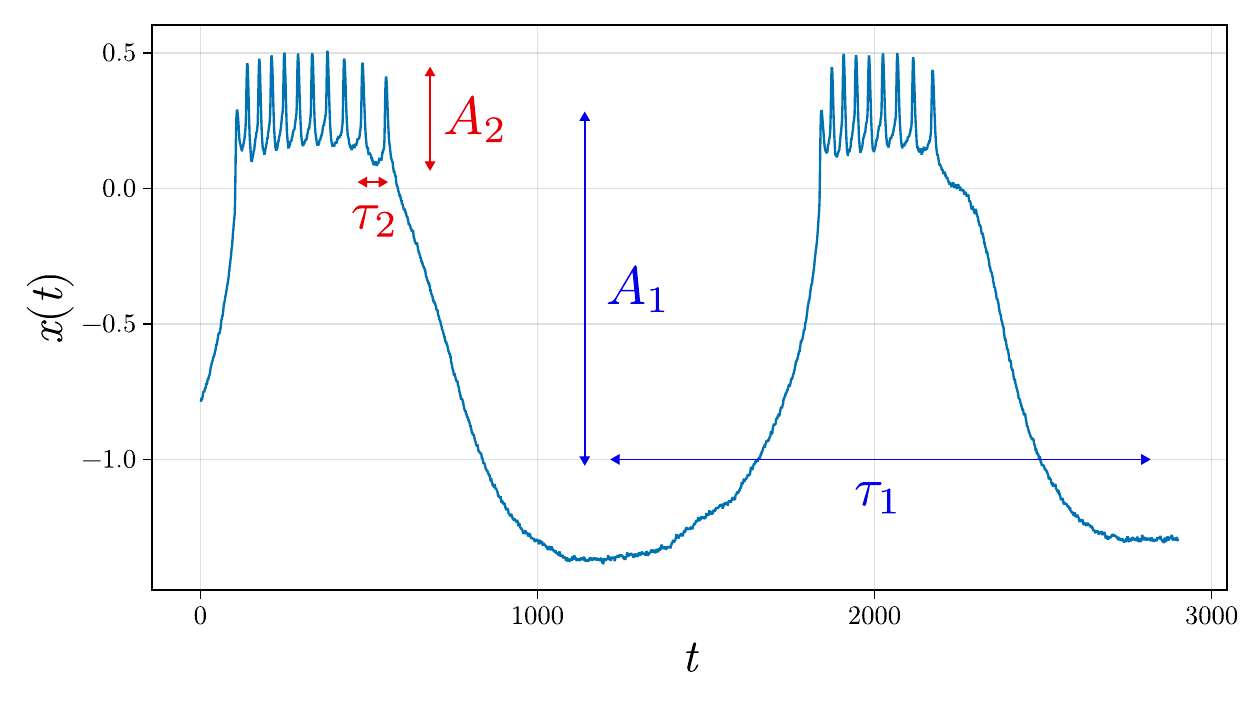}
    \caption{Illustrative example of time series two different time scales and magnitudes. Data taken from experimental measurements of lobster lateral pyloric (LP) neuron courtesy for Abarbanel et al. \cite{abarbanel1996synchronized}.}
    \label{fig:lobster}
\end{figure}

Secondly, reconstruction from a uniform delay embedding that is sufficient is not necessarily optimal. Here, we must clarify that the definition of optimal presumes some criterion or notion of quality. Casdagli noted that the quality of an embedding, defined as the reconstruction's robustness to noisy data for prediction, can vary locally throughout different regions of the attractor. Additionally, we should also consider that invariant measures such as the Lyapunov exponent also vary locally. Hence, the selection of a single embedding lag implies that all these variations may be averaged. \\

One obvious way to address these problems is to include multiple non-uniform time lags when constructing the delay vector:
\begin{equation}
    \vec{x}(t) = (x(t), x(t-\tau_1), ..., x(t-\tau_{m-1})).
\end{equation}

The selection of delay lags represents a combinatorially hard problem that grows with increasing embedding dimension. The methods proposed for constructing non-uniform delay embedding often involve the iterative selection of time lags to gradually construct a delay vector until the required embedding dimension is reached.

\subsection{Garcia and Almeida}
One of the earliest methods of choosing non-uniform delays was proposed by \cite{garcia2005multivariate}. They proposed a variation of the global false nearest neighbours methods applied to the problem of selecting time delays. This method also recursively selected lags using a proposed $N$-statistic over multiple embedding cycles. At the end of each cycle, the false nearest neighbours algorithm is used to assess the quality of the newly constructed embedding. This process is repeated until the false nearest neighbour statistic $F$ decreases below a critical threshold. 

The general algorithm goes as follows:

\begin{enumerate}
    \item For the selection of the first time lag $\tau_1$, a 2D delay embedding $\vec{x}(t)$ is first done with respect to some prospective time lag $\tau^*$ to be tested,
    $$\vec{x}(t) = (x(t), x(t-\tau^*)).$$

    \item Identify the closest neighbour $\vec{x}(t_j)$ for each point $\vec{x}(t_i)$ in the embedding reconstruction taking into account a Theiler window.

    \item Calculate the Euclidean distances $d_{1,\tau^*}(t_i), d_{2,\tau^*}(t_i)$ between any given two points,
    \begin{align}
        d_{1,\tau^*}(t_i) &= || \vec{x}(t_i) - \vec{x}(t_j) || \notag \\
        d_{2,\tau^*}(t_i) &= || \vec{x}(t_i+\delta t) - \vec{x}(t_j+\delta t) ||,
    \end{align}

    where $\delta t$ is the sampling time of the data. Simply put $d_{1,\tau^*}$ is the spatial separation between pairs of nearest neighbours in the reconstructed state space, and $d_{2,\tau^*}$ is the resulting separation one step forward in time.

    \item Calculate the $F$-statistic taken as the proportion of points whose distances ratio $d_{2,\tau^*}/d_{2,\tau^*}>10$,
    \begin{equation}
        F(\tau^*)= \frac{1}{N}\sum_{i=1}^N \mathbb{1} \left( \frac{d_{2,\tau^*}(t_i)}{d_{1,\tau^*}(t_i)}>10 \right),
    \end{equation}

    where $N$ is the length of the time series and $\mathbb{1}$ is the indicator function. The threshold of 10 is  heuristically selected by the authors based on the numerical calculations of Kennel et al. \cite{kennel1992determining}. The time lag corresponding to the first minimum in $N(\tau^*)$ is taken to be the embedding lag.

    \item Repeat steps 1-4 but with the dimension increased by one and including the newly selected time lag. Therefore, the selection of the $m^{th}$ embedding lag in a non-uniform embedding procedure will require neighbours and distances $d_{2,\tau^*},d_{2,\tau^*}$ to be calculated using the embedding with $m-1$ lags that have already been chosen and the new candidate lag $\tau^*$,

    $$\vec{x}(t) = (x(t),x(t-\tau_1),...,x(t-\tau_{m-1}),x(t-\tau^*)).$$
\end{enumerate}

\subsection{Continuity Statistic}

The continuity statistic was first proposed by \cite{pecora2007unified}. as a way to procedurally construct non-uniform delay vectors based on the idea of functional independence between vector coordinates. Takens' and Sauer both discussed the requirement that an embedding reconstruction requires vectors whose coordinates are independent. Pecora et al. proposed using a test for calculating the functional dependence between the components of a delay vector's components in order to assess the quality of an embedding. A functional dependence between vector coordinates implies,
\begin{equation}
    x(t-\tau_{m}) = F(x(t), x(t-\tau_1),...,x(t-\tau_{m-1})),
\end{equation}
where $F$ is some arbitrary function. Constructing a non-uniform delay embedding requires iteratively building of a collection of time lags $\vec{\tau} = \{ \tau_1, ...\tau_{m-1}\}$ that minimises the likelihood of a functional dependence between components. In each iteration, a prospective lag $\tau_i$ is tested for functional dependence with the existing lagged components corresponding to $\vec{\tau}$. If there is no significant functional dependence, then $\tau_i$ may be added to the collection of lags. To test this, the authors assume that $F$ is smooth and use the property of continuity to quantify functional dependence. 

\begin{enumerate}
    \item Consider an existing $m+1$ dimensional embedding $\vec{x}_m(t) \in \mathbb{R}^m$ constructed from lag $\tau=\{ \tau_1, ... \tau_{m}\}$ and a potential new embedding lag to be tested $\tau_{m+1}$. To test the functional dependence of a new lag, select a reference reference point $\vec{x}_m(t_0)$ in embedded space. If a smooth functional dependence exists, then the continuity condition states that points $\vec{x}_{m,i}$ nearby the reference point ($|| \vec{x}_{m,i} - \vec{x}_m(t_0) ||<\delta$) in reconstructed space $\mathbb{R}^m$ should have lagged $m+1^{th}$ components that are also close by to each other ($|x_i(t-\tau_{m+1}) - x(t_0-\tau_{m+1})|<\epsilon$).

    \item Calculate the proportion $p$ of points $\vec{x}_{m,i}$ whose lagged components lie within $\epsilon$ of the reference point's lag component. 

    \item Compare $p$ against a null hypothesis; i.e. that correspondence between these sets is purely by chance. Large values of $p$ suggest a strong relationship between the $m$-dimensional reconstruction and the new $\tau_{m+1}$ lagged component. Pecora et al. suggest the usage of a binomial distribution with a critical value of  $p^*=0.5$ in order to decide if a functional dependence exists with respect to some chosen $\epsilon$ due to its simplicity and robustness to noise.

    \item Repeat the continuity test for decreasing values of $\epsilon$ until the null hypothesis fails to be rejected. The smallest possible value for rejecting the null hypothesis is given as $\epsilon^*$. This value is averaged over a collection of reference points sampled from the data to calculate the continuity statistic $\langle \epsilon^*\rangle(\tau_{m+1})$.

    \item The new lag $\tau_{m+1}$ is taken as the lag corresponding to the relative maxima of the continuity statistic profile. This is repeated until the desired embedding dimension (as per Whitney's theorem) is reached. Pecora et al. also propose an undersampling statistic that can be used as a termination criterion for iterative selection of time delays. Further details can be found in the original paper. 
\end{enumerate}

\begin{figure}[h]
    \centering
    \includegraphics[width = 0.8\textwidth]{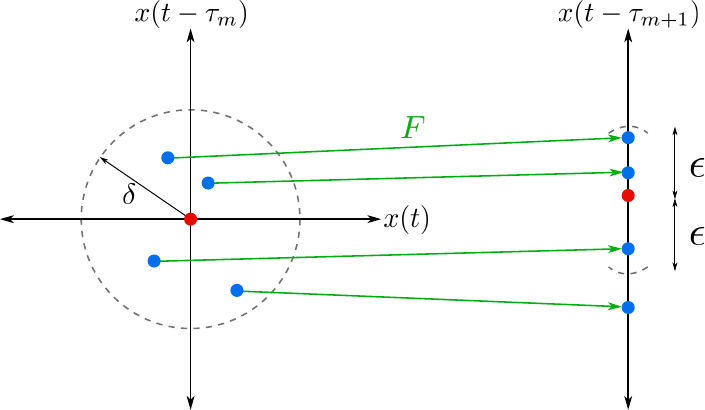}
    \caption{Calculation of the continuity statistic}
    \label{fig:continuity_statistic}
\end{figure}

\subsection{PECUZAL}

A criticism of the continuity statistic method is the ambiguity in selecting the optimal lag $\tau$ at each embedding iteration . In the original paper of Pecora et al., the definition of \textit{relative maxima} is unclear and there is no objective criterion for selecting the best lag between multiple prospective local maxima. Additionally, the method also does not consider the effects of selecting different distances $\delta$ used to define nearby points in the reconstruction. Finally, the undersampling statistic originally proposed as a breaking condition for the embedding algorithm is computationally intensive, and does not inform on which of the prospective lags should be selected.\\

Kr\"{a}mer et al. suggested that the continuity statistics approach could be combined with the L-statistic method to provide a fully automated method of constructing non-uniform embedding delays \cite{kramer2021unified}. 

\begin{enumerate}
    \item Calculate the continuity statistic profile as per Pecora et al.
    \item Identify the collection of lags $\mathcal{T} = \{\tau_{1}^{*},...\}$ corresponding to all local minima of $\langle \epsilon^*\rangle(\tau_{m})$.
    \item Calculate $L$-statistic of the current embedding arrangement (this will initially be $m=1$). (i.e $L_m$)
    \item Calculate the improvement in the $L$-statistic as a result of including each of the potential lags in the $\mathcal{T}$
    $$\Delta L_i = L_{m}-(L_{m+1} |\tau^*_i)$$

    \item Select the lag $\tau_i^*$ corresponding to the biggest decrease in the $L$-statistic (i.e. find the biggest value of $\Delta L_i$) and add this lag to the list of lags.

    \item Repeat step 1-5 until $|L_{m+1}-L_{m}|>0$
\end{enumerate}

\subsection{Other Methods}

We include for the interested reader a selection of alternative methods that have also been proposed for tackling the parameter selection problem of non-uniform embedding.

\begin{itemize}
    \item Maximum Derivatives on Projection (Nichkawde, 2013)- quantifies the unfolding of reconstructed attractors baseed on derivatives of some inferred functional dependence $F$.
    \item Reduced autoregressive models (Judd \& Mees, 1998)- Progressively select lags that give good prediction performance when used  in a pseudo-linear autoregressive model. Termination is determined using measure of model description length.
    \item Search Optimisation Algorithms - using various optimisation algorithms to search state space with an objective function that can be chosen by the user. Examples include Ant Colony Opimisation (ACO) \cite{dorigo1997ant, shen2013optimal} and Monte Carlo Decision Tree Search (MCDTS) \cite{kraemer2022optimal}.
    \item Significant Times on Persistent Strands (SToPS) - a persistent homology approach to detect and quantify the significance of pseudo-periodic structures in time series to inform and rank embedding lag selection \cite{tan2023selecting}.
\end{itemize}

\chapter{Invariants}\label{chap:Invariants}

From embedding theory, the transformation given by

\begin{equation}
    y_{t}= (x_{t},x_{t-\tau}, ..., x_{t-(m-1)\tau})
\end{equation}
yields an embedding that preserves the dynamics of a system up to a diffeomorphism. In the simplest uniform embedding case, this requires the selection of embedding dimension $m$ and embedding lag $\tau$. For state prediction, the next step would be build a model function $F$ such that,
\begin{equation}
    y_{t+1} = F(y_{t}).
\end{equation}
The process of constructing $F$ is non-trivial and will be discussed in later sections. But what can we learn about the underlying dynamical system just from examining the observed data (i.e. system characterisation)? Since $y_t$ and $x_t$ are related by way of a (complicated) transformation map, we would naturally seek identify measures, quantities and features that remain unchanged with respect to this map.

A quantity $s(\{ y_{t} \}_t)$ is said to be invariant if it doesn't depend on the way in which measurements were made. In the ideal case, we want the measurement function $h:\mathbb{R}\to \mathbb{R}$ to be sufficiently smooth as it bestows several useful analytical properties.

\begin{definition} \textbf{-- Invariant quantities}\\
    Let $s:T \times X \to \mathbb{R}$ be a function that measures some quantity. The quantity $s$ is invariant with respect to measurement function $h:S \to \mathbb{R}, S \subset X$ if
    
    $$s(\{ x_{t} \}_{t})=s( \{ h(x_{t} \}_t).$$
\end{definition}

In more mathematical terms, we can refine further to provide a definition of an invariant measure:

\begin{definition} \textbf{-- Invariant measaure}\\
    Let $X$ be a compact metrisable space, and $\mathcal{B}$ be a Borel $\sigma$-algebra (i.e. smallest subset containing all +open sets of $X$ that is closed under complement, countable unions and countable intersections. i.e this just gives nice properties). Let $T:X\to X$ be a measurable map (i.e. pre-images and images of $T$ are both measurable$\implies$ $T^{-1}\mathcal{B}\subset \mathcal{B}$ ). A probability measure $\mu:\mathcal{B}\to \mathbb{R}^{+}$ is invariant under $T$ if $\mu({T^{-1}}B) = \mu(B),\,\, \forall B\in \mathcal{B}$.
\end{definition}

Invariant measures are nice because they don't depend on how things have been measured. In practice, real world measurements are not smooth due to digitisation. But given fine enough resolution, it is sufficiently smooth. Furthermore, many invariant quantities measure some form of ``relevant'' quantity to the dynamics of the system. We can estimate a value of an invariant (and sometimes provide confidence intervals ) from a time series by calculating time and/or spatial averages. There is a huge amount of advanced technical work and theory justifying the validity of taking such averages. However, it is beyond the scope and aim of this text to delve deep into this theory. Nevertheless, a brief discussion is provided to give a flavour and a starting point for interested readers.

\section{Ergodicity and the Natural Measure}

In studying invariants where spatial and time averages must be taken, it is useful to consider the property of ergodicity. Simply put, a system is ergodic if the the trajectory of almost every point in phase space eventually passes arbitrarily close to every other points. Recall that this property is guaranteed in measure preserving map on bounded domains by the Poincar\'{e} Recurrence Theorem.

\begin{theorem} \textbf{ - Ergodic trajectories}\\
    A time evolution in a set $S$ is ergodic if and only if all ergodic components $R\in S$ either have $\mu(R) =0$ or $\mu (R) = \mu (S)$. Alternatively, an invariant measure $\mu \in \mathcal{M}(\mathcal{X}, \mathcal{T})$ is ergodic if whenever $T^{-1}B=B$ for some $B\subset \mathcal{B}$, we have either $\mu(B)=0$ or $\mu(B)=1$.
\end{theorem}

As it turns out, the existence of ergodic trajectories provides some very useful properties when trying to estimate invariants. Ergodicity guarantees that trajectories will eventually explore the entire attractor manifold (invariant set) of the system. Therefore, the calculation of a time average can be equated to the calculation of a spatial average. In essence, ergodicity provides us with two different ways to perform the averaging required for evaluating invariant measures:
\begin{enumerate}
    \item Averaging over time (i.e. across multiple iterates, or a long trajectory)
    \item Taking a weighted integral over an invariant density $\mu (x)$ (i.e. the long term probabilistic distribution of observed states at any given time.)
\end{enumerate}

The invariant density,
$$d\mu = \mu(x) dx,$$
is also called the \textbf{natural measure} on the attractor is the average time a typical trajectory spends in a infinitely small phase space element $dx$. Ergodicity is guaranteed when this measure is independent of the initial condition since all trajectories visit everywhere eventually. From this, we can describe long time averages of some quantity $\mathcal{O}$ in terms of the invariant density,

\begin{equation}
    \langle \mathcal{O} \rangle_{t}= \int \mathcal{O}(x)\mu(x)\,dx.
\end{equation}

The natural measure $\mu(x)$ can be calculated numerically if one is given a very long trajectory and constructing a histogram of observed positions $x$ in state space:
\begin{equation}
    \mu(x) = \lim _{n\to \infty} \frac{1}{n} \sum_{k=1}^{n} \delta^{d}(x-y_{k}),\quad y_k\in \mathbb{R}^d
\end{equation}
where $\delta^{d}$ is the Dirac delta function. The natural measure $\mu(x)$ is invariant because it is left the same under the forward mapping $F$. When iterating a long trajectory of length $N$ forward by one step, this is equivalent to appending a new value to the trajectory $x_{N+1}$ and removing $x_{1}$.

\begin{theorem} \textbf{ - Invariance of the natural measure}\\
    Let $\mu (x)$ be the natural measure a dynamical system in state space $S\subset \mathbb{R}^d$ characterised by the forward map $F$ (i.e. $y_{k+1} = F(y_k)$). Let $f:S\to \mathbb{R}$ be a measurable quantity. The spatial average of $f$ across $S$ is given by
    \begin{align*}
        \bar{f} &= \int_{S}\mu (x) f(x) \, dx^{d} \qquad \text{-- (spatial average)}\\
        &\approx \frac{1}{N}\sum _{i=1}^{N} f(y_k)\\
        &= \frac{1}{N} \sum_{i =1}^{N} f(F^{k-1}(y_{1))} \qquad \text{-- (time average)}\\
    \end{align*}
    To show invariance under forward mapping $F$,
    \begin{align*}
        \int_{S}\mu (x) f(F(x)) \, dx^{d} &= \bar{f} + \frac{1}{N} \left[ f(y_{N+1}-f(y_{1})) \right]\\
        &\approx \bar{f}\quad  \text{for}\,\,\, N\to\infty
    \end{align*}
    Hence, the natural measure is invariant under the forward mapping for an infinitely long trajectory.
\end{theorem}

Given the convenient properties of the natural measure, calculating interesting invariant measures now relies on the selection of ``interesting/useful'' functions of $f$.

\section{Lyapunov Exponents}

Recall that an important characteristic of chaotic dynamical systems is the sensitivity dependence on initial conditions (SDIC). Namely, trajectories with very close initial conditions diverge exponentially quickly. In simple terms, this is partially described by the following description,

\begin{equation}
    |\delta(t)| \approx \delta_{0}e^{\lambda t},
\end{equation}
where $\delta$ is the separation distance between a pair of trajectories and $\lambda$ is a quantity called the maximal Lyapunov exponent. One indicator of chaotic behaviour (SDIC) is the presence of a positive maximal Lyapunov exponent $\lambda$. In fact, as we will see, Lyapunov exponents are invariant for ergodic systems. We will look at two different ways of looking at Lyapunov exponents. First from a more theoretical approach in Section \ref{sec:lyapunov_theory}, and then from a computational approach in Section  \ref{sec:lyapunov_comp}.

\subsection{Lyapunov Exponents (the painful way)}\label{sec:lyapunov_theory}

Let us consider a $C^r$ ($r\geq 1$) vector field,
\begin{equation}
    \dot{x}=f(x), \quad \in \mathbb{R}^{n}.
\end{equation}
Let $x(t,x_0)$ be a solution trajectory of this system with initial condition $x(0, x_0)=x_0$. As Lyapunov exponents are concerned with the divergence of initially close orbits, we want to know the geometry associated with the attraction and/or repulsion of the orbits near and relative to $x(t,x_0)$. One can draw similarities with the stability analysis of fixed points. When studying fixed points, we linearise about the fixed point and observe the behaviour of states that start within some neighbourhood of the fixed point (i.e. attraction/repulsion). Analysing Lyapunov exponents is akin to replacing fixed points with whole orbits.\\

\begin{figure}
    \centering
    \includegraphics[width = 0.8\textwidth]{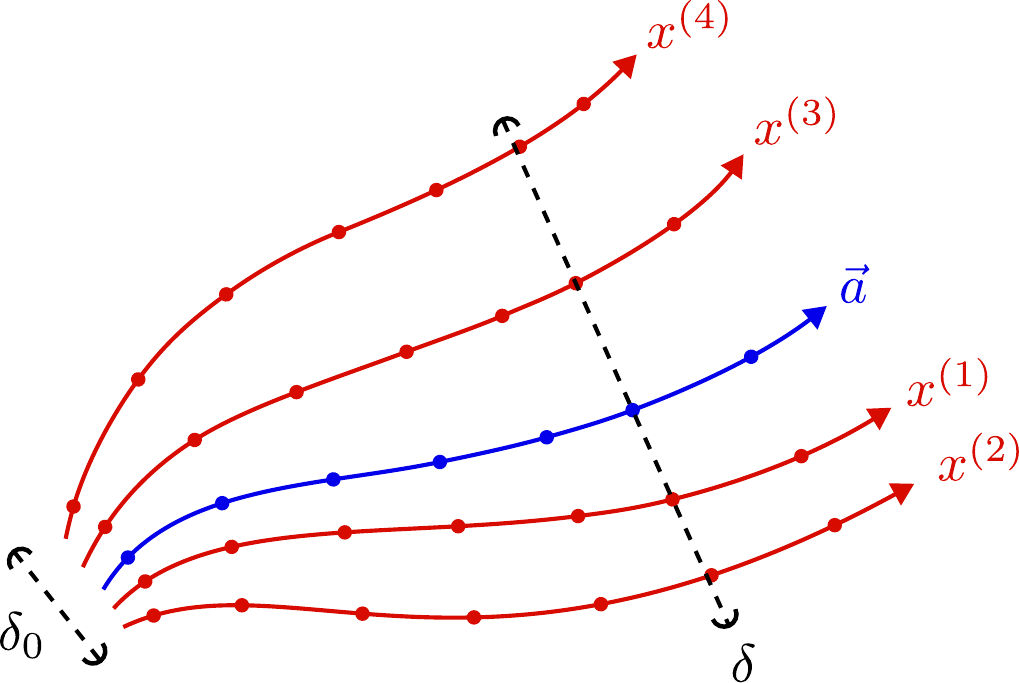}
    \caption{Linearisation with respect to a trajectory $a$}
    \label{fig:lyapunov_exponents}
\end{figure}

How does one linearise the above system about an orbit? For this, we consider replacing the orbit $x(t,x_0)$ with a sequence of $N$ points. 
\begin{align}
    \vec{a}&= [a_0,a_1,...,a_{N}]\notag\\
    &=[x(0,x_{0}),x(t_1,x_0),...,x(t_{N-1},x_{0})]
\end{align}

We then proceed with linearising the system about each point on the orbit. To do this, let there be a change of coordinates $\xi_{i}=x-a_i$ and treat $a_i$ like a fixed point. From the linearisation theorem, we can derive the following,
\begin{subequations}
    \begin{align}
        x_i &= \xi_{i}+a_i\\
        \dot{x}_i &= \dot{\xi}\\
        f(x_i) &= f(a_{i}) + f'(a_{i})(x-a_i)+...\\
        &\approx 0 + f'(a_i)\xi_i\notag
    \end{align}
\end{subequations}

Therefore, linearising about the point $a_i$, we can write down a linearly equivalent dynamical system given by
\begin{equation}
    \dot{\xi}_i = f'(a_{i})\xi_i.
\end{equation}

Repeating the above procedure for each point $a_i$ on the orbit we are linearising about, we can write down a system of ODE the describe the linearised dynamics about the orbit $x(t,x_0)$,
\begin{align}
    \dot{\xi}_{0}&= Df(x(0,x_0))\xi_{0} \notag\\
    \dot{\xi}_{1}&= Df(x(t_1,x_0))\xi_{1}\\
    \vdots\notag\\
    \dot{\xi}_{N-1}&= Df(x(t_{N-1},x_0))\xi_{N-1},\notag\\
\end{align}
where each $\xi_i$ is actually a vector defined in $\mathbb{R}^n$. We can simplify this notation by writing the following matrix ODE,
\begin{equation}
    \dot{\boldsymbol{\xi}}=Df(x(t,x_0))\boldsymbol{\xi}
\end{equation}

This system is an appropriate representation for $N\to\infty$. As this system is a matrix ODE, solutions are also $n\times N$ matrices. Let $X(t;, x(t,x_0))$ be the fundamental solution matrix of the linearised system about the orbit $x(t,x_0)$. Suppose we want to measure the amount of expansion/contraction of the neighbourhood around $x(t,x_{0})$ in the direction of $\vec{v}$. We can define the coefficient of expansion in the direction $v\in \mathbb{R}^n$ along the trajectory through $x_{0}$,
\begin{equation}
    \chi_{t}(x_{0},v)=\frac{||X(t;x(t,x_0))v||}{||v||}.
\end{equation}

Note that $\chi_{t}(x_{0},v)$ is a time-independent quantity that also depends on a particular orbit. Assuming a form of exponential divergence (i.e. $\chi_t$ scales exponentially), we can define the Lyapunov exponent in the direction $v$ along the trajectory passing through $x_{0}$ as
\begin{equation}
    \lambda(X(t;x(t,x_0)),x_{0},v)=\sup \lim_{t\to\infty} \log \chi_{t}(x_{0},v).
\end{equation}

Some key remarks

\begin{enumerate}
    \item $\lambda$ is an asymptotic quantity and thus assumes that $x(t,x_{0})$ exists for all $t>0$. This requires the phase space to be compact, boundaryless manifold or $x_{0}$ lies in a positively invariant region (there are also interesting results for negatively invariance).
    
    \item$\lambda(X(t;x(t,x_0)),x_{0},\vec{0}) = -\infty$
    
    \item Because the fundamental solution $X$ is trajectory specific, technically the results are also specific to the trajectory. However, if the system is ergodic, this does not matter as the trajectory will always visit all parts of phase space given a sufficiently large $t\to \infty$.

    \item $\lambda(X(t;x(t,x_{0})),x_{0},v)=\lambda(X(t;x(t,x_{0})),x_{1},v)$ if $x_{1}=x(t_1,x_0)$ for $t_{1}>0$. Therefore, we may simplify our notation

    $$\lambda(X(t;x(t,x_{0})),x_{0},v) = \lambda(v)$$
\end{enumerate}

\begin{theorem} \textbf{-- Geometrical properties of Lyapunov exponents}\\
    For any vectors $f,g\in \mathbb{R}^n$, and nonzero constant $c\in\mathbb{R}$,
    \begin{align}
        \lambda(f+g) &\leq \max\{ \lambda(f),\lambda(g) \}\\
        \lambda(cg) &= \lambda(g)
    \end{align}
\end{theorem}

\begin{theorem} \textbf{-- Lyapunov directions and subspaces}\\
    For any $r \in \mathbb{R}$,
    $$\{ g\in \mathbb{R}^{n} | \lambda(g) \leq r \}$$
    is a vector subspace of $\mathbb{R}^n$.
\end{theorem}

This result shows that for a state space of dimension $n$, there exists at most $n$ distinct Lyapunov exponents associated with the trajectory. More precisely, we can order these Lyapunov exponents in a specific manner to form a Lyapunov spectrum.

\begin{proposition} \textbf{-- Lyapunov spectra and nested subspaces}\\
    The set of numbers
    $$\{ \lambda(g) \}_{g\in \mathbb{R}^{n},g\neq0}$$
    takes at most $n=dim(\mathbb{R}^n)$ values $\nu_i$ which can be ordered such that
    $$\nu_{1}> \nu_{2}>...>\nu_{s},\quad 1\leq s\leq n.$$
    Furthermore, let $L_{i}= \{ g\in \mathbb{R}^{n}| \lambda(g) \leq \nu_{i} \}$. Then 
    $$\{ 0 \} = K_n \subset ... \subset L_{s+1} \subset L_{s} \subset ... \subset L_{1} = \mathbb{R}^n,$$
    such that $L_{i+1} \neq L_{i}$, and $\lambda(g)=\nu_i$ if and only if $g\in L_{i}\setminus L_{i+1}$, $1 \leq i\leq s$. The multiplicity of a Lyapunov exponent $\nu_{i}$ denoted by $k_{i}$ is given by,
    $$k_{i}=dim(L_{i})-dim(L_{i+1})$$
\end{proposition}

The above proposition on Lyapunov spectra establishes that the strictly decreasing ordering of Lyapunov exponents can only be achieved if the corresponding subspace $L_{i}$ shrinks. Therefore, the associated directions for the Lyapunov spectrum form a basis that procedurally builds up the state space $\mathbb{R}^n$.

Similar to the analysis of fixed points in flows and maps where it is required to establish the existence and uniqueness of trajectories. We then ask, under what conditions does the supremum limit (and so the Lyapunov exponent of the associated trajectory) exist? To answer this, we will require two more definitions.

\begin{definition} \textbf{-- Normal basis}\\
    A basis $\{ e_{1,}..., e_n \}$ of $\mathbb{R}^n$ is said to be a normal basis if 
    $$\sum\limits_{i=1}^{n} \lambda (e_{i})\leq \sum _{i=1}^{n} \lambda(f_i),$$

    where $\{ f_{1}, ..., f_{n} \}$ is any other basis of $\mathbb{R}^n$. In other words, a normal basis is one that is associated with a Lyapunov spectrum whose sum is minimal.
\end{definition}

\begin{definition} \textbf{-- Regular family}\\
    The fundamental solution matrix $X(t; x(t,x_{0}))$ is called regular as $t\to \infty$ if
    \begin{enumerate}
        \item $\lim _{t\to \infty} \frac{1}{t} \log | \det X(t;x(t,x_0)) |$ exists and is finite, and
        \item for each normal basis $\{ e_{1},..., e_{n} \}$
        $$\sum_{i=1}^{n}\lambda (e_{i}) = \lim_{t \to \infty} \frac{1}{t} \log | \det X(t;x(t,x_0)) | $$
    \end{enumerate}
\end{definition}

\begin{theorem} \textbf{-- Existence of Lyapunov exponents}\\
    If $X(t;x(t,x_0))$ is regular as $t\to \infty$, then
    $$\lambda (e) = \lim_{t \to \infty} \frac{1}{t} \log \chi_{t}(x_{0,}e)$$

    exists and is finite for any vector $e \in \mathbb{R}^n$
\end{theorem}
Lyapunov exponents exist and are well defined as long as there exists a fundamental solution matrix $X(t, x(t,x_0))$ that is regular. All that remains is to check for this property. This problem is addressed be the \textbf{multiplicative ergodic theorem}. This theorem shows that with respect to some invariant measure, almost all trajectories (i.e. except a set of measure zero) give rise to regular fundamental solution matrices.

\subsection{Lyapunov exponents (the easy way)}\label{sec:lyapunov_comp}

Now that we have convinced ourselves on the general theory and ideas of Lyapunov exponents, we discuss a more practical approach to understanding and eventually computing Lyapunov exponents and spectra. As previously discussed, chaotic systems have the property of \textbf{sensitivity to initial conditions (SDIC)}, which  means infinitesimally close vectors in space give rise to trajectories that diverge exponentially fast. This divergence can be completely described with the tools of linearised dynamics and Lyapunov exponents.

Let us consider two different trajectories $\mathbf{x}_n$ and $\mathbf{y}_n$ in $m$-dimensional state space that are initially nearby each other for a dynamical system characterised by a map $F$. We can use this formulation to write a linear map centred around trajectory $x_n$ that describes the separation of this pair of trajectories,

\begin{align}
\mathbf{y}_{n+1}-\mathbf{x}_{n+1} &= F(\mathbf{y}_{n}) - F(\mathbf{x}_{n}) \notag \\
&= J_{n}(\mathbf{y}_{n}-\mathbf{x}_{n})+ \mathcal{O}(||\mathbf{y}_n-\mathbf{x}_n||^2),
\end{align}
where $J_{n}=J(x_n)$ is the $m \times m$ Jacobian matrix of $F$ at $x$. For simplicity, we express the pertubation as $\boldsymbol{\delta}_{n} = \mathbf{y}_{n}-\mathbf{x}_n$ and write down the equivalent system
\begin{equation}
    \boldsymbol{\delta}_{n+1} \approx J_{n}\boldsymbol{\delta}_{n}.
\end{equation}

Focusing on the time step $n$, let $\mathbf{e}_i$ be an eigenvector of $J$ with associated eigenvalue $\Lambda_{i}$. Thus, we can express $\boldsymbol{\delta}_n$ in terms of the eigenvectors of $J_n$ with basis coefficients $\beta_i$ where,
\begin{equation}
    \boldsymbol{\delta}_{n}=\sum_{i=1}^{m} \beta_{i} \mathbf{e}_i.
\end{equation}

It follows then that,
\begin{equation}
    \boldsymbol{\delta}_{n+1}=\sum_{i=1}^{m} \beta_{i} \Lambda_i\mathbf{e}_i.
\end{equation}

Simply put, we are expressing the initial perturbation in terms of the $m$ linear eigenspaces and estimating their contraction/expansion in each eigenvector direction.

Up to now, the analysis has been local in time and space as we are linearising with respect to a specific point $\mathbf{x}_n$. Clearly, as $n$ changes and the trajectory moves to a different point in state space, the values of $J_{n}, \mathbf{e}_{i}$ and $\Lambda_i$ will change. In order to characterise the system as a whole, we need to define global objects and quantities. This can be done calculating an average over the different stretching factors and directions. To do this, we observe that forward iterations on the trajectory of length $N-1$ is equivalent to multiplying the $N$ linearised matrices $J_n$ together. Following this logic, we may then write down the following eigenvector equation with respect to a direction vector $\mathbf{u}_{i}^{(N)}$,
\begin{equation}
    \left( \prod_{n=1}^{N} J_{n}\right) \mathbf{u}_{i}^{(N)}=\Lambda_{i}^{(N)}\mathbf{u}_{i}^{(N)}.
\end{equation}

We may then define the $i^{th}$ global Lyapunov exponent $\lambda_{i}$ to be and average over the trajectory,

\begin{equation}
    \lambda_{i} = \lim _{N\to \infty} \frac{1}{N} \ln |\Lambda_{i}^{(N)}|.
\end{equation}
Similar to the conclusions in its theoretical formulations, the existence and uniqueness of the above limit is given by the \textbf{multiplicative ergodic theorem}. This is a highly nontrivial result, in particular since the multiplication of matrices is noncommutative and the logarithm cannot be exchanged with the formation of eigenvalues. For one-dimensional maps the Jacobian is a real number and the above definition reduces to,
\begin{equation}
    \lambda=\lim_{N \to \infty} \frac{1}{N} \sum_{n=1}^{N} \ln |f'(x_n)|.
\end{equation}

For a state space in $\mathbb{R}^m$, the set of $m$ different exponents is called the Lyapunov spectrum. For all initial conditions except the set of measure zero which does not lead to the natural invariant measure, the spectrum is the same. Furthermore, Lyapunov exponents are invariant under smooth transformations of the state space. If the dynamics given by $F$ is invertible, every invariant measure of $F$ is also invariant under the time reversed dynamics $F^{-1}$. Computing this values will result in a Lyapunov spectrum that is equal in magnitude, but with opposite signs. The eigenvectors corresponding to the Lyapunov exponents correspond to the stretching and contracting directions of trajectories. These directions are tangential to the global invariant manifolds (i.e stable and unstable manifolds).

\subsection{Wolf's Algorithm}

In practice, we would like to calculate Lyapunov exponents from data. How does one go about doing this?  One method of doing this is Wolf's algorithm \cite{wolf1985determining, rosenstein1993practical}, which is used to calculate the entire Lyapunov spectrum and is given as follows:

\begin{enumerate}
    \item Given a time series $x(t)$, perform a phase space reconstruction using a time delay embedding with delay vectors:
    $$\mathbf{y}(t) = (x(t), x(t-\tau), x(t-(d-1)\tau)).$$

    \item Take an initial point $\mathbf{y}(t_{0})$ and find its nearest neighbour $\mathbf{y}^{NN}(t_{0}) = \mathbf{y}(t_{0}^{NN})$ making sure to exclude temporal nearest neighbours (i.e. Theiler window)

    \item Calculate the distance between the pair of points,
    \begin{align*}
        L(t_{0}) &= |\mathbf{y}(t_{0})-\mathbf{y}^{NN}(t_{0})|\\
        &= |\mathbf{y}(t_{0})-\mathbf{y}(t_{0}^{NN})|
    \end{align*}

    \item At a later time $t_{1} = t_{0} + n$, where $n$ is the number of steps, calculate the the new distance,
    \begin{align*}
        L'(t_{1}) &= |\mathbf{y}(t_{0}+n)-\mathbf{y}(t_{0}^{NN}+n)|
    \end{align*}

    Note that $n$ should be chosen such that $L(t_{0}^{*})$ does not exceed the length scale of the attractor (e.g. to avoid the folding effects of the attractor).

    % \item Calculate the unit vector of the separation direction

    % $$\mathbf{\hat{v}}(t_{0})=\frac{1}{L(t_{0})} \left[ \mathbf{y}(t_{0}^{*})-\mathbf{y}^{NN}(t_{0}^{*}) \right]$$

    \item Let $t_{1}= t_{0}+n$ and repeat steps 2-4.
    % but with another randomly sampled. However, re-select the nearest neighbour $\mathbf{y}^{NN}(t_{1})$ this time with the additional constraint that the initial separation direction is as close as possible to $\mathbf{\hat{v}}(t_{0})$.

    \item Continue the above for as many fiducial points $M$ as required ($t_{M}$). This is then used calculate the first (maximal) Lyapunov exponent,
    $$\lambda_{1} = \frac{1}{t_M-t_{0}} \sum_{k=1}^{M} \log \frac{L'(t_{k})}{L(t_{k-1})}$$

    \item For the next largest positive Lyapunov exponent (if it exists) $\lambda_{2}$, repeat the above procedure with respect to a fiducial trajectory. But this time instead of selecting a single nearest neighbour, select two nearest neighbours and calculate the are $A(t_{0})$ formed by the points rather than the length $L(t_{0})$. Similarly, when $A(t_{0}^{*})$ exceeds the size scale of the attractor, resample another two points near the fiducial trajectory such that the resulting triangle plane's orientation is as close to that of $A(t_{0}^{*})$. Following the same formula as step 8 gives the value of $\lambda_1+\lambda_2$, which may be subsequently used to calculate $\lambda_{2}$.

    \item Repeat the above for as many positive Lyapunov exponents as required. 
    
\end{enumerate}

\begin{figure}
    \centering
    \includegraphics[width = 0.85\textwidth]{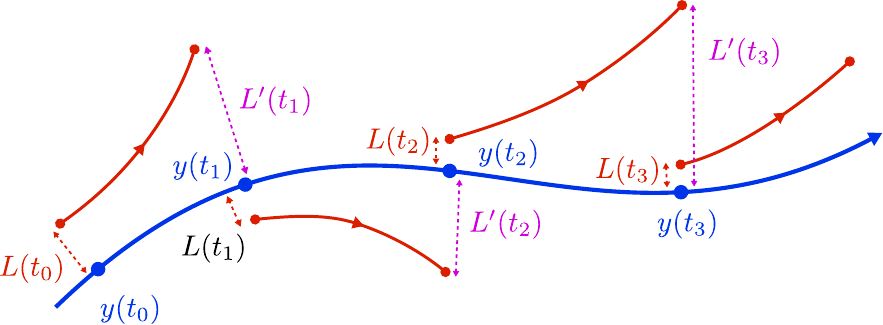}
    \caption{Wolf algorithm for calculating the maximal Lyapunov exponent $\lambda_1$}
    \label{fig:wolf_1}
\end{figure}

\begin{figure}
    \centering
    \includegraphics[width = 0.85\textwidth]{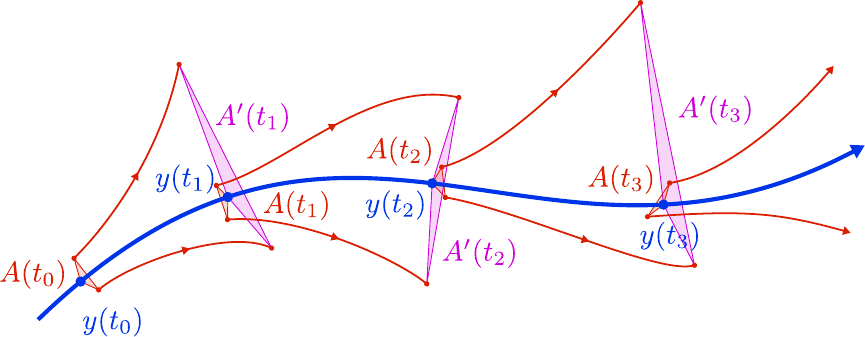}
    \caption{Wolf algorithm for calculating the first two Lyapunov exponents $\lambda_1 + \lambda_2$}
    \label{fig:wolf_2}
\end{figure}

In practice, the algorithm is typically used to only obtain the first (maximal) Lyapunov exponent. This truncated algorithm is also referred to as Rosenstein's algorithm \cite{rosenstein1993practical}. These algorithms are relatively simple to understand and implement. However, they possess a few flaws:
\begin{itemize}
    \item Continual resetting by replacement doesn't overcome inaccuracies associated with chaotic trajectories
    \item Only a single trajectory is followed.
\end{itemize}
\vspace{1em}
There are seveeral ways to address this:
\begin{itemize}
    \item Compute an average over a collection of distance in a neighbourhood of the reference point (not just the nearest neighbour).
    \item Choose a time span
    \item Choose many initial points $n_0$ rather than a single trajectory.
\end{itemize}
\vspace{1em}
Of these solutions, the last one produces the most significant improvement. To do so, the average interpoint distances of balls of nearest neighbours are tracked instead to represent the divergence of their trajectories,
\begin{equation}
    S(n) = \frac{1}{N}\sum_{n_{0}=1}^{N}\ln\big{(}\frac{1}{|{\cal U}(B_{n_{0}})|}\sum_{y_{s}\in{\cal U}(B_{n_{0}})}|y_{n_{0}+n}-y_{s+n}|\big{)}.
\end{equation}
Here, ${\cal U}(B_{n_{0}})$ is a ball of near-neighbours whose inter-point distances are within within $\varepsilon$.  $S(n)$ is a profile calculated for a fixed value of $\varepsilon$ and embedding.  It is essential to draw many curves over a range of $\varepsilon$ and $m$ -- embedding dimension.  This is so a linear scaling region can be identified within which an estimated slope is the estimate for the largest Lyapunov exponent.\\

Consider $\Delta_{0}=S_{n}-S_{n'}$ and $\Delta_{l}=S_{n+l}-S_{n'+l}$. If $|\Delta_{l}\approx|\Delta_{0}|e^{\lambda l}$ then the maximal Lyapunov  exponent $\lambda_{1}\approx\frac{1}{l}\ln\frac{|\Delta_{l}|}{|\Delta_{0}|}$.  We average over the trajectory $\{s_{n}\}$ across different embedding dimensions $m\ge m_{0}$, so that
\begin{equation}
    S(\varepsilon, m, t) = \langle\ln\frac{1}{|{\cal U}_{n}|}\sum_{s_{n'}\in{\cal U}_{n}}|s_{n+t}-s_{n'+t}|\rangle_{n}.
\end{equation}

\begin{figure}
    \centering
    \includegraphics[width = 0.95\textwidth]{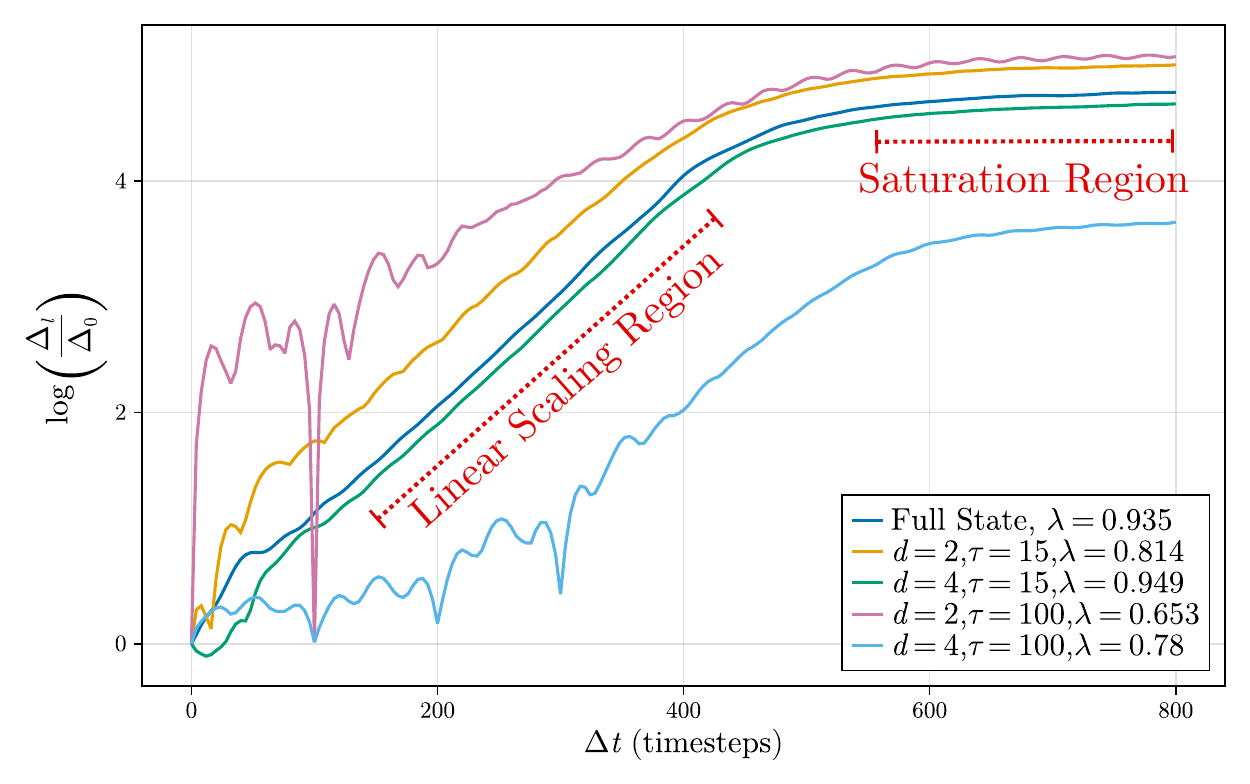}
    \caption{Extracting the Lyapunov exponent as a linear scaling region. Different curves show effects of embedding dimension that is too small, or delay lag that is too large.}
    \label{fig:lyapunov_scaling}
\end{figure}

Determine numerous ${\cal U}_{n}$ and $t=1, 2, \dots, T$ but not too long into the future.  A plot of $S$ (on log-scale) against $t$ can reveal a scaling region (straight line part) where the slope gives $\lambda_{1}$ (see Figure \ref{fig:lyapunov_scaling}). \\

If we have reliable estimates of the values of the Lyapunov spectrum then there are interesting/useful relationships between Lyapunov exponents and other dynamical invariants. \\

\textbf{Pesin's Identity} relates Kolmogorov-Sinai entropy (look it up) and the positive Lyapunov exponents via,

\begin{equation}
    h_{KS} = \sum_{i:\lambda_{i}>0}\lambda_{i}.
\end{equation}

The \textbf{Kaplan-Yorke Conjecture} relates information dimension to Lyapunov exponents. Namely,
\begin{equation}
    D_{KY} = k + \frac{\sum_{i=1}^{k}\lambda_{i}}{|\lambda_{k+1}|}
\end{equation}
where $\sum_{i=1}^{k}\lambda_{i}\ge 0$ and $\sum_{i=1}^{k+1}\lambda_{i}<0$.  This dimension is also called Lyapunov dimension and is conjectured to be the information dimension.

\subsection{Spurious Lyapunov Exponents}

Consider a dynamical system in state space $\mathbb{R}^d$. Recall that Takens' embedding theorem provides a state space reconstruction whose dynamics and manifold are diffeomorphic to the original dynamics as long as $m \geq 2d+1$. From the above details on Lyapunov exponents, a dynamical system in $\mathbb{R}^m$ will have at most $m$ unique Lyapunov exponents in its spectrum. 

Now, in most cases it is inevitable that we will need to embed observed time series data and all our analyses will be conducted in $\mathbb{R}^m$. Therefore, na\"{i}ve numerical calculation of Lyapunov exponents would provide a spectrum of up to $m$ exponents. Yet, we know from that such a system really has at most $d<m$ Lyapunov exponents, since the original dynamics are defined in $\mathbb{R}^d$. This leaves us with a problem: what do we do with the \textit{extra} Lyapunov exponents, and are they important?

One conclusion we can make is that the original system cannot have more than $d$ exponents. Therefore, there exists $d-m$ \textit{spurious exponents} that have no relation to the dynamics. How then does one go about identifying which exponents are spurious? As proposed by Parlitz \cite{parlitz1992identification}, one can consider the time-reversible properties of Lyapunov exponents. If we reverse the time ordering of the trajectories and recalculate the the Lyapunov spectrum, in theory only the real Lyapunov exponents will be preserved albeit with a negative sign. All others would be classified as spurious. However, Kantz and Schreiber \cite{kantz2003nonlinear} note that this does not work as well in practice and hence most calculations focus only on the maximal exponent $\lambda_{1}$.

\section{Fractal Dimension and Scaling}

Let us first consider the density distribution of rationals as an example. There are more formal definitions, but in the case of the number line this means that for any pair of rational numbers $i$ and $j$, it is possible to identify another rational number $k$ such that $i<k<j$. Now, given an $\epsilon \in \mathbb{R}$ that is small, how many open balls of with radius $\epsilon$ would be required to cover all rational numbers? Let this number be $N(\epsilon)$. Intuitively, this question is interested at looking at the scaling properties of some geometric property (e.g. the density of points). The case of rational numbers is rather trivial as they are are uniformly dense, and thus we expect $N(\epsilon) \propto \frac{1}{\epsilon}$. The number of open balls scales linearly with granularity. But is this true for all mathematical objects?

An alternative to the uniformly distributed, dense rational numbers, we consider an illustrative example from geography. The main land mass of Norway consists of a rugged terrain bordered by numerous fjords and irregularly shaped coastal boundaries. How long is the coastline of Norway? As one would expect, answering this question would depend on the resolution used to calculate the edges of the coastline. Coarse resolutions ignore the detailed border of the fjords and provide smooth boundary. Increasingly fine resolutions expectedly cause increases in the calculated length of the coastline. Given the (approximately) self-similar nature of this boundary, one can ask how does the coastline length scale with resolution? How can one generalise this when studying fractal structures?\\

\subsection{Box-counting Dimension}

The box counting dimension (fractal) dimension provides a way for us to address this scaling question. Consider a surface in $S\in \mathbb{R}^n$ and grid boxes with side length $\epsilon$. Let $N(\epsilon)$ be the number of boxes required to cover $S$. In the simplest case, a line segment of length $L$ will scale as 
\begin{equation}
    N(\epsilon) \approx \frac{L}{\epsilon^{1}}.
\end{equation}
\noindent Similarly, for area,
\begin{equation}
    N(\epsilon) \approx \frac{A}{\epsilon^{2}},
\end{equation}
and volume,
\begin{equation}
    N(\epsilon) \approx \frac{V}{\epsilon^{3}}.
\end{equation}

Following this line of logic, generally for a manifold of some notion of dimension $d$,
\begin{equation}
    N(\epsilon) \approx \frac{C}{\epsilon^{d}},
\end{equation}
\noindent where $C$ is some constant that is a measure of the manifold. We may use this definition to derive an expression for $d$
\begin{align}
    \log N(\epsilon) &\approx \log C -d \log \epsilon \notag\\
    d &= \frac{\log N(\epsilon) - \log C}{-\log \epsilon}.
\end{align}

With this, we can define the box-counting dimension

\begin{definition} \textbf{-- Box-counting dimension}\\
    The box-counting (fractal) dimension $d$ of a manifold is defined as 

    $$d = \lim_{\varepsilon\rightarrow 0}\frac{\log N(\varepsilon)}{-\log\varepsilon}$$
\end{definition}

This quantity can be exactly calculated for simple mathematical objects such as fractal. As an example, we consider the Cantor set. Recall that the Cantor set is constructed by the following algorithm

\begin{enumerate}
    \item Let the $\mathcal{C}_{0}$ be the unit interval $[0,1]$.
    \item Remove the middle third line segment of $\mathcal{C}_{0}$. Let the union of the two resulting disjoint intervals be $\mathcal{C}_1$.
    \item For each disjoint interval in $\mathcal{C}_{n}$, remove the the middle third and let the union of the all the resulting disjoint intervals be $\mathcal{C}_{n+1}$
    \item The Cantor set is defined as 
    $$\mathcal{C} = \lim_{ n\to \infty } \mathcal{C}_{n}$$
\end{enumerate}

A notable property of $\cal{C}$ is that it is an example of a set of measure zero. Consider the $1^{st}$ iterate of the Cantor set ($\cal{C}_{1}$). This set can be most efficiently covered by a single ball with diameter $\epsilon_1 = 1$. Similarly for $\cal{C}_2$, we require two open balls each with diameter $\epsilon_{2} = \frac{1}{3}$. In general the $n^{th}$ iterate can be covered by $2^{n}$ balls, each with diameter $\frac{1}{3^{n}}$.Thus, in the limit as $n \to \infty$, the box counting dimension of the Cantor set is calculated as,
\begin{align}
d &= \lim_{\epsilon \to 0} \frac{\log 2^{n}}{-\log \frac{1}{3^{n}}} \notag \\
&= \frac{\log 2}{\log 3}\\
&= 0.6309 \notag 
\end{align}

The box-counting dimension can be understood as a rough upper bound to the the Hausdorff dimension (which also measures smoothness). These two quantities generally differ only for some constructed examples \cite{kantz2003nonlinear}.

\begin{figure}
    \centering
    \includegraphics[width = 0.7\textwidth]{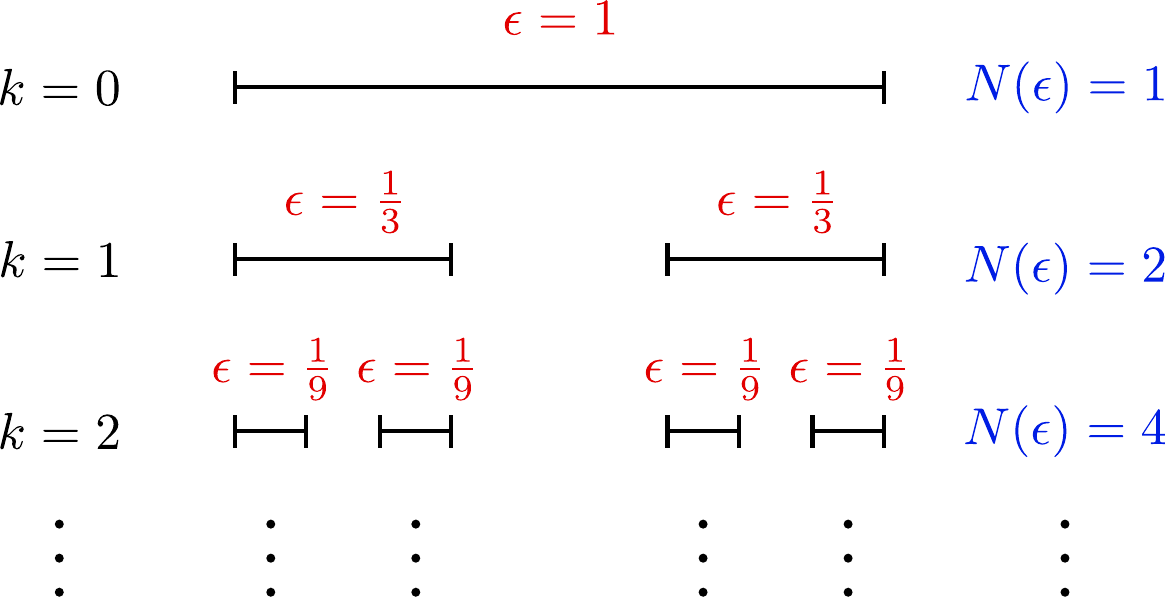}
    \caption{Cantor set construction}
    \label{fig:cantor}
\end{figure}

\subsection{Correlation and Generalised Dimension}

Recall previously that attractors and bounded dynamical systems possess a form of natural measure in phase space that describe the density distribution of points. Furthermore, natural measures are usually not homogeneous. One can see this from the distribution of orbits on the wings of the Lorenz attractor, where trajectories cluster in along fractal orbital bands. Therefore, it would not be very accurate to equate the scaling dimension of an attractor merely by the dimension of its support (i.e. the set/manifold it lives on). Ideally, we would need to consider the density of points and weight different regions of the attractor appropriately.

To account for the natural measure of the system, let $p_{\epsilon}(x) = \int_{\mathcal{U}_{\epsilon}} d\mu (x)$ be the probability of finding a typical trajectory in a ball $\mathcal{U}_{\epsilon}(x)$. We can thus define the following quantity, termed the correlation integral that describes some form of self similarity,
\begin{equation}
    C_{q}(\epsilon) = \int_{x} (p(x)_{\epsilon})^{q-1}d\mu(x) = \langle (p_{\epsilon})^{q-1} \rangle,
\end{equation}
where $q$ is an integer parameter. Just like the box counting dimension, we can define an exponential scaling relationship that describes some form of \textit{dimension},
\begin{equation}
    C_{q} \propto \epsilon^{(q-1)D_{q}}, \quad \epsilon \to 0
\end{equation}
Taking equality and limits, we arrive at a form for the generalised dimension

\begin{definition} \textbf{-- Generalised dimension}\\
    The generalised dimension of order $q$ is defined as,
    $$D_{q}=\lim_{\epsilon \to 0} \frac{1}{q-1} \frac{\ln C_{q}(\epsilon)}{\ln \epsilon}$$
    where $C_{q}(\epsilon)$ is the correlation integral,
    $$C_{q}(\epsilon) = \int_{x} (p(x)_{\epsilon})^{q-1}d\mu(x) = \langle (p_{\epsilon})^{q-1} \rangle$$
\end{definition}

There are special names given for the first few orders of generalised dimension:
\begin{itemize}
    \item $q=0$: Capacity/Box-counting dimension
    \item $q=1$: Information dimension
    \item $q=2$: Correlation dimension (measure of mutual occurrence)
\end{itemize}
Generally, if $q_1 > q_2$, then $D_{q_1}<D_{q_2}$.

\subsection{Grassberger-Procaccia Algorithm}
The correlation dimension contains useful information about the fractal and scaling properties of a system's attractor and so we would like to have a method to calculate it. A common method for calculation is the Grassberger-Procaccia algorithm \cite{grassberger1983measuring}. To do so, we first embed in $d_{E}>2m+1$ where $m$ is the true dimension.  Because we don't know $m$ we do this for a number of $d_{E}$.  Consider the correlation sum,
\begin{equation}
    \mathcal{C}(\varepsilon, N) = \frac{2}{N(N-1)}\sum_{i=1}^{N}\sum_{j=i+1}^{N}\Theta(\varepsilon-\|x_{i}-x_{j}\|),
\end{equation}
where $\Theta(x)$ is the Heaviside step function and $C(\varepsilon, N)$ is evaluated across the pairs of points within a distance $\varepsilon$ of each other.

\begin{definition} \textbf{-- Correlation dimension}\\
    The correlation dimension is defined to be
    $$D_{2} = \lim_{\varepsilon\rightarrow 0}\lim_{N\rightarrow\infty}\frac{\log \mathcal{C}(\varepsilon, N)}{\log\varepsilon}.$$
\end{definition}

That is, assume the correlation sum scales like a power law, $\mathcal{C}(\varepsilon)\approx\varepsilon^{D_{2}}$ for small $\varepsilon$ and a high enough $d_{E}$.  If we take the lkimit of $N\to\infty$, the correlation sum approaches the correlation integral,
\begin{equation}
    \mathcal{C} \rightarrow \mathcal{C}(\epsilon,N):=\int\int \Theta(\|x-y\|<\varepsilon)d\mu(x)d\mu(y).
\end{equation}

The Grassberger-Procaccia algorithm directly calculates the correlation sum and plots $\log \mathcal{C}(\varepsilon, N)$ against $\log\varepsilon$.  One then finds a linear scaling region and approximates the slope to get $D_{2}$.  As with most estimation algorithms for estimating invariant measures, best practice requires generating multiple trajectories with different initial conditions, and the usage of a Theiler window when identifying near neighbours.

\begin{figure}
    \centering
    \includegraphics[width = \textwidth]{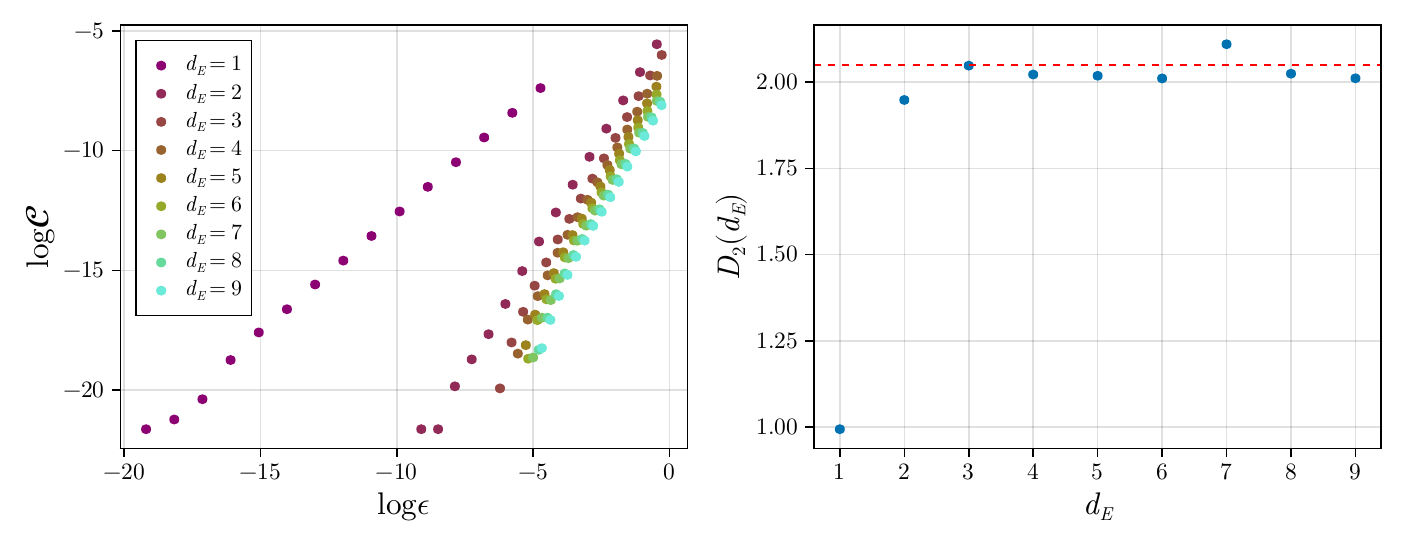}
    \caption{Using the linear scaling region of the correlation sum to infer the correlation dimension from the Lorenz data with $\tau=15$ timesteps and varying embedding dimension $d_E$. Red dashed lined on right figure corresponds to the true correlation dimension value of $2.05$}
    \label{fig:GP_algorithm}
\end{figure}

Note, the error estimate of the fit are not confidence intervals for the correlation dimension estimate.  Noise will fill the available space so in the case of noise or high noise contaminated signals expect estimates of $D_{2}$ to be the embedding dimension.  Always look at (and show) the Correlation Sum plots and then decide, given the data, if extracting a number of $D_{2}$ is useful, i.e., show the scaling plots and the region used to estimate $D_{2}$.

One of the problems with the Grassberger-Proccacia algorithm is that it assumes that a scaling region exists and then makes a best guess at what it should be.  Kantz and Schreiber \cite{kantz2003nonlinear} discuss various palliative measures to overcome this. We note in passing that other dimensions have also been described---most common of these is the box-counting dimension $d_0$, the Hausdorff dimension, and the information dimension. These deviate from the correlation dimension in the choice of metric used to measure closeness in the correlation integral.

In what follows we describe pragmatic alternatives to the correlation dimension that are better suited to application to finite data sets.  We first begin by listing several challenges with the vanilla Grassberger-Procaccio scheme.

\begin{itemize}
    \item In the G-P implementation, the scaling region is first assumed to exist---then found.
    \item Scaling in the correlation integral is bounded at large scales by the size of the attractor. This is indicated by the saturation of profile.
    \item At small length scales, even with noiseless data, quantisation effects produce artifacts that alter the estimates of $\log(\mathcal{C})$ (which is another problem), quantization screws things up.
    \item At small length scales the distribution of inter-point distances is biased because the points (coming from a trajectory) are not independent.
    \item There are less points for small length scales, and many points at large scales, this introduces a statistical correlation as the sample size from which the correlation integral is estimated varies with scale.
    \item As a result of all of the above, the scaling region is always over a finite range of length scales---this represents deterministic structural properties of the attractor rather than the asymptotic.
    \item Estimates of variance (error) in the estimate of correlation dimension are not forthcoming.
\end{itemize}

\subsection{Judd's Algorithm}

One refinement of the Grassberger-Procaccia algorithm was proposed by Judd \cite{judd1992improved} by expressing the formula for the correlation integral $\mathcal{C}$ to be in terms of a probability distribution. This gives the following alternative formula for the correlation dimension
\begin{equation}
    d_{c}(V, \mu) = \limsup_{\epsilon \to 0} \frac{\log \text{Pr}(||X-Y||<\epsilon| X,Y \in V) }{\log \epsilon},
\end{equation}
where $\mu$ is the natural measure, and $V$ is the set in phase space (e.g. the attractor). Here, the correlation sum $C(\epsilon)$ can be thought of as an estimator of
\begin{equation}
    \rho_V(\epsilon) = \text{Pr}(||X-Y||<\epsilon| X,Y \in V).
\end{equation}

The second refinement of Judd's algorithm is to model the attractor as a composition of two compoennts $V=U\times Z$ where $Z$ is a connected subset of a smooth manifold that is topologically equivalent to a closed subset of Euclidean space, and $U$ is a ``Cantor set-like'' fractal part. This is used to show the following
\begin{equation}
    \rho_{V} \approx \epsilon^{d_{C}}p(\epsilon)
\end{equation}
where $\epsilon^{d_C}$ relates to the fractal self-similarity due to $U$, and $p(\epsilon)$ corresponds to the smooth curvature of $Z$. Furthermore, $p(\epsilon)$ is modelled as a polynomial of degree $k$, and typically $k=1,2,3$. The motivation of this approach is that the separation of random points on the set $V$ are broadly determined by the self similarity structure of $U$ (the exponential term), and the bending and twisting of the smooth manifold (polynomial term) $Z$. The polynomial term will depend on the underlying set we are trying to estimate and so needs to be estimated from the data as well. Interested readers are encouraged to look at the original paper for details.

\subsection{Gaussian Kernel Algorithm}

Both the Grassberger-Procaccia and Judd algorithm do not account for the effect of noise. This consideration should be quite important as noise can effectively destroy the self-similarity fractal structure at the resolution equal to magnitude of the noise, which can result in overestimates of the correlation dimension. One alternative approach is to assume that all observations $X$ and $Y$ contain observational Gaussian noise. To do so, Diks et al. \cite{diks2003correlation} and Yu et al. \cite{yu2000efficient} modify the original correlation integral by replacing the Heaviside function with a Gaussian kernel.

\begin{definition} \textbf{-- Gaussian kernel correlation sum}\\
    The correlation sum assumming observations with Gaussian noise corruption is given as,
    $$T(h,m) = \int \int w \left( \frac{||x-y||-\epsilon}{h} \right)d\mu (x) d\mu (y)$$
    where $h$ is the bandwidth parameter, $m$ is the embedding dimension of the time delay embedding and
    $$w(x) = e^{-\frac{x^2}{4}}$$

    is the Gaussian kernel.
\end{definition}

Let $\rho_m$ be the natural measure associated with the $m$-dimensional time delay embedding, and assume the self-similarity relationship
\begin{equation}
    T(h,m)\propto h^{d_C}.
\end{equation}
\noindent The discrete formulation of the the Gaussian kernel correlation sum is given by
\begin{equation}
    T(h,m) = \frac{1}{N(N-1)} \sum_{i=1}^{N} \sum_{j\neq i}^{N} e^{-||x_{i}- y_{j}||^{2}/4h^2}.  
\end{equation}
Here, all observations must be normalised to zero mean, unit standard deviation. Following this, proceed to estimate $d_{C}$ using the linear scaling region as per the Grassberger-Procaccia algorithm.

\subsection{Information Dimension}

Note that the generalised dimension $D_{q}$ exists for different values of $q$. The case of $D_{1}$ is often called the \textbf{information dimension} and can be derived from the original generalised dimension definition.

\begin{definition} \textbf{-- Information dimension}\\
    $$D_{1}=\lim_{\epsilon \to 0} \frac{\langle \ln p_{\epsilon} \rangle_{\mu}}{\ln \epsilon}$$
\end{definition}

The information dimension is named after the fact that $\langle p_{\epsilon} \rangle_{\mu}$ is the average Shannon information needed to specify a point $x$ with accuracy $\epsilon$. The information dimension specifies how this amount of information required scales with resolution $\epsilon$.

\section{Entropies and Information Theory}

The entropy of a dynamical system is another measure that is invariant under topological conjugacy, and comes in several different forms. It has several anologues from thermodynamics and is a part of \textbf{information theory}. This theory has developed since the 1940 whose main contributions came from Shannon, Renyi, Kolomogorov and Sinai. Information theory attempts to mathematically describe the measurement and flow of information within dynamical processes, and thus provides an important approach to time series analysis.\\

We begin by first considering a static distribution of a discrete random variable $X$ with $M$ possible outcomes. Let $p_{X}$ be the probability that $X=M$. From this, we define information and Shannon entropy.

\begin{definition} \textbf{-- Shannon information and entropy}\\
    The information contained by an observation or outcome $X=i$ of a random process is given by
    $$I_{i} = -\log(p_X(i)).$$
    where $p_X$ is the probability mass function of random variable $X$. The Shannon entropy of a discrete probability distribution $p_{X}$ is defined as the average information contained in the distribution
    $$H(X) = -\sum_{i} p_{X}(i)\log (p_{X}(i)).$$
    Similarly, for continuous random variable $X$,
    $$H(X) = -\int p_{X}(x)\log p_{X}(x) dx$$
\end{definition}

Similar to the correlation integral, we can defined a generalised $q$-order form of entropy, called the $q$-order Renyi entropy

\begin{definition} \textbf{-- R\'{e}nyi entropy}\\
    Let observations $x$ be defined on a set of disjoint boxes $\mathcal {P}_j$ of side length $\epsilon$. Let $p_{j}= \int_{\mathcal{P}_{j}}d\mu(x)$ be the fraction of the measure contained in the $j^{th}$ box. The $q$-order R\'{e}nyi entropy with resolution $\epsilon$ is defined as

    $$\tilde{H}_{q}(\epsilon) = \sup_{\mathcal{P}_{\epsilon}}\frac{1}{1-q}\ln \left( \sum\limits_{j} p_{j}^{q} \right)$$
\end{definition}

In this form, the Shannon entropy is a specific case of the R\'{e}nyi entropy with $q = 1$. For any given probability distribution, the entropy increases as $p_X$ tends to a uniform distribution.

Another entropy that is interesting to look at is the Kolmogorov-Sinai entropy. This quantity tries to measure the predictability and rate of if information loss in a dynamical system.

\begin{definition} \textbf{-- Kolomogorov-Sinai entropy}\\
    Let $\cal P _\epsilon$ be a division of phase space of some observable into $m$ partitions $I_{1}, I_{2},...,I_{m}$ with size scale $\epsilon$. Consider the joint probability $p_{i_{1}, i_{2},...,i_{m}}$ that at an arbitrary time $n$ the observable falls into partition $I_{i_1}$, and at time $n+1$ it falls into interval $I_{i_2}$ etc. Consider the following entropy,

    $$K_{n}=-\sum\limits_{i_{1}, i_{2},...,i_{m}} p_{i_{1}, i_{2},...,i_{n}} \log p_{i_{1}, i_{2},...,i_{n}}.$$

    The Kolmogorov-Sinai entropy is defined as

    $$h_{KS} = \lim_{T \to 0} \lim_{\epsilon \to 0^{+}} \lim_{N \to \infty} \frac{1}{NT} \sum\limits_{n=0}^{N-1} K_{n+1}-K_{n}$$
\end{definition}

Similar to how Shannon entropy describes the spread of a probability distribution, Kolmogorov-Sinai does a similar job but within the context of dynamical systems. Intuitively, $h_{KS}$ describes the rate at which information is lost in the forward evolution of some trajectory. This is done by replacing continuous dynamics with a discrete Markov chain. Another way to look at $h_{KS}$ is that given the knowledge of the entire historical trajectory, how much further can one reliably predict in the future before the historical information becomes irrelevant?

\chapter{Prediction and Function Approximation}

\section{Modelling with Dynamical Systems}

George E.P. Box was a statistician who coined the phrase in his 1979 report, saying that ``All models are wrong, but some are useful - now it would be very remarkable if any system existing in the real world could be exactly represented by any simple model" \cite{box1979robustness}. Generally speaking, the task of modelling falls right in the centre of the work scope typically asked of any applied mathematician. Often times, it is difficult to analyse and describe what we observe in the real world fully and exactly. To handle this, one needs to make simplifications such that subsequent analyses are tractable and workable. For the applied mathematician, this means carefully choosing the correct assumptions. Once a model is constructed, it may then be used for a variety of tasks such as prediction and characterisation.

As a general rule, we should aim to construct models that are sufficient for the analyses and conclusions that we need to conduct. This means accounting for the amount and quality of data, and the needs of the tasks to be performed. However in the context of dynamical systems, it would be pertinent to consider these favourable model qualities:
\begin{itemize}
    \item[\textbf{Q1.}] \textbf{Does the model predict well?} - Prediction, if desired, can come in many forms. These can include future prediction (next $n$ time steps), prediction of the onset of some phenomenon, detection of change points, classification etc. Each of these cases have their own unique measures of performance. For example: prediction horizon, F1 recall, AUROC, Matthews correlation coefficient.
    \item[\textbf{Q2.}] \textbf{Does the model reproduce the observed dynamics?} - Good prediction performance does not necessarily imply that the model accurately describes the dynamics. For example, it is not difficult to produce short-term forecasts of chaotic time series. In contrast, long term forecasts are difficult because any small errors in the model result in exponentially increasing prediction errors over time. If we come to terms with this natural limitation, one alternative way to assess models is that it should at least preserve the overall dynamics of observed system. Models that produce long term forecasts that are incorrect but behave similarly to the system (e.g. reproduces invariant measures, creates similar attractors) are more informative than models that fail to do so (e.g. by drifting away from the domain of interest)
    \item[\textbf{Q3.}] \textbf{Is the model simple and generalisable?} - It is important not to conflate model simplicity with generalisability. Simple models may have a fewer number of parameters that need to be adjusted, and the structure of the model may also be relatively straightforward. This is often a favourable trait because large numbers of parameters can act as sources for model misspecification, which in turn can impact model performance. Generalisability is the ability for a model to perform well on tasks that are outside its training/fit data. There are two kinds of generalisability that should be considered: generalisation within the domain (interpolation) and generalisation outside the domain (extrapolation). Of these, the latter is much harder to achieve, but also much more useful. Overall, simple models do not imply generalisbility (and vice versa), but there is certainly a case to be made that there is a tendency for the two properties to be correlated. There are exceptions to this such as neural networks, of which there is active research to formulate a generalisation theory of model performance in these cases.
\end{itemize}

\vspace{1em}
These discussions are quite broad and may be applied to any field of mathematics that involves the task of modelling and prediction. The dynamical systems approach, and more specifically time series analysis, refines this by assuming that any given system may be generally represented as an arbitrary autonomous dynamical system:
\begin{equation}
    \dot{x}=f(x), \quad x(t+\delta t) = \phi_{\delta t}(x(t)).
\end{equation}
We may also consider a discrete representation of this as well,
\begin{equation}
    x_{n+1} = f(x_n).
\end{equation}
In many cases, the aim is to approximate or model the function $f$ (or evolution operator $\phi$). If this can be accurately done, it is possible to perform prediction and system characterisation. However, it is important to note that this may not be possible in all cases. It is quite possible that we must settle for alternative, more inferior outcomes.

\section{Autoregressive Processes}

One of the most common and simplest approaches for time series prediction is the fitting of an autoregressive model to observed data. In it's most general form, an autoregressive model $AR(p)$ can be defined as
\begin{equation}
    X_{t}= \sum\limits_{i=1}^{p} \varphi_{i}X_{t-i} + \epsilon_{t}.
\end{equation}
One can also consider larger and nonuniformly spaced lag terms in the model. For stochastic time series, it is possible to include moving average and noise terms as well,
\begin{equation}
    X_{t}= \epsilon_{t} +\sum\limits_{i=1}^{p} \varphi_{i} x_{t-i} + \sum\limits_{i=1}^{q} \theta_{i}\epsilon_{t-i}, \quad \epsilon_{i}\sim \mathcal{N}(0,\sigma^2).
\end{equation}
In general, the fitting of an autoregressive model $AR(p)$ is equivalent to fitting a linear model from a $p-1$ dimensional delay embedding. Care must be taken in the selection of $\varphi_{i}$ in order to ensure that iterative applications of the model do not result in unbounded growth in values. To do so, we can write the autoregressive process in terms of the backshift operator $BX_{t}= X_{t-1}$:
\begin{subequations}
    \begin{align}
        \Phi(B) &= 1 - \sum\limits_{j=1}^{p} \varphi_{j}B^{j}= 1 -\varphi_{1}B - \varphi_{2}B^2-...-\varphi_{p}B^{p}\\
        \Phi(B) X_{t} &= \epsilon_t 
    \end{align}.
\end{subequations}

\begin{theorem} \textbf{-- Stability of autoregressive processes}\\
    An autoregressive process $AR(p)$ is stationary if and only if all the roots $\beta_i$ of the characteristic polynomial $\Phi (B)=0$ are greater than one.
\end{theorem}

\section{Nearest Neighbours (the method of analogues)}

A simple approach to next time step prediction is via the method of analogues also known as nearest neighbour prediction, which is based on a simple weather prediction model. Suppose we wish to predict tomorrow's $(n+1)$ weather. One straightforward way is to look back in history and find another day $(m<n)$ whose weather is the closest match for today $(n)$. We can then conclude that tomorrow's weather $(n+1)$ will be the same as day $(m+1)$, up to some assumed degree of error.\\

In terms of time series (embedded) $\{ x(t) \}$ data, the nearest neighbour prediction algorithm proceeds as follows:
\begin{enumerate}
    \item Let the current state be $x(t)$.
    \item Find the nearest neighbour $z^{NN}(t) = x(k)$ where $t-k>l_{Th}$, and $l_{Th}$ is a Theiler window.
    \item The predicted value is given by
    $$x(t+1) = x(k+1)$$
\end{enumerate}

Clearly this approach is very limited and has numerous potential problems. For example, it may be that the observed point does not possess any near neighbours that are sufficiently close. One will likely encounter problems as well when trying to perform predictions for states that are near the boundaries of invariant regions (e.g. near a separatrix). We can do better by using a local average of $N$ closest neighbours instead:
\begin{enumerate}
    \item Let the current state be $x(t)$.
    \item Find the $N$ nearest neighbours $z_{i}^{NN}(t)=x(k_i)$ within some Theiler window.
    \item The predicted value is given by
    $$x(t+1) = \frac{1}{N}\sum\limits_{i=1}x(k_i+1)$$
\end{enumerate}

The $N$ nearest neighbour approach is more robust to noise. However, it is prone to error when trying to predict points on the attractor that have low density or from sparse data as the nearest neighbour may be quite far away. This can be addressed by choosing neighbours based on some upper bound of distance $\epsilon$ from the reference instead (i.e. find all $|| z_{i}^{NN}(t) x^{NN}(t) || < \epsilon$). Additionally, one can also weight the prediction contributions of each neighbour based on the distance from the reference point $x(t)$ using a kernel function. For example, we can use an exponential kernel function with size parameter $S$,
\begin{equation}
    w(d) = e^{-Sd/\bar{d}},
\end{equation}
where $d = || z(t) - x(t)||$ and $\bar{d}$ is the average of these distances. Thus, we may make the prediction as
\begin{equation}
    x(t+1) = \frac{\sum x^{NN}(k_{i}+1) w(d_{i})}{\sum w(d_{i})}
\end{equation}

\section{Neural Networks (who cares about number of parameters anyway?)}

\subsection{The perceptron}

Neural networks are a wide class of mathematical models that have increasingly gained attraction. However, their existence is no new and has existed within the literature for more than 70 years, and draws inspiration from historical perceptions of how biological neurons communicate. Sequence of neurons pass messages with each other, with signals eventually building pass a threshold and subsequently ``fire''. This also led to the Hebbian description "neurons that fire together, wire together". We begin our discussino on neural networks by first looking at the universal approximation theorem \cite{hornik1989multilayer}.

\begin{theorem} \textbf{-- Universal Approximation Theorem}\\
    Let $\mathcal{C}(X, \mathbb{R}^m)$ denote the set of continuous functions from a subset $X \subset \mathbb{R}^n$ to $\mathbb{R}^m$. Let $\sigma \in C(\mathbb{R}, \mathbb{R})$. Note that $(\sigma \circ x)_{i} = \sigma (x_{i})$ is applied element wise. Then the following holds:\\
    
    $\sigma$ is non-polynomial if and only if for every $n,m \in \mathbb{N}$, compact $K\subset \mathbb{R}^n$, $f\in \mathcal{C}(K, \mathbb{R}^{m})$, $\epsilon>0$, there exists $k\in \mathbb{N}$, $A\in \mathbb{R}^{k\times n}$, $b\in \mathbb{R}^k$, $C\in \mathbb{R}^{m\times k}$ such that
    $$\sup_{x\in K}|| f(x)-g(x) ||<\epsilon,$$
    where $g(x)=C\cdot(\sigma\circ (A\cdot x+b))$
\end{theorem}

\begin{figure}
    \centering
    \includegraphics[width = 0.85\textwidth]{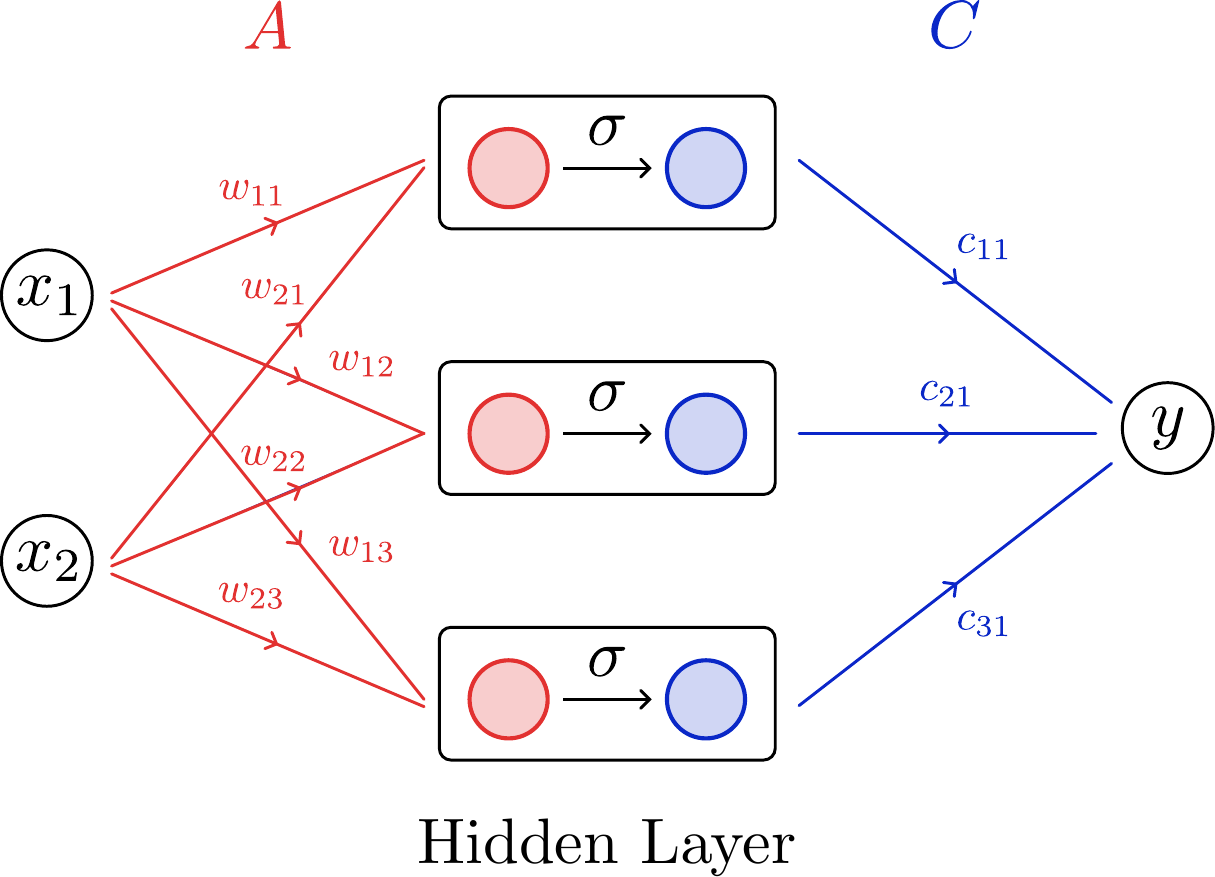}
    \caption{Feedforward neural network (perceptron) with one hidden layer.}
    \label{fig:MLP}
\end{figure}

The implications of the universal approximation theorem when one considers the model structure of a basic neural network, also termed the perceptron. In essence, the theorem states that as long as $\sigma$ is a non-linear function that is non-polynomial (e.g. exponential, logarithmic, sigmoidal), a neural network with a single hidden layer can approximate any continuous function given that its hidden layer is wide enough. This extreme flexibility means that one can easily use the neural network architecture to train the often complicated and non-analytical evolution operators for time-series prediction:
\begin{align}
    x(t+1) &= f(x(t)) \notag\\
    &= \phi^{NN}(x(t))\\
    &= C \cdot (\sigma \circ (A\cdot x +b)) \notag
\end{align}

Here, $C$ and $A$ are the output and input vectors to the perceptron, and $b$ is a bias vector. There are numerous options for the choice of nonlinearity $\sigma$ (e.g. $\tanh$, ELU, ReLU, etc.). The choice of which can impact the performance of the model.

\subsection{Gradient Descent Algorithm}

Fitting the perceptron model requires the selection (training) of the constant elements (weights) in $C, A$ and $b$ that minimises function (prediction) error, also termed the loss $\cal L$. More generally, the loss is understood as
\begin{equation}
    \mathcal L(x) = d(f(x), \phi^{NN}(x)),
\end{equation}
where $d$ is a distance metric/kernel. This is done via the gradient descent algorithm as outlined below:
\begin{enumerate}
    \item Let the collection of initial model weights be given as $\lambda_{0} = \{ \alpha_{0,1}, \alpha_{0,2},... \}$ and the evaluation of the neural network using these weights as $\phi^{NN}_{\lambda_0}$.
    \item Define input data $\{ x_{n}\}$ and output data $\{ y_{n}\}$ between which we want to learn a mapping.
    \item Calculate predicted output data by inputting $\{ x_{n}\}$ into the network:
    $$\hat{y}_{n}= \phi^{NN}_{\lambda_{0}}(x_n)$$
    
    \item Calculate the loss given by
    $$\mathcal{L}(x_{n})= d(y_n,\hat{y}_n)$$

    \item Calculate the gradient vector of the loss with respect to the model weights
    $$\vec{\nabla} = \frac{\partial \mathcal{L}(x_{n})}{\partial \lambda_{0}} =  \left( \frac{\partial \mathcal{L} (x_n)}{\partial \alpha_{0,1}}, \frac{\partial \mathcal{L} (x_n)}{\partial \alpha_{0,2}},... \right)$$

    \item Update model weights to reduce loss
    $$\lambda_{n+1} = \lambda_{n} - \beta \cdot \frac{\vec{\nabla}}{|| \vec{\nabla} ||},$$
    
    where $\beta>0$ is a small value called the learning rate.  
\end{enumerate}

\vspace{1em}
The gradient descent algorithm does not guarantee that the global optimum will be reached and will often settle on local optima instead. Furthermore, due to the complexity and dimension of the optimisation landscape, it is usually impossible to find the global minima. However, some local minima may still be sufficiently useful for prediction. There are several refinements that can be made to improve the gradient descent algorithm and partially address these flaws:
\begin{itemize}
    \item \textbf{Batch gradient descent -} calculate the average of the loss across all input data before evaluating the gradient and updating weights. This produces a loss surface that is smoother and easier to optimise along.
    \item \textbf{Stochastic gradient descent -} Addresses the same problem as batch gradient descent, but by selecting random input points rather than averaging across the whole batch in order to minimise the effect input temporal correlations between successive weight updates.
    \item \textbf{Mini-batch gradient descent -} Do both batch and stochastic gradient descent together by randomly choosing small subsamples of inputs to smooth across at each update step.
    \item \textbf{Momentum -} Includes some notion of momentum in the calculation of the $\vec{\nabla}$
    
    $$ \vec{\nabla}_{n+1} = \gamma\vec{\nabla}_{n} + (1-\gamma)\vec{\nabla}^*,$$

    where $\vec{\nabla}^*$ is the newest estimate of the gradient. This implementation avoids sudden changes in the gradient direction and causes the algorithm to overshoot local optima, thus encouraging better exploration of the loss landscape.

    \item Other optimisation algorithms (e.g. ADAM, AdaGrad, AdaDelta)
\end{itemize}

\subsection{Variations of Neural Networks}
\label{sec:rnn}

\textbf{Feedforward Neural Networks (FNN) - }This is an extension of the classical perceptron to include multiple hidden layers.\\

\noindent \textbf{Recurrent Neural Networks (RNN) and LSTMs - }Architecture designed for sequence-sequence prediction. Typically only one hidden layer, but can have more. Each hidden layer node possesses a state value that evolves according to the equation\\
\begin{equation}
    \vec{s}(t+1) = \sigma \circ (C_{in} \vec{x}(t)+C_{rec}\vec{s}(t) +\vec{b}),
\end{equation}

and next step prediction is given as
\begin{equation}
    \vec{x}(t+1) = C_{out}\vec{s}(t+1).
\end{equation}

This structure builds in recurrence within the model and allows the network to essentially be a universal autoregressive process. However, training RNNs required the usage of modified gradient descent algorithm to account for errors propagating through time. This is done via the backpropagation through time (BPTT) algorithm that tries to unfold the recurrent network with $T$ time step recurrence as a feed forward network with $T$ hidden layers. However, this algorithm suffers from the problem of exploding and vanishing gradients due to the successive multiplication of very large/small numbers. Long short term memory (LSTM) networks are a refinement of the RNN that encodes a memory buffer in the nodes to mitigate this problem \cite{hochreiter1997long}.\\

\begin{figure}
    \centering
    \includegraphics[width = 0.85\textwidth]{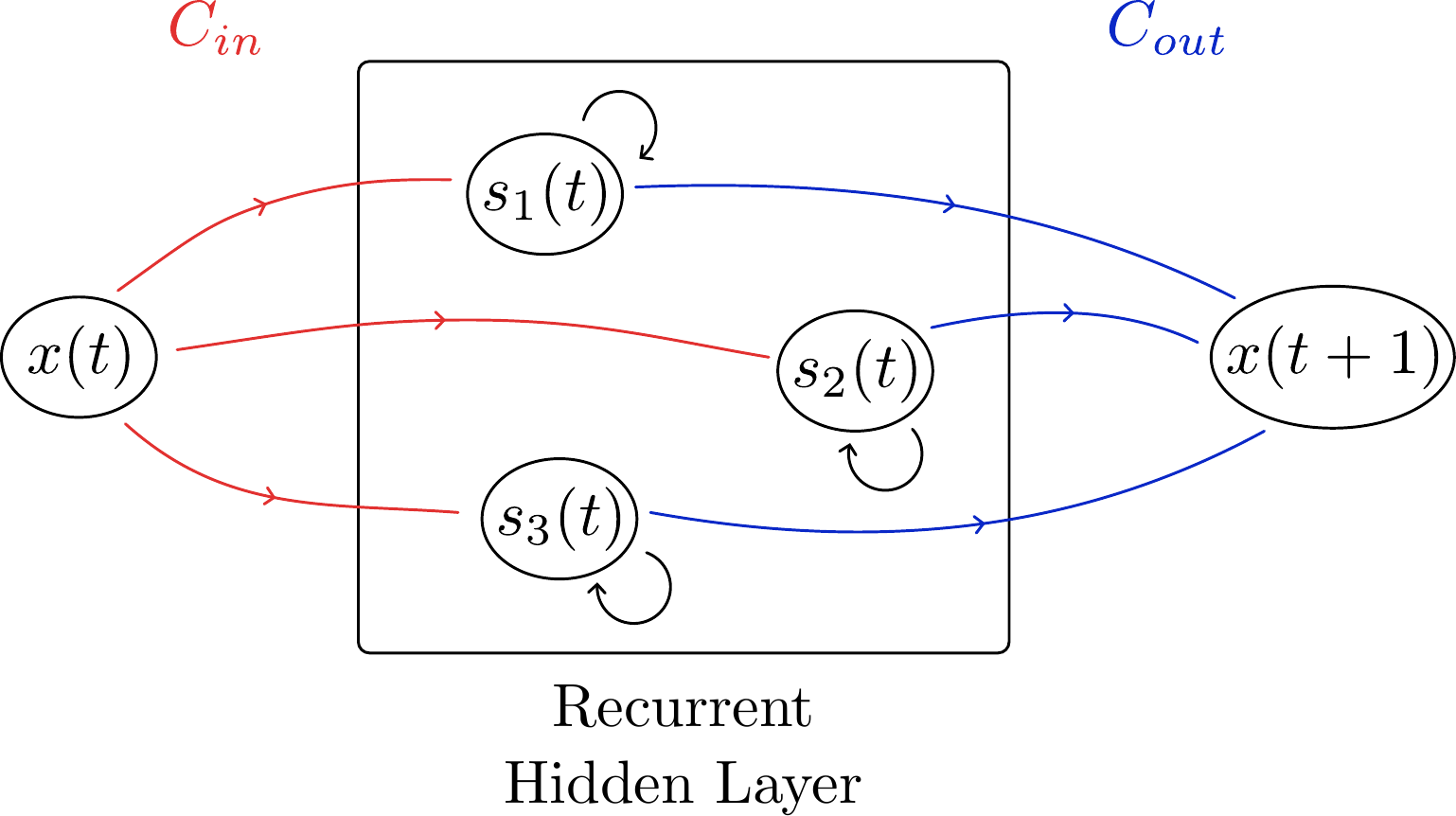}
    \caption{Structure of a single hidden layer recurrent neural network with three hidden nodes.}
    \label{fig:RNN}
\end{figure}

\noindent \textbf{Radial Basis Networks -} Similar formulation as a perceptron, but with different nonlinearity functions. Specifically, this class of models use radial basis functions 
\begin{equation}
    \phi^{RBN}_{j}(x(t)) = C \cdot (\sigma \circ ( \frac{|| x(t) - c_j||}{r_j} ))
\end{equation}

\noindent \textbf{Reservoir Computers - } A variation of RNNs that avoid the direct optimisation of hidden node weights and instead only trains the readout weights, reducing the training algorithm into a regularised least squares regression problem. Tailored more for time series analysis, sequence of time series are fed into the RNN and allowed to echo within the hidden layer to produce basis functions, which are subsequently used for future prediction. Essentially functions as an in build infinite dimensional delay embedding.  Reservoir computers essentially function as an infinite-dimensional delay embedding (see Chapter \ref{chap:reservoir_comp})

\section{Symbolic Regression (i.e. Giving neural networks to an algebraist/functional analyst)}

Neural networks are universal function approximators. This is achieved by approximating the desired function with a nested collection of nonlinearities of a predefined structure. Accuracy in reproducing the function is then dependent on the optimisation of constants. More generally, any feedforward neural network is essentially a function with the mathematical form
\begin{equation}
    \phi^{NN}(x) = \sigma(A...\sigma(\sigma(Ax+b)+b)...+b)
\end{equation}
Whilst useful for numerical analysis, this form suffers from several problems
\begin{itemize}
    \item \textbf{Interpretability - }The functional representation does not intuitively describe the underlying behaviour. Think replacing an exponential with a Taylor expansion.
    \item \textbf{High complexity - }There are large numbers of model parameters that require training, resulting in very high model complexity. A majority of which are not used and the model is not parsimonious.
\end{itemize}
\vspace{1em}
One alternative to neural networks is symbolic regression. This approach tries to analytically reconstruct the original function $f$ without imposing a pre-defined set mathematical form. This is the done with the aim of hopefully producing a more mathematically interpretable and parsimonious solution. The sparse identification of nonlinear dynamics (SINDy) algorithm is one current algorithm that aims to do this directly from data \cite{brunton2016discovering}. The algorithm is as follows:
\begin{enumerate}
    \item Let there be an $n$-dimensional multivariate time series given by $\mathbf{x}(t) = (x_{1}(t), x_{2}(t),...,x_{n}(t))$ measured at $m+1$ time steps. Construct two matrices $\mathbf{X},\mathbf{\dot{X}} \in \mathbb{R}^{m\times n}$ with elements containing the observed states and numerical derivatives respectively.

    \begin{subequations}
        \begin{align}
            \mathbf{X} = \begin{bmatrix} \mathbf{x}^{T}(t_{1}) \\ \mathbf{x}^{T}(t_{2})  \\ \vdots  \\ \mathbf{x}^{T}(t_{m}) \end{bmatrix} = 
            \begin{bmatrix} x_{1}(t_{1}) & x_{2}(t_{1}) & \dots & x_{n}(t_{1}) \\ x_{1}(t_{2}) & x_{2}(t_{2}) & \dots & x_{n}(t_{2}) \\ \vdots & \vdots & \ddots & \vdots \\ x_{1}(t_{m}) & x_{2}(t_{m}) & \dots & x_{n}(t_{m}) \end{bmatrix}\\ 
            \dot{\mathbf{X}} = \begin{bmatrix} \dot{\mathbf{x}}^{T}(t_{1}) \\ \dot{\mathbf{x}}^{T}(t_{2})  \\ \vdots  \\ \dot{\mathbf{x}}^{T}(t_{m}) \end{bmatrix} = 
            \begin{bmatrix} \dot{x}_{1}(t_{1}) & \dot{x}_{2}(t_{1}) & \dots & \dot{x}_{n}(t_{1}) \\ \dot{x}_{1}(t_{2}) & \dot{x}_{2}(t_{2}) & \dots & \dot{x}_{n}(t_{2}) \\ \vdots & \vdots & \ddots & \vdots \\ \dot{x}_{1}(t_{m}) & \dot{x}_{2}(t_{m}) & \dots & \dot{x}_{n}(t_{m}) \end{bmatrix}
        \end{align}
    \end{subequations}
    
    \item Let $F = \{f_{1}, f_{2},..., f_{K} \}$ be a library of candidate nonlinear functions (e.g. 1, $x$, $x^2$, $\sin x$ ...). Construct a matrix $\Theta (\mathbf{X}) \in \mathbb{R}^{m \times nK}$ consisting of element evaluations of the basis functions $F$ on the values of $\mathbf{X}$.

    \item Let $\Xi = (\xi_{1}, \xi_{2},...,\xi_{n})$ be a set of sparse vectors of coefficients that determine which nonlinearities are active:

    $$ \mathbf{\dot{X}}= \Theta (\mathbf{X})\Xi$$

    Perform sparse optimisation for each column in the above equation to identify the require coefficients. The result is an approximate analytical expression for the vector field equations.
\end{enumerate}

\vspace{1em}

A more recent alternative to SINDy is the Equation Learner Network (EQL) \cite{martius2016extrapolation} that modifies the basic feedforward neural network to consist of base operations and attempts to construct the desired function $f$ from the ground up. In general however, the field of symbolic regression is rapidly evolving and interested reader may consider further research in the area. \\

For the interested reader, there is the approach of physics-informed machine learning recently developed for solving PDEs. This approach generally refers to the case where a structured physics model is combined with a neural network to solve predefined PDE problem. The physics informed neural network (PINN) algorithm \cite{raissi2019physics} consists of an architecture in two halves: a feedforward neural network for the unknown function $u$, and a partial derivative PDE component that is evaluated numerically. These halves are used in conjunction to calculate a loss that can be subsequently backpropagated to update the neural network weights.

\chapter{Model Selection}

In comparison to other tasks in time series analysis, models are relatively easy to construct and can be optimised. However, model selection is not so straightforward and poses a bigger challenge. To illustrate this, consider a model used to optimise the transmission of a text message of a given length $N$. We may consider two possible approaches for sending the message:

\begin{enumerate}
    \item Directly encode the message verbatim and send it along. This approach has high precision and accuracy, but is computationally expensive and does not scale well.

    \item Encode the general semantics of the message (i.e. the simple model) alongwith sufficient accompanying detail (i.e. the error distribution) for the receiver to reconstruct the message. This approach has lower precision and is prone to errors, but is much cheaper to send and is more parsimonious.
\end{enumerate}

\vspace{1em}

The above illustration describes one key aspect of model selection: balancing between complexity and accuracy. High complexity models (although direct message passing is hardly a model) are powerful and flexible, but also may be prone to error, difficult to interpret, computationally expensive and may not generalise well. On the other hand, low complexity models may be more generalisable (owing to their simplicity), cheap to calibrate and run, but may not be as accurate or flexible. So how does one go about selecting models?\\

Model selection may be done with the use of several information criteria. These criteria are similar to a cost function that balances the weight of model complexity and accuracy. This reduces the task of model selection to an optimisation problem where the information criterion score needs to be minimised. Here, we consider three different information criteria with increasing levels of selectiveness. These are the Akaike Information Criterion (AIC), Bayes (Schwarz) Information Criterion (BIC), and the model description length (MDL).

\section{Akaike Information Criterion (AIC)}

Suppose we have a collection of observed data points that are generated from an underlying true probability distribution $g(y)$. In the modelling problem, we would like to construct a parametric model $f(y|\theta)$ with model parameters $\theta$ that approximates $g(y)$. Specifically, we assume that we have already selected the class of models each with $k$ parameters and distributions given by
\begin{equation}
    \mathcal{F}(k) = \{ f(y|\theta_k, \, \theta_k \in \Theta (k) \},
\end{equation}
where $\Theta (k)$ is the collection of all possible $k$ dimensional parameter vectors.\\

Intuitively, good models $f(y|\theta_k)$ should correspond to a closer matching with the real distribution. Given this, we would like to assess the goodness of fit for our model.  Our objective is to search amongst the collection of all fitted models that provides the best approximation for $g(y)$. To do this, we will first need to define a method for quantifying similarity.

A \textbf{divergence} is a function commonly used to describe the similarity between two distributions. The usage of the term ``divergence'' is not to be confused with its other definition  $\nabla \cdot\mathbf{F}$. Divergences are very similar to metrics, but with looser conditions. Most notably, a divergence does not need to be symmetric (i.e. $\mathcal{D}(f,g) = \mathcal{D}(g,f)$).

For our purposes, we consider a divergence derived from information theory -- the Kullback-Leibler (KL) divergence

\begin{theorem} \textbf{-- Kullback-Leibler divergence}\\
    $$\mathcal{D}_{KL}(g||f) = \mathbb{E} \left[ \log \left( \frac{g(y)}{f(y | \theta_k)} \right) \right]$$
    Furthermore,
    $$ \mathcal{D}_{KL}(g||f) = 0 \iff f = g$$
\end{theorem}
The KL divergence $\mathcal{D}_{KL}(g||f)$ measures the amount of information required to encode $g(y)$ using a code optimised for encoding $f$. Whilst $\mathcal{D}_{KL}$ is not a formal metric, its value grows to reflect the disparity between $f$ and $g$. Therefore, $\mathcal{D}_{KL}(g||f)$ measures how well the model distribution $f$ approximates the true distribution $g$. To proceed further, we define
\begin{equation}
    d(\theta_k) = \mathbb{E} [ -2\log f(y| \theta_k) ].
\end{equation}
Using this, it follows that
\begin{align}
    2 \mathcal{D}_{KL}(g||f) &= \mathbb{E} \left[ -2\log \left( \frac{g(y)}{f(y| \theta_k)} \right) \right] \notag\\
    &= \mathbb{E}[ -2 \log g(y) + 2 \log f(y | \theta_k) ]\\
    &= d(\theta_k) - \mathbb{E} [ -2 \log g(y) ]  \notag.
\end{align}
As we have no control over $g(y)$, we aim to perform all our optimisations with resepct to $\theta_k$. Ranking models based on $\mathcal{D}_{KL}$ would be equivalent to ranking based onn just $d(\theta_k)$. Therefore, it should be possible to rank a model with parameters $\hat{\theta}_k$ based on the evaluation of the value:
\begin{equation}
    d(\theta_k) = \mathbb{E} [ -2\log f(y| \hat{\theta}_k) ].
\end{equation}
Unfortunately, it is not possible to estimate $d(\hat{\theta}_k/0$ directly as it depends on observations $y$, which subsequently depends on the original $g(y)$. (Akaike, 1974) proposed an biased estimator by using the term $-2 \log f(y | \hat{\theta}_k)$  where $y$ are the observed data so far. In this formulation $-2 \log f(y | \hat{\theta}_k)$ is a negatively biased estimator. This can be accounted for by using a bias adjustment $\delta>0$ given by
\begin{align}
    \delta &= \mathbb{E} [2d(\hat{\theta}_k)] - \mathbb{E} [-2 \log f(y | \hat{\theta}_k]\\
    &\approx k. \notag
\end{align}
Combining the above components yields an expression for the Akaike Information Criterion (AIC)
\begin{definition} \textbf{-- Akaike Information Criterion (AIC)}\\
    \begin{align*}
        \text{AIC} &= -2 \log f(y | \hat{\theta}_k) + 2k\\
        &= 2k - 2\ln \hat{\mathcal{L}}
    \end{align*}
    where $\hat{L}= \max ( \mathcal{L} ( \hat{\theta}_k | y ))$ is the maximised value of the likelihood function.  Detailed derivations are provided in \cite{cavanaugh2019akaike}.
\end{definition}

In simple terms, the AIC score penalises the goodness of fit of a model by the number of required parameters. The $2k$ term ensures that models are not overfit with large numbers of parameters.

\section{Bayesian Information Criterion (BIC)}

The BIC was proposed by Schwarz as a Bayesian approach to tackling the model selection problem \cite{schwarz1978estimating}. Again, suppose that there is a true underlying model distribution $g(y)$, and we aim to select a model $M_k$ from a group of parametric class of distributions. In other words, given a model class $M_K$, there exists a $k$-dimensional parameterisation set of probabilty density functions
\begin{equation}
    \mathcal{F}(k) = \{ f(y|\theta_k | \theta_k \in \Theta (k) \},
\end{equation}
Once we select a model class $M_k$, we can calculate the likelihood function $\mathcal{L}(\theta_k |y)$ conditioned on observed data $y$. Let $\hat{\theta}_k$ denote a vector of parameter estimates obtained by maximising $\mathcal{L} (\theta_k |y)$ over $\Theta (k)$ with respect to $M_k$.

Suppose we can formulate a collection of models $M_{k_1}, M_{k_2}, \dots, M_{k_L}$ (e.g. different form, statistics, covariance structure etc). We ask the question: which $M_k$ corresponds to the best approximation of $g(y)$?

As is common in all Bayesian approaches, we first begin with a statement of Bayes theorem.

\begin{definition} \textbf{-- Bayes' Theorem}\\
    Let $A$ and $B$ be to random events, and the $P(A)$ and $P(B)$ are the probability of each event occurring respectively. Then,
    $$P(A|B) = \frac{P(B|A) P(A)}{P(B)}$$

    $P(A)$ is called the prior distribution of $A$, and $P(B)$ is the marginal distribution of $B$.
\end{definition}

Let $\pi(k)$ be the discrete prior distribution of model classes $M_{k_1}, M_{k_2}, \dots, M_{k_L}$, and $g (\theta_k | k)$ denote the prior on $\theta_k$ given the model $M_k$. Applying Bayes theorem, we get:

\begin{equation}
    h( (k, \theta_k) | y) = \frac{\mathcal{L}( \theta_k (y) g (\theta_k |k) \pi (k) )}{m(y)}
\end{equation}
where $m(y)$ is the marginal distribution of observations $y$. A Bayesian formulation of model selection tries to maximise the posterior likelihood $h( (k, \theta_k) | y)$ with respect to $k$. We may also write this as a posterior probability,
\begin{equation}
    P(k|y) = \frac{1}{m(y)} \pi (k) \int_{\Theta_k} \mathcal{L} [\theta_k|y] g(\theta_k |k) d\theta_k
\end{equation}
Alternatively, this is equivalent to minimising the log posterior,
\begin{equation}
    -2\ln P(k|y) = 2 \ln m(y) -2\ln \pi(k) - 2\ln \left[ \int \mathcal{L} [\theta_k|y] g(\theta_k |k) d\theta_k \right]
\end{equation}
As we are only interested in model selection, we can rewrite the above expression in a simpler form by discarding the first term

\begin{equation}
    -2\ln P(k|y) = -2\ln \pi(k) - 2\ln \left[ \int \mathcal{L} [\theta_k|y] g(\theta_k |k) d\theta_k \right]
\end{equation}

This expression can be further simplified and terms can be discarded to yield the following simplified expression for the BIC.

\begin{definition} \textbf{-- Bayesian (Schwarz) Information Criterion (BIC)}\\
    \begin{align*}
        \text{BIC} &= k \ln (n) -2 \ln \mathcal{L} [\theta_k|y]\\
        &= k \ln (n) -2 \hat{\mathcal{L}}
    \end{align*}
    where $\hat{L}= \max ( \mathcal{L} ( \hat{\theta}_k | y ))$ is the maximised value of the likelihood function. Detailed derivations are provided in \cite{neath2012bayesian}.
\end{definition}

Observe that there are similarities between BIC and AIC. The two scores only differ by their penalty term. BIC applies a penalty on the number of parameters scaled by the number of observations. Generally, BIC imposes stricted considerations for model selection. In the Bayesian framework, one expects that in the presence of more observed data, it should be possible to construct a better model. Thus, the BIC uses the quality and quantity of observed data as a regulariser when performing model selection.

\section{Model Selection and Coding Theory}

\subsection{The message passing problem revisited}

Suppose we have a time series $x_t$, $t=1,2,3,\ldots N$ ($x=(x_1,\ldots,x_N)$). By Takens' Theorem, for sufficiently large $d$ we can \textit{reconstruct} the underlying dynamics via a time delay embedding,
\begin{equation}
    v_t=(x_t,x_{t-1},\ldots x_{t-d+1})
\end{equation}
where the dynamic evolution of $v_{t+1} = \phi(v_t)$ is equivalent to the evolution operator of the underlying dynamical system. We want to estimate $\phi$. We can do this via a model $f:\mathbb{R}^d \rightarrow \mathbb{R}$ of the form
\begin{equation}
    \label{model}
    x_{t+1} \approx  f(x_t,x_{t-1},\ldots x_{t-d+1};\lambda)
\end{equation}
where $\lambda\in\mathbb{R}^k$ parameterises the set of possible models.  But how do we choose the best model?  For example, if one were to restrict to a class of models, say radial basis models, there remains the question of how many centres should be included and to what precision? We require an objective principle for choosing between alternate models. We examine this problem using ideas from \textit{information theory}.

Consider a message passing problem. We wish to transmit a data stream between two parties.  For example, the stream $1, 1, 2, 3, 5, 8, 13, 21, 34, 55, \dots$.  There are two ways of achieving this:
\begin{enumerate}
    \item Simply send the raw data verbatim    
    \item Send a model of the data, initial conditions and the model's \textit{prediction} errors so that the recipient can reconstruct the data stream.
\end{enumerate}

\vspace{1em}
The data stream in the example above can be recognized as the first few terms of the Fibonacci sequence. To recreate the stream we can use a \textbf{Perfect Model}, i.e.  $x_{1}$ and $x_{2}$ and the rule $x_{t}=x_{t-1}+x_{t-2}, ~~t\ge 3$ are sent.  We should also send the number of terms so that the recipient knows when to stop.

A \textbf{Perfect Model} may not always be available -- it is unknown, we made an error transcribing the model and so on -- but it may still be useful to send an \textbf{Imperfect Model} with errors rather than sending the raw data stream itself.  For example, rather than the perfect model of the Fibonacci sequence we instead send the rule $x_{t} = 2x_{t-1}, ~~t\ge 2$ with $x_{1}=1$ and also send the errors $e_{t}$ so that the data stream is $\hat{x}_{t} = x_{t} + e_{t}$.  If the model is ``good'' (small model and small errors) then the cost of sending the model plus the errors may still be less than simply sending the raw data. This is especially true if one limits the range of predictions expected of the recipient (e.g. limit to the first 100 terms).

We choose to define \textit{best} to mean the model that offers the most compact description of the data -- up to some finite specified precision. That is, if we were to transmit the original data down a finite-capacity communication channel, what model would offer the most compact description of the data? We are thus presented with two choices:
\begin{enumerate}
    \item \textbf{Perfect model - }the model captures the dynamics exactly.  Solve for $\lambda$, then send $\lambda$ and the initial conditions $x_0$
    \item \textbf{Imperfect model - }estimate $\lambda$, and optimise a trade-off
    \begin{itemize}
        \item The model is bad: it is quicker (cheaper) to transmit the data $x$
        \item The model is good: sending the model parameters $\lambda$, plus $x_0$ and the model prediction errors, is better.
    \end{itemize}
\end{enumerate} 

Therefore, we should aim to construct models such that
\begin{equation}
    L(model) + L(errors) < L(data),
\end{equation}
where $L$ is the cost of sending numbers. How does one quantify $L$? To do so, we draw on ideas from \textbf{Coding Theory}.\\

\subsection{Coding Theory}

Let $A$ represent a \textbf{finite alphabet} of symbols. Communication on a \textbf{digital} channel requires the representation of this alphabet via a separate set of symbols $B$. Let $B$ denote the set of code symbols permissible on the channel, and denote by $B^*$ the set of words over that alphabet.  A word $(b_1,b_2,\ldots,b_\ell)$ of $\ell$ symbols from $B$ is an ordered sequence of members of $B$.  $B^*$ contains all such permissible sequences.  Hence, we have a map ${\cal L}:A\rightarrow B^*$ with $a\mapsto (b_1,b_2,\ldots,b_\ell)$ with $b_i\in B$.\\

Define the {\em code length} $L:A\rightarrow \mathbb{Z}^+$ such that $a\mapsto \ell$, i.e. $L(a)=\ell$ is the length of the encoding of $A$ with the alphabet $B$.  Ideally, we want codes that are both {\em unambiguous} and {\em efficient}.\\

\begin{example} \textbf{-- Morse code}\\
    In Morse Code, $|B|=3$ (dot, dash and space).  Morse code is unambiguous---spaces separate codewords composed of dots and dashes and there is a one-to-one correspondence between $B^*$ and $A$.  Morse code is also efficient -- more frequently occurring symbols in $A$ have shorter codewords in $B^*$.  That is, if $L(a_1)<L(a_2)$ then $a_1$ is more common than $a_2$.   Consider a binary encoding of Morse code, i.e., let $B=\{0,1\}$ and $A=\{{\rm dot},{\rm dash},{\rm space}\}$. Consider the four possible codings in Table~\ref{tab:morse}.
    
\end{example}

\begin{table}
    \centering
    \begin{tabular}{*{10}{r|ccc|l}}
    & space & dot & dash &\\
    & (\textbf{\textvisiblespace}) & ($\bullet$) & (\textbf{--}) &\\
    1. & $0$ & $1$ & $11$ & ambiguous\\
    2. & $00$ & $01$ & $11$ & unambiguous (but incomplete)\\
    3. & $0$ & $10$ & $11$ & unambiguous and more efficient (complete)\\
    4. & $00$ & $1$ & $01$ & unambiguous and still \\
    &&&& more efficient (w.r.t. frequency)
    \end{tabular}
    \caption{Codings of Morse code}
    \label{tab:morse}
\end{table}

\subsection{Coding Trees}
Coding trees are a useful tool to help to ascertain if a coding is ambiguous or unambiguous.  For an unambiguous coding, codewords used are the ``leaves'' of the tree.  That is, no codeword is a prefix subset of any other codeword.  This is also called a \textit{prefix-free} coding.  A \textit{complete} code has no unused leaves.  Consider each of the binary encodings of Morse code symbols $\{\mbox{dot}, \mbox{dash}, \mbox{space}\}=A$ given in Table~\ref{tab:morse}, where $B=\{0, 1\}$.

\begin{enumerate}
    \item
    \begin{tabular}{cc}
    \begin{tabular}{ccc}
    space & dot & dash \\
    0 & 1 & 11
    \end{tabular}
    &
    \includegraphics[scale=0.40]{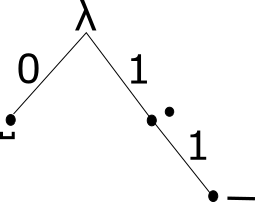}
    \end{tabular}\\
    The problem with this encoding is ambiguity.  If one transmits the message $11$, it is ambuguous if this refers to a code of ``$\bullet\, \bullet$'' or ``\textbf{--}''
    \item
    \begin{tabular}{cc}
    \begin{tabular}{ccc}
    space & dot & dash \\
    00 & 01 & 11
    \end{tabular}
    &
    \includegraphics[scale=0.40]{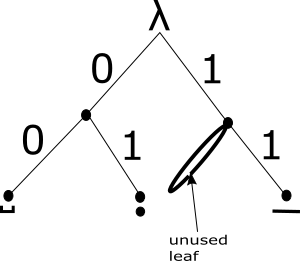}
    \end{tabular}\\
    This encoding is unambiguous but it is incomplete, as there is an unused leaf for the code 10.
    \item
    \begin{tabular}{cc}
    \begin{tabular}{ccc}
    space & dot & dash \\
    0 & 10 & 11
    \end{tabular}
    &
    \includegraphics[scale=0.40]{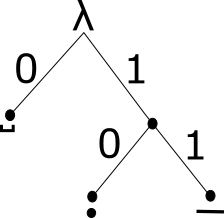}
    \end{tabular}\\
    This encoding is unambiguous, complete (no unused leaves) but it is not fully effficient as it does not account for the relative occurrence frequency of each symbol. Ideally, more frequent symbols should utilise shorter encodings.
    \item
    \begin{tabular}{cc}
    \begin{tabular}{ccc}
    space & dot & dash \\
    00 & 1 & 01
    \end{tabular}
    &
    \includegraphics[scale=0.40]{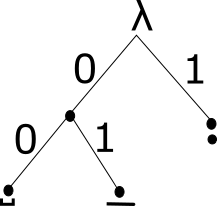}
    \end{tabular}\\
    This encoding is unambiguous, complete and it is also more efficient.
\end{enumerate}

\subsection{Optimal encodings}
We have used the term in brief, but what does it mean to have an \textbf{efficient encoding}?  We aim to construct a code such that short codewords are used for more frequently occurring (more probable) symbols. The length of the codewords should depend on the expected frequency of occurrence in such a way that the expected message length is minimal. Let 

\begin{equation}
    p(a)=\text{frequency/probability of encountering symbol }a\in A.
\end{equation}

The \textbf{efficiency} of the encoding can now be defined by the average (expected) codeword length

\begin{equation}
    \label{efficiency}
    \text{Efficiency} = \sum_{a\in A}p(a)L(a).
\end{equation}
We have the following result known as the \textbf{Kraft Inequality}.

\begin{corollary} \textbf{-- Kraft's Inequality}\\
    For a prefix-free code, $\sum_{a\in A}2^{-L(a)}\leq 1$, when $|B|=2$.  Furthermore, if the code is complete, then $\sum_{a\in A}2^{-L(a)}=1$. This result can be proved using a tree diagram and is left as an exercise for the reader.
\end{corollary}

Suppose $A=\{a_1,a_2,\dots, a_n\}$ and $p_i= \text{frequency  of }a_i$ ($\sum p_i=1$). Define $\ell_i=\lceil -\log_2 p_i\rceil$ (where $\lceil x\rceil$ is the smallest integer not less than $x$). Then,
\begin{align}
    \sum_i 2^{-\ell_i} & = \sum_i 2^{-\lceil-\log_2 p_i\rceil} \notag \\
    &\leq \sum_i 2^{\log_2 p_i} \notag \\
    &= \sum_i p_i = 1
\end{align}
 
The following result allows us to consider the definition for an \textbf{optimal encoding}.  
\begin{theorem} \textbf{-- Optimal encoding}\\
    Let $p=(p_1,p_2,\ldots,p_n)$ and $q=(q_1,\ldots q_n)$ where $0\leq p_i\leq1$, $0\leq q_i\leq 1$ and $\sum_{i=1}^n p_i=\sum_{i=1}^n q_i=1$. Then 
    $$-\sum_{i=1}^n p_i\log{p_i}\leq-\sum_{i=1}^n p_i\log{q_i}$$
    
    with equality iff $p_i=q_i$, $\forall i$.
\end{theorem}

The quantity $H(p)=-\sum_{i=1}^np_i\log_2 p_i$ is called the \textbf{entropy} of the information source.  The entropy  $H(p)$ is a lower bound on the average codeword length:

\begin{equation}
    \sum_{a\in A}p(a)L(a)\geq H(p)
\end{equation}
and, if the $p_i$ are all powers of $2$, then $\lceil-\log_2 p_i\rceil$ are integers and we can achieve $H(p)$.  Thus, given an alphabet $A$ and frequencies $p=(p_1,p_2,\ldots)$, then there exists a prefix-free code over a binary alphabet $B=\{0,1\}$, such that 
\begin{equation}
    \sum_{a\in A}p(a)L(a)-H(p) \leq 1
\end{equation}
That is, $H(p)\leq \sum_{a\in A}p(a)L(a) \leq H(p)+1$.

\noindent Recall Kraft's inequality, if a code is uniquely decodable its lengths must satisfy
\begin{equation}
    \sum_{i}2^{-l_{i}} \le 1.
\end{equation}
Therefore for any lengths satisfying the Kraft Inequality there exists a prefix-free code with these lengths.\\

\subsubsection{Code lengths for ensembles}
Optimal source code lengths for an ensemble are equal to the Shannon information content
\begin{equation}
    l_{i} = \log_{2}\frac{1}{p_{i}},
\end{equation}
Conversely, any choice of codelengths defines implicit probabilities,
\begin{equation}
    q_{i} = \frac{2^{-l_{i}}}{Z}
\end{equation}
where $Z$ is normalized to be one if the code is complete.\\

\subsubsection{Source coding theorem for symbol codes}
For an ensemble $X$ there exists a prefix-free code whose expected length satisfies
\begin{equation}
    H(X) \le L(C,X) < H(X) + 1,
\end{equation}
where
\begin{equation}
    H(X) = \sum_{i}p_{i}\log\frac{1}{p_{i}}~~\mbox{-- entropy}
\end{equation}
and
\begin{equation}
    L(C,X) = \sum_{i}p_{i}l_{i}~~\mbox{-- expected length of a code}
\end{equation}

\begin{corollary}\textbf{-- Upper bound of data code length}\\
    Let data $D=\{\alpha_1,\alpha_2,\ldots \alpha_N\}$ where $\alpha_i\in A$ and assume symbols in $D$ occur with probability $p$.  Let $L(D)$ denote the total number of bits required to send this data. Then,

    \begin{align}
        L(D) &= \sum_{i=1}^N L(\alpha_i) \notag \\
        &= \sum_{i=1}^N\lceil-\log_2 p(\alpha_i)\rceil \notag \\
        &\approx N \sum_{a\in A}p(a)\lceil-\log_2 p(a)\rceil \notag \\
        &\leq N(H(p)+1)
    \end{align}
\end{corollary}

\subsection{Huffman Coding}
How, then, do we construct optimal codes?  A \textbf{Huffman Coding} is guaranteed to achieve the optimum coding.  Suppose we write $(w,p)$ for the (codeword, frequency) pair for a symbol $a \in A$. Then the following algorithm achieves an optimum code in the sense that the average codeword length is minimal.
\begin{enumerate}
    \item Find the two pairs $(w,p)$ and $(w',p')$ with the smallest probabilities (frequencies)---break ties arbitrarily.
    \item Replace these two pairs with one pair $(\hat{w},\hat{p})$ where $\hat{p}=p+p'$, $w=\hat{w}_0$ and $w'=\hat{w}_1$.
    \item Repeat, until all that is left is one pair $(w,p)=(\lambda,1)$.
    \item Back-substitute.
\end{enumerate}

\vspace{1em}
The codewords will not be uniquely determined -- we break ties arbitrarily -- but their lengths are. 

\begin{example} \textbf{-- Huffman coding}\\
    Compute the Huffman coding for the following set of codeword/frequency pairings: $\{(a,0.1), (b,0.2), (c,0.4), (d,0.3)\}$.  
\end{example}

\begin{table}
    \centering
    \begin{tabular}{cc}
    Symbol & Freq./Prob. \\
    \hline
    a & 0.1 \\
    b & 0.2 \\
    c & 0.4 \\
    d & 0.3 \\
    A & 0.3 \\
    B & 0.6 \\
    \hline
    \end{tabular}
    \caption{Huffman coding (non-unique)}
    \label{tab:huffman}
\end{table}

\begin{figure}[ht!]
    %\begin{tcolorbox}[colback=white]
    \centering
    \includegraphics[scale=0.95]{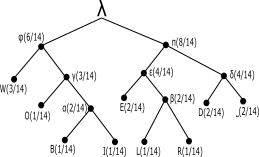}
    \caption{Coding tree for a Huffman coding of WWW.}
    \label{fig:www}
    %\end{tcolorbox}
\end{figure}

\begin{example} \textbf{-- Huffman coding 2 }\\
    Compute the Huffman coding for the string: WORLD\textvisiblespace WIDE\textvisiblespace WEB\\

    Including the space character there are $14$ symbols.  A table of the symbols, their frequenices, entropies, Huffman code and code length are given in Table~\ref{tab:huffmanwww}.  A particular coding tree is shown in Figure~\ref{fig:www}.
\end{example} 

%%%%%%%%%%%%%%%%%%
%%%%%%%%%%%%%%%%%%%

\begin{table}
    \centering
    \begin{tabular}{ccccc}
    Symbol $a_{i}$ & $p_{i}$ & $h_{i}$ & Code $C(a_{i})$ & $l_{i}$ \\
    \hline
    W & 3/14 & 2.2 & 00 & 2 \\
    O & 1/14 & 3.8 & 010 & 3 \\
    R & 1/14 & 3.8 & 1011 & 4 \\
    L & 1/14 & 3.8 & 1010 & 4 \\
    D & 2/14 & 2.8 & 110 & 3 \\
    \textvisiblespace & 2/14 & 2.8 & 111 & 3 \\
    I & 1/14 & 3.8 & 0111 & 4 \\
    E & 2/14 & 2.8 & 100 & 3 \\
    B & 1/14 & 3.8 & 0110 & 4 \\
    $\alpha$ & 2/14 & & & \\
    $\beta$ & 2/14 & & & \\
    $\gamma$ & 3/14 & & & \\
    $\delta$ & 4/14 & & & \\
    $\varepsilon$ & 4/14 & & & \\
    $\phi$ & 6/14 & & & \\
    $\pi$ & 8/14 & & & \\
    \hline
    \end{tabular}
    \caption{Huffman coding of the WWW.}
    \label{tab:huffmanwww}
\end{table}

\textbf{Disadvantages of Huffman coding}
\begin{itemize}
    \item depends on the ensemble (the ``text'')
    \item ignores context (q followed by the letter u)
    \item potential extra-bit overhead
\end{itemize}

\vspace{1em}

The scheme described above provides a method to optimally encode a given (single) number -- from $\mathbb{N}$, $\mathbb{Z}$ or $\mathbb{Q}$; exactly how to extend the scheme to each of these cases is left as an exercise for the reader.  Although Huffman codes are optimal symbol codes. But for practical purposes symbol codes may not always be ideal. This is because the encoding does not include a method for indicating the end of a sequence of code.  We can either introduce a new symbol (``comma'') to delimit the distinction between codewords (ala. ``space'' in Morse code); or, design an encoding which is self-delimiting.  That is, 
we build int a method to indicate when a codeword ends an a new one begins. In essence, we communicate the length of the integer $l_{b}(n)$ then communicate the original integer.  We must also communicate the length of the length of the integer and so on.\\

\subsection{Self-delimiting Decoders}
Consider a sequence of binary bits.
\begin{enumerate}
    \item If the first digit is a zero, then the encoded integer is $1$. Stop.
    \item Otherwise, consider the first three digits (the first of which is $1$) as the current code-block. If the third digit is $0$, then the first two digits encode the desired integer (i.e. $10=2$, $11=3$), and stop.
    \item Otherwise consider, the first two digits encode the length of the next block to be read---if $10$ then let $w=3$, if $11$ then let $w=4$.
    \item Discard the current code-block and read the next $w$ bits as the next (new) code-block. 
    \item \label{repeat} If the $w+1$-th bit (i.e. the next bit after the current code-block is $0$, then the current $w$-bit code block encodes the desired number in standard binary---stop.
    \item Otherwise, the current $w$-bit code block encodes the length $u$ of the next block to be read. Discard the current code block, the first $w$ bits, read the next $u$ bits as the (new) current code-block, let $w=u$ and go back to Step. \ref{repeat}.
\end{enumerate}

\begin{example} \textbf{-- Self-delimiting codes of integers}\\
    Realizations of the binary self-delimiting code for some integers are shown in Table~\ref{tab:selfdelimit}.
\end{example}

%%%%%%%%%%%%
\begin{table}[h]
    \centering
    \begin{tabular}{cl}
    integer & codeword\\
    1 & 0 \\
    2 & 10\ 0\\
    3 & 11\ 0\\
    4 & 10\ 100\ 0\\
    7 & 10\ 111\ 0\\
     14 & 11\ 1110\ 0\\
     15 & 11\ 1111\ 0\\
     16 & 10\ 101\ 10000\ 0
     \end{tabular}
     \caption{Realizations of a binary self-delimiting code for integers.}
     \label{tab:selfdelimit}
\end{table}
%%%%%%%%%%%%%%%%

\begin{example} \textbf{-- Self-delimiting codes of integers 2}\\
    Encode last year $2017$ using a binary self-delimiting code.

    \begin{align*}
        2017_{2} & = & 11111100001 \mbox{---11 bits} \\
        11_{2} & = & 1011 \mbox{---4 bits} \\
        4_{2} & = & 100 \mbox{---3 bits} \\
        3_{2} & = & 11 \mbox{---2 bits}
    \end{align*}

    Thus, send $10~100~1011~11111100001~0$.  Recall, step 3.\\

    \textbf{How to decode: }Check the first three digits.  If the third digit is a one then the first two digits denote the number of digits to be read next, say $r_{n}$.  Check the $r_{n}+1$ digit.  If it is zero then the $r_{n}$ digits is the number else it's the number of digits to be read next.

\end{example}

\textbf{Elias Encoder}\\
An operation encoder is given by \textbf{Elias's encoder} (taken from \cite{mackay2003information})
\begin{itemize}
    \item[] Write '0'
    \item[] Loop $\{$
    \begin{itemize}
        \item[] If $\lfloor \log n\rfloor = 0$ halt.
        \item[] Prepend $C_{2}(n)$ to the written string
        \item[] $n:=\lfloor\log n\rfloor$
    \end{itemize}
    \item[] $\}$
\end{itemize}

\vspace{1em}
One can calculate the number of bits required for a self-delimiting code to encode $n\in\mathbb{N}$ as,
\begin{equation}
    L^{*}(n) = \lceil \log_2 n\rceil +\lceil \log_2 \lceil \log_2 n\rceil \rceil + \lceil \log_2 \lceil \log_2 \lceil \log_2 n\rceil\rceil\rceil+ \ldots +1 
\end{equation}
For encoding integers $n\in\mathbb{Z}$, we define $n' = 2n, n\ge 0$ and $n=-(2n-1), n<0$ then for $n\in\mathbb{Z}$, $L^{*}(2n)$. That is, {$n\in\mathbb{Z}$} if $n>0$ then $n\mapsto 2n$, otherwise $n\mapsto -2n+1$ and hence
\begin{equation}
    L^*(n) = \lceil \log_2 2n\rceil +\lceil \log_2 \lceil \log_2 2n\rceil \rceil + \lceil \log_2 \lceil \log_2 \lceil \log_2 2n\rceil\rceil\rceil+ \ldots +1 
\end{equation}
To extend the encoding to rational floating point numbers $\lambda \in \mathbb{Q}$, we first consider the following representation,
\begin{equation}
    \lambda=a\times 2^b
\end{equation}
where $a\in[\frac 12,1)$ and $b\in\mathbb{Z}$. Due to limited computational memory, a truncated approximation to $\lambda$ can be written as,
\begin{equation}
    \overline\lambda=\pm 0.a_1a_2a_3a_4\ldots a_n\times 2^{\pm m}
\end{equation}
where we are only ever interested in finite-bit representations ($n<\infty$) and two extra bits will be required to encode the sign of the mantissa and the exponent.  Let $\delta$ denote the {\em precision} of our finite representation $\overline\lambda$ of $\lambda$
\begin{equation}
    \delta>|\lambda-\overline\lambda|
\end{equation}
and hence truncating (round-off) at $n$-bits gives $\delta=2^{-n}$ ($\delta=2^{-(n+1)}$, respectively).

Therefore, if one chooses to encode a real number, this can be achieved by instead encoding a truncated $\bar{\lambda}$ via an $n$-bit floating point with a require code length of
\begin{align}
    L^*(\lambda)&= L^*(\pm a_1a_2a_3a_4\ldots a_n)+L^*(2m) \notag \\
    &= L^*(2\times 2^{(n-1)})+ L^*(2m) \notag \\
    &= L^*(\lceil\delta^{-1}\rceil)+L^*(2m)
\end{align}

\section{Minimum Model Description Length (MDL)}
We wish to evaluate the performance (``goodness'') of our model by computing the cost of using that model to describe (encode) the data.
\begin{enumerate}
    \item Huffman coding for large amounts of data with \textbf{known} distribution (ensemble dependent)---we can use this for model prediction errors.
    \item Self-delimiting code for small amounts of data with \textbf{no specific distribution}, $\mathbb{Z}^{+}, \mathbb{Z}$ or $\mathbb{R}$ with specified precision---we can use this for model parameters.
\end{enumerate}
\vspace{1em}

Recall, we have data $x_{1}, x_{2}, \dots, x_{n} \in \mathbb{R}$ to finite precision and a model
\begin{equation}
    x_{t+1} = f(x_{t}, x_{t-1}, \dots, x_{t-d};\lambda) + e_{t+1}, ~~\lambda\in\mathbb{R}^{k}
\end{equation}
with errors $e_{t}$ that have a known (or assumed) distribution $P_{e}$.\\

\textbf{N.B.} The $e_{t}$ are also of finite precision and hopefully on average less than $x_{t}$ to encode.  We also need initial conditions $x_{-(d-1)}, x_{-(d-2)}, \dots, x_{0}$ and W.L.O.G. we include these with $\lambda$.  So, altogether we have some parameter vector, $\theta$ which includes $\lambda$, initial conditions, \textbf{and} precisions $\delta$.  Both $e$ and $x$ are finite precision vectors, hopefully, $L(e)\ll L(x)$. Therefore, we require a two-part code: parameters first, then $e_{t}$.  The parameters define $P_{e}$ and $f$ and allow decoding of $e_{t}$.\\

So, let $\theta$ denote the complete parameter vector ($\lambda$, initial conditions and everything---and from hereon in we let $|\theta|=k$ and the $j$\textsuperscript{th} component of $\theta$ we denote by $\lambda_{j}$ and say that it has precision $\delta_{j}$). Then we can evaluate the model by sending a two part code: first $\theta$, then $e$ with total code length
\begin{align}
    L(x,\theta)&= L(x|\theta)+L(\theta)\notag\\
    \nonumber&= L(e)+L(\theta)\notag\\
    \nonumber&\leq n(H(P_e)+1)+L(\theta)\\
     &\approx  n(H(P_e)+1)+ \sum_{j=1}^k \left(L^*(\lceil\delta_j^{-1}\rceil)+L^*(2m_j)\right)\notag
\end{align}
where $m_j$ is the exponent of the corresponding $\overline\lambda_j$. What we do in practice, however, is
\begin{align}
    L(x,\theta) & = L(\theta) + L(x|\theta) \\
    & \approx L(\theta) + L(x|\hat{\lambda}) + \frac{1}{2}\delta^{T}Q\delta \\
    &\approx \sum_{j=1}^{|\theta|}(L^{*}(2m_{j})+L^{*}(\lceil\delta^{-1}_{j}\rceil)) + L(x|\hat{\lambda}) + \frac{1}{2}\delta^{T}Q\delta
\end{align}

So, for some chosen model class (e.g., radial basis, neural networks, AR process...) the general model description length approach for model selection is as follows: 
\begin{enumerate}
    \item Find $\hat{\lambda}$ that minimizes $H(P_{e})$ but really log-likelihood of an assumed error distribution. Typically, we assume additive Gaussian for the model errors
    \item Given $\hat{\lambda}$ choose optimal values of $\delta$
    \item Use the total code length to select best model (the one with lowest total code length)
\end{enumerate}

\noindent \textbf{Roadmap:}\\
It is convenient to first make some assumptions and approximations.
\begin{itemize}
    \item Bound the exponents of $\overline\lambda_{j}$ by a constant, say $\gamma (=32)$.
    \item Approximate $L^{*}$ by the leading term so that with $\gamma$ the code length of the parameters becomes

    $$L(\theta) \approx \sum_{j=1}^{k}\log(\frac{\gamma}{\delta_{j}})$$
    
    \item Find $\hat\lambda$ that minimizes $L(e)$ which is $-\log(P_{e})$
\end{itemize}

\begin{figure}
    %\begin{tcolorbox}[colback=white]
    \centering
    \includegraphics[width=0.8\textwidth]{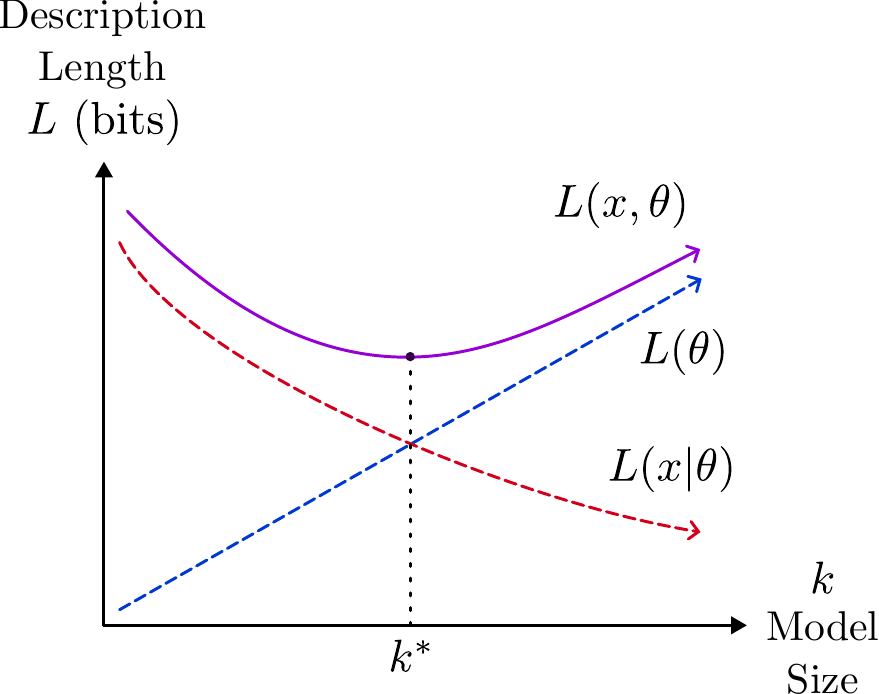}
    \caption{Typical profile of the total code length of a generalized linear model with respect to model size $k$.  There is a trade-off between the cost of the model prediction errors and the cost of the model parameters.}
    \label{fig:glm}
    %\end{tcolorbox}
\end{figure}

\noindent \textbf{Roadmap: continued}
\begin{itemize}
    \item Given $\hat\lambda$ choose optimal $\delta$
    \item Compute ($L(e)+L(\theta)$) and use this to select the best model (i.e. optimal $k$) of the general form (\ref{model}).
    \item \textbf{N.B.} $P_{e}$ depends on $\sigma^{2}$ which we must also find a solution to and send along with its precision so,
    $$L(\theta) \approx \sum_{j=0}^{k}\log(\frac{\gamma}{\delta_{j}})$$

    with $\delta_{0}$ the precision of $\hat{\sigma}^{2}$.
    \item We will also work in nats rather than bits (i.e., use $\log$ instead of $\log_{2}$).
\end{itemize}
Consider a generalized linear model of the form
\begin{equation}
    f(x_{t},x_{t-1},\dots,x_{t-d+1};\lambda) = \sum_{j=1}^{k}\lambda_{j}\phi_{j}(x_{t}, x_{t-1}, \dots, x_{t-d+1})
\end{equation}
where $k$ non-zero, $\lambda_{j}$ at some precision $\delta_{j}$.   For a linear model
\begin{equation}
f(x_t,x_{t-1},\dots x_{t-d+1};\lambda) =\sum_{j=1}^d\lambda_j x_{t-j+1}.
\end{equation}
A typical profile for $L(x,\theta)$ with respect to model size $k$ is sketched in Figure~\ref{fig:glm}. 

%%%%%%%%%%%%%%%
%%%%%%%%%%%%%%%%

Approximately, $L(\theta)\approx\sum_{j=1}^{k}L^{*}(\lceil\delta^{-1}\rceil)+\mbox{constant}$.  We bound the mantissa $m_{j}$ by $\gamma$ which will require its own precision.  To find the precisions, i.e., $\delta$'s, consider
\begin{equation}
    L(x) \approx L(x|\hat{\lambda}) + \frac{1}{2}\delta^{T}Q\delta + O(\delta^{3}) + \sum_{j=1}^{k}L^{*}(\lceil\delta_{j}^{-1}\rceil) + \sum_{j=1}^{k}\log\gamma .
\end{equation}

We work in nats (natural log bits) and really have $\sum_{j=0}^{k}\log\frac{\gamma}{\delta_{j}}$ since $j=0$ corresponds to $\sigma^{2}$ and its precision $\eta$, see below.  We want to minimize $L(x)$ with respect to $\delta_{i}$.  Therefore,

\begin{equation}
    \frac{\partial}{\partial\delta_{i}}L(x) = 0
\end{equation}
which implies
\begin{equation}
    \frac{1}{\delta_{i}} = [Q\delta]_{i}
\end{equation}
and $[\cdot]_i$  denotes the $i^{}th$ component of the vector.  We solve this for the $\delta$'s and so need to find $Q$.  To do this we consider a specific model.  In particular, for a generalized linear model, specify the basis functions.  Let $x^{(t)}=x_{t-1},\ldots x_{t-d}$ then:
\begin{itemize}
    \item $\phi_j = \phi(a_jx^{(t)}-b_j)$ yields a feedforward neural network;
    \item $\phi_j = \phi(\frac{\|x^{(t)}-c_j\|}{r_j})$ is a radial basis network.
\end{itemize}

Let $V$ be the $(n-d)\times k$ matrix where the columns and rows are the evalution of each basis function on each data point respectively.  Then the model can be expressed in matrix notation as
\begin{equation}
    e = x - V\lambda
\end{equation}
Assume (for Huffman coding) the errors $e$ have a Gaussian distribution $e\sim N(0,\sigma^{2})$ so that
\begin{equation}
    P_{e_{i}} = \frac{1}{\sqrt{2\pi\sigma^{2}}}e^{-\frac{e_{i}^{2}}{2\sigma^{2}}}
\end{equation}
So,
\begin{align}
    L(x,\theta) & = L(x|\theta) + L(\theta) \notag\\
    & = -\sum_{i=1}^{n}\log P_{e_{i}} + L(\theta) \notag\\
    & \approx \frac{n}{2}\log 2\pi\sigma^{2} + \frac{1}{2\sigma^{2}}(x-V\lambda)^{T}(x-V\lambda) + L(\theta)
\end{align}
We seek the optimal value of $\lambda=\hat{\lambda}$ and $\sigma^{2}=\hat{\sigma}^{2}$.  Thus, solve $D_{\lambda}L(x,\theta)=0$ and $D_{\sigma^{2}}L(x,\theta)=0$ simultaneously.\\

\noindent For $D_{\lambda}L=0$:
\begin{equation}
    \big{[}-\frac{1}{\sigma^{2}}V^{T}(x-V\lambda)\big{]}=0
\end{equation}
implies $\hat{\lambda} = (V^{T}V)^{-1}V^{T}x$, the usual ``least squares'' solution.\\

\noindent For $D_{\sigma^{2}}L=0$:
\begin{equation}
    \big{[}\frac{n}{2}\frac{1}{\sigma^{2}}-\frac{1}{2(\sigma^{2})^{2}}(x-V\lambda)^{T}(x-V\lambda)\big{]}=0
\end{equation}
implies $\hat{\sigma}^{2}=\big{[}(x-V\lambda)^{T}(x-V\lambda)\frac{1}{n}\big{]}=\frac{e^{T}e}{n}$.  Now,
\begin{equation}
    Q = \left [\begin{array}{cc} D_{\lambda\lambda}L(x,\theta)|_{\hat{\lambda},\hat{\sigma}^{2}} & D_{\lambda\sigma^{2}}L(x,\theta)|_{\hat{\lambda},\hat{\sigma}^{2}} \\
    D_{\sigma^{2}\lambda}L(x,\theta)|_{\hat{\lambda},\hat{\sigma}^{2}} & D_{\sigma^{2}\sigma^{2}}L(x,\theta)|_{\hat{\lambda},\hat{\sigma}^{2}}
    \end{array}\right ]
\end{equation}

\noindent We find $D_{\lambda\lambda}L|_{\hat{\lambda},\hat{\sigma}^{2}} = \frac{1}{\hat{\sigma}^{2}}V^{T}V$, $D_{\lambda\sigma^{2}}L|_{\hat{\lambda},\hat{\sigma}^{2}}=0$ and
\begin{align}
    D_{\hat{\sigma}^{2}\hat{\sigma}^{2}}L|_{\hat{\lambda},\hat{\sigma}^{2}} & = -\frac{n}{2}\frac{1}{(\sigma^{2})^{2}}+\frac{2}{2(\sigma^{2})^{3}}e^{T}e|_{\hat{\lambda},\hat{\sigma}^{2}} \notag\\
    & = -\frac{n}{2}\frac{1}{(\hat{\sigma}^{2})^{2}}+\frac{n}{(\hat{\sigma}^{2})^{2}} \notag\\
    & = \frac{n}{2}\frac{1}{(\hat{\sigma}^{2})^{2}}
\end{align}

\noindent Therefore,
\begin{equation}
    Q = \left [\begin{array}{cc}
    \frac{1}{\hat{\sigma}^{2}}V^{T}V & 0 \\
    0 & \frac{n}{2}\frac{1}{(\hat{\sigma}^{2})^{2}}
    \end{array}\right ]
\end{equation}
The precisions $\delta$ and $\eta$ can thus be obtained from
\begin{equation}
    \big{[}\frac{1}{\hat{\sigma}^{2}}V^{T}V\delta\big{]}_{j} = \frac{1}{\delta_{j}}
\end{equation}
and
\begin{equation}
    \frac{n}{2}\frac{1}{(\hat{\sigma}^{2})^{2}}\eta = \frac{1}{\eta},
\end{equation}
The solution of which is given by $\eta=\sqrt{\frac{2}{n}}\hat{\sigma}^{2}$. This is used to provide the following expression for model description length (MDL):
\begin{align}
    L(x,\theta) & \approx \frac{n}{2}\log(2\pi\hat{\sigma}^{2}) + \frac{1}{2\hat{\sigma}^{2}}e^{T}e + \frac{1}{2}\delta^{T}Q\delta + \sum_{j=0}^{k}\log(\frac{\gamma}{\delta_{j}}) \notag \\
     & \approx \frac{n}{2}\log(2\pi\hat{\sigma}^{2}) + \frac{1}{2\hat{\sigma}^{2}}e^{T}e +\frac{1}{2}(k+1)+(k+1)\log\gamma  \notag \\
     & - \sum_{j=1}^{k}\log\hat{\delta}_{j} - \log\eta \notag \\
     &=\frac{n}{2}\log 2\pi + \frac{n}{2}\log\frac{e^{T}e}{n} + \frac{n}{2} + (k+1)(\frac{1}{2}+\log\gamma)  \notag \\
     & - \sum_{j=1}^{k}\log\hat{\delta}_{j} - \log\sqrt{\frac{2}{n}}\hat{\sigma}^{2}  \notag \\
     & =  (\frac{n}{2}-1)\log\frac{e^{T}e}{n} + (k+1)(\frac{1}{2}+\log\gamma) - \sum_{j=1}^{k}\log\hat{\delta}_{j}  \notag \\
     & + \frac{n}{2}(1+\log 2\pi)+\frac{1}{2}\log\frac{n}{2}
\end{align}

\begin{definition} \textbf{-- Model Description Length (MDL)}\\
    For a model with Gaussian errors, the model description length is given as
    $$DL(k) =(\frac{n}{2}-1)\log\frac{e^{T}e}{n} + (k+1)(\frac{1}{2}+\log\gamma) - \sum_{j=1}^{k}\log\hat{\delta}_{j} \notag \\
    + \frac{n}{2}(1+\log 2\pi)+\frac{1}{2}\log\frac{n}{2}$$
\end{definition}

\textbf{N.B.} For linear models BIC (SIC) is an asymptotic approximation to the description length \cite{rissanen1998stochastic,judd1992improved}.\\

Recall the rationale.  We have a time series $y_{t}$ which we assume comes from some measurement of the state of a dynamical system evolving on some low-dimensional attractor.  We embed using the method of time delay embedding to reconstruct an equivalent state space with states $x_{t}$, the dynamics of which are in some sense equivalent to the unknown ``true'' dynamics. So, quantities invariant under coordinate transformations calculated directly from the (embedded) data or indirectly from models of the evolution help to characterize the original system.\\

If we construct a global model, how do we select a useful model where
\begin{align}
    y_{t+1} & = f(x_{t}) \notag \\
    & = \sum_{i=1}^{k}\lambda_{i}\phi_{i}(x_{t}) + e_{t}
\end{align}
where $\phi_{i}$ are basis functions so that $e=y-V\lambda$?  Note, the model is linear in the (unknown) parameters $\lambda$.  We send the estimated parameters, precisions and prediction errors.  To compare models we use the model description length (MDL)
\begin{equation}
    DL_{k}(\theta) = (\frac{1}{2}n-1)\ln\frac{e^{T}e}{n} + (k+1)(\frac{1}{2}+\ln\gamma) - \sum_{j=1}^{k}\ln\hat{\delta}_{j}+\frac{1}{2}n(1+\ln 2\pi)+\frac{1}{2}\ln(\frac{1}{2}n)
\end{equation}
where $n$ is the number of data, $\gamma$ is a constant, $\hat{\delta}$ solves $[Q\delta]_{j}=\frac{1}{\delta_{j}}$ and $\theta$ represents the data, parameters, precisions etc.\\

\subsection{Model Selection with MDL}

\textbf{General Plan}\\
The model with minimum description length is chosen as the ``best'' or most useful model. Thus we have the following roadmap for model selection.
\begin{enumerate}
    \item Select the model class (i.e. the basis functions in $y_{t}=\sum_{i=1}^{k}\lambda_{i}\phi_{i}(x_{t})$ so that $y=V\lambda$)
    \item Find $\hat\lambda$ that minimizes $e = y - V\lambda$.
    \item Given $\hat\lambda$ choose optimal $\delta$ (i.e., solve $[Q\delta]_{i}=\frac{1}{\delta_{i}}$).
    \item Compute $DL_{k}$ and use this to select the best model (i.e. optimal $k$ such that $DL_{k}$ is a minimum).
    \item Alternatively, or in addition, compute $AIC(k)$ or $SIC(k)$ and select the best model (i.e. minimum $AIC(k)$ or $SIC(k)$).
\end{enumerate}

\noindent\textbf{Subset selection}\\
Model selection essentially reduces to solving the optimisation problem
\begin{center}
minimize $|y - V\lambda |$ subject to $N(\lambda) = k$.
\end{center}
where $N$ is the number of parameters. This is difficult.  To tackle this, we start with a fixed $k$  with Gaussian i.i.d. noise assumptions for the errors and solve the corresponding optimization problem,
\begin{center}
    \begin{tabular}{ll}
        minimize & $\frac{1}{2}e^{T}e$ \\
        over & $e$, $\lambda$ \\
        subject to & $V\lambda + e = y$, $N(\lambda)=k$.
    \end{tabular}
\end{center}
In this form, $N(\lambda)=k$ is the complicating factor.  The idea is to start with a model of size $k$ and then determine from within a set of candidate basis functions the ``best'' function to add to the model to increase its size to $k+1$.  Then determine which of the $k+1$ basis functions is ``best'' culled from the model to return the size to $k$ basis functions.  This is done in an iterative fashion to find the ``optimal'' model, where optimal means the model with minimum description length.\\

\noindent\textbf{Growing the model}\\
If we denote by $B$ the set of basic (model terms) indices, namely,
\begin{equation}
    B = \{ j:\lambda_{j}\ne 0\},
\end{equation}
so $N(\lambda) = |B|$ then we can investigate via sensitivity analysis the effect of changing the size of $B$.  It is convenient to use auxillary variables and apply Lagrangian optimization theory.  We can express the constraint $N(\lambda)=k$ in terms of non-basic variables as
\begin{equation}
    \lambda_{j} = u_{j},~~j\notin B,
\end{equation}
where the parameters $u$ are equal to zero.  The Lagrangian of the optimization problem is
\begin{equation}
    L = \frac{1}{2}e^{T}e + \nu^{T}(y-V\lambda-e) + \mu^{T}(u-\lambda),
\end{equation}
where $\mu$ and $\nu$ are dual variables.  Optimizing over $e$ and $\lambda$ without constraints ($\frac{\partial L}{\partial e} = 0 = \frac{\partial L}{\partial \lambda}$) gives
\begin{equation}
    \nu = e
\end{equation}
and
\begin{equation}
    \mu = -V^{T}e.
\end{equation}

Now, $\mu$ is the dual variable corresponding to the non-basic constraint and so is the sensitivity to changes in $u$ at optimality.  Thus the largest element of $\mu$ in absolute value should be added to the model, i.e., made basic, to give the greatest marginal payoff.\\

\noindent\textbf{Culling the model}\\
We want to find the basic variable to remove from the model which does the least damage to the mean square error.  We consider the dual optimization problem by sustituting the optimal values of $\lambda$ and $e$ into the Lagrangian, i.e.,
\begin{center}
    \begin{tabular}{ll}
        maximize & $-\frac{1}{2}\nu^{T}\nu + \nu^{T}y + \mu^{T}u$ \\
        over & $\mu$, $\nu$ \\
        subject to & $(V^{T}\nu)_{j} = w_{j}$, for $j\in B$.
    \end{tabular}
\end{center}
In a manner similar to before, $w=0$ but is kept as a parameter but we set $u=0$ immediately.  Note, the constraint involving $w$ affects the basis variables and so the sensitivity to $w$ tells about the variables already in the basis.  Since the problem is convex the primal and dual solutions are the same.  Hence, the dual vairable to $\mu$ is $\lambda$ and 
\begin{equation}
    \lambda_{j} = \frac{\partial \psi}{\partial w_{j}},
\end{equation}
where $\psi$ is the optimal solution.  Thus, select the smallest element of $\lambda_{j}$ in absolute value to make non-basic, i.e., cull from the model, to do the least marginal damage to the payoff.

This motivates the following algorithm for determining and assessing a ``good'' model.  Let $V_{B}$ be the n x k matrix formed from the columns of $V$ with indices in $B$.  Let $\lambda_{B}$ be the least squares solution to $y=V_{B}\lambda$ and let $e_{B}=y - V_{B}\lambda_{B}$.
\begin{enumerate}
    \item Let $S_{0}=(\frac{1}{2}n-1)\ln(\frac{y^{T}y}{n})+\frac{1}{2}+\ln\gamma$.
    \item Let $B=\{j\}$ where $V_{j}$ is the column of $V$ such that $|V_{j}^{T}y|$ is maximum.
    \item Let $\mu=V^{T}e_{B}$ and $i$ be the index of the component of $\mu$ with maximum absolute value.  (Index $i$ is the index coming \textit{in} to the basis.)  Let $B'=B\cup\{i\}$.
    \item Calculate $\lambda_{B'}$.  Let $o$ be the index in $B'$ corresponding to the component of $\lambda_{B'}$ with smallest absolute value.  (Index $o$ is the index going \textit{out} of the basis.)
    \item If $i\ne o$ then put $B=B'\setminus\{o\}$ and go to step 3).
    \item Define $B_{k}=B$ where $k=|B|$.  Find $\delta$ such that $(V_{B}^{T}V_{B}\delta)_{j}=\frac{1}{\delta_{j}}$ for each $j=\{1, 2, \dots, k\}$ and calculate $S_{k}=(\frac{1}{2}n-1)\ln (\hat{e}^{T}\hat{e}/n)+(k+1)(\frac{1}{2}+\ln\gamma)-\sum\limits_{j=1}^{k}\ln\hat\delta_{j}$.
    \item If $S_{k}<S_{k-1}$ then go to step 3).
    \item Take basis $B_{k}$ with $S_{k}$ a minimum as the optimal model.
\end{enumerate}

\chapter{Surrogates}

In the study of time series and dynamical systems, we often solve the associated forward problem. That is, given a set of equations (ODE), solve the ODE to generate the relevant time series. In reality, the inverse problem is more common, and also more difficult. Given an observed time series, what can we learn about the dynamical process. There are numerous challenges that are encountered when tackling the inverse problem:
\begin{itemize}
    \item Learning a global model may be useful for prediction, but is not necessarily informative
    \item Noise and nonlinearity can affect performance of linear autocorrelation.
    \item Estimating invariant measures (e.g. correlation dimension, Lyapunov exponent) is prone to errors and misinterpretation unless large amounts of good quality data is available.
\end{itemize}
We can also consider the following related questions as well:
\begin{itemize}
    \item[\textbf{Q1.}] What is the cause of the variability (dynamics) in a time series?
    \item[\textbf{Q2.}] How can we come up with simple ways to test and reject "simple" explanations
\end{itemize}

\vspace{1em}
Why is this interesting? Having knowledge of the underlying form of the data generating process can guide the selection of appropriate analysis methods to be used. For example, if dynamics are nonlinear, we would expect autocorrelation to not be useful. If it is stochastic, then Lyapunov exponents are not informative. Based on these, we aim to find out which models are best suited for analyses. One way to tackle this is by constructing surrogates for hypothesis testing. This is done by addressing two simpler questions:
\begin{enumerate}
    \item[\textbf{Q3.}] Is the process nonlinear?
    \item[\textbf{Q4.}] If so, what class of nonlinear model is warranted, and can produce results consistent with the data?
\end{enumerate}

\vspace{1em}
Hypothesis testing in statistics is the idea that one can compute some statistical value from observed data, compare that value to some theoretical distribution and conclude that an underlying null hypothesis is either rejected or that we fail to reject said hypothesis. Surrogate hypothesis testing follows the general procedure:
\begin{enumerate}
    \item Define the null hypothesis $H_{0}$ regarding the type of data.
    \item Define a data generating process that is consistent with $H_{0}$ (i.e. the surrogate process), but preserves other aspects consistent with the observed data (e.g. mean, Fourier spectrum etc.)
    \item Generate an ensemble of surrogate data $s_{i}(t)$
    \item Define a test statistic, $d$ (Lyapunov exponent, correlation dimension) and calculate the distribution of $d(s_{i}(t))\sim \mathcal{D}$
    \item Calculate $d_{o}=d(x(t))$, the test statistic of the original observed data.
    \item Does $d_{0}$ adhere to the distribution $\mathcal{D}$?
        \begin{itemize}
            \item \textbf{Yes} $\implies$ $H_{0}$ cannot be rejected. Either the data is prescribed by the surrogate process, or $d()$ was not an adequate choice to discriminate between the two processes
            \item \textbf{N}o $\implies$ Reject $H_{0}$ and conclude the $x(t)$ is a result of a more ``complicated'' process.
        \end{itemize}
\end{enumerate}

\vspace{1em}
The surrogate hypothesis testing procedure requires two things:
\begin{enumerate}
    \item A null hypothesis with a corresponding algorithm with which to generate surrogate time series.
    \item Discriminating statistics for decision making via $p$-value testing.
\end{enumerate}

\section{The Hierarchy of Surrogate Algorithms}
Basic surrogate testing involves a collection of three main hypothesis, with each containing their own surrogate generation algorithm. For simplicity, we name these Algoirithm 0, 1 and 2. If a time series rejects the null hypotheses of all algorithms, we may conclude that the observations are the result of a more complicated processes.

\subsection{Algorithm 0: Random Shuffle Surrogates}
\textbf{Hypothesis: }The observed data is the product of an i.i.d noise process. 

To test this hypothesis, surrogates are generated to preserve the probability distribution of the data, but destroy any temporal structure. This can be achieved by shuffling the time-order of the data. The steps of algorithm 0 are given as follows:
\begin{enumerate}
    \item Let $\{ x_{t}\}_{t=1}^{N}$ be the original observed time series
    \item Let $\{ s_{t}\}_{t=1}^{N}$ be the newly generated time series
    \item $s_t$ is calculated by shuffling the time ordering of ${x_t}$ (sampling without replacement
    \begin{enumerate}
        \item $u=\text{randn}(N,1)$
        \item $idu = \text{sort}(u)$
	\item $s_{t}= x\left[ idu \right]$
    \end{enumerate}
\end{enumerate}
\vspace{1em}
The result of the shuffling removes any autocorrelation in the data and hence $s_{t}\perp s_{t-\tau} \forall \,\tau$

\subsection{Algorithm 1: Fourier Transform Surrogates}

\textbf{Hypothesis: }The observed data is essentially linearly filtered Gaussian noise (equivalent to an AR process).\\

Given a scalar time series $x_{t}$, a linear filter is a transformation that produces another time series $y_{t}$ of the following form:
\begin{equation}
    y_{t} = \sum\limits_{j=1}^{N} \beta_{j} x_{t-j},
\end{equation}
where $\beta_{j}$ are constants. It is also possible to express a linear filter in the frequency domain as well:
\begin{equation}
    Y_{\omega} = \beta_{\omega}X_{\omega}
\end{equation}
There $\vec{\beta}$ (i.e. the collection of $\beta_{\omega}$) is known as the transfer function. \\

Algorithm 1 tries to produce a time series that is partially random \cite{theiler1992testing}. The underlying data values are randomly generated, but the linear applies smoothing process such that temporally close values are correlated. Therefore, all interesting features of the data are purely due to the autocorrelation function. Generally, the application of linear filters merely alters the power spectrum.

In practice, we aim to construct a surrogate time series $s_t$ that preserves the same autocorrelation function (i.e. the same power spectrum), but is random otherwise. This is achieved by shuffling the phases of the data, but keeping the amplitude structure:
\begin{enumerate}
    \item Let $\{ x_{t} \}_{t=1}^{N}$ be the original time series observed.
    \item Apply a discrete Fourier transform $X_{n} = \mathcal{F}(x_{t}) = \{ R_{n}e^{-i\pi \theta_{n}+\phi_{n}} \}$
    \item Randomise the complex phases:
    \begin{enumerate}
        \item Let $\psi_{n}\sim \mathcal{U}(0, 2\pi)$
	\item Multiply each amplitude by $e^{i\psi}$ to produce $S_{n}$. Make sure that $\psi(f) = -\psi(-f)$ so that the inverse Fourier transform is real valued
    \end{enumerate}
    \item Apply the inverse Fourier transform $\{ s_{t} \}_{t=1}^{N} = \mathcal{F}^{-1}(S_n)$
\end{enumerate}
\vspace{1em}
For a discriminating statistic, linear autocorrelation is insufficient as the algorithm intentionally preserves this quantity. Therefore, a nonlinear form of autocorrelation such as average mutual information is an appropriate statistic. There are several problems associated with Algorithm 0:
\begin{itemize}
    \item The numerical Fourier transform is imperfect as it does not fully preserve the autocorrelation structure. 
    \item Fourier transforms can also be problematic when working with aperiodic, or non matching ends.
\end{itemize}

\subsection{Algorithm 2: Amplitude Adjusted Fourier Surrogates}

\textbf{Hypothesis: }The observed data is described by linearly filtered noise, but then transformed again through a smooth,  observation function $h$, where $h: \mathbb{R} \to \mathbb{R}$ is a monotonic nonlinear function (i.e. $h(a) > h(b) \iff a > b$).\\

The null hypothesis for Algorithm 2 proposes that all observed nonlinearity is due to the observation function and not inherent dynamics. This is useful when the data is clearly not Gaussian, but one suspects the underlying cause to be a linear stochastic process \cite{theiler1996constrained}. The process is as follows:

\begin{enumerate}
    \item Let $\{ x_{t} \}_{t=1}^{N}$ be the observed data. Generate a sequence of $N$ Gaussian random numbers $\{ z_{t} \} \sim \mathcal{N}( 0 , 1)$ and reorder $z_{t}$ according to the same rank distribution as $x_{t}$ to copy the power spectrum.
    \item Apply the Fourier transform Algorithm 1 to $\{ z_{t} \}$ to produce a sequence $\{ \hat{s}_{t} \}$ that has the same power spectrum (autocorrelation structure).
    \item Reorder $\{ x_{t} \}$ according to the rank ordering $\{ \hat{s}_{t} \}$ to produce a surrogate $\{ s_{t} \}$ that has the same spectrum, autocorrelation, mean and standard deviation as the original observed data.
\end{enumerate}

In practice, the AAFT (Algorithm 2) does not perfectly preserve the autocorrelation and power spectrum. This is because the rank ordering (rescaling) process leads to an altered spectrum. An attempt to better preserve the power spectrum is by an iterative procedure of the above aptly named the Iterative Amplitude Adjusted Fourier Transform \cite{schreiber1996improved},
\begin{enumerate}
    \item Calculate the squared amplitudes of the Fourier transform $\{ x_{t}^{(0)} \}$ given by $\{ S_{k}^{2} \}$.
    \item Apply random shuffle Algorithm 0 to observed data $\{ x_{t} \}$ and produce $\{ \hat{s}_{t}^{(i)} \}$.
    \item Take the Fourier transform of $\{ \hat{s}_{t}^{(i)} \}$ and replace its amplitude with $\{ S_{k}^{2} \}$ to preserve the desired power spectrum.
    \item Take the inverse Fourier transform, keeping the phases of the sequence. This enforces the correct spectrum, but incorrect scalar distribution.
    \item Reorder $\{ x_t \}$ according to the untransformed $\{ \hat{s}_{t}^{(0)} \}$ to produce the surrogate $\{ \hat{s}_{t}^{(i+1)} \}$
    \item Check power spectrum and repeat if not close enough.
\end{enumerate}

A discriminating statistic suitable for these surrogates is correlation dimension.

\section{Other Surrogate Algorithms}

The standard battery of Algorithm 0, 1 and 2 surrogates are simple to implement and useful for determining whether a nonlinear analysis is warranted. However, there are other algorithms which are refinements to these three tests that aim more general (or specific) null hypothesis. Due to their large number and variations, we will only briefly discuss a small selection.

\subsection{Truncated Fourier Transform Surrogates}

The TFTS algorithm is a variation of the the standard Fourier transform surrogate (Algorithm 1). This algorithm aims to test for a certain class of non-stationary processes, namely, data with long-term trends (periodicity) but otherwise possess irregular fluctuations. Specifically, we want to investigate whether there is nonlinearity in these irregular fluctuations or if they attributed to stochastic processes. The idea is to ``destroy'' nonlinearity in the irregular fluctuations but preserve the ``global'' behaviours (i.e. trends and periodicity).

\begin{figure}[!ht]
    %\begin{tcolorbox}[colback=white]
    \centering
    \includegraphics[width = 0.7\textwidth]{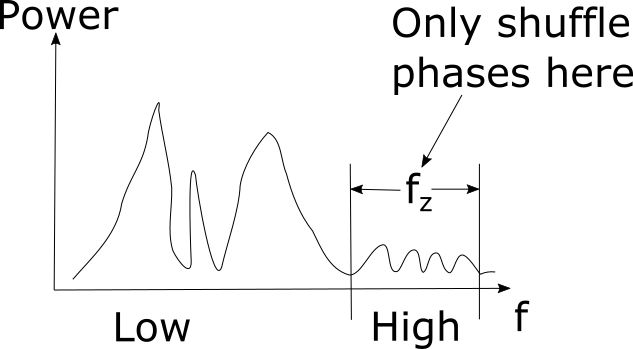}
    \caption{The idea underlying TFTS.}
    \label{fig:TFTS}
    %\end{tcolorbox}
\end{figure}

To achieve this, take the Fourier Transform of the data and examine the power spectrum.  Rather than shuffling all phases, only randomly shuffle phases at high frequency, i.e., those frequencies that lie in some frequency band $f_{z}$ (see Figure \ref{fig:TFTS}).  This parameter requires tuning.  The surrogates are the inverse Fourier transform of this truncated shuffled fourier transform.  Since only the high freqency components are shuffled the long-term trends (low periods) are preserved and only the irregular fluctuations are randomized.  Any nonlinearity within the short-term fluctuations should be destroyed.  A discriminating statistic for these surrogates can be average mutual information over a range of lags.

%%%%%%%%%%%%%%%%%%%%%%%%%%%%%%%
%%%%%%%%%%%%%%%%%%%%%%%%

\subsection{Small Shuffle Surrogates (SSS)}

The TFTS method is a modification of Algorithm 1 surrogate generation to investigate evidence for nonlinearity in irregular fluctuations of data with long-term trends.  It works by shuffling phases in a restricted band of high frequencies in the frequency domain.  A surrogate generation method that works in the time domain and can be considered an analogous modification of Algorithm 0 surrogate methods is the small shuffle surrogates (SSS) \cite{nakamura2005small}.  The null hypothesis for this test is that irregular fluctuations are iid random variables, i.e. no short-term dynamics. As mentioned the method is essentially Algorithm 0 but on a smaller scale.
\begin{enumerate}
    \item Start with the time series $x(t)$ and $i(t)$ be the index of $x(t)$ so that $i(t)=t$.
    \item Let $g(t)$ be a Gaussian random number then for $A>0$, $i'(t)=i(t)+Ag(t)$.
    \item Rank order $i'(t)$ so that the index of $i'(t)$ is $\hat{i}(t)$.
    \item Set $s(t)=x(\hat{i}(t))$.
\end{enumerate}

\subsection{Pseudo-periodic Surrogates (PPS)}

These are applicable to data with strong periodicities.  The null hypothesis tests against the data being consistent with a periodic orbit plus uncorrelated noise, hence Pseudo-Periodic Surrogates (PPS) \cite{small2002applying}.  Take the time series $\{x_{t}\}$ and embed using a time delay embedding with lag $\tau$ and dimension $d$ to get reconstructed states $\{z_{t}\}$.
\begin{enumerate}
    \item Randomly choose $s_{1}\in\{z_{t}\} (i=1)$.
    \item Select a neighbour $s_{j}$ randomly with probability $P(s_{j}=z_{t})\approx \exp (\frac{-\|z_{t}-s_{j}\|}{\rho})$.
    \item Set $s_{i+1}=s_{j+1}$, increment $i$.
    \item Repeat until $i>N$.
\end{enumerate}

\vspace{1em}
The method is a mixture of nearest neighbour prediction and the S-map predictor.  Now, choose as the scalar time series $y_{t}$ (the surrogate) the first component of the vectors $s_{t}$.  Observe, $s_{t}$ and $z_{t}$ are approximately the same ``attractor'' but an embedding of $\{y_{t}\}$ is only an approximation of the attractor.  The surrogates follow (approximately) the same vector field as the original data but fine scale dynamics are obliterated.  The parameter $\rho$ need tuning: for large $\rho$ we have effectively Algorithm 0 surrogates; for small $\rho$ there will be no variation (we keep selecting $s_{i}$ and so follow the original data.  Somewhere in the middle there is a $\rho$ where intercycle dynamics are destroyed but intracycle dynamics are preserved.  (See the paper for how to select a good $\rho$).  cf. with cycle shuffled surrogates.\\

\subsection{Attractor Network Surrogates}

The PPS algorithm generates surrogates by constructing trajectories that randomly switch between neighbouring orbits. This destroys the fractal structure of the trajectory, but approximately preserves the vector field. One can also consider a discretised form of the dynamics. Similar to symbolic dynamics, we can reduce flows in continuous phase space into discrete jumps between non-overlapping (not necessarily equally sized) regions on phase space. In the traditional symbolic dynamics discussed earlier in the unit, any set of dynamics can be expressed by an $n^{th}$ order Markov chain for infinite $n$. (e.g. for Lorenz we can consider the $L$ and $R$ wings as symbols). One can consider a refinement of this idea by increasing the number of discretised regions $N$, as $N\to \infty$, the dynamics may be arbitarily approximated by a $1^{st}$ order Markov chain (i.e. the memoryless property).

Attractor network surrogates (ANS) utilise the above principle by constructing a discrete Markov chain representation of system dynamics consisting of regions whose bounds are defined by a set of training data \cite{tan2023network}. Once constructed, one can perform random walks along the Markov chain to produce a surrogate time series. Like PPS, ANS broadly preserves the vector field (or in this case the phase space transition probabilities) of a target system, and hence replicating the ergodic measures (Lyapunov exponent and correlation dimension) of the observed data, but is stochastic rather than deterministic. This presents an even stronger surrogate test for assessing classification algorithms.

\section{Constraining Surrogates}

In the general Algorithms 0, 1 and 2, the aims are to preserve a given property of the data such as the randomness or linear properties. In this case, we may describe these as ``constrained realisations'' for the generated surrogate data because we surrogates preserve certain properties of the original data. However, one may also utilise \textbf{unconstrained realisations}. Further discussion warrants a proper definition of two concepts.

\begin{definition} \textbf{-- Constrained surrogates}\\
    Let $H_{0}$ be a null hypothesis for a given surrogate test and $\mathcal{P}$ be the set of all data generating processes that are consistent with $H_{0}$. Suppose all element processes of $\mathcal{P}$ can be parameterised by a set of parameters $\alpha$. Therefore, any data generating process consistent with $H_{0}$, $f \in \mathcal{P}$, can be uniquely defined by its associated parameters $\alpha(f)$. Let $x_{t}$ be observed data and $s_{t}$ be a set of generated surrogate data. A surrogate generation algorithm generates constrained realisations if the estimate of parameters $\alpha(x_{t})$ and $\alpha(s_{t})$ follow the same distribution, but are not necessarily equal. That is, there is no dependency on a fitted model. In contrast, unconstrained surrogates generated from $\alpha(x_{t})$. This can be achieved by fitting some model process (e.g. maximum likelihood) to the observe data, and using the the resulting model to generate new data.
\end{definition}

\begin{definition} \textbf{-- Pivotal statistics}\\
    A statistic is pivotal if the distribution of statistic values obtained from $f\in \mathcal{P}$ is independent of the selection of $f$. In other words, statistics of data sets generated by all the processes $f$ consistent with the null hypothesis $H_{0}$ follow the same distribution.
\end{definition}

 A pivotal statistic means that one does not have to ensure that the surrogate generation scheme is constrained when performing hypothesis testing. Otherwise, one does. Unbiased estimates of correlation dimension provide a pivotal statistic when the hypothesis being tested is the general class of dynamics consistent with the correlation dimension (i.e. all the nonlinear dynamics and chaos we have discussed so far). Hence, with the correlation dimension (and other dynamical invariants) one may use model based methods (which produce constrained surrogates) to test for membership of a broad class of dynamical systems.

\chapter{Reservoir Computing}\label{chap:reservoir_comp}

Lets recall our definition of the RNN from section~\ref{sec:rnn} together with the evolution function (with some slight modifications to the time indexing of the activation states),
\begin{equation}
    \vec{s}(t) = \sigma \circ (C_{in} \vec{x}(t)+C_{rec}\vec{s}(t-1) +\vec{b}),
\end{equation}
and next step prediction
\begin{equation}
    \vec{x}(t+1) = C_{out}\vec{s}(t)\,.
\end{equation}
Under a typical RNN setup, each of the parameters contained in $C_{in}\in\mathbb{R}^{n \times k}$, $C_{rec}\in\mathbb{R}^{k\times k}$, $\vec{b}\in\mathbb{R}^{k}$ and $C_{out}\in\mathbb{R}^{m\times k}$ require tuning, which for even a modest network size $k$ can quickly become computationally intensive.
Instead, what if we can leverage our understanding of dynamical systems and design our network in some sophisticated way initially in order to reduce the computational requirements later on?

This was the focus of the work by Herbert Jaeger in 2001~\cite{jaeger2001echo}.
Jaeger observed that under certain parameter setups, the states of a RNN can be seen as an ``echo'' of the input history.
Under these conditions there are a number of interesting properties that emerge that allow us to simplify the training required from the reservoir computer.
While often the development of reservoir computing is attributed to Jaeger and his echo state networks (ESN), Mass et al.~\cite{maass2002real} introduced the idea independently in 2002 with their liquid state machines (LSM).

\section{The Echo State Property}

\begin{definition} \textbf{-- Echo States}\\
    Consider a left-infinite input sequence $\vec{x}(-\infty),\dots,\vec{x}(T-1),\vec{x}(T)$ with resulting network states $\vec{s}(-\infty),\dots,\vec{s}(T-1),\vec{s}(T)$ and $\vec{s}\,'(-\infty),\dots,\vec{s}\,'(T-1),\vec{s}\,'(T)$ where $\vec{s}(t) = f_{RC}(\vec{s}(t-1),\vec{x}(t))$ and $\vec{s}\,'(t) = f_{RC}(\vec{s}\,'(t-1),\vec{x}(t))$.
    A network is said to possess \textit{echo states} if the states are uniquely determined by the input sequence $\vec{x}$ (ie. it holds that $\vec{s}(T)=\vec{s}\,'(T)$).
\end{definition}

The ability of reservoir computers hinges upon the echo state property, although more recently the term \textit{consistency} property has begun to see use~\cite{lymburn2019consistency}.
The idea is that, in networks which possess echo states, the initialisation of the network does not matter and each state is a function of the inputs alone.
In practice, checking for echo states involves assessing three characteristics of the underlying system.

\begin{theorem} \textbf{-- Practical Echo State Property}\\
    A network is said to possess \textit{echo states} if the network:
    \begin{itemize}
        \item is uniformly state contracting,
        \item is state forgetting, and
        \item possesses fading memory.
    \end{itemize}
\end{theorem}

\begin{definition} \textbf{-- Uniformly State Contracting}\\
    A network is called \textit{uniformly state contracting} if there exists a null sequence $\delta(h)$ (ie. a sequence that tends to $0$) such that for all right-inﬁnite input sequences $\vec{x}(T),\vec{x}(T+1),\dots,\vec{x}(\infty)$ with states $\vec{s},\vec{s}\,'$ it holds that $||f_{RC}(\vec{s}, \bar{x}(h) - f_{RC}(\vec{s}\,', \bar{x}(h))||_2 < \delta(h)$ for all input sequence prefixes $\bar{x}(h) = \vec{x}(T),\vec{x}(T+1),\dots,\vec{x}(h)$ with $h\geq0$.
\end{definition}

\begin{definition} \textbf{-- State Forgetting}\\
    A network is called \textit{State Forgetting} if there exists a null sequence $\delta(h)$ such that for all left-inﬁnite input sequences $\vec{x}(-\infty),\dots,\vec{x}(T-1),\vec{x}(T)$ with states $\vec{s},\vec{s}\,'$ it holds that $||f_{RC}(\vec{s}, \bar{x}(h) - f_{RC}(\vec{s}\,', \bar{x}(h))||_2 < \delta(h)$ for all input sequence suffixes $\bar{x}(h) = \vec{x}(T-h),\dots,\vec{x}(T-1),\vec{x}(T)$ with $h\geq0$.
\end{definition}

The first two of these properties relate to the contraction of network states over time.
These ensure that states from ESNs with different initial conditions will converge over time with respect to the same input sequence, and will not diverge so long as the same input sequences continues to be fed in.

\begin{definition} \textbf{-- Fading Memory}\\
    A network is said to have \textit{fading memory} when the outputs associated to inputs that are close in the recent past are close, even when those inputs may be very different in the distant past.
    
    Formally, a network has fading memory if for all left-inﬁnite input sequences $\vec{x}(-\infty),\dots,\vec{x}(T-1),\vec{x}(T)$ there exists a null sequence $\delta(h)$ such that for all input sequence suffixes $\bar{x}(h) = \vec{x}(T-h),\dots,\vec{x}(T-1),\vec{s}(T)$ with $h\geq0$, for all left-infinite input sequences $\bar{w}\cdot\bar{x}(h)$ and $\bar{w}\,'\cdot\bar{x}(h)$ (ie. a left-infinite sequence whose suffix is $\bar{x}(h)$), for all states $\vec{s}$ \textit{end-compatible} with $\bar{w}\cdot\bar{x}(h)$ and states $\vec{s}\,'$ end-compatible with $\bar{w}\,'\cdot\bar{x}(h)$ it holds that $||\vec{s}-\vec{s}\,'||_2<\delta(h)$.
\end{definition}

\begin{figure}[h]
    \centering
    \includegraphics[width=0.8\textwidth]{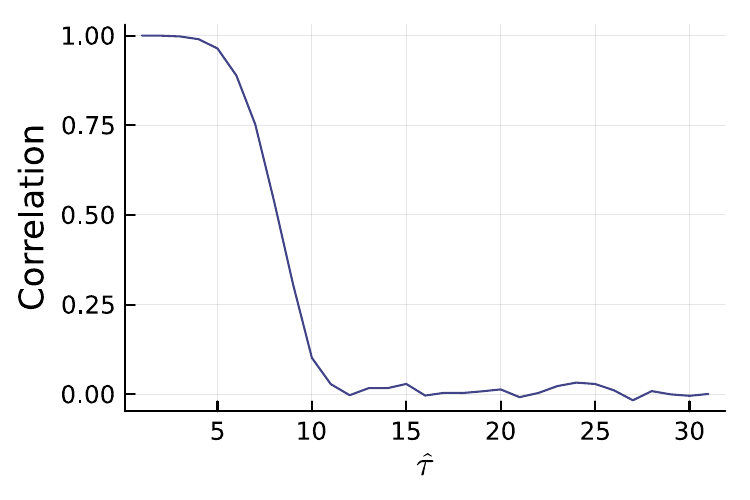}
    \caption{Demonstration of the fading memory property via recollection performance.}
    \label{fig:fading_memory}
\end{figure}

\begin{figure}[h]
     \centering
     \includegraphics[width=0.8\textwidth]{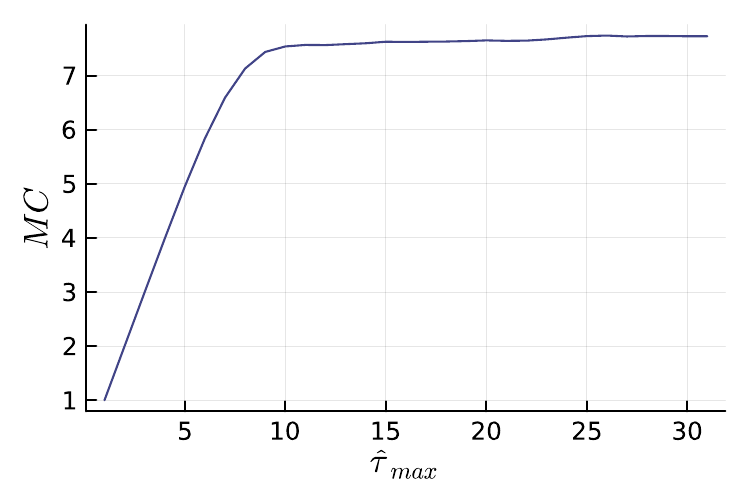}
     \caption{Convergence of the memory capacity.}
    \label{fig:memory_capacity}
\end{figure}

The fading memory property states that the impact of past inputs on the current state of the ESN decreases the further (temporally) those past inputs are from the current state.
Practically, this means that every state of the network is a function of only a finite number of previous inputs of the network such that
\begin{equation}
\begin{split}
    \vec{s}(t) &= \mathcal{E}_\infty\Big(\vec{s}(t-1),\vec{s}(t-2),\dots\Big)\\
    &\approx\mathcal{E}\Big(\vec{x}(t-1),\vec{x}(t-2),\dots,\vec{x}(t-\tau_{MC})\Big)\,.
\end{split}
\end{equation}
where $\tau_{MC}$ is a measure of how long the reservoir's memory is, and is in some way related to the memory capacity ($MC$).
The memory capacity for a particular reservoir computer can be quantified in a number of ways, however the original method as proposed by Jaeger~\cite{jaeger2001echo} is the most common
\begin{equation}
    MC = \sum_{\hat{\tau}=0}^{\hat{\tau}_{max}} \text{correlation}\bigg( \{x(t)\}_{t=1}^{T-{\hat{\tau}}}, C_{out}\{\vec{s}(t)\}_{t=\hat{\tau}+1}^{T} \bigg)
\end{equation}
where $x$ is a signal of standard Gaussian noise and $C_{out}$ are the learned parameters for predicting the $\hat{\tau}$th previous value of the signal.
The memory capacity converges for sufficiently large $\hat{\tau}_{max}$, with a demonstration presented in Figures \ref{fig:fading_memory} and \ref{fig:memory_capacity}.

At this point one may note similarities between the reservoir state representation of the input signal $\vec{x}$ and the representation generated by an embedding.
Indeed, the earliest mention of reservoir computers acting as embedding machines is to the original work by Jaeger~\cite{jaeger2001echo}.
The real breakthrough however came in the work of Hart et al.~\cite{hart2020embedding} who were able to put forward a proof that ESNs do indeed generate an embedding of the underlying signal with positive probability.
Hart~\cite{hart2024generalised} would later go on to prove that continuous time linear reservoir computers also generate an embedding almost surely.

\section{Reservoir Design}

Although thus far we have looked at RNNs, at its simplest a reservoir computer is an input-output automaton that utilises a dynamical system (called the \emph{reservoir}) with echo states as the central unit for information processing.
Time series $\vec{x}(t)$ are fed into the dynamical system as a form of driving signal and the system is allowed to evolve naturally according to
\begin{equation}
    \frac{d\vec{s}}{dt}(t) = f_{RC}(\vec{s}(t), \vec{x}(t), V)
\end{equation}
where $f$ describes the dynamics and $V$ collects the \emph{fixed} parameters of the reservoir.
The states $\vec{s}(t)$ (called the \emph{activation states}) of this system can then be recorded and used for training purposes.

There are two corollaries from the use of a fixed dynamical system:
\begin{enumerate}
    \item The only training is from $\vec{s}(t)$ to our desired output $\vec{y}(t)$.
    \item The burden of developing a good model simplifies to choosing a useful dynamical system.
\end{enumerate}
While this second point may appear daunting at face value, in practice there are a vast number of systems that may be considered.

\subsection{Echo State Networks (ESN)}

\begin{figure}
    \centering
    \includegraphics[width = 0.95\textwidth]{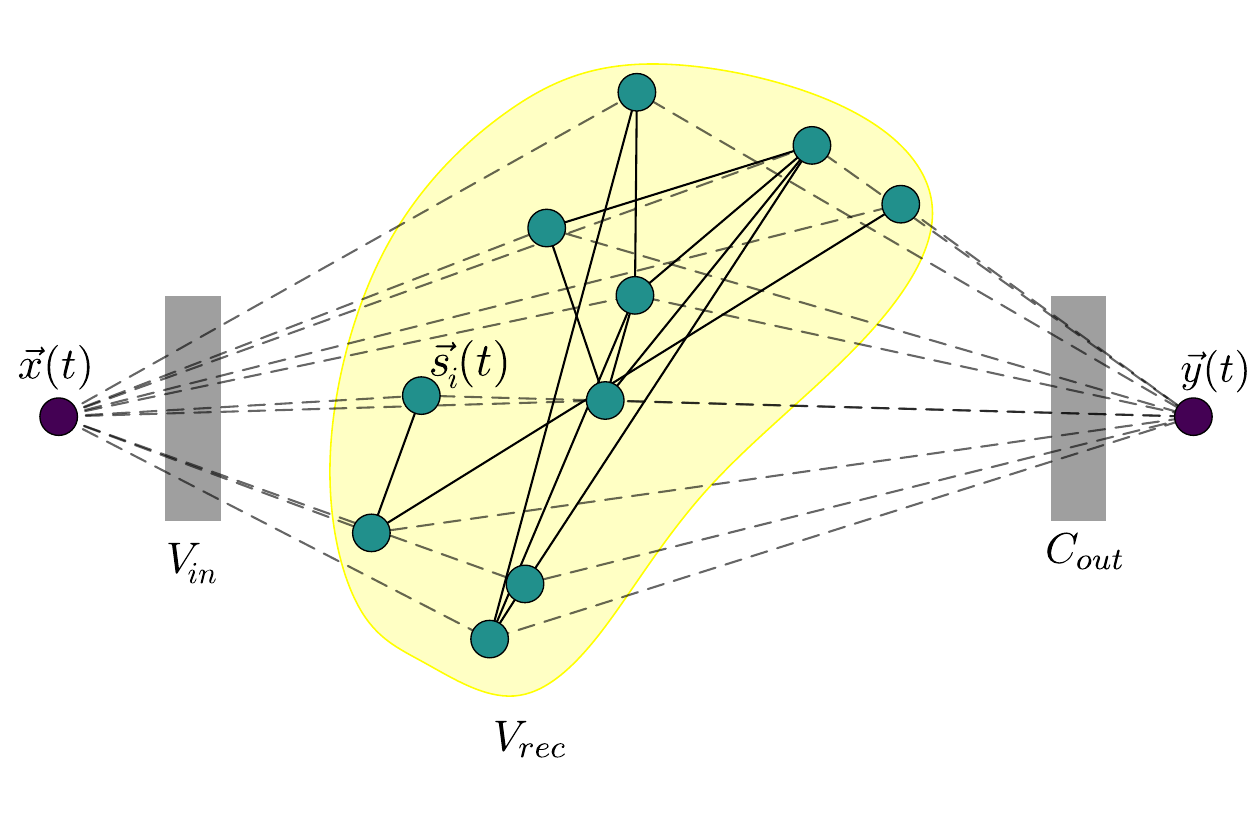}
    \caption{Visualisation of an echo state network.}
    \label{fig:ESN_diagram}
\end{figure}

The most well-studied reservoir computer is the ESN, which we may think of as a fixed RNN.
In this case we are now working in discrete time so may rewrite our evolution function as
\begin{equation}
    \vec{s}(t+1) = \tanh(V_{in} \vec{x}(t)+V_{rec}\vec{s}(t) +V_{bias})
\end{equation}
noting the choice of a sigmoid activation function which is typical for ESNs.

The fixed parameters $V$ in the case of an ESN are generated randomly:
\begin{itemize}
    \item The input matrix $V_{in}\in\mathbb{R}^{n\times k}$ typically draws elements randomly from an appropriate distribution (typically uniform or Gaussian).
    \item The reservoir adjacency matrix $V_{rec}\in\mathbb{R}^{k \times k}$ is initially generated as the adjacency matrix of an Erd\H{o}s-R\'{e}nyi random graph with an average degree $d$ then scaled to have a maximum absolute eigenvalue (or \emph{spectral radius}) of $\rho$.
    \item The bias matrix $V_{bias}\in\mathbb{R}^{k}$ is generated in a similar way to $V_{in}$.
\end{itemize}

The control of the dynamics and resulting impact on performance therefore falls to our selection of the \emph{hyper}-parameters $k$, $d$ and $\rho$.
Since we do not look to train these hyper-parameters, we rely on a series of heuristic approaches to their selection.,
What follows is a summary of these heuristics, the work of Luko{\v{s}}evi{\v{c}}ius~\cite{lukovsevivcius2012practical} provides a more thorough introduction to implementing an ESN.\\

\noindent \textbf{Reservoir Size $k$}\\
The general intuition for the reservoir size is that bigger is better.
Unlike traditional RNNs where utilising a large number of nodes leads to significantly increased training loads, this is mitigated by the fixed nature of the ESN.
As such it is common to set $k$ as large as reasonably possible for the hardware we are using.
It is common to see ESNs with $k$ in the order of hundreds or thousands.\\

\noindent \textbf{Average Degree $d$}\\
Of the hyper-parameters of relevance, $d$ is of the least importance.
General agreement amongst the RC community is that $d$ should be sufficiently small so as to keep the network sparse.
Typical values vary around $5\%$ of the network size $k$, and rarely are seen in excess of $10\%$.
More generally, work in RC has generally shown that the choice of edges is relatively unimportant.
Other network structures (eg. lattice or ring networks) have been shown to function similarly to random networks, and so often little consideration is given to these structural concerns and instead more focus is given to the critical hyper-parameter $\rho$.\\

\noindent \textbf{Spectral Radius $\rho$}\\
The spectral radius, also called the internal gain, is the mechanism by which we are able to control the dynamics of the reservoir and ensure the presence of echo states.
It is often useful to think of $\rho$ as controlling the \emph{memory} of the reservoir and thus consider the following extremes:
\begin{itemize}
    \item $\rho=0$ removes all memory in the ESN (each state is a funciton of the input only)
    \item $\rho\gg1$ the ESN memory dominates the input signal and the states become uninformative.
    \item $\rho\approx 1$ There is a balance between the forcing from the input and the memory of the ESN.
\end{itemize}

As we increase $\rho>1$ we approach a tipping point past which the ESN memory dominates.
This tipping point is often referred to as the edge of stability or the edge of chaos due to the change in dynamics observed.
As such a common heuristic is to set $\rho=1$, ensuring memory while preventing deconstructive dynamics, however research has shown that setting $\rho$ closer to the edge of chaos tends to benefit performance.
Methods for selecting $\rho$ remains an area of active research.\\

\noindent \textbf{Other Parameters}\\
It is not uncommon to see a more hyper-parameter-dense variation on the ESN evolution equation given by 
\begin{equation}
    \vec{s}(t+1) = (1-\alpha)\vec{s}(t) + \alpha\tanh(\eta V_{in} \vec{x}(t)+V_{rec}\vec{s}(t) +V_{bias})\,.
\end{equation}
This introduces the following additional hyper-parameters:
\begin{itemize}
    \item $\alpha$ is an additional mechanism that can be used to artificially increase the amount of memory in the reservoir.
    For signals with large time scales this can allow the reservoir to better extract dynamics, although similar behaviour can be induced through other hyper-parameters.
    \item $\eta$ is a parameter used to control the magnitude of the input.
    Due to the sigmoid activation function $\tanh(\cdot)$, signal values outside the range $[-1,1]$ are compressed and information is lost.
    Similar results to introducing the hyper-parameter $\eta$ can be obtained by normalising the input signal $\vec{x}$ prior to input.
\end{itemize}

\subsection{Other Reservoir Computers}

\noindent \textbf{Liquid State Machines (LSM)}\\
Another type of RNN based reservoir computer, LSM utilise a spiking neural network as opposed to the more simple structure of the ESN.
Many of the relevant intuitions for the ESN are also applicable to the LSM in terms of structural design; randomly generated matrices, high dimensional networks, control of hyper-parameters to ensure dynamics.\\

\noindent \textbf{Single Delay Node}\\
Rather than the spatio-temporal ESN, Appeltant et al.~\cite{appeltant2011information} considered a temporal only approach by time-multiplexing the input and feeding it into a single node with a delay.
The common choice of delay system is the Mackey-Glass oscillator such that
\begin{equation}
    \frac{ds}{dt}(t) = -s(t) + \frac{
        \eta \left( s(t-\tau) + \gamma x(t)\right)
    }{
        1 + \left( s(t-\tau) + \gamma x(t)\right)^p
    }
\end{equation}
where $\eta, \tau, \gamma$ and $p$ are all hyper-parameters.
Note that due to the time-multiplexing the input and output are both scalar series, so the high dimensional states are taken by reading $s(t)$ at various delays.\\

\noindent \textbf{Agent Based Reservoir Computers}\\
Instead of nodes we may consider a large number of agents and measure their behaviour to an external drive.
This was the focus of work by Lymburn et al.~\cite{lymburn2021reservoir} which utilised a agent based model for a swarm under the influence of an external predator as a reservoir computer.
These have seen relatively little exploration, but previous honours projects have also considered a territorial predator model with a single prey as a driving signal.\\

\noindent \textbf{Physical Reservoir Computing}\\
Part of the appeal of RC is that rather than digital systems we may instead use physical dynamical systems to greatly increase the speed at which states are generated.
There are numerous physical models that have arisen in response to this:
\begin{enumerate}
    \item Actual ``liquid'' state machines
    \item Photonic circuits
    \item Spintronics
    \item Soft robotics
\end{enumerate}
We don't go into details of physical RC here, but the book by Nakajima.~\cite{nakajima2021reservoir} covers many of the current research areas.

\section{Training a Reservoir Computer}

Since the dynamical system is fixed, the only parameters that require tuning are those connecting the activation states to the desired output ($C_{out}$).

\subsection{Time Series Prediction}

The most common learning task for reservoir computers is time series prediction, where we are attempting to map our input signal $\vec{x}(t)$ to some desired output signal $\vec{y}(t)$.
Since there is no feedback between $C_{out}$ and the dynamics of the reservoir, the training process can be handled through linear regression such that
\begin{equation}
    \hat{y}(t) = C_{out}\vec{s}(t)
\end{equation}
where $\hat{y}(t)$ is the resulting approximation the target signal $\vec{y}(t)$.
Typically the determination of $C_{out}$ is handled by ridge regression (or Tikhonov regularisation) which determines the relevant parameter values by solving the least squares problem
\begin{equation}
    \min_{C_{out}} \left( \sum_t \big( \vec{z}(t) - C_{out}\vec{s}(t) \big)^2 + \beta\sum_i \big(C_{out,i}\big)^2 \right)
    \label{eqn:ridge_regression}
\end{equation}
where $\vec{z}(t)$ is our training target, used to distinguish the signal used in training to the signal $\vec{y}(t)$ used in testing.
The parameter $\beta$ is called the regularisation parameter and is used to prevent over-fitting in the training step.
A closed form solution to the above minimisation problem can be written as
\begin{equation}
    C_{out} = (\vec{s}\vec{s}^\top + \beta I)^{-1} \vec{s} \vec{z}\,.
\end{equation}

\subsection{System Simulation}

Alternatively, we may be interested in building a model which is able to simulate and forecast the original system.
Such a task is simply an extension of the time series prediction task, beginning with solving the regression problem in equation \ref{eqn:ridge_regression} with a target signal corresponding to a one-step forecast of the input signal (ie. $\vec{z}(t)=\vec{x}(t+1)$) such that
\begin{equation*}
    \vec{x}(t+1) = C_{out}\vec{s}(t)
\end{equation*}
as in the traditional RNN case.
This approximation can then be fed back into the reservoir computer in place of the input signal to facilitate free-run forecasts of the system.
In the case of an ESN the new evolution equation reads
\begin{equation}
    \vec{s}(t) = \tanh(V_{in}C_{out}\vec{s}(t-1) + V_{rec}\vec{s}(t-1) + V_{bias})\,.
\end{equation}

An issue remains here as to how we should choose the initial state of the reservoir for the free-run prediction.
Ideally, we would have a sufficiently long series leading up to the start of our free-run to allow the reservoir states to converge.
In the absence of this we run into issues with not falling onto the attractor, however methods to overcome this have been considered such as training another reservoir computer to map from the input space to the output space (ie. $\vec{z}(t)=\vec{x}(t)$).

\subsection{Time Series Classification \& Reservoir Time Series Analysis}

Time series classification is another task where reservoir computing has seen extensive use.
Rather than relating points in time, classification looks to relate an entire series to a single classification or label.
The primary difference here is that rather than training from activation states to an output series, the activation states are summarised into some \textit{feature} $H$.
A large number of features can then be generated for different series and then a classifier can be trained.
Typical classifiers include support vector machines and random forests.

There are a number of different features that have been considered, however some of the most common are:
\begin{align}
    H_{RC}^{last} &= \vec{s}(T) \\
    H_{RC}^{mean} &= \frac{1}{T}\sum_t \vec{s}(t) \\
    H_{RC}^{R} &= C_{out} \text{ with } \vec{z}(t) = \vec{x}(t+1)
\end{align}
These features all essentially aim to quantify some aspect of the underlying system that allows us to distinguish between signals from differing systems.

While developed with supervised machine learning in mind, the idea of quantifying aspects of dynamical systems in some characteristic feature is not unfamiliar to you.
This was the driving idea behind invariant measures, and indeed some such measures have been approximated with reservoir computers for the task of time series classification.
While \textit{classification} focuses on the supervised learning task, \textit{reservoir time series analysis} is the broader field that looks at methods for extracting information from the embedding generated by a reservoir computer (the reservoir space embedding)~\cite{thorne2022reservoir}.\\

\begin{figure}
    \centering
    \includegraphics[width=0.325\textwidth]{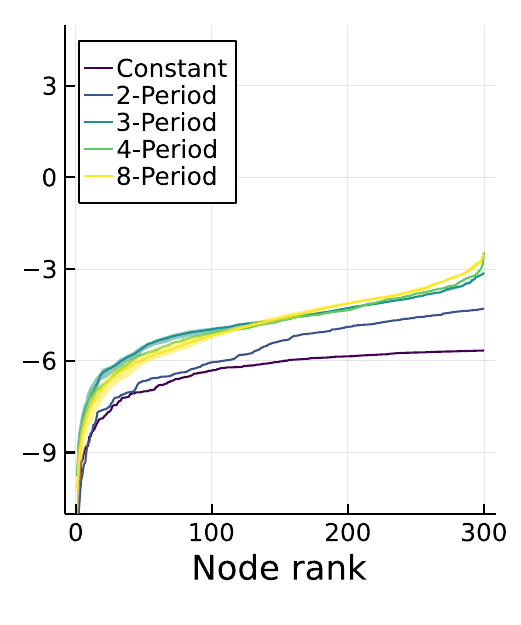}
    \includegraphics[width=0.325\textwidth]{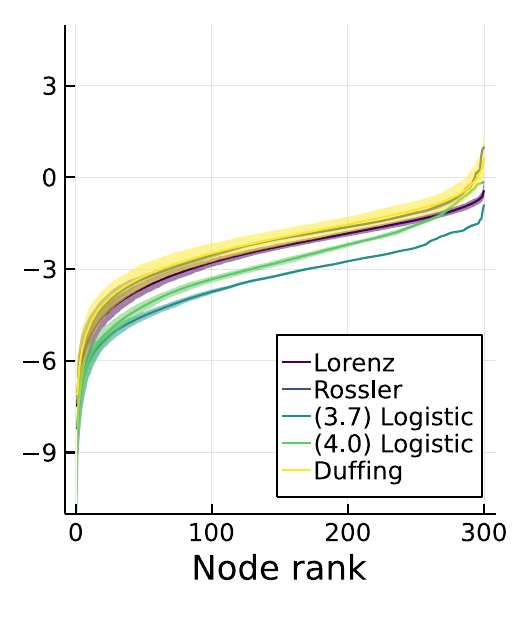}
    \includegraphics[width=0.325\textwidth]{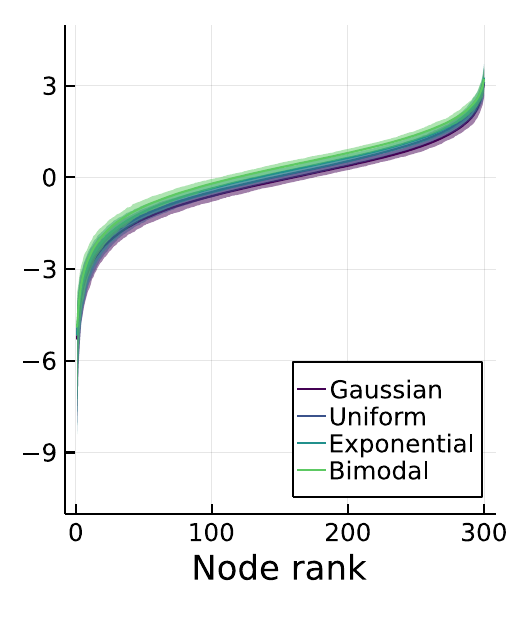}
    \caption{Comparison of readout magnitude curves generated from a 1-step prediction for periodic, chaotic and random signals (left ot right).}
    \label{fig:rcc_curves}
\end{figure}

Reservoir time series analysis is still novel and many of the methods currently overlap with those from classification, however some efforts have been made to develop specialised invariant measures.
One example is the reservoir computing complexity (RCC)~\cite{thorne2022novel}.
The idea behind the RCC is that reservoir use a greater proportion of nodes when forecasting complex systems compared to simple systems, and this greater proportion is quantified in the magnitude of the readout vector $C_{out}$ in the presence of regularisation.
Visualisation of this phenomenon is provided in Fig.~\ref{fig:rcc_curves}.
The behaviour in the curves can be summarised by taking an approximation of the area under the curves
\begin{equation}
    H_{RC}^{cplx} = \log\left(\sum_{i=1}^{k-1} \frac{\hat{C}_{out,i}+\hat{C}_{out,i+1}}{2}\right)
\end{equation}
where $\hat{C}_{out}$ is the rank ordered variation of $C_{out}$.

\chapter{Recurrence Plots}

Recall that we can describe dynamics in terms of the recurrent behaviour.
In principle, there should be equivalent information encoded in analysing the recurrence times and states as in the continuous flow.
Recurrence (once we define some kind of neighbourhood) is simpler to analyse than the full state dynamics and can be useful way of quickly getting a basic understanding of the system.
So how does one construct an intuitive representation of recurrence?

\begin{figure}[h]
    \centering
    \includegraphics[width=\textwidth]{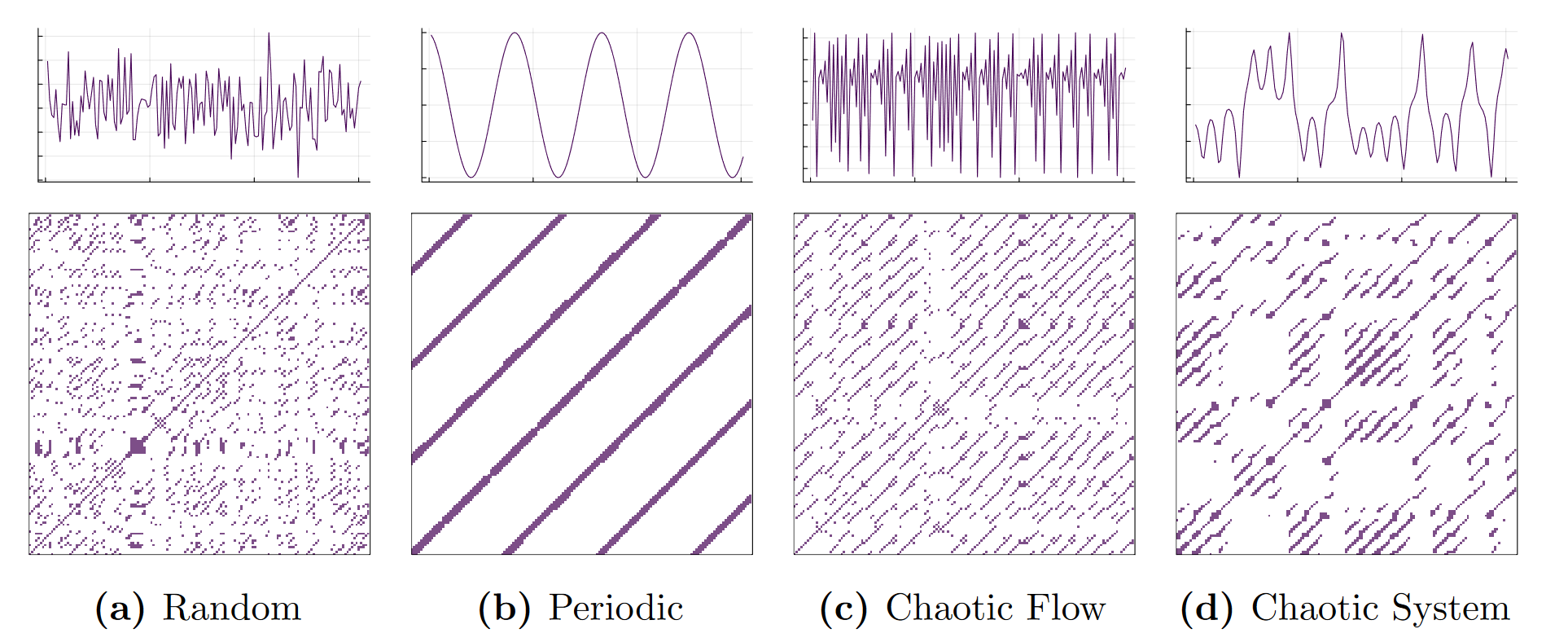}
    \caption{Example $RP$s for signals with varying complexity. Signals correspond to random noise, a periodic sine wave, and observations from a chaotic Lorenz system and a chaotic logistic map (sub-figures (a)-(d) respectively).}
    \label{fig:recurrence_plots}
\end{figure}

One way to visualise and quantify recurrence is through the usage of a recurrence plot $RP\in\mathbb{R}^{T \times T}$ such that
\begin{equation*}
    R_{ij} = \begin{cases}
        1, & d(\vec{x}(i), \vec{x}(j))<\epsilon\\
        0, & \, \text{otherwise}
    \end{cases}
\end{equation*}
where $d(\cdot)$ is some (not necessarily Euclidean) distance metric and $\epsilon$ is the recurrence threshold.
Recurrence plots (in theory) contain all the ``interesting" details of a time series' dynamics.
Thus, certain patterns and structures in $RP$ can be related to dynamical properties:
\begin{itemize}
    \item \textit{Diagonal lines} encode the duration that recurrent states align, with the length of the line equal to the length of time. The distribution of lengths also helps characterise the mixing behaviour of the system. Notably the lengths of these diagonal lines are proportional to the Lyapunov time ($\frac{1}{\lambda_{max}}$).
    \item \textit{Vertical lines} can highlight stationary/steady states, or states which are highly laminar.
    \item Vertical voids or \textit{recurrence times} provide information about when states recur, with regularity suggesting periodic dynamics.
\end{itemize}
By plotting $RP$ (eg. as a heatmap) one can qualitatively analyse these above structures to understand the underlying system, with an example of this provided in Figure \ref{fig:recurrence_plots}.

\section{Construction Considerations}

At this point it is apparent there are a number of choices that need to be made before one can generate a recurrence plot, and making these decisions is generally the most difficult part of such analysis.\\

\noindent \textbf{Distance Metric $||\cdot||$}\\
The most common choices for distance metrics when working with numerical time series data are the Manhattan distance
\begin{equation}
    d_1(\vec{a}, \vec{b}) = \sum_{i} |a_i - b_i|\,,
\end{equation}
the Euclidean distance 
\begin{equation}
    d_2(\vec{a}, \vec{b}) = \sqrt{\sum_{i} (a_i - b_i)^2}
\end{equation}
and the supremum (or maximum) distance
\begin{equation}
    d_\infty(\vec{a}, \vec{b}) = \max_i(|a_i - b_i|)\,.
\end{equation}
Of these methods the supremum norm holds specific advantages in decreasing computational costs when generating $RP$s as well as allowing certain results such as the correlation dimension to be derived analytically.

This freedom over the choice of norm however does allow the ideas of recurrence to be transitioned to series with different structures.
For example, a Kronecker delta function as a distance metric can construct $RP$s on symbolic series~\cite{faure2010recurrence}.
More recently the edit distance has been used for defining $RP$s on series of event data, a method that continues to see ongoing application and research~\cite{marwan2023challenges}.\\

\noindent \textbf{Recurrence Threshold $\epsilon$}\\
Most commonly the threshold $\epsilon$ is chosen and fixed based on some heuristic regarding the distances expected; large enough that recurrences occur but not so large that everything recurs.
Rather than a fixed distance, we can instead choose $\epsilon$ based on a fixed \textit{recurrence rate} by choosing a percentage of points that we expect should recur.
Other heuristic and numerical approaches also exist for estimating this fixed threshold with different trade-offs with regards to resolution and noise robustness.\\

\noindent \textbf{Embedding $\vec{x}$}\\
Ideally one would construct the $RP$ on the full state space trajectory, however as discussed throughout this unit rarely do we actually have all this information.
Instead it is common for embedding to be constructed and the embedded space used for determining which states recur.
Interestingly there is ongoing debate as to the importance of this step, as research has shown that $RP$s generated from observations without embedding are qualitatively similar to those with embedding and many important results (such as the derivation of the correlation dimension) can be done without embedding.
Despite this, the general consensus amongst the community is that if one is able to perform an embedding before analysing a $RP$ one should do so!

\section{Recurrence Quantification Analysis}

Where recurrence plots provide a qualitative assessment of system behaviour, recurrence quantification analaysis (RQA) formalises and quantifies the geometrical properties from recurrence plots.
This information can then be used for tasks like change point detection, determining complexity measures or periodic-chaotic transition detection.
Broadly speaking RQA approaches follow a similar idea; summarise the empirical distribution of some qualitative feature in the $RP$, then calculate some quantitative measure based on that distribution.
Measures that look at the distribution of diagonal $\mathcal{H}_{diag}$ lines, vertical $\mathcal{H}_{vert}$ lines and recurrence times $\mathcal{H}_{recc}$ are the most common of these qualitative features.

\subsection{Point Density}

The first couple of features we look at here concern simply the distribution of points in the $RP$ rather than the line structures.
For a fixed distance threshold, the first of these is the recurrence rate
\begin{equation}
    RR = \frac{1}{T^2}\sum_{i=1}^{T}\sum_{j=1}^{T} RP_{i,j},
\end{equation}
which measures the density of recurrence points.
If we exclude the line of identity, this corresponds to our definition of the correlation sum.
A similar value
\begin{equation}
    N_{avg} = \frac{1}{T}\sum_{i=1}^{T}\sum_{j=1}^{T} RP_{i,j},
\end{equation}
which measures the average number of recurrences (or neighbours) across points in the phase space can also be calculated.

With a fixed recurrence rate, we may instead work backwards to determine the fixed distance threshold
\begin{equation}
    \epsilon_{RR} = \max_{i,j}\Big(d(\vec{x}(i), \vec{x}(j)) \Big| RP_{i,j}=1\Big)\,.
\end{equation}

\subsection{Diagonal Lines}

The histogram of diagonal lines is given by
\begin{equation}
    \mathcal{H}_{diag}(l) = \sum_{i,j=1}^{T} (1-RP_{i-1,j-1})(1-RP_{i+l,j+l})\prod_{k=0}^{l-1} RP_{i+k,j+k},
\end{equation}
with corresponding experimental distribution
\begin{equation}
    \mathcal{P}_{diag}(l) = \frac{\mathcal{H}_{diag}(l)}{\sum_{\hat{l}=l_{\min}}^{T} \mathcal{H}_{diag}(\hat{l})},
\end{equation}
noting that $l_{min}$ should agree with the minimum length of lines considered for the respective features.
Diagonal line features are used as indicators of chaos-order transitions in dynamical systems
Some of the common quantitative features derived from $\mathcal{H}_{diag}$ are detailed below.\\

\noindent \textbf{Determinism}\\
Processes with uncorrelated, weakly correlated, stochastic or chaotic behaviour cause none or short diagonals.
Deterministic processes cause longer diagonals and less isolated recurrence points.
A measure of determinism (or predictability) therefore is the ratio of recurrence points which form long ($l>l_{min}$) diagonals to all recurrence points (excluding tangential motion and the line of identity),
\begin{equation}
    DET = \frac{
            \sum_{l=l_{min}}^{T} l \mathcal{H}_{diag}(l)
        }{
            \sum_{l=1}^{T} l \mathcal{H}_{diag}(l)
        }.
\end{equation}

\noindent \textbf{Average Diagonal Line Length}\\
Average diagonal line length is related to divergence of trajectories and so is a measure of mean prediction time,
\begin{equation}
    L_{avg} = \frac{
            \sum_{l=l_{min}}^{T} l \mathcal{H}_{diag}(l)
        }{
            \sum_{l=l_{min}}^{T} \mathcal{H}_{diag}(l)
        }\,.
\end{equation}

\noindent \textbf{Longest Diagonal}\\
The longest diagonal line is related to the exponential divergence of the phase space trajectory, and so is proportional to the Lyapunov time of the system.
A related measure is the inverse in the \textit{divergence}, which is therefore proportional to the largest Lyapunov exponent.
The faster the trajectory segments diverge, the shorter the diagonal lines and the higher the divergence,
\begin{align}
    L_{max} &= \max(l_{i} | \mathcal{H}_{diag}(l_i) > 0) \propto \frac{1}{\lambda_{max}}\,, \\
    DIV &= \frac{1}{L_{max}}.
\end{align}

\noindent \textbf{Recurrence Entropy}\\
The Shannon entropy of the distribution of diagonal lines, 
\begin{equation}
    ENTR = -\sum_{l=l_{\min}}^{T} \mathcal{P}_{diag}(l) \log_2\left(\mathcal{P}_{diag}(l)\right)\,,
\end{equation}
reflects the complexity of the recurrence plot with respect of diagonal lines.
For example, both uncorrelated noise and periodic signals tend to have low $ENTR$ values.

\subsection{Vertical Lines}

The histogram of vertical lines is given by,
\begin{equation}
    \mathcal{H}_{vert}(l) = \sum_{i,j=1}^{T} (1-RP_{i,j-1})(1-RP_{i,j+l})\prod_{k=0}^{l-1} RP_{i,j+k},
\end{equation}
with corresponding experimental distribution $\mathcal{P}_{diag}(l)$ calculated as in the diagonal line case.
Vertical lines are predominantly a sign of trajectories that pass sufficiently close to the neighbourhood at some other time (ie. tangential motion).
However, vertical lines have also been shown to indicate the presence of laminar states, with corresponding measures being able to detect chaos-chaos transitions in addition to chaos-order transitions
Some of the common quantitative features derived from $\mathcal{H}_{vert}$ are detailed below.\\

\noindent \textbf{Laminarity}\\
Calculated analogously to $DET$, the laminarity given by, 
\begin{equation}
    LAM= \frac{
            \sum_{l=l_{min}}^{T} l \mathcal{H}_{vert}(l)
        }{
            \sum_{l=1}^{T} l \mathcal{H}_{vert}(l)
        },
\end{equation}
tells us the occurrence of certain states that do not change or change slowly for some time, known as \textit{laminar states}.
These states are typical of intermittency.\\

\noindent \textbf{Trapping time}\\
The average length of vertical structures or \textit{trapping time},
\begin{equation}
    TT = \frac{
            \sum_{l=l_{min}}^{T} l \mathcal{H}_{vert}(l)
        }{
            \sum_{l=l_{min}}^{T} \mathcal{H}_{vert}(l)
        }\,,
\end{equation}
estimates the mean length a system will reside within a particular state, and can be useful for detecting chaos-chaos transitions and so can be used to detect intermittency.\\

\noindent \textbf{Longest Vertical \& Entropy}\\
Analogous measure for the maximal vertical line $V_{max}$ and entropy of vertical line structures $VENTR$ also exist, although are proportionally less used.
Interpretations are similar as in the diagonal line case, with $V_{max}$ reflecting the maximal time trapped in a particular state while $VENTR$ is again a measure of complexity.

\subsection{Recurrence Times}

The distribution of recurrence times is given as the distribution of vertical distances between recurrent structures
\begin{equation}
    \mathcal{H}_{recc}(l) = \sum_{i,j=1}^{T} RP_{i,j-1}RP_{i,j+l}\prod_{k=0}^{l-1} (1-RP_{i,j+k})\,,
\end{equation}
again with a corresponding $\mathcal{P}_{recc}$.
The recurrence times relate to the \textit{point-wise dimension}, which is an invariant of the system.
In particular the point-wise dimension $D_{p}$ is related to the mean recurrence time,
\begin{equation}
    r_{avg} = \frac{
            \sum_{l=l_{min}}^{T} l \mathcal{H}_{recc}(l)
        }{
            \sum_{l=l_{min}}^{T} \mathcal{H}_{recc}(l)
        }\,,
\end{equation}
by the relation $r_{avg}\sim \epsilon^{D_{p}}$.\\

Analogous features such as ratios of lengths, maximal length and entropy can be calculated for the recurrence times as with the diagonal and vertical lines.
In general recurrence time measures have been shown to detect a number of periodic and chaotic transitions in dynamical systems, including transitions involving strange non-chaotic attractors.

% \subsection{Other Features}

% Dynamical invariants can be derived from recurrence plots (see Eckmann and Marwan papers) but it is difficult and not as immediate as suggested by RQAs.  $K_{2}$ can be estimated from a recurrence plot (e.g. length distribution of diagonal lines).  Recall the Kaplan-Yorke conjecture relates $K_{2}$ and the sum of the positive Lyapunov exponents and so for systems with one positive maximum Lyapunov exponent recurrence plots provide an alternative estimate.  We can also estimate $D_{2}$ the correlation dimension and JRP's can be used to estimate joint Renyi entropies.  (See, (Marwan, 2007))\\

% RQA was developed and matured using symmetric recurrence plots (metric threshold).  The original recurrence plot paper led to asymmetric recurrence plots (topological threshold).  Recently quantitative measures and techniques to analyze these original recurrence plots have been proposed.

\chapter{Dynamical Networks}
Up until now, we have extensively studied dynamical systems in isolation of the canonical form,
$$ \dot{x} = f(x), x\in \mathbb{R}.$$
However, what happens if two or more individual dynamical systems are allowed to interact with each other? In the extreme case, one can consider a network (not necessarily fully connected) of interacting dynamical systems. Such a system is often referred to as a \textbf{dynamical network} as it refers to dynamics (temporal behaviour) occurring along some topological structure (the network), the latter of which encodes some form on information flow between individual dynamical systems.

\section{Coupling and Synchronisation}
\subsection{Coupling Functions}
Consider two dynamical systems of the following form
\begin{align*}
    \dot{x} &= f_{1}(x),\\
    \dot{y} &= f_{2}(y) + G(x,y),
\end{align*}
where $x,y \in \mathbb{R}^m$ and $G: \mathbb{R} \times \mathbb{R} \to \mathbb{R}$. These two systems adhere to the canonical form ($\dot{x}=f(x)$) with the main difference being the addition of an interaction term in the second system. The term $G$ is called the \textbf{coupling function} and describes the mechanism in which one system affects the other. The above case corresponds to unidirectional coupling where $x$ influences $y$, but not reverse (sometimes called the driving and driven system respectively). The case where coupling is applied to direct state variables $x,y$ is known as amplitude coupling. If one limits the dynamical systems to the study of oscillators whose phase can be defined, coupling in these terms is similarly referred to as phase coupling (see \ref{sec: Sync})

The coupling function can be described as containing two parts: its (1) strength, and (2) form,
$$G(x,y) = Kq(x,y)$$
The former is usually mediated by way of a constant multiplier $K$, where larger values indicate stronger coupling. The form of the coupling function $q(x,y)$ describes the nature of the interaction and can heavily affect the dynamics of the coupled system in unpredictable ways \cite{stankovski2017coupling}. Two of the most common forms of coupling functions are the direct coupling $q(x,y) = q(y)$ representing unidirectional influence, and diffusive coupling $q(x,y) = q(y-x)$.

The case for diffusive coupling popularised by Kuramoto in the study of phase coupled oscillators is a particularly interesting case \cite{kuramoto1984chemical}. The mathematical form of diffusive coupling represents the case where the strength of interaction is dependent on the difference between two states. Typically, this is set such that a larger difference corresponds to stronger coupling. The case where $q'(0)<0$ is also called dissipative coupling, referring to coupling forces pushing the system to converge towards the same state. Conversely, the $q'(0)>0$ is called repulsive coupling for the same reason. \cite{kuramoto1984chemical, stankovski2017coupling}

\subsection{Two-body systems and synchronisation}\label{sec: Sync}
The presence of sufficiently strong diffusive coupling between two systems can result in the convergence system states, a phenomenon known as synchronisation. This results in both systems having the same state that evolve together in time. To demonstrate this, we consider the following bi-directionally coupled pair of systems,
\begin{align*}
    \dot{x} &= f(x) + Kg(y-x),\\
    \dot{y} &= f(y) + Kg(x-y).
\end{align*}
To better analyse the extent of synchrony between both systems, we restrict our analysis to an auxiliary variable $z = x-y$, which describes the deviation between the two systems. Thus we can write the following equation for the evolution of $z$,
\begin{align*}
    \dot{z} &= \dot{x} - \dot{y}\\
    &= f(x) - f(y) + K[g(-z) - g(z)].
\end{align*}
If we assume that $g(z)>0$ for $z>0$, and $g(-z) = -g(z)$ (or at least $g(-z) \propto -g(z)$) then the evolution of the state difference simplifies to,
\begin{equation*}
    \dot{z} = f(x)-f(y) - 2Kg(z).
\end{equation*}
Analysing this reduced system, we can observe that there is a fixed point at $z^{*}=0$ (i.e. $x=y$). Furthermore, whilst we cannot guarantee the sign of $f(x)-f(y)$ for all time, it is clear that for sufficiently large values of $K$ and dissipative coupling ($g'(0) < 0$), this fixed point $z^{*}$ is stable.

To illustrate the onset of synchronisation, consider the following pair of dissipatively coupled R\"{o}ssler oscillators \cite{zou2012geometric},

\begin{equation*}
\begin{aligned}[c]
\dot{x}_{1} &= -y_{1} - z_{1} + K(x_{2} - x_{1})\\
    \dot{y}_{1} &= x_{1} + ay_{1}\\
    \dot{z}_{1} &= 0.4 + z_{1}(x_{1} - 8.5)\\
\end{aligned}
\qquad\qquad
\begin{aligned}[c]
\dot{x}_{2} &= -y_{2} - z_{2} + K(x_{1} - x_{2})\\
    \dot{y}_{2} &= x_{2} + ay_{2}\\
    \dot{z}_{2} &= 0.4 + z_{2}(x_{2} - 8.5)\\
\end{aligned}
\end{equation*}

\begin{figure}
    \centering
    \includegraphics[width=0.95\linewidth]{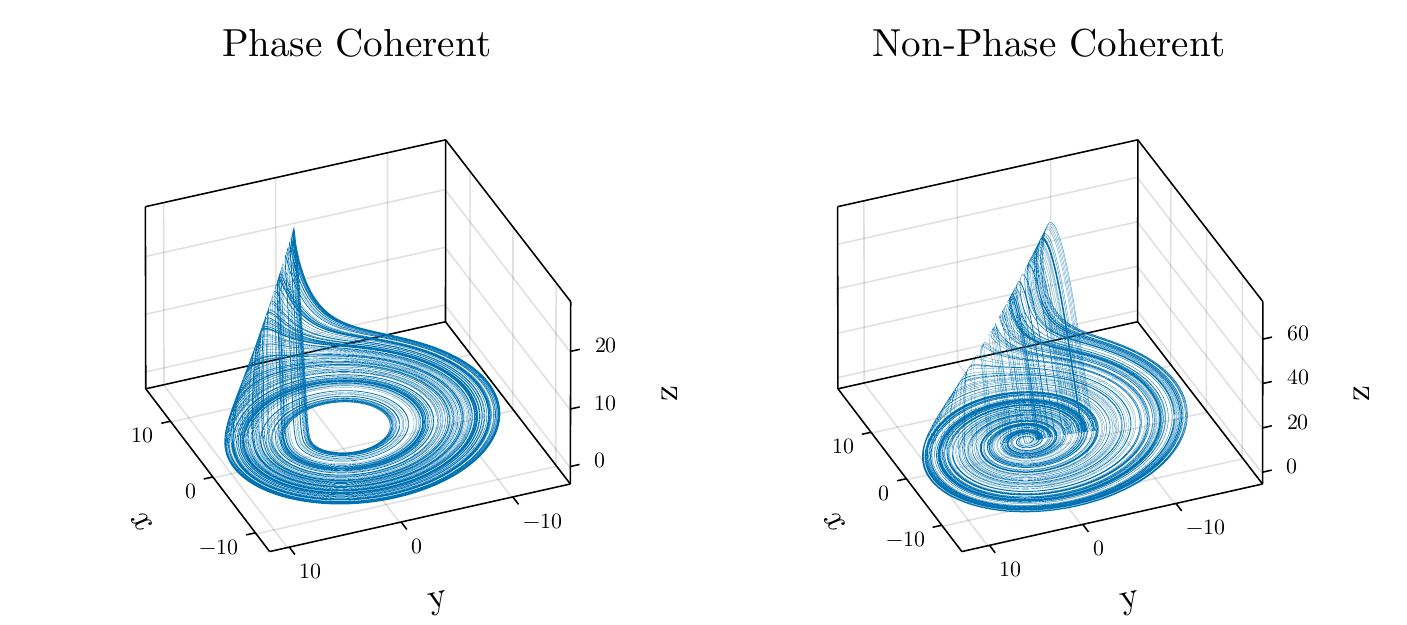}
    \caption{R\"{o}ssler atractor in the phase coherent $a = 0.165$ (left) and non-phase coherent $a = 0.265$ regime.}
    \label{fig:rossler}
\end{figure}

The regular R\"{o}ssler attractor consists of orbits rotating about an unstable fixed point in the origin alongside a folding on one size of the manifold in the $z$-direction (see Figure \ref{fig:rossler}). This folding results in chaotic behaviour. For the $a = 0.165$, the system dynamics lie within a phase coherent regime. The term ``phase-coherent'' refers to the rotational frequency of trajectories being roughly constant throughout all parts of the attractor. Hence, any two trajectories should only differ in their phases $\phi$ and radii (with respect to the $xy$ plane). For simplicity, we may estimate the phase at any point in a trajectory as $\phi = atan(\frac{y}{x})$.

\begin{figure}
    \centering
    \includegraphics[width=\linewidth]{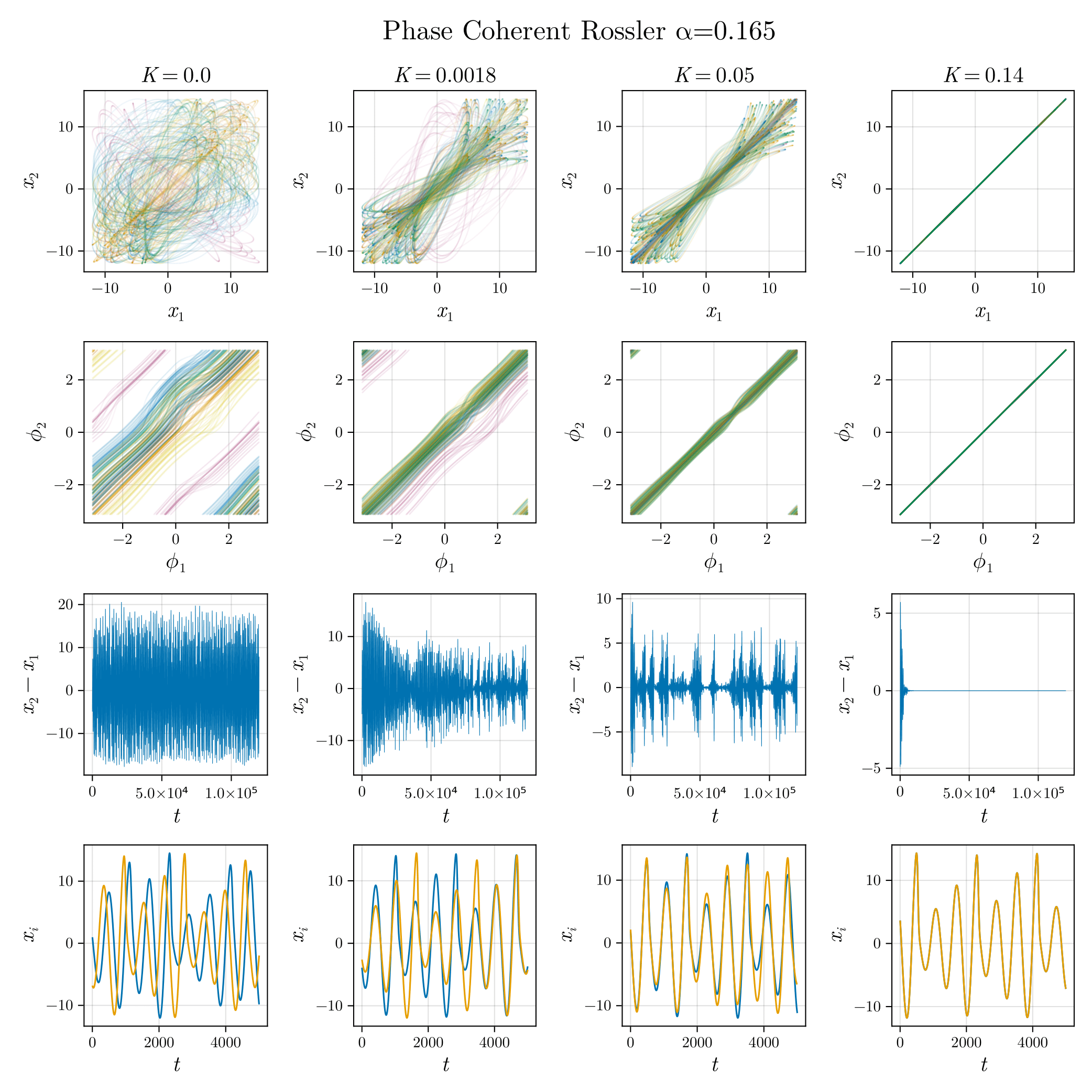}
    \caption{Two coupled phase coherent $(\alpha = 0.165)$ R\"{o}ssler oscillators with increasing coupling strength $K$ from left to right: asynchronous (no coupling), phase synchronisation, intermittent synchronisation, followed by complete synchronisation at $K=0.14$. Note the $\phi_1 - \phi_2$ plots showing diagonal slopes characteristic of phase coherence where all dynamics naturally adhere to a single a dominant frequency.}
    \label{fig:pc_rossler_sync}
\end{figure}

\begin{figure}
    \centering
    \includegraphics[width=1\linewidth]{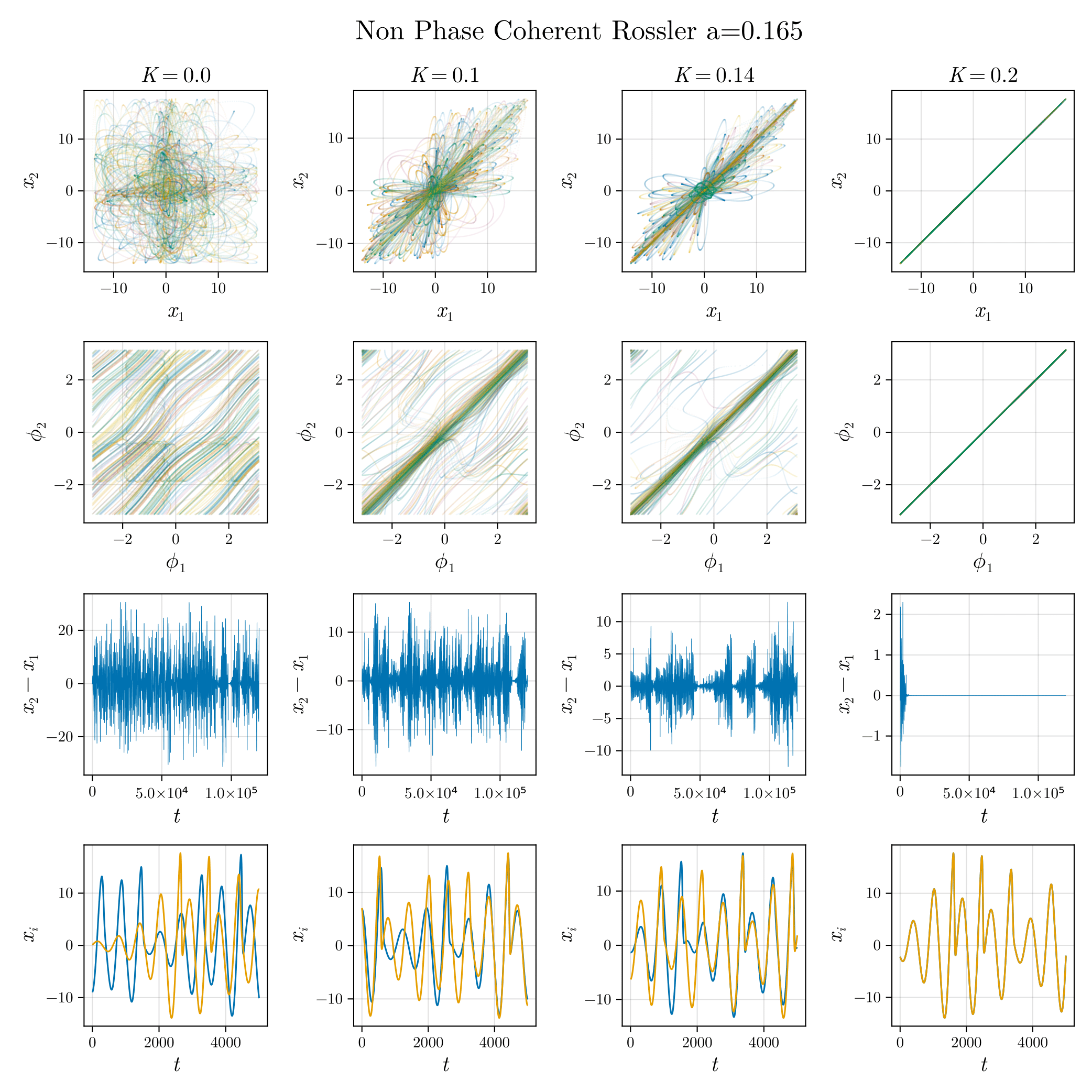}
    \caption{Two coupled R\"{o}ssler oscillators in the non-phase coherent regime $(\alpha = 0.265)$. Lack of phase coherence results in a larger coupling $K$ required to achieve similar synchronisation behaviour to the phase coherent case.}
    \label{fig:npc_rossler_sync}
\end{figure}

At $K = 0$, both oscillators do not interact at all and their phases differ according to their initial conditions. However, as $K$ is gradually increased the coupled dynamics of the oscillators undergo several different stages of change (see Figures \ref{fig:pc_rossler_sync} \& \ref{fig:npc_rossler_sync}). For low to moderate coupling, phase synchronisation occurs where the phases of the oscillators begin to converge, but differ in the amplitude. Increasing $K$ slightly, the system begins to exhibit intermittent synchronisation (not to be confused with intermittency) where the pair of oscillators converge to one another and synchronise for short periods at a time. This is shown in episodic collapses in the deviation of the coupled $x$-components. Finally, full synchronisation occurs for sufficiently large $K$ corresponding to a permanent convergence of both trajectories. 

\section{Dynamics on Networks}
\subsection{Oscillator Networks and Chimeras}
The coupling of oscillators can be easily extended to include an arbitrary $N$ number of (usually identical) oscillators. The exact wiring of the couplings may also be altered directly and need not be all-to-all connected. This can generally be represented in network form with the following equation:
\begin{equation*}
    \dot{x}_{i} = f(x_i) + K \sum_{j = 1}^{N} A_{ij}G(x_{i},x_{j}),
\end{equation*}
where $A$ is the adjacency matrix of the underlying network where nodes represent individual oscillators, and edges represent coupling between a pair of oscillators. If $A$ is symmetric, the coupling is bidirectional.

Previously, it was assumed that the behaviour of coupled identical phase-oscillators were relatively straightforward with phase synchronisation occurring for sufficiently strong coupling. These results were presumed to extend to networks of oscillators where sufficiently large coupling result in total synchronisation of the entire network. However, Kuramoto et al. \cite{kuramoto1984chemical} demonstrated that for certain network configurations and initial conditions it is possible for the coexistence of both coherent and incoherent clusters to occur. Therefore, a portion of the oscillators in the network are phase locked, whilst the remainder continue to be asynchronous. The coexistence of these two incongruous states led to these phenomena being labelled as ``chimera states'' after the mythological Greek beast.

The existence and stability of chimera states with a network of globally coupled oscillators has been studied from several perspectives. Focusing on the prototypical network of Kuramoto oscillators containing two clusters $\sigma, \sigma'$,
\begin{equation*}
    \frac{d\theta_{i}^{\sigma}}{dt} = \omega - \mu \langle \sin (\theta^{\sigma}_{i} - \theta^{\sigma}_{j} + \alpha) \rangle_{j\in \sigma} - \nu \langle \sin (\theta^{\sigma}_{i} - \theta^{\sigma}_{j} + \alpha) \rangle_{j\in \sigma'},
\end{equation*}
where $\mu$ and $\nu$ correspond to the intra- and inter-cluster coupling strengths, Ott \& Antonsen demonstrated that macroscopic dynamics of an inifinitely large network can be reduced to a set finite ODEs \cite{ott2008low}. For a system of two clusters, this can be further refined to 2 pairs of equations (case for 1 cluster shown below):
\begin{align*}
    0 &= \dot{\rho}_{X} + \frac{\rho^{2}_{X}-1}{2}[\mu \rho_{X} \cos\alpha + \nu\rho_{Y}\cos(\phi_{Y}-\phi_{X}-\alpha)],\\
    0 &= -\rho_{X}\dot{\phi_{X}} + \rho_{X}\omega - \frac{1 + \rho_{X}^{2}}{2}[\mu \rho_{X} \sin\alpha + \nu\rho_{Y} \sin (\phi_{X} - \phi_{Y} + \alpha)],
\end{align*}
where $X,Y$ identify the two clusters, $\rho_{X},\rho_{Y}$ are order parameters (degree of synchrony) of each cluster, and $\phi_{X},\phi_{Y}$ are their phases \cite{abrams2008solvable}. A dynamics of the chimera can be studied by restricting the above equations to the manifold $\rho_{X} = 1$ (i.e. cluster $X$ is fully synchronised), and applying the change of coordinates $r = \rho_{Y}$, $\psi = \phi_{X} - \phi_{Y}$,
\begin{align*}
    \dot{r} &= \frac{1-r^{2}}{2}[ \mu r\cos\alpha + \nu\cos (\psi - \alpha) ],\\
    \dot{\psi} &= \frac{1+r^2}{2r}[\mu r \sin\alpha - \nu \sin (\psi - \alpha)] - \mu\sin\alpha - \nu r \sin (\psi + \alpha).
\end{align*}
This reduced set of equations has a trivial fixed point at $r = 1$ (the fully synchronised state), as well as a non trivial fixed points for $r<1$ (chimera state). Interested readers should refer to \cite{panaggio2016chimera} for a more detailed discussion.

\subsection{FitzHugh-Nagumo Chimera model}
To numerically demonstrate a system exhibiting chimera states, we refer to the FitzHugh-Nagumo oscillator example provided by \cite{scholl2016synchronization, omelchenko2015robustness}. Briefly, the FitzHugh-Nagumo oscillator is a simplified form of the popular Hodgkin-Huxley model used to describe the neuron spiking dynamics with equations given by,
\begin{align*}
    \dot{u} &= u - \frac{u^3}{3}-v + I,\\
    \epsilon\dot{v} &= u + a - bv,
\end{align*}
where $I, a$ and $b$ are constants, and state variables $u$ and $v$ describe the activation and inhibition mechanisms in a neuron. For chimeras, we consider a network of $N$ nonlocally coupled FitzHugh-Nagumo (FHN) oscillators,
\begin{align*}
    \epsilon \dot{u_k} &= u_k - \frac{u^3_k}{3} - v_k + \sigma \sum{j = 1}^N A_{kj} [b_{uu} (u_j - u_k) + b_{uv} (v_j - u_k)  ],\\
    \dot{v_k} &= u_k + a + \sigma \sum_{j = 1}^{N} A_{kj} [ b_{vu}(u_j - v_k) + b_{vv}(v_j - v_k)],
\end{align*}
where oscillators are indexed by $k$, $\epsilon > 0$ is a small parameter that provides a time scale separation and $\sigma$ is the coupling strength. Additionally, the coupling of oscillators contain both direct $u-u$ and $v-v$, and cross-couplings $u-v$.

The term nonlocal coupling refers to interactions between oscillators that are non direct spatial neighbours. In the simplest case, we consider oscillators arranged in a one-dimensional ring topology with coupling terms of the following form,
\begin{equation*}
    \sigma \sum_{j = 1}^{N} A_{kj} = \frac{\sigma}{2R} \sum_{j= k - R}^{k+ R}1.
\end{equation*}
where $R$ is the number of neighbours each node is connected to. A convenient value to track is the coupling $r = \frac{R}{N}$. A value of $r=\frac{1}{N}$ corresponds to nearest-neighbour (local) coupling, and $r = 0.5$ is all-to-all coupling/global (i.e. fully connected network). Intermediate values of $r$ corresponds to nonlocal coupling.

\begin{figure}
    \centering
    \includegraphics[width=0.97\linewidth]{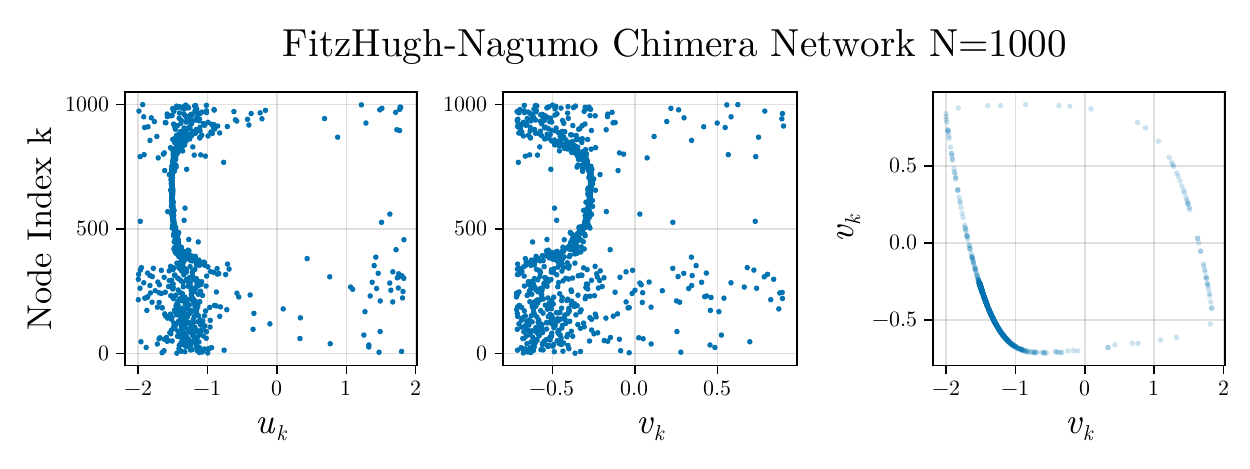}
    \caption{Simulation of the FHN chimera network with 1000 nodes. Parameters are from \cite{scholl2016synchronization}: $\epsilon = 0.05, a = 0.05, r = 0.35, \sigma=0.1, \phi = \pi /2 -0.1$. Initial conditions are randomly distributed uniformly on the circle $u_k^2 + v_k^2 = 4$. Nodes $500\sim700$ appear as a synchronous cluster, and other nodes are asynchronous. }
    \label{fig:FHN_chimera}
\end{figure}

Numerical simulation of the model yields the the results in Figure \ref{fig:FHN_chimera}. Notably, for carefully chosen parameter values, FHN oscillator ring networks with nonlocal coupling yields a contiguous cluster of oscillators that are fully synchronised coexisting with a separate group of asynchronous (low order) oscillators.

% \subsection{Information Flow and Network Inference}
% SECTION TO ADD IN THE FUTURE

\section{Opinion Dynamics}
\subsection{A dynamical networks approach to social physics}
An alternative but nonetheless illustrative example of a dynamical network model is in social physics and opinion dynamics. The study of opinion dynamics and consensus formation is an active area of research that aims to model and understand the driving forces and mechanisms of human interaction in  group behaviour. Broadly speaking, the general network formulation of opinion dynamics is given as follows. Let an opinion network be described by a group of $N$ agents whose individual opinions at a given time $t$ are represented by a $k$-dimensional vector $\mathbf{x}_i(t) \in \mathbb{R}^k$ , and an evolution operator $\phi$ that maps current opinions to future states. Consequently, the opinion state of the social network at time $t$ can be written as a matrix $X_t \in \mathbb{R}^{N\times k}$ where,
\begin{equation}
    X_{t+1} = \phi(X_t).
\end{equation}
Given the above formulation, the classic opinion dynamics problem aims to understand the effects of $\phi$ on opinion formation in the network as the system approaches steady state. The generality of the above formulation, primarily in the definition of $\phi$, allows for the extreme flexibility in the range of dynamics that may be produced. The classical formulations of $\phi$ commonly assumes the case of a fully connected structure (i.e. nodes interact with all other nodes). However, these ideas may be generalised to include non-fully connected random networks, often at the expense of analytical guarantees on convergence and consensus formation.

There are several characteristics that are of interest when studying opinion formation. A common feature that is often studied is the number of clusters in the steady state which may be categorised as one of three different outcomes: (1) consensus - where all nodes possess the same opinion, (2) polarisation - the formation of two dominant opinion clusters, and (3) fragmentation - more than two opinion clusters. Other characteristics that may be of interest may include convergence time, cluster sizes and distribution of opinions at the end state.

Numerous models have been proposed to model opinion dynamics based on different mechanisms. Some classical models include the voter model for discrete dynamics \cite{holme2006nonequilibrium, horstmeyer2020adaptive, sood2008voter}, and the DeGroot \cite{degroot1974reaching} and Friedkin-Johnsen \cite{friedkin1990social} models when considering probability distributions of opinions. An alternative model that considers continuous valued opinions are bounded confidence models proposed by Hegselman-Krause \cite{hegselmann2002opinion} and Deffuant \cite{deffuant2000mixing, weisbuch2004bounded}. These models account for the concept of \textit{selective exposure} \cite{hickok2022bounded, sears1967selective}, where agents are unlikely to interact with neighbours who hold vastly differing opinions to themselves. This is enforced by defining a threshold (radius/confidence bound) $c$ which determines if at any given time, whether a pair of interacting neighbours will compromise in their held opinions by some amount $\mu$. Generally, analytical results are largely intractable for even modest complications in the network such as random graphs (scale free, small world etc), heterogenous agent thresholds $c$ and asymmetrical confidence thresholds. As a result, analysis of these cases often rely on Monte Carlo numerical simulations to obtain results. However, recent work by Meng et al. have provided extensive empirical results on the effect of various network topologies and the resulting opinion dynamics \cite{meng2018opinion}. Game theoretic approaches to understanding agent interactions have also been proposed by Di Mare et al. \cite{di2007opinion}.

\subsection{Asynchronous vs. Synchronous Simulations}
Similar to networks of oscillators, much of approaches to studying opinion dynamics require the usage of numerical simulations due to their mathematical intractability. However, unlike oscillators, opinion dynamic models differ in that interactions occur at discrete time steps (difference equation) rather than in a continuous manner (ODEs). Simulating these networks often employs an agent-based approach. This episodic nature of opinion updates presents two different ways in which simulations may be done.

Synchronous updates are a computationally simple and cheap way to simulate opinion dynamics. In this approach, all agents update their state simultaneously based entirely on the knowledge of the previous time step. For linear models such as the DeGroot (DG) and Friedkin-Johnsen (FJ) models, a synchronous approach can greatly simplify the coding required to implement (i.e. iterate a matrix equation). Synchronous update routines tend to exhibit more rapid convergence behaviour as every time step is guaranteed a state update. However, one criticism is that such an approach lacks realism as individuals do not (normally) coordinate to update their opinions at discrete time steps.

Asynchronous updates are a more ``continuous'' approach to simulating agent-based/opinion models, and is particuarly applicable when studying opinion dynamics occurring on top of non-regular network topologies (i.e. anything that is not a ring, lattice, fully connected or cannot be described by a continuous spatial approximation). They are implemented by randomly selecting an agent (or edge) and performing a local state update to that age whilst keeping other agents constant. This is process is repeated until convergence is achieved. Unsurprisingly, such an approach is far more computationally expensive than synchronous updates often with slower convergence as time steps correspond to smaller units of time and scales inversely with the number of agents. Observed dynamics can differ greatly between synchronous and asynchronous simulations, even for the same class of model and after accounting for the different in time step length.

\subsection{DeGroot Model and Friedkin-Johnsen model}
One of the earliest models for opinion formation was proposed by DeGroot et a.. \cite{degroot1974reaching}. The model consists of a linear matrix model given as follows,
\begin{equation*}
    \mathbf{o}^{(t+1)} = \mathbf{A}\mathbf{o}^{(t)} = \dots = \mathbf{A}^{(t+1)} \mathbf{o}^{(0)},
\end{equation*}
where $\mathbf{o}^{(t)}\in \mathbb{R}^{n}$ is a vector of real values opinion and $\mathbf{A} \in \mathbb{R}^{n \times n}$ is a doubly stochastic matrix corresponding to the level of influence that any given agent/individual receives from others. Therefore, an agent updates their opinion as a weighted sum of their current opinion and connected neighbours'. This model is relatively simple an allows for analytical tractability with dynamics characterised by $\mathbf{A}$.

An extension to the DeGroot model is the Friedkin-Johnsen (FJ) model. The FJ model alters the evolution equation by including terms to describe the behaviour of stubborn agents,
\begin{equation*}
    \mathbf{o}^{(t+1)} = \mathbf{DA}\mathbf{o}^{(t)} + \mathbf{(I-D)}\mathbf{o}^{(0)},
\end{equation*}
where $\mathbf{D}$ is a diagonal matrix with non-zero entries corresponding to an individual's susceptibility to changing their opinions. The addition of the second term describes the attachment of an individual to their initial opinions and thus indicate a level of stubbornness in agents. In both DG and FJ models, updates are applied synchronously and thus individuals only have knowledge of their neighbours' most recently observed opinion state. 

\subsection{Voter Model}
In contrast to the the deterministic DG and FJ models, the voter model is a stochastic opinion model that may be simulated using either a synchronous or asynchronous scheme. For a given network of $N$ nodes, each node may possess an binary opinion $o_i^{(n)} \in \{0,1\}$. In the classic formulation, agents follow three rules:
\vspace{0.5mm}
\begin{enumerate}
    \item If asynchronous, pick a voter at random
    \item Voter adopts the state of a random neighbour
    \item Repeat until consensus is reached
\end{enumerate}
\vspace{2mm}
More generally, voter models can be used to describe agents whose binary state flips according to some transition function $c(x,\beta)$ where $x$ is a location/node and $\beta$ is the current overall state/configuration of the network (i.e. the state at $x$ is given by looking up $\beta(x)$).

\subsection{Boundary Confidence Model (BCM)}

Bounded confidence models account for the phenomenon where agents may be more inclined to interact and be influenced by neighbours who are of a more similar opinion to them \cite{noorazar2020classical}. This is enforced by defining a threshold (radius/confidence bound) $c$ which determines if at any given time, whether a pair of interacting neighbours will compromise in their existing opinions by some amount $\mu$. There are two main classes of bounded confidence models that are widely studied. The original proposed approach, termed the Deffuant-Weisbuch (DW) model \cite{deffuant2000mixing}, utilises an asynchronous update rule where only one pair of agents' opinions are updated in each time step. Namely, at each time step a random edge corresponding to a pair of neighbouring nodes $i,j$ are selected and their opinions change according to the following equation,

\begin{equation}
    x_i(t+1) = f(x_i(t), x_j(t)),
\end{equation}

\begin{equation}
    f(x_i(t), x_j(t)) = \begin{cases}
                            x_i(t) + \mu(x_j(t) - x_i(t)),\quad & |x_i(t) - x_j(t)|<c\\
                            x_i(t), \quad & \text{otherwise}
                        \end{cases}.
    \label{eq:5_BCM}
\end{equation}

An alternative formulation is the Hegselman-Krause (HK) model where all nodes' opinions are updated synchronously at each time step based on the opinions in the previous time step \cite{hegselmann2002opinion}. In this formulation, each agent considers the weighted average opinion $\bar{x}_i(t)$ of all its neighbours and updates its opinion accordingly,

\begin{equation}
    f(x_i(t), \bar{x}_j(t)) = \begin{cases}
                            x_i(t) + \mu(\bar{x}_i(t) - x_i(t)),\quad & |x_i(t) - \bar{x}_i(t)|<c\\
                            x_i(t), \quad &\text{otherwise}
                        \end{cases}.
\end{equation}

This alternative formulation allows for a more inclusive interaction strategy as agents with more extreme opinions may be still partially considered, provided that $\bar{x}_i$ lies within the confidence threshold \cite{noorazar2020classical}. The BCM may also be altered to account for heterogeneous agents by allowing for different threshold values for each agent.

Apart from the simplest network structures (e.g. lattice, fully connected, ring) where simplifications and mean field approaches may be applied, analytical models for BCM behaviour are often difficult or intractable. Such simple network structures are rare in large social networks as individuals are unlikely to have direct connections with all other members. Instead, other non-trivial network structures such as rich clubs, cliques and hubs may be present. These features can skew the influence of individuals and thus change the resulting convergence behaviour. Therefore, most experiments rely on computer simulations to provide insights.

The functional form of BCM inherently results in agents either forming consensus or fragmentation as any two interacting agents will either be static or converge towards each other. Generally, smaller thresholds $c$ and larger compromise rates $\mu$ tend to result in quicker convergence. In the intermediate case, network topology affects the characteristic features such as the convergence time and resultant number of clusters \cite{meng2018opinion}.

\chapter{Ordinal Time Series Analysis}\label{chap:OrdinalTSA}

\section{Symbolic Dynamics}
Many analytical tools and phase space analysis approaches are useful for describing the behaviour of well-defined dynamical systems. However, these approaches because become more difficult to apply as system dynamics increase in complexity. In real world applications where only numerical observation data (e.g. time-series) is present, analytical approaches that rely on the presence of mathematical structure become all but easily applicable.

The traditional approach to tackle this conundrum is to rely on the principles of embedding theorem in the hopes of reconstructing the trajectory in a state space that is diffeomorphic to the real state space (if it exists) vis-à-vis delay embedding methods. If one wishes to learn the dynamics for tasks such as prediction, one could employ a host of model construction methods to approximate the evolution operator. However, the delay embedding approach augments data into the space of $\mathbb{R}^n$ which presents several numerical difficulties if the observed is noisy or limited in resolution.

One workaround to dealing with noisy and imperfect data is to instead analyse a reduced approximation of the dynamics (i.e. the symbolic dynamics). Recall in Chapter \ref{chap:DynamicalSystemsTheory}, symbolic dynamics refers to a reduction of dynamics on continuous-values (e.g. $\mathbb{R}^n$) to one on a finite set of values such that much of the interesting dynamical features (e.g. periodic orbits, homoclinic/heteroclinic orbits, stability) are approximately preserved. Namely, periodic sequences in the original system should correspond to periodic sequences in the symbolic dynamics. For example, one can choose to reduce the dynamics of the tent map and Bernouulli map defined on the unit interval into a sequence of two symbols $(L,R)$ (left as an exercise for the reader) and show that periodic orbits are preserved.

The dimensional reduction achieved through symbolic dynamics can afford significant computational efficiency and noise robustness. What is less clear is the requisite transformation to convert time series into symbolic dynamics. Generally, this transformation to symbolic dynamics may be done by employing an appropriate partition of state space with each symbol corresponding to a unique partition. However, not all partitions are created equal and the extent to which symbolic dynamics preserve the true dynamics can differ between choice of partitions. One therefore must take care when defining partitions.

\begin{definition} \textbf{-- Partitions}\\
    Let $\mathcal{S}$ be a set. A partition $\mathcal{P}(X)$ is defined as a collection of disjoint sets $\{ X_i \}_i$ such that
    
    $$\mathcal{S} = \bigcup_{i} X_i$$

\end{definition}

Some examples of constructing partitions of state space include the following:
\begin{enumerate}
    \item Voxel partitioning \cite{st2022analysis} (i.e. split state space into tiny boxes)
    \item Neasrest neighbour partitions \cite{tan2023network}
    \item Hyperplanes (e.g. sign of $x$-component in Lorenz, ordinal patterns)
\end{enumerate}

\section{Ordinal Methods}
\subsection{Ordinal representations of time series}
\begin{figure}
    \centering
    \includegraphics[width=\linewidth]{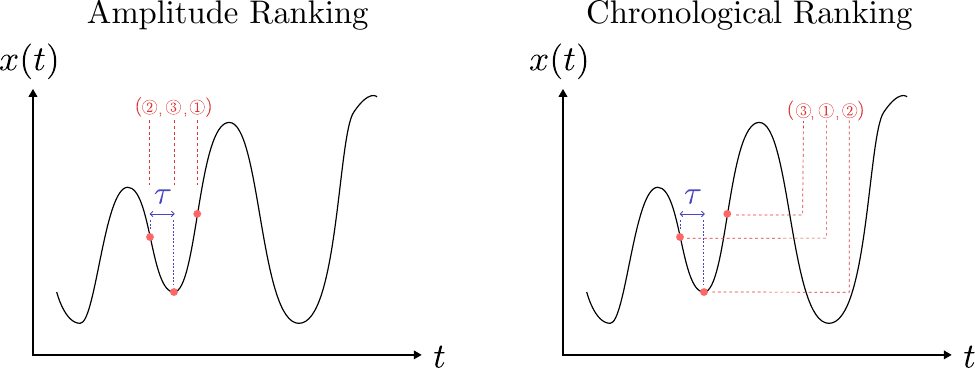}
    \caption{Two different ways to define ordinal symbols from a sequence of values. Amplitude ranking vs. Chronological ranking}
    \label{fig:OrdinalSymbols}
\end{figure}

Ordinal methods aim to utilise the rank ordering of observed values in time series to construct a discrete representation. To do this, we consider a sequence of $m$ points in a time series each separated by a lag $\tau$. This sequence of points is then mapped to a vector of integers representing an ordinal symbol, denoted as $\pi_i$ that describes the ordinal ranking of the magnitudes in the $m$ points (see Figure \ref{fig:OrdinalSymbols}).

There are two commonly employed but functionally equivalent methods for defining ordinal symbols. Amplitude ranking constructs a ranking such that the $m$\textsuperscript{th} integer in the symbol corresponds to the magnitude ranking of the $m$\textsuperscript{th} element in the sequence of observations. Chronological ranking uses the time ordering of the sequence to define the symbol where the $m$\textsuperscript{th} of the ordinal symbol is time index of the $m$\textsuperscript{th} largest element in the sequence. In both methods, applying the transformation procedurally across moving windows of length $\tau(m-1)$ results in a sequence of ordinal symbols $\pi(n)$, termed the ordinal time series, where $\pi(n) \in \mathcal{X}$ and $\mathcal{X}$ is the collection of all possible ordinal symbols. Notably, $\mathcal{X}$ is finite with $|\mathcal{X}|=m!$.

\subsection{Ordinal partitions}
The astute reader will notice numerous similarities in the construction of the ordinal symbol with that of the delay vector in uniform delay embedding. Simply put, an $m$ order ordinal symbol set corresponds to a partitioning of delay reconstructed phase space into $m!$ regions using a set hyperplanes that pass all pass through the origin. As a special case, order $m=2$ ordinal symbols encode information about whether a given subsequence of the time series is increasing or decreasing. Furthermore, it is possible to show that construction of ordinal symbols yields a proper complete partition if rank ordering ties are broken systematically\\

To show this, let $x(t)$ be a time series and $\vec{x}(t)$ be an $m$-dimensional delay vector constructed in the conventional way,
\begin{equation*}
    \vec{x}(t) = (x(t), x(t-\tau), \dots, x(t-(m-1)\tau)) = (x_{1}(t), x_{1}(2), \dots , x_{m}(t)),
\end{equation*}
where $m \in \mathbb{N}, m>2$ and $\vec{x}(t) \in \mathcal{M} \subseteq \mathbb{R}^m$. We denote $m$ as the order of the ordinal construction, and $\mathcal{M}$ as the reconstructed attractor manifold. Furthermore, let $\mathcal{X} = \{ \pi_{i} \}_{i=1}^{m!}$ be the set of ordinal partitions. By definition, each ordinal symbol $\pi_i$ represents a particular rank ordering of delay vector elements that fulfill the following set of nested inequalities,
\begin{equation}
    x_{i_1} \geq x_{i_2} \geq \dots \geq x_{i_m},
    \label{eq:ordinal}
\end{equation}
where $\{i_1,\dots, i_m \} \in \{1,\dots,m\}$ are index of rank ordered elements in the delay vector. We may then define the partition corresponding to $\pi_i$ as,
\begin{equation*}
    W_i = \{ \vec{x} \in \mathbb{R}^m \, |\,  x_{i_1} \geq x_{i_2} \geq \dots \geq x_{i_m}\}.
\end{equation*}
As \ref{eq:ordinal} consists of $m-1$ separate linear inequalities, $W_i$ is thus a wedge region with boundaries defined by $m-1$ hyperplanes all of which pass through the origin. Furthermore, 
\begin{equation*}
    W_i \cap W_j \neq \emptyset \quad \iff \quad \exists \, i,j \;\; \text{s.t.} \;\;x_i = x_j,
\end{equation*}
and $\mathcal{P}=\bigcup_{i} W_i = \mathbb{R}^m$. I.e. ordinal symbols describe a partition of $\mathbb{R}^m$ up to the boundaries of each partition, where these boundaries correspond to equally valued element in the delay vector component. To construct a proper partition such that $W_i \cap W_j = \emptyset$ for all $i\neq j$, we can systematically break ties according to the element index such that if $x_i = x_j$ where $i>j$, we assign $\vec{x}$ to the closest partition corresponding to condition $x_i > x_j$.

Because the manifold on which dynamics occur $\mathcal{M}$ lies on a subset of $\mathbb{R}^m$, it is not guaranteed that trajectories visit all possible partitions. In practice, not all of the $m!$ ordinal symbols are required to fully describe the dynamics. These unused symbols are also described as forbidden symbols. Therefore, for a noiseless system it is sufficient to consider instead a smaller set of ordinal symbols $\hat{\mathcal{P}}\subset \mathcal{P}$ whose partitions span $\mathcal{M}$ instead. As noise is introduced, observed trajectories may appear to enter previously forbidden partitions that are not part of $\mathcal{M}$ and thus increase the number of utilised ordinal symbols.

\section{Network Approaches}
\subsection{Ordinal Partition Transition Networks}
\begin{figure}
    \centering
    \includegraphics[width=0.65\linewidth]{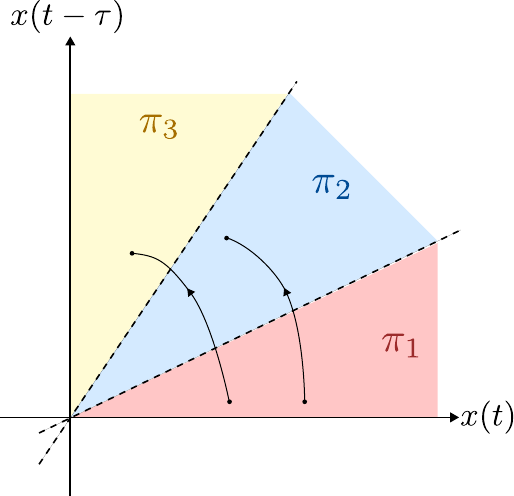}
    \caption{Forward mapping or ordinal symbols with finite number of partitions is not guaranteed to be unique. Example trajectory shown in figure. Third axis for $x(t-\tau)$ not shown for illustrative purposes.}
    \label{fig:ordinalmap}
\end{figure}

Similar to the ideas of the Poincar\'{e} map, dynamics in the original system is reduced to a map with the added constraint of having a discrete finite domain,
\begin{equation*}
    \pi (n+1) = f_{m}(\pi (n)),\qquad \pi (n) \in \mathcal{X}_{(m)},
\end{equation*}
where $\mathcal{X}_{(m)}$ is the set of all possible ordinal symbols generated from an $m$-order partition. For discrete time observations (i.e. most real world data) the dynamics of the map $f_{m}$ are approximate to the true dynamics as $m\to \infty$. This corresponds to infinitely granular wedge shaped partitions of the reconstructed state space. For finite $m$ with discrete observations, $f_{m}$ is neither guaranteed to reproduce the dynamics nor be a function. It is possible for two different trajectories that originate from the same partition to map into different partitions (see Figure \ref{fig:ordinalmap}).
Accepting this limitation, one may choose to abandon all hopes for a deterministic description and instead employ a more stochastic approach to understanding the dynamics. To do this, we can consider using the transition probabilities between symbols,
\begin{equation*}
    p(\pi(n+1) = \pi_i\, | \, \pi(n) = \pi_{i_n},\pi(n-1)=\pi_{i_{n-1}}\dots).
\end{equation*}
Because the next step probabilities are dependent on the historical observed trajectory of symbols, the behaviour of the transitions do not necessarily have the Markov property (i.e. next state transition probabilities are purely dependent on my current state). However, the Markov property can be assumed if the underlying dynamics are ergodic. Therefore, we may fully describe the stochastic approximation of the dynamics using only the one step transition probabilities,
\begin{equation*}
    p_{ij} = p(\pi (n+1) = \pi_j \, | \, \pi(n) = \pi_i).
\end{equation*}
In practice, we can numerically approximate values of $p_{ij}$ by enumerating the frequency of observed transitions between any two ordinal symbols $\pi_i$ and $\pi_j$. Once calculated, these probabilities may be represented in the form of an ordinal transition network whose connectivity is a weighted matrix $A \in \mathbb{R}^{|\mathcal{X}| \times |\mathcal{X}|}$ with entries $A_{ij} = p_{ij}$. Here each node represents a single unique ordinal symbol, and directed edges correspond to the transition probabilities between pairs of symbols.

Similar to forbidden symbols, some transitions between particular pairs of symbols may also never occur for a given dynamical system. These transitions are thus called forbidden transitions and can be used as a preliminary measure of the dynamical complexity of given time series. Simpler dynamics such as periodic orbits only utilise a small fraction of the possible ordinal symbols, as well as the number of possible transitions. These fractions increase for more ``complex'' systems such as chaotic dynamics, and fully saturate in the case of purely stochastic noise where all possible symbols and transitions may be observed. 

\section{Complexity Measures}
Suppose we require a way to distinguish between periodic, chaotic and random data. If we consider these three different regimes as lying on a spectrum, one can describe each case as increasing in dynamical complexity. Periodic signals are least complex, and random i.i.d noise are the most complex. The complexity of a time series can also be related to the predictability of the time series.

As their name suggests, complexity measures are a way to quantify and compare the complexity of different time series. A simple example that have been discussed thus far are Lyapunov spectra, where positive values indicate SDIC conditions. More sophistacated examples include the various suite of recurrence plot measures, and reservoir computing complexity measures.

Apart from comparing two different systems, one of them main practical uses of complexity measures is in regime classification and change point detection. Time series observations of systems in different dynamical regimes should result in characteristically different values of complexity (e.g. modulation of a bifurcation parameter between chaotic and periodic regimes, or change point from sinus rhythm to ventricular tachycardia in ECG). Complexity measures are a way to transform the time series into a form that is more amenable to the traditional classification and change point detection methods by extracting and summarising key dynamical features of the time series. In practice, a suite of complexity measures are used in tandem to classify and detect dynamical changes.

\subsection{Permutation Entropy}
One can use the number and frequency distribution of ordinal symbols calculate the complexity of a given time series. More regular and periodic time series utilise a smaller subset of the possible $m!$ collection of ordinal symbols. In contrast, higher complexity time series may venture into a wider variety of partitions and thus utilise a larger number of ordinal symbols. Bandt \& Pompe \cite{bandt2002permutation} proposed the idea of permutation entropy, that leverages on the frequency distribution of ordinal symbols to quantify complexity. Permutation entropy can be calculated as,
\begin{equation*}
    h^{PE} = -\sum_{i=1}^{m!} p_{i} \log_2 p_{i},
\end{equation*}
where $p_{i}$ is the frequency of observing the $i$\textsuperscript{th} ordinal symbol in the time series. This expression is an application of the Shannon entropy calculation to probability distribution of ordinal symbols. High values of permutation entropy indicate a more diverse set of ordinal symbols, and thus implies a time series that is more complex. One notable property of the permutation entropy is that it is bounded below by the Kolmogorov-Sinai entropy $h_{KS}$ for the general case where $m\to \infty$. Therefore, permutation entropy can also be interpreted as as a proxy measure of the predictability in a time series.\\

\subsection{Conditional Permutation Entropy}
An extension of the permutation entropy is the conditional permutation entropy that was proposed by Unakafov \cite{unakafov2014conditional}. Instead of just calculating the frequency distributions of ordinal symbols, the conditional permutation entropy focuses on the transition probabilities between pair of ordinal symbols,
\begin{equation*}
    h_{CPE} = -\sum_{i=1}^{m!}p_{i} \sum_{j=1}^{m!} p_{ij} \log_2 (p_{ij}),
\end{equation*}
where $p_{ij}$ is the probability of moving from symbol $\pi_i$ to $\pi_j$. This entropy measure was proposed as a way to address the inadequacies for the standard permutation entropy $h_{PE}$ to estimate the Kolmogorov-Sinai entropy $h_{KS}$ for finite $m!$. Whilst having the computational complexity of permutation entropy, $h_{CPE}$ is a better estimator for $h_{KS}$ with equality being true for periodic dynamical systems $(h_{KS} = 0)$ given a sufficiently large $m!$.

\section{Ordinal Poincar\'{e} Sections}
\subsection{Poincar\'{e} sections and first return maps}
\begin{figure}
    \centering
    \includegraphics[width=\linewidth]{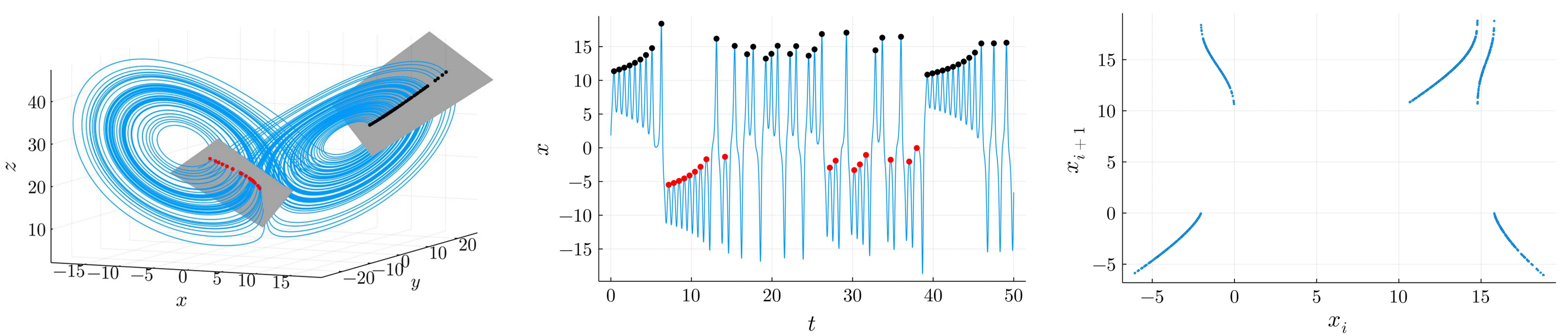}
    \caption{Lorenz Poincar\'{e} map taken from \cite{shahriari2023ordinal}. From left to right: (1) Poincar\'{e} section from the successive maxima method, (2) maxima on the $x$-component of the timeseries, (3) successive maxima first return map}
    \label{fig:lorenz_pmap}
\end{figure}

As briefly discussed, first return maps (Poincar\'{e} maps) provide a convenient way to reduce the dimensionality by one when analysing a continuous dynamical system. However, defining a Poincar\'{e} section required to construct the first return map is not necessarily easy nor equally performant for all possible choices of section. One conventional method for constructing a section is to take the successive values of maxima in the time series as values in the map sequence. This approach was first used by Lorenz in 1963 to define a Poincar\'{e} section that cuts across the lobes of the Lorenz attractor (see Figure \ref{fig:lorenz_pmap}). This technique can also be used to construct first return maps for other systems such as the R\"{o}ssler system to show that its dynamics belong to a family of one-hump maps.

\subsection{Ordinal methods for Poincar\'{e} sections}
An alternative way to define Poincar\'{e} section using ideas from ordinal time series was presented by Shahriari et. al. \cite{shahriari2023ordinal}. The method of ordinal Poincar\'{e} sections relies on the the hyperplane boundaries of ordinal partitions to define a large family of Poincar\'{e} sections. The algorithm proceeds as follows:
\begin{enumerate}
    \item Let $x(t)$ be a time series. For a given order $m$ and delay lag $\tau$, calculate the corresponding ordinal partition, set of ordinal symbols $\pi_1,\pi_2,\dots, \pi_n$, and occurrence frequency $K_1,K_2,\dots,K_n$.
    \item Separate the full time series $x(t)$ according to their corresponding ordinal sequences to construct $n$ different time series $\hat{x}_1(t),\hat{x}_2(t),\dots,\hat{x}_n(t)$.
    \item For each component time series $\hat{x}_i(t)$, construct another ordinal partition with ordinal symbols $\hat{\pi}_1,\hat{\pi}_2,\dots,\hat{\pi}_m$ and occurrence frequency $\hat{p}_{i_1},\hat{p}_{i_2},\dots,\hat{p}_{i_m}$.
    \item Calculate the weighted permutation entropy $h_{W_{i}}$ for component time series $\hat{x}_i(t)$ with weights given by the occurrence frequency of $K_i$ in the original complete time series $x(t)$,
    \begin{equation*}
        h_{W_{i}} = -\sum_{j=1}^{m}K_{i}p_{i_j} \log _2 (K_{i} p_{i_j} )
    \end{equation*}
    \item One can also calculate the entrance weighted permutation entropy using weights $\hat{K}_1,\hat{K}_2,\dots,\hat{K}_n$ calculated on the frequency of first entrances into a given ordinal partition (i.e. if there is a string of time series observations that have the same ordinal symbol $\pi_i$, this only contributes a count of one). Thus the entrance weighted permutation entropy for a given ordinal time series component $x_i(t)$ is calculated as,
    \begin{equation*}
        h_{EW_{i}} = -\sum_{j=1}^{m}\hat{K}_{i}p_{i_j} \log _2 (\hat{K}_{i} p_{i_j} )
    \end{equation*}
    \item Rank each ordinal symbol $\pi_i$ according to its weighted (or entrance weighted) permutation entropy. For each symbol $\pi_i$, find the series of time indices $\{ t^{i}_{1}, t^{i}_{2},\dots  \}$ corresponding to the first entrances into the ordinal partition for $\pi_i$. The first return map corresponding to the ordinal symbol $\pi_i$ is given as $(x(t^{i}_{n}, x(t^{i}_{n+1}))$
\end{enumerate}

\begin{figure}
    \centering
    \includegraphics[width=\linewidth]{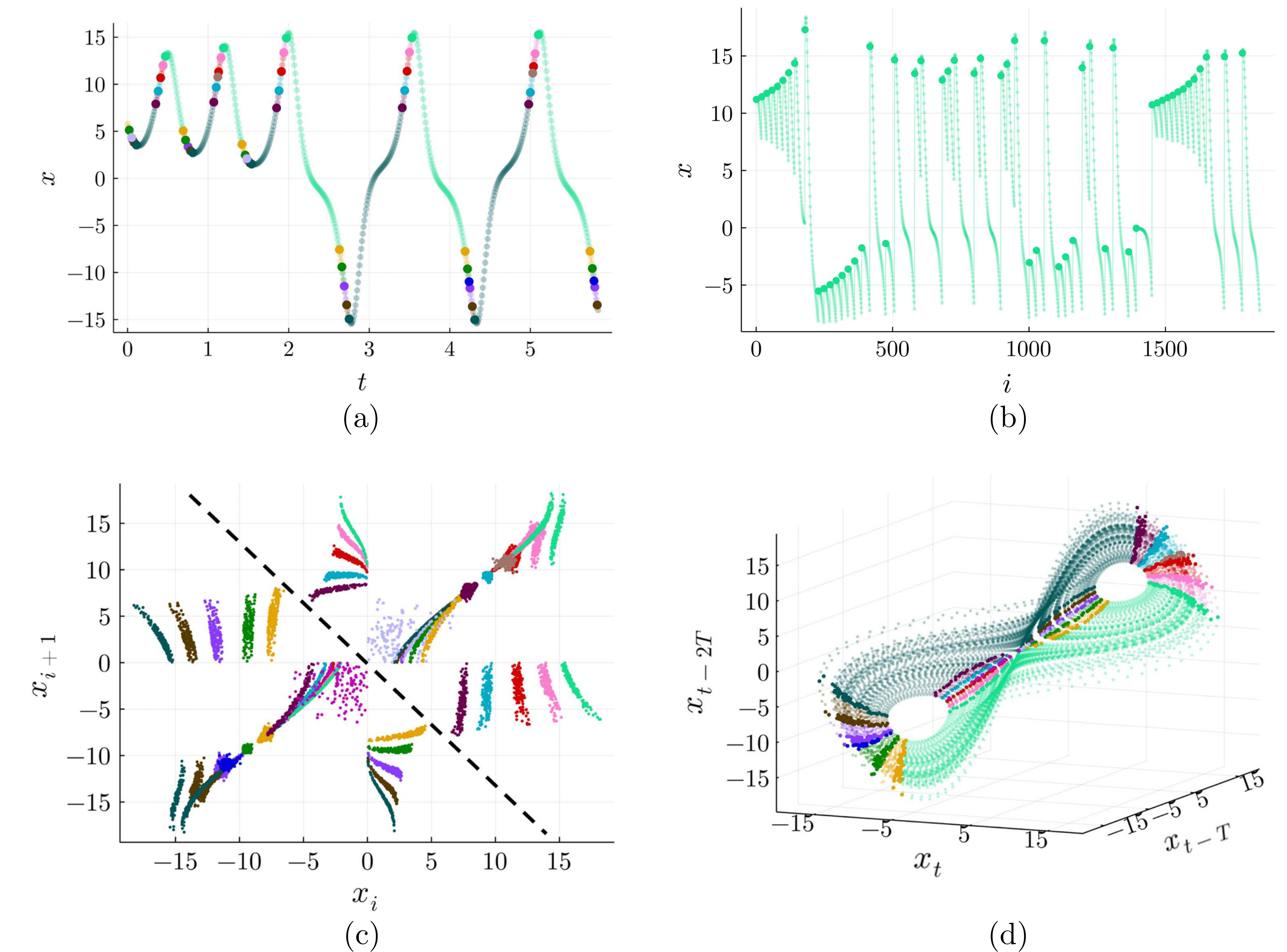}
    \caption{Construction of the ordinal Poincar\'{e} sections: (a) Time series coloured according to ordinal symbol with entrance points bolded. (b) Extracted component timeseries corresponding to the monotonical decreasing ordinal sequence $\pi_1 = (4,3,2,1)$. (c) family or first return maps coloured according to ordinal symbol. (d) Lorenz time delay attractor with entrance point bolded and coloured according to ordinal symbol.}
    \label{fig:ordinalmap_construction}
\end{figure}

\begin{figure}
    \centering
    \includegraphics[width=\linewidth]{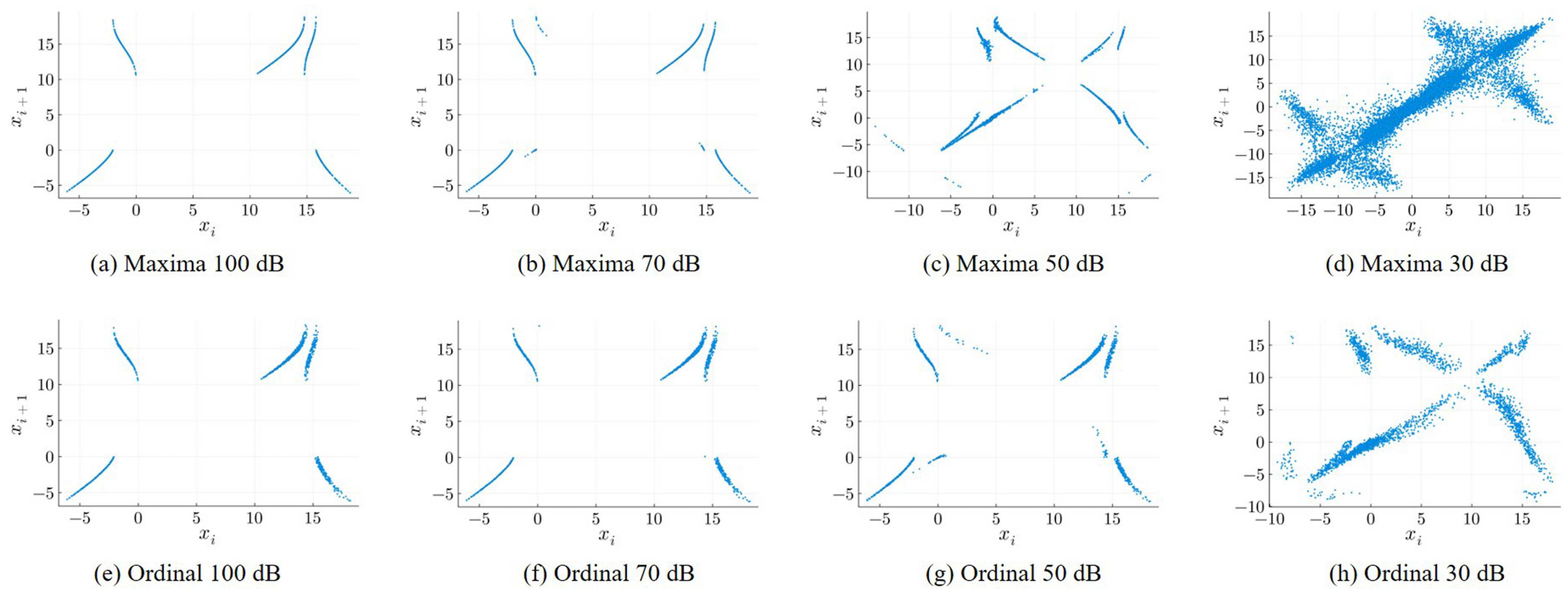}
    \caption{Constructed first return maps for the Lorenz system at various signal to noise ratios. Noise applied is additive Gaussian. Top row (a-d) are max constructed using the successive maxima approach, bottom row (e-g) are those using the ordinal Poincar\'{e} section.}
    \label{fig:ordinalmap_robust}
\end{figure}

The ordinal Poincar\'{e} sections method produces a whole family of candidate Poincar\'{e} sections that can be used to construct first return maps. The weighted and entrance weighted entropies are used as a selection criteria where the maps corresponding to ordinal symbols with the highest entropy values are preferred. This choice results in Poincar\'{e} sections which cut across attractors in a smooth and well-defined way. The first return maps constructed using ordinal Poincar\'{e} sections are also more robust to noise when compared against the conventional successive maxima first return maps (see Figure \ref{fig:ordinalmap_robust}).

\chapter{Exercises and Computational Experiments}

\begin{enumerate}[itemsep=1em]

\item (Linearisation) Consider the following 2D dynamical system with real constant parameters $\delta, \mu$, where $\delta \neq 0$. How does the stability of the fixed point at the origin vary with $\delta$ and $\mu$?
\begin{eqnarray*}
\dot{x} &=& -\delta x -\mu y + xy\\
\dot{y} &=& \mu x -\delta y + \frac{1}{2}(x^2-y^2)
\end{eqnarray*}

\item (Centre manifolds) Recall that the characteristic equation for centre manifolds of flows is given by 
\begin{eqnarray*}
    \mathcal{N}(h(x)) = Dh(x)[ Ax + f(x,h(x))]-Bh(x)-g(x,h(x)) =  0
\end{eqnarray*}

Derive a similar characteristic equation for calculating centre manifolds on discrete maps.

\item (Stability) Let $\bar{x}(t)$ be a solution to an autonomous ODE system with $\bar{x}(0) = x_0$. The positive orbit through $x_0$ for $t > t_0$ is defined as:

$$
O^{+}(x_{0}, t_{0}) = \{ x \in \mathbb{R}^{n} |\, x = \bar{x}(t), \,t \geq t_{0}, \, \bar{x}(t_{0}) = x_{0} \}
$$

Following this, consider the following stability definitions:\\

\textit{\textbf{Orbital Stability} -- $\bar{x}(t)$ is orbitally stable if given $\epsilon > 0$, $\exists \delta > 0$ such that for any other solution $y(t)$ where $|x(t_{0}) - y(t_{0})| < \delta$, then $d(y(t), O^{+}(x_{0}, t_{0})) <\epsilon$ for $t>t_{0}$\\}

\textit{\textbf{Asymptotic Orbital Stability} -- $\bar{x}(t)$ is asymptotically orbitally stable, if it is orbitally stable and $\exists b > 0 $ such that, if $|x(t_{0}) - y(t_{0})| < b$, then}
$$\lim_{t \rightarrow \infty} d(y(t), O^{+}(x_{0}, t_{0})) = 0.$$

Show that all trajectories of the following vector field are asymptotically orbitally stable:
\begin{eqnarray*}
\dot{\theta} &=& -\sin^{2} \theta + (1-r^{2})\\
\dot{r} &=& r(1-r),
\end{eqnarray*}
where $(\theta,r) \in S^{1}, \mathbb{R}$. 

\item (Centre manifolds) Consider the following dynamical system with bifurcation parameter $\epsilon$, where $\epsilon$ is small,
\begin{eqnarray*}
\dot{x} &=& 2x+2y\\
\dot{y} &=& x+y+x^4 + \epsilon y^2
\end{eqnarray*}
\begin{enumerate}
    \item Show that the system has one fixed point for $\epsilon \geq 0$, and three when $\epsilon<0$
    \item Calculate the linearised dynamics about the fixed point $(0,0)$ and show that it is non-hyperbolic and that linearisation is insufficient. Hence, conclude that it must either be an unstable fixed point, or a saddle.

    \item Use the eigenvalues and eigenvectors found in part (b) to perform a linear change of coordinates to $(u,v)$. Show that the resulting system is given by

    \begin{eqnarray*}
    \begin{bmatrix}
           \dot{u} \\
           \dot{v}
     \end{bmatrix}
     = 
     \begin{bmatrix}
       0 & 0 \\
       0 & 3
     \end{bmatrix}
     \begin{bmatrix}
           u \\
           v
     \end{bmatrix}
     + \frac{1}{3}
     \begin{bmatrix}
           -2 [(u+2v)^4 + \epsilon(-u+v)^2]  \\
           (u+2v)^4 + \epsilon(-u+v)^2
     \end{bmatrix}
    \end{eqnarray*}
    where $u$ and $v$ are the centre and and unstable components respectively.

    \item Find a $4^{\text{th}}$ order approximation for the centre manifold.
    \begin{eqnarray*}
    h(u) &=& a_1 u^2 + a_2 u \epsilon + a_3 \epsilon^2 + a_4 u^3 + a_5 u^2 \epsilon + a_6 u \epsilon ^2 \\
    &+& a_7 \epsilon^3 + a_8 u^4 + a_9 u^3 \epsilon + a_{10} u^2 \epsilon^2 + a_{11} u \epsilon^3 + a_{12} \epsilon^4 
    \end{eqnarray*}

    Note: For systems with constant parameters $\epsilon$. The $n^{\text{th}}$ order approximation includes all terms where the sum of the powers of $u$ and $\epsilon$ are less than or equal to $n$. It is a good idea to proactively remove terms that exceed order $4$ whenever they arise to reduce the number of terms. 

    \item Using the centre manifold, describe the stability of the fixed point at the origin. How does it change with varying $\epsilon$?

    \item Draw the global flow of the dynamical system in the $(u,v)$ coordinates. How does the behaviour change with $\epsilon$?
\end{enumerate}

\item (Recurrence) Show that the tripling map $f(x_n) = 3x_n \text{\,mod\,} 1$ is a length preserving function on the unit interval. Verify that this does not imply $\mu(U) = \mu (f(U))$, where $\mu$ is the Lebesgue measure (i.e. length) and $U$ is a subinterval of $[0,1]$.

\item (Non-autonomous systems) Consider the following non-autonomous dynamical system
\begin{eqnarray*}
\dot{x} &=& x(x+a)(x-a) + q(t)
\end{eqnarray*}

\begin{enumerate}
    \item Let $q(t) = 0$. Find all the fixed points. For those that are stable, provide the basin of attraction.
    
    \item Find the conditions for $q(t)$ required for some stable basin of attraction to exist.
    
    \item Assuming $q(t)$ adheres to the conditions in part (b). Let $\Phi_{\infty}(x_{0}|q(0) = q_{0})$ long time evolution of a point with initial condition $x(t) = x_{0}$ and initial perturbation $q(t) = p_0$. Does $\Phi_{\infty}(x_{0}|q(0) = q_{0}) \implies \Phi_{\infty}(x_{0}|q(0) = p_{0})$ for all $q_{0} \neq p_{0}$?

    \item Let $q(t) = \frac{-2a^{3}}{3\sqrt{3}}  \cos (\omega t + \phi )$. Generate the time series and omega limit sets for various values of $\omega$ and $phi$ across multiple initial conditions. Do the findings align with your expectations from part (c)?
\end{enumerate}

\item \label{q:lorenz_embed}(Embedding) Perform a numerical integration of the Lorenz equations (see below) with the standard parameter values ($\sigma=10$, $\beta=\frac 83$, $\rho=28$).
    \begin{eqnarray*}
    \dot{x} &=& \sigma (y-x)\\
    \dot{y} &=&x(\rho-z)-y\\
    \dot{z} &=&xy-\beta z
    \end{eqnarray*}
    demonstrate that (for ``suitably chosen'' embedding parameters) time delay reconstruction from the $x$ and $y$ component alone lead to an object topologically equivalent to the original attractor. However, show that the $z$ component alone of the Lorenz equations does not lead to a valid embedding. Provide an explanation for this. Furthermore, propose a way for reducing the Lorenz time series into symbolic dynamics.

\item (Embedding) Consider the $x$ component of the chaotic Lorenz system (see Q\ref{q:lorenz_embed}).
    
    \begin{enumerate}
        \item Implement an algorithm produces a global principal value embedding.
        \item Using the Julia {\tt DynamicalSystems} package or your own implementation, calculate the ideal embedding lag based on the first minimum of the auto mutual information. Compare the reconstructed attractor from PCA and uniform delay embedding.
    
        \item Consider a $m$-dimensional PCA embedding constructed with maximum lag $\tau_{max}$. Would the components be equivalent to a uniform $m$-dimensional embedding with lag $\tau = \frac{\tau_{max}}{m}$?
    \end{enumerate}

\item (SDIC) Consider the Bernoulli map,
    
    \begin{eqnarray*}
    x_{n+1} = \begin{cases}
    2x_n, \quad &x_n < 0.5\\
    2x_n -1, \quad &x_n \geq 0.5
    \end{cases}
    \end{eqnarray*}

    \begin{enumerate}
        \item Simulate the Bernoulli map with various rational and irrational initial conditions. Plot the trajectories and verify that rational iniital conditions always yield periodic behaviour.
        
        \item Let $a_{0} = \frac{\pi}{20}$ and $b_{0}$ be equal to $a_0$ but rounded to 10 decimal places. Show that the separation $a_{n}$ and $b_{n}$ scales approximately exponentially with base $2$. (You might want to plot on a logarithmic scale).
    
    \end{enumerate}

\item (Information Theory, Algebra) Prove the following identity for mutual information in terms of the joint and marginal probabilities

\begin{align*}
I(x,y) &= \sum_{i,j}p_{i,j}^{(xy)} \ln p_{i,j}^{(xy)} - \sum_{i}  p_{i}^{(x)} \ln p_{i}^{(x)} -\sum_{j} p_{j}^{(y)} \ln p_{j}^{(y)} \\
    &= \sum_{i,j} p_{i,j}^{(xy)} \ln \frac{p_{i,j}^{(xy)}}{p_{i}^{(x)}p_{j}^{(y)}}
\end{align*}

\item (Model Selection, Algebra) Let there be observed data $x_y$ and $y_t$, and a class of models $\phi_\theta$, where $\theta$ are a set of $k$ model parameters. Suppose the optimum set of parameters $\hat{\theta}$ has already been found such that
\begin{eqnarray*}
    y_t = \phi_{\hat{\theta}}(x_t) + \epsilon_t, \, \epsilon_t \overset{i.i.d}{\sim} \mathcal{N}(0, \sigma^2).
\end{eqnarray*}
Show that the AIC for this model is given (up to a constant) by the form

\begin{eqnarray*}
    AIC(\phi_{\hat{\theta}}) = 2k + n \log \left( \frac{SSE}{n} \right)
\end{eqnarray*}

\item (Invariants) Calculate the fractal dimension of the Menger sponge (googling will help.)

\item (Invariants) Show that any affine transformation of a bounded set $\mathcal{B} \in \mathbb{R}^n$ for $n \in \mathbb{Z}^+$ does not alter the correlation dimension of that set.

\item (Correlation Integrals) Show that the generalised dimension of order $q=1$ (i.e. in the information dimension) is given by the following expression.

\begin{eqnarray*}
    D_1 = \lim_{\epsilon \to 0} \frac{\langle \ln p_\epsilon \rangle_\mu}{\ln \epsilon}
\end{eqnarray*}

\item (Entropy) Show that $H(p)$ satisfies the property that $H(p) < -\sum_{i =1}^n p_1 \log q_i$ for probability distributions $p = (p_1,...,p_n)$ and $q = (q_1,...,q_n)$ if $p \neq q$.

\item (Neural Networks) Consider a general feedforward neural network that maps $\mathbb{R} \to \mathbb{R}$. Let $n_l$ and $n_h$ be the number hidden layers and nodes per hidden layer respectively. Assuming that the model errors after training are Gaussian, and that the precision of all the parameters are constant and not optimised. Provide an expression for the model description length in terms of $n_l$ and $n_h$. Suppose the mean error scales exponentially with the number of hidden nodes. How would this change the final expression?

\item (Information density) English and Indonesian are three languages that contain very different written linguistic structures. For example, Indonesian, a Malayic language is contains common usage of prefixes and suffixes in order to construct new words and semantics.\\

    Consider the following translations of Article 1 in the Universal Declaration of Human Rights:\\

    \textbf{English}\\
    All human beings are born free and equal in dignity and rights.
    They are endowed with reason and conscience and should act towards one another in a spirit of brotherhood.\\

    \textbf{Indonesian}\\
    Semua orang dilahirkan merdeka dan mempunyai martabat dan hak-hak yang sama.\\
    Mereka dikaruniai akal dan hati nurani dan hendaknya bergaul satu sama lain dalam semangat persaudaraan.\\

    % \textbf{Japanese (just Hiragana)}

    % すべて　の　にんげん　は、うまれながら　に　して　じゆう　で　あり、　かつ、
    % そんげん　と　けんり　と　に　ついて　びょうどう　で　ある。\\
    % にんげん　は、　りせい　と　りょうしん　と　を　さずけられて　おり、\\
    % たがい　に　どうほう　の　せいひン　を　もって　こうどう　しなければ\\
    % ならない。\\

    % \textbf{Japanese (with Kanji)}\\
    % すべて　の　人間　は、　生まれながら　に　して　自由　で　あり、　かつ、\\
    % 尊厳 　と　権利　と　に　ついて　平等　で ある。　人間　は、　理性　と　良心　と　を　授けられて　おり、　互い　に　同胞　の　精神　を　もって　行動　しなければ ならない。\\

    % Write a program automatically construct a Huffman code from a text and a set of codeword frequencies. Calculate the average code-word length (in bits) and the encoding length of each of the above passages. How do the lengths of each of the passages compare against each other? How would using Kanji in Japanese, or the frequent usage of prefixes and suffixes affect the number of codewords their relative frequencies? How would this affect the information density of a text of given length?

    Write a program automatically construct a Huffman code from a text and a set of codeword frequencies. Calculate the average code-word length (in bits) and the encoding length of each of the above passages. How do the lengths of each of the passages compare against each other? How would the frequent usage of prefixes and suffixes affect the number of codewords their relative frequencies? How would this affect the information density of a text of given length?

    \item (Prediction) Integrate a trajectory of the Lorenz equations with random initial conditions. With the resulting time series, implement some code to perform a nearest neighbour prediction and apply to the integrated time series with (a) no embedding, (b) 2D delay embedding, (c) 3D delay embedding for various embedding lags. Comment on your findings. What happens when you decrease the amount of training data and/or add observational Gaussian noise?

    \item (ECG, experimental) Using your nearest neighbour predictor, try and perform one-step and multi-step freerun predictions for the \texttt{ECG.csv} dataset. How well does it perform? Are there any patterns to when the predictor fails? What are the potential causes for a nearest neighbour predictor to fail. If you had the time, what are some ways that you could construct a better predictor.

\item (Prediction) Integrate a trajectory of the Lorenz equations with random initial conditions. With the resulting time series, implement some code to perform a nearest neighbour prediction and apply to the integrated time series with (a) no embedding, (b) 2D delay embedding, (c) 3D delay embedding for various embedding lags. Comment on your findings. Are there particular areas where predictions tend to fail and what is their cause? What happens when you decrease the amount of training data and/or add observational Gaussian noise? Provide some suggestions on altering the algorithm that could improve predictions.

\item (Challenge - Lyapunov Exponents) Implement Wolf's algorithm (without the assistance of pre-built dynamical systems packages) and numerically estimate the maximal Lyapunov exponent for 2 continuous chaotic systems of your choice. Construct two versions of the algorithm: (a) the standard approach with pairs of trajectories, and (b) by considering neighbourhoods. Make estimations for two cases: (i) using data from the full state vector, (ii) using a delay reconstruction. Compare your results against the established Lyapunov exponents from the literature (or using existing pre-built tools) and discuss your findings. What are the challenges associated with the algorithm?

\item (ECG, experimental) Using your nearest neighbour predictor, try and perform one-step and multi-step freerun predictions for the \texttt{ECG.csv} dataset. How well does it perform? Are there any patterns to when the predictor fails? What are the potential causes for a nearest neighbour predictor to fail. If you had the time, what are some ways that you could construct a better predictor. 

\item (Ordinal Time Series Analysis) Recall the Lorenz dynamical system given by
    \begin{eqnarray*}
    \dot{x} &=& \sigma (y-x)\\
    \dot{y} &=&x(\rho-z)-y\\
    \dot{z} &=&xy-\beta z
    \end{eqnarray*}

    \begin{enumerate}
        \item Numerically integrate the equations for different time step sizes\\$\Delta t = (0.001, 0.05, 0.01)$

        \item Generate a plot of the first return map (FRM) using a single component of the time series.

        \item Using the time series you generate in part (a), apply an ordinal partition with dimension $m=4$. What are the occurrence probabilities of each ordinal partition?

        \item Choose one of the highest weights of the ordinal partitions, and plot the FRM for each time series. What does the FRM for the different time series and their components look like?

        \item Construct a delay embedding of the $x$-component time series into 3 dimensions with an appropriately chosen lag. How can you explain the difference between the FRMs based on the location of the corresponding points to this specific ordinal partition on the embedded attractor?
    \end{enumerate}

\item (Surrogates) Consider the R\"{o}ssler chaotic dynamical system operating in the phase given by
    \begin{eqnarray*}
        \dot{x} &=& -y-z\\
        \dot{y} &=& x+ay\\
        \dot{z} &=& b + z(x-c)
    \end{eqnarray*}
    \begin{enumerate}
        \item The R\"{o}ssler system exhibits two unique regimes: phase coherent $(a,b,c)=(0.165,0.4,8.5)$ and non-phase coherent $(0.265,0.4,8.5)$. Integrate with these values and describe the geometric and dynamical differences that can be observed from the time series and its reconstructed attractor. How do the Fourier spectra of the signals from these two systems differ? You may wish to search and read some articles to find out more and discuss your findings.

        \item Write an algorithm to construct AAFT surrogates, and pseudo-periodic surrogates from the $x$-component of the non-phase coherent system. Describe the differences in between the two surrogate time series.
    \end{enumerate}

\item (Reservoir Computing) Implement a working reservoir computer and use it to simulate the Lorenz system.
    
    \begin{enumerate}
        \item Write a function \texttt{create\_ESN} which takes in hyper-parameters and returns the relevant matrices for an ESN. Assume the input to the ESN is \textit{scalar}.
        \item Write a function \texttt{run\_ESN} which takes in a scalar input signal $x$ and the relevant matrices to generate the activation states of an ESN.
        \item Train your ESN to take in a trajectory of the $x$-component of the Lorenz system and output the same component one time step in the future. Use a time step of $\Delta t=0.05$ when generating your Lorenz trajectory. You must choose and clearly state sensible hyper-parameters, however they do not need to be optimised. Plot the prediction on an independent testing series and report the error as a correlation.
        \item Reproduce Figures \textbf{31} and \textbf{32} in Chapter \ref{chap:reservoir_comp} and hence state the memory capacity of your reservoir computer.
        \item Vary the spectral radius of your reservoir computer and describe the impact on both the memory capacity and performance error in the prediction task.
        \item Use your reservoir computer to generate a free-run prediction of the $x$-component of the Lorenz system.
    \end{enumerate}
    
\item (Recurrence Plots) Below is a bifurcation diagram generated from the $x$-component of the Ikeda map,
    \begin{align*}
        x_{n+1} &= 1 + r(x_n \cos(t_n) - y_n \sin(t_n) ) \\
        y_{n+1} &= r(x_n \sin(t_n) + y_n \cos(t_n) )
    \end{align*}
    with $t_n = 0.4 - \frac{6}{1+x_n^2+y_n^2}$.
    
    Traces of the mean, standard deviation and minimum value across the parameter range are shown in red, blue and green respectively.
    \begin{center}
        \includegraphics[width=0.7\textwidth]{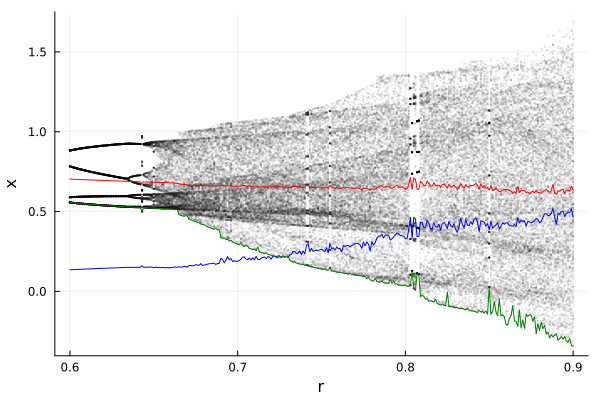}
    \end{center}
    \begin{enumerate}
        \item Generate recurrence plots from the $x$-component of the Ikeda map for parameter values $r\in \{0.60, 0.65, 0.70\}$ without embedding and comment on your observations.
        Repeat this process with embedding.
        \item Reproduce the bifurcation diagram of the system and plot traces of a scalar diagonal line, vertical line and recurrence time feature from RQA.
        Comment on the ability of the features to detect regime changes.
        
        \textit{NOTE: The final figure should be a bifurcation diagram with three traces, however the features may need to be standardised to align the different scales.} 
        
        \item Repeat the above process with three scalar features from other nonlinear time-series analysis fields (e.g. reservoir time series analysis, symbolic analysis, persistent homology). They may all be from the same analysis field or three different ones.
    \end{enumerate}
    \textit{NOTE: To make plotting more manageable, use only 200 points for the bifurcation diagrams but 1000 points (or more) for calculating features. Don't forget to remove transients from your trajectories.}

    \item (MDL) Suppose we have a linear model with parameters $\theta = (\theta_1,\theta_2,...\theta_k)$. The model output is $X\theta$, and the desired model output is $Y$ where
    \begin{equation*}
        Y = (y(t_1), y(t_2), ... y(t_n))^T \text{ , } X = 
        \begin{bmatrix}
            x_1(t_1) & x_2(t_1) & ...  & x_k(t_1) \\
            x_1(t_2) & x_2(t_2) & ...  & x_k(t_2) \\
            ... & ... & ...  & ... \\
            x_1(t_n) & x_2(t_n) & ...  & x_k(t_n) \\
        \end{bmatrix}
        \text{ .}
    \end{equation*}

    In order to find the description length of this model we must approximate the description length of its errors ($E = Y-X\theta$) and its parameters ($\theta$). We approximate the length of the errors using the differential entropy of the error probability distribution. Assuming that the errors are normally distributed with zero mean and variance $\sigma^2$, the description length of $n$ error points is given by 
    \begin{equation*}
        L(E|\theta)\approx\frac{n}{2}\log(2\pi e\sigma^2)\text{ where }\sigma^2 = \frac{E^TE}{n}\text{ .}
    \end{equation*}
    If the parameters are encoded as floating-point numbers to precision $\gamma$ (I.e. $\theta_i = 0.a_1a_2...a_{\log(\gamma)}\times e^{2m}$) and we ignore their sign, any self-delimiting prefixes and the cost of encoding $2m$ then the length of the parameters ($L(\theta)$) is just $-k\log(\gamma)$. By adding a precise perturbation** $\delta\in\mathbb{R}^k$ we can write $\theta_i$ as $0.a_1a_2...a_{(\log(\delta_i)-2m)}\times e^{2m_i}$ and reduce $L(\theta)$ by $\sum_i(\log(\gamma)-\log(|\delta_i|)+2m_i)$, I.e
    \begin{equation*}
        L(\theta + \delta) = L(\theta)+k\log(\gamma)-\sum_i(\log(|\delta_i|)-2m_i)
    \end{equation*}
    The optimal rounding perturbation minimises the sum of $L(\theta+\delta)$ and $L(E|\theta+\delta)$, i.e.:
    \begin{align*}
        \frac{\partial}{\partial\delta_i}\left(L(E|\theta)+\delta^T\frac{\partial L(E|\theta)}{\partial\delta}+\frac{1}{2}\delta^T\frac{\partial^2 L(E|\theta)}{\partial\delta^2}\delta\right)\\
            +\frac{\partial}{\partial\delta_i}\left( L(\theta)+k\log(\gamma)-\sum_i(\log(|\delta_i|)-2m_i) \right) =0
    \end{align*}
    The minimum total description length of the model (which occurs when $\delta$ satisfies eq 3) is given by
    \begin{align*}
        L_{tot} &= L(E|\theta+\delta_{opt})+L(\theta+\delta_{opt})\\
        &= \frac{n}{2}\log(2\pi e\sigma^2)+\delta^T\frac{\partial L(E|\theta)}{\partial\delta}+\frac{1}{2}\delta^T\frac{\partial^2 L(E|\theta)}{\partial\delta^2}\delta\\
        &- {\partial\delta^2}\delta-\sum_{i=1}^k(\log(|\delta_i|)\text{  (ignoring the 2m) .}
    \end{align*}
    If linear regression is used to find the parameters $\theta$, equation 3 can be simplified, and the total description length will be similar to the expression in the unit reader.Consider the description length if the model is instead trained using ridge regression with penalty $\alpha$. In the case of ridge regression, $\theta$ satisfy
    \begin{equation*}
        \frac{\partial}{\partial\theta}\left(E^TE+\alpha\theta^T\theta \right)=0\text{ .}
    \end{equation*}

    **This is only really the case if $\delta_i$ happens to align exactly with the difference between the true parameter values and their values rounded to lower precision. However, assume that it is true for all $\delta\in\mathbb{R}^k$ including when $|\delta_i|>|\theta_i|$.  

    \begin{enumerate}
        \item What is the solution for $\theta$ to the ridge regression equation (eq 4)?
        \item What are $\frac{\partial L(E)}{\partial\theta}$ and $\frac{\partial^2 L(E)}{\partial\theta^2}$? (you will need eqs 1 and 4)
        \item Find an expression for $\delta$ by simplifying Eq 3 (Hint: most terms in Eq 3 have no $\delta$ dependence)
        \item Find an expression for $L_{tot}$ and simplify it as much as possible.
        \item What happens in the limits $\alpha \rightarrow 0$ and $n\rightarrow\infty$?
    \end{enumerate}

    \item (Dynamical Networks) Consider an dynamical network model for continuous opinion dynamics. Recall that the 2 opinion BCM model is given as follows:
    $$
    x_{i}(t+1) = \begin{cases}
        x_{i}(t) + \alpha(x_j(t)-x_{i}(t)), &\, |x_{i}(t) - x_{j}(t)|< \rho\\
        x_{i}(t), &\, |x_{i}(t) - x_{j}(t)|\geq \rho\\
    \end{cases}
    $$
    \begin{enumerate}
        \item Propose a way to extend the above opinion update model to account for 3 different opinion directions. How would you generalise this for the case of $n$ opinions?
        
        \item Consider a 3 opinion dynamical network to represent voter sentiment in an election, corresponding to three political groups: Labor, LNP and Independent. Implement an agent based simulation of the BCM model with your proposed model in section (a) with two different configurations, (i) synchronous and (ii) asynchronous updates. Present simulation for varying network topologies and confidence thresholds $\alpha$.
        
        \item Assume that all individuals subsequently choose to vote in an election at each time step and all votes are equally valued. How do election outcomes differ if using first-past-the-post vs. preferential voting, and for various values of $\alpha$? (You will need to choose a sensible way to determine the vote preference of an individual based on their opinion. How do these results change as the number of different opinions $n$ increase?
    \end{enumerate}
\end{enumerate}

\backmatter%%%%%%%%%%%%%%%%%%%%%%%%%%%%%%%%%%%%%%%%%%%%%%%%%%%%%%%
% \include{author/glossary}
% \include{author/solutions}
% \printindex

\printbibliography[heading=bibintoc,title={\textbf{References}}]

%%%%%%%%%%%%%%%%%%%%%%%%%%%%%%%%%%%%%%%%%%%%%%%%%%%%%%%%%%%%%%%%%%%%%%

\end{document}